%% file: main_expert.tex
\documentclass[11pt,preprint]{imsart}

\input{./style.tex}

\alglanguage{pseudocode}

\begin{document}

\begin{frontmatter}

	\title{Optimal Permutation Estimation in Crowd-Sourcing problems}
	\runtitle{Optimal ranking}
	\begin{aug}
		\author{Emmanuel Pilliat, Alexandra Carpentier, and Nicolas Verzelen}
		\runauthor{Pilliat et al.}
	\end{aug}

	\begin{abstract}
	Motivated by crowd-sourcing applications, we consider a model where we have partial observations from a bivariate isotonic $n\times d$ matrix with an unknown permutation $\pi^*$ acting on its rows. Focusing on the twin problems of recovering the permutation $\pi^*$ and estimating the unknown matrix, we introduce a polynomial-time procedure achieving the minimax risk for these two problems, this for all possible values of $n$, $d$, and all possible sampling efforts. Along the way, we establish that, in some regimes, recovering the unknown permutation $\pi^*$ is considerably simpler than estimating the matrix.  
	\end{abstract}

\end{frontmatter}

\maketitle

\section{Introduction}

We consider a crowd-sourcing problem with $n$ experts and $d$ questions. For an unknown matrix $M$, $M_{i,j}\in [0,1]$ stands for the ability of expert $i$ at question $j$. For the purpose of calibrating the model, we receive partial and noisy observations of the matrix $M$ and our goal is to rank the experts according to their ability. Earlier models in crowd-labelling problems or in the related problems of pairwise comparisons typically assumed that the matrix $M$ belongs to a parametric model~\cite{bradley1952rank,luce2012individual,thurstone1927law,dawid1979maximum,braverman2008noisy}, a prominent example being Bradley-Luce-Terry model. While there has been significant progress in this direction,  such models do not tend to fit well real-world data~\cite{mclaughlin1965stochastic,ballinger1997decisions}. 

To address this issue, there has been a recent interest in the class of permutation-based models~\cite{chatterjee2015matrix,shah2015estimation,shah2016stochastically,mao2020towards,chatterjee2019estimation,flammarion2019optimal,pananjady2022isotonic,shah2020permutation} where it is only assumed that the matrix $M$ satisfies some shape-constrained conditions before one (or two) permutations acts on the rows (and possibly on the columns) of $M$. Quite surprisingly, it has been established in~\cite{shah2016stochastically} that, at least in some settings, the matrix $M$ can be estimated at the same rate in those non-parametric models as in classical parametric models by relying on the least-square estimator on the class of permuted bi-isotonic matrices.
Unfortunately, the corresponding  class of matrices is highly non-convex and no polynomial-time algorithm is known for computing this least-square estimator. 
 Furthermore, known computationally efficient procedures such as spectral estimators~\cite{chatterjee2015matrix,chatterjee2018matrix} only achieve sub-optimal convergence rates. This has led several authors to conjecture the existence of computational-statistical trade-offs~\cite{flammarion2019optimal,shah2019feeling}. Despite recent progress in this direction~\cite{mao2020towards,liu2020better}, the fundamental limits of polynomial-time algorithms for this class of problems remain largely unknown.

Arguably, for most applications, the primary objective is to recover the underlying permutation $\pi^*$ acting on the rows or equivalently to rank the experts accordingly. While estimation of the full matrix $M$ is closely related to ranking, it is also of a quite different nature as argued below. In this work, we investigate the estimation of the permutation $\pi^*$ by characterizing the minimax risk for estimating $\pi^*$ in a permuted shape-constrained model, introducing a polynomial-time procedure nearly achieving this risk bound. As a byproduct, we also disprove the existence of a computational-statistical gap for the reconstruction of the matrix $M$.

\subsection{Problem formulation}

A bounded matrix $B\in [0,1]^{n\times d}$ is said to be bi-isotonic if it satisfies $B_{i,j}\leq B_{i+1,j}$ and $B_{i,j}\leq B_{i,j+1}$ for any $i\in [n-1]$ and $j\in [d-1]$. Henceforth, we write $\mathbb{C}_{\text{BISO}}$ for the collection of such $n\times d$ bounded  bi-isotonic matrices.  

In this work, we assume that the matrix $M$ is a row-permuted bi-isotonic matrix as in~\cite{mao2020towards,liu2020better}. In other words, up to a single permutation $\pi^*$ of $[n]$, the matrix $M_{\pi^{*-1}}$ defined by $(M_{\pi^{*-1}})_{i,j}= (M_{\pi^{*-1}(i),j})$ is bi-isotonic. From a modeling viewpoint, this amounts to assuming that the $d$ questions are ordered from the most difficult question to the most simple question. The permutation $\pi^*$ is not necessarily unique, but the corresponding permuted matrix $M_{\pi^{*-1}}$ is unique. 
Despite that, we refer, with a slight abuse of terminology, to $\pi^*$ as the oracle permutation. With this definition, $\pi^{*-1}(i)$ corresponds to any $i$-th smallest row (or equivalently expert to use the crowd-sourcing terminology) in the matrix $M$. In the following, the $i^{\mathrm{th}}$ row of $M$ is referred to as \emph{expert} $i$, whereas the $k^{\mathrm{th}}$ column is referred to as  \emph{question} $k$.

We consider an observation-scheme where the statistician has partial access to noisy observations $Y$ of $M$ such that 
\beq\label{eq:model_0}
Y= M + E\ , 
\eeq
where the entries of $E$ are centered, independent, subGaussian - see definition 2.2 of \cite{wainwright2019high} - with Orlicz norm at most $\zeta$, but are not necessarily identically distributed. In particular, this model encompasses binary observations $Y_{i,k}\sim Ber(M_{i,k})$ which arise in crowd-labelling problems, in which case we have $\zeta=1$. 
In the following, we refer to $\zeta$ as the noise level.

As usual in the literature --e.g.~\cite{mao2020towards}, we use the Poissonization trick to model the partial observations. Given some $\lambda>0$, which is henceforth referred as the sampling effort, we have $N= Poi(\lambda nd)$ observations of the form
\beq\label{eq:model_partial}
(x_t,y_t),\quad t=,1\ldots, N,
\eeq
where the position $x_t$ is sampled uniformly in $[n]\times [d]$, and $y_{t}=M_{x_t}+ E_{x_t}$ is
an independent observation of matrix $Y$ of~\eqref{eq:model_0} at position $x_t$. Conditionally to $N$, this scheme is equivalent to a uniform sampling scheme with replacement~\cite{mao2018minimax}. If $\lambda < 1$, then a specific entry of $Y$ is sampled at least once with probability $1-e^{-\lambda}$ which is close to $\lambda$. More generally, $\lambda$ corresponds to the expected number of times a specific entry of $Y$ is observed, so that $\lambda >1$ would correspond to the situation where entries are sampled multiple times.

Since our aim is to recover the permutation $\pi^*$ from the partial observations $(x_t,y_t)$, we consider, for some estimator $\hat \pi$, the following error metric 
\beq
\label{eq:permloss}
l(\hat{\pi}; \pi^*)= \|M_{\hat \pi^{-1}} - M_{\pi^{*-1}} \|_F^2 \enspace , 
\eeq
where $\|.\|_F$ stands for the Frobenius norm. This loss quantifies the distance between the matrix $M$ ordered according to the oracle permutation $\pi^*$ and the matrix $M$ ordered according to the estimated permutation. When $\pi^*$ is not unique, the error $l(\pi^*,\pi^{*'})$ between any two oracle permutations is zero. If $\hat{\pi}$ and $\pi^*$ only differ by a transposition or equivalently if the ranking $\hat{\pi}$ and $\pi^*$ only differ on two experts, then $l(\hat{\pi}; \pi^*)$ is twice the square Euclidean distance between the corresponding rows of $M$. More generally, $l(\hat{\pi}; \pi^*)$ interprets as the sum over all $i=1,\ldots,n$ of the square Euclidean distance between the $i$-th smallest row of $M$ according to $\hat{\pi}$ and according to the oracle ranking $\pi^*$.

\medskip 

The loss~\eqref{eq:permloss} is ubiquitous when one aims at estimating the matrix $M$ in Frobenius norm, that is building an estimator $\widehat{M}$ such that $\|\widehat{M}-M\|_F^2$ is as small as possible --see e.g.~\cite{shah2016stochastically,mao2020towards,liu2020better}. Indeed, estimating $\pi^*$ by $\hat{\pi}$ is a first step towards building an estimator of $M$ by doing as if $M_{\hat{\pi}^{-1}}$ was bi-isotonic. It turns out that the error in $\|\widehat{M}-M\|_F^2$ decomposes as the sum of two terms, one of them being $l(\hat{\pi},\pi^*)$ while the other one does not really depend on $\hat{\pi}$. Conversely, an estimator $\widehat{M}$ can be easily transformed into an estimator $\hat{\pi}$ whose loss $l(\hat{\pi},\pi^*)$ is controlled by $\|\widehat{M}-M\|_F^2$. See~\cite{shah2016stochastically,mao2020towards} for further discussions.
In summary, controlling $l(\hat{\pi}; \pi^*)$ is important in order to evaluate to what extent $\pi^*$ is well estimated, but it is also the key stepping stone towards a good estimation of the matrix $M$.

In some works, the authors directly consider distances on the symmetric group of permutations. Examples include the Kendall tau distance $d_{KT}(\pi,\pi')= \sum_{(i,j): \pi(i)<\pi(j)}\1\{\pi'(i)>\pi'(j)\}$ or the $l_{\infty}$ distance $\|\pi-\pi'\|_{\infty}= \max_{i\in[n]}|\pi(i)-\pi'(i)|$ --see~\cite{braverman2008noisy,mao2018minimax} in the noisy sorting model. However, those distances are not well suited to handle the non-parametric class of bi-isotonic matrices, because to control them we would need to make assumptions on the separation between the rows of the matrix $M$ --see Appendix A of~\cite{shah2016stochastically}.

Equipped with this notation, we consider the minimax risk of permutation recovery as a function of the number $n$ of experts, the number $d$ of question, the sampling effort $\lambda$, and the noise level $\zeta$. 
\beq\label{eq:minimax_risk_perm}
\cR^*[n,d,\lambda, \zeta]= \inf_{\hat{\pi}}\sup_{\pi^*\in \Pi_n}\sup_{M:\,  M_{\pi^{*-1}}\in \mathbb{C}_{\text{BISO}}} \E\left[\|M_{\hat{\pi}^{-1}}- M_{\pi^{*-1}}\|_F^2\right]\ ,
\eeq
where $\Pi_n$ stands for the collection of all permutations of $[n]$.
In particular, our general aim is to tightly control this minimax risk and, if possible, to provide a computationally efficient procedure achieving this minimax risk. Although our primary interest lies in the estimation of $\pi^*$, we also consider the minimax estimation risk of $M$
\beq\label{eq:minimax_reconstruction}
\cR^*_{est}[n,d,\lambda, \zeta]	= \inf_{\hat{M}}\sup_{\pi^*\in \Pi_n}\sup_{M:\,  M_{\pi^{*-1}}\in \mathbb{C}_{\text{BISO}}} \E\left[\|\hat{M}- M\|_F^2\right]\ . 
\eeq
as  studied in~\cite{shah2016stochastically,liu2020better,mao2020towards,pananjady2020worst} in order to assess the performances of our computationally efficient procedures.

\subsection{Related work and open questions}

The most relevant body of work to the current paper is that on estimating square matrices $M$ satisfying the so-called strong stochastic transitivity class (SST)~\cite{chatterjee2015matrix,shah2016stochastically}. A matrix $M$ belongs to the SST class if (i) $M$ is skew-symmetric that is $M+M^T= e e^T$ where $e$ is the constant vector of size $n$ and (ii) there exists a common permutation $\pi^*$ of $[n]$ such that row and column-permuted matrix $M_{\pi^{*-1}\pi^{*-1}}$ is bi-isotonic. This class is suited for considering pairwise comparisons problems. Shah et al.~\cite{shah2016stochastically} consider the full observation setting, namely a setting where each entry of the matrix $M$ is observed once in noise - which is in some sense akin to $\lambda = 1$ in our Poissonian scheme\footnote{In the Poissonian scheme, each entry is observed at least once with probability $1-e^{-\lambda}$.}. They proved that the minimax risk for estimating $M$ in square Frobenius distance is, up to logarithmic terms, of the order of $n$ and is achieved by the corresponding least-square estimator over the SST class. Unfortunately, this estimator cannot be efficiently computed. They also analyzed an efficient spectral estimator achieving the rate $n^{3/2}$. This rate is also achieved~\cite{shah2016stochastically} by the near-linear time Borda count algorithm $CRL$ that simply ranks the individuals according to the row sums of the observations and then plugs the corresponding permutation to estimate $M$. See also~\cite{chatterjee2019estimation} for related results. 
This led some authors~\cite{flammarion2019optimal,shah2019feeling} to conjecture the existence of a $\sqrt{n}$ computational gap for SST matrices and for other shape-constrained matrices with unknown permutation.

In crowd-sourcing problems where $M\in [0,1]^{n\times d}$, non-parametric models~\cite{mao2020towards} assume that the matrix $M$ is bi-isotonic up to a permutation $\pi^*$ of the rows (experts) - and sometimes also up to a permutation $\tau^*$ of the columns (questions)\footnote{This would correspond to the situation where the corresponding ordering of the questions is also unknown.}.
In this paper as in this literature review, we focus however solely on the case where $M$ is bi-isotonic up to a permutation $\pi^*$ of the rows (experts). Mao et al.~\cite{mao2020towards} have established the minimax risk $\cR^*_{est}[n,d,\lambda,1]$ for estimating $M$ in the specific case where $n\geq d$. In the arguably most interesting regime of partial observations $\lambda\leq 1$, they prove that this minimax risk is of the order of $n/\lambda \wedge (nd)$. This rate is achieved by the inefficient least-square estimator. Furthermore, Mao et al.~\cite{mao2020towards} were the first to narrow the conjectured computational gap by introducing a new efficient procedure called one-dimensional sorting. In the square case $n=d$ with full observations, these procedures achieve (up to log terms) the rate $n^{5/4}$ for estimating the matrix $M$, thereby improving over the previous $n^{3/2}$ barrier.

Recently, this rate was improved by Liu and Moitra~\cite{liu2020better} in a specific instance of the problem  where $n=d$ and one has access to a sub-polynomial number of noisy independent samples of the complete matrix $M$ from~\eqref{eq:model_0} -- which is akin to our Poissonian scheme for $\lambda$ being sub-polynomial in $n,d$. They introduce a polynomial-time procedure achieving the rate $n^{1+o(1)}$ for permutation recovery and  matrix estimation which, up to the factor $n^{o(1)}$, turns out to be minimax optimal for both problems. As a consequence, in this very specific instance, the computational gap turns out to be nonexistent.

\medskip

There remain important open problems to characterize the  estimation of $\pi^*$ and $M$ in crowdsourcing problems. 
\begin{itemize}
    \item Beyond the case $n\geq d$ handled by Mao et al.~\cite{mao2020towards}, the minimax risk $\cR^*[n,d,\lambda, 1]$ of estimation of the permutation $\pi^*$ - as well as the minimax risk $\cR^*_{est}[n,d,\lambda, 1]$ of estimation of the matrix $M$ -  are unknown. In particular, in the rectangular case where $n\ll d$, the number of questions exceeds the number of experts is both relevant from a practical~\cite{shah2020permutation} and a conceptual perspective. Indeed, the analysis of the least-square estimator of~\cite{mao2020towards} and related works is based on entropy calculation of the class of permuted bi-isotonic matrices. While the minimax risk $\cR^*_{est}[n,d,\lambda, 1]$ turns out to be (up to logarithm terms), characterized by this entropy, this is not always the case for the estimation of $\pi^*$ as many matrices $M$ share the same permutation $\pi^*$. As a consequence, even if we leave aside computational constraints, pinpointing the optimal risk $\cR^*[n,d,\lambda, 1]$ for estimating $\pi^*$ requires quite different arguments.

    \item 	Beyond the toy "over-complete" observation model in the square case $n=d$ of Liu and Moitra~\cite{liu2020better}, it remains unclear whether there is a computational gap for general rectangular settings with partial observations. 
\end{itemize}.

\subsection{Our Contributions}\label{ss:contributions}

Echoing with these open problems, we make the following contributions in this work:
\begin{itemize}
    \item First, we  characterize (up to polylogarithmic multiplicative terms) the minimax risk $\cR^*[n,d,\lambda, \zeta]$ of permutation recovery, this, for all possible number of experts  $n\geq 1$, number of questions  $d\geq 1$, noise level $\zeta\geq 0$, and almost all sampling efforts $\lambda >0$.  When $n \ll d$, we prove in particular that $\cR^*[n,d,\lambda, 1] \ll \cR^*_{est}[n,d,\lambda, 1]$ in all non-trivial regimes, highlighting that when $n \ll d$, the problem of permutation recovery is statistically easier than the problem of matrix estimation.
	   \item  Moreover, we introduce a polynomial-time procedure achieving this risk bound, thereby establishing that there does not exist any significant computational-statistical trade-off for the problem of recovering a single permutation $\pi^*$. While our procedure borrows some of the ingredients of Liu and Moitra~\cite{liu2020better}, we need to introduce several new ideas to deal with the significantly more involved case $n \ll d$. Since an estimator $\hat{\pi}$ of $\pi^*$ can be easily combined with a least-square estimator of a bi-isotonic matrix to estimate the matrix $M$ --see e.g.~\cite{shah2016stochastically,mao2020towards} -- we also deduce a  polynomial time estimator $\widehat{M}$ which nearly achieves the minimax estimation risk $\cR^*_{est}[n,d,\lambda, 1]$, thereby proving that this problem does not either exhibit any computational-statistical trade-off, thereby answering the open problem of~\cite{mao2020towards}.
\end{itemize}

\medskip 
To provide a glimpse of our results, let us describe the minimax risks on the arguably most interesting case where the noise level $\zeta$ is of order $1$ as in the Bernoulli observation setting and where $\lambda < 1$ which corresponds to a partially observed matrix. In Section~\ref{sec:partial}, we establish that the minimax risk 
 $\cR^*[n,d,\lambda,1]$ of permutation recovery is (up to polylogarithmic multiplicative terms) of the order of 
\beq\label{eq:minimax_intro} 
\left[\frac{nd^{1/6}}{\lambda^{5/6}}\bigwedge \frac{n^{3/4}d^{1/4}}{\lambda^{3/4}}+\frac{n}{\lambda}\right]\bigwedge nd,
\eeq
whereas the  minimax reconstruction risk $\cR^*_{est}[n,d,\lambda, 1]$ is of the order
\beq\label{eq:minimax_intro_reconstruction} 
\left[\sqrt{\frac{nd}{\lambda}}\bigwedge \frac{nd}{\lambda^{2/3}(n\vee d)^{2/3}}+\frac{n}{\lambda}\right]\bigwedge  nd.
\eeq
We display in Figure~\ref{tab:sum} a summary of our results in the specific case where we also have $\lambda = 1$ on top of $\zeta = 1$, and will discuss this case more in details, as it highlights one of our main findings.

\begin{figure}
        \begin{tabular}{ c || c | c | c }
          
          ~ & $n \leq d^{1/3}$ & $d^{1/3} \leq n \leq d$ & $n \geq d$\\ \hline\hline
         Permutation estimation: $\cR^*[n,d,1,1]$ &  $nd^{1/6}$ & $n^{3/4}d^{1/4}$& $n$ \\ \hline
          Matrix estimation: $\cR^*_{est}[n,d,1,1]$ & $nd^{1/3}$ & $\sqrt{nd}$& $n$ \\ 
        \end{tabular}
	\caption{Summary of the minimax risks  (up to poly-logarithmic terms) for permutation estimation ($\cR^*[n,d,1,1]$) and matrix estimation ($\cR^*_{est}[n,d,1,1]$) in the specific cases where $\lambda,\zeta = 1$.}
	\label{tab:sum}
\end{figure}

A first comment is that the minimax risk of matrix estimation $\cR^*_{est}[n,d,1, 1]$ can be interpreted through the covering numbers of the space of permuted bi-isotonic matrices as in~\cite{mao2020towards}.  For $n\geq d$ both minimax risks - $\cR^*[n,d,1, 1]$, $\cR^*_{est}[n,d,1, 1]$ - are of the order of $n$ so that recovering the permutation $\pi^*$ is as hard as estimating the matrix $M$ (up to logarithmic factors). This is the regime studied in the literature, see~\cite{mao2020towards, liu2020better}.  When the number  $d$ of questions is large - $n \ll d$ - then the regimes are more tricky. There are two of them, depending on whether $n$ is larger than $d^{1/3}$ or not, and in both regimes $\cR^*_{est}[n,d,1, 1]$ is significantly larger than $\cR^*[n,d,1,1]$. More regimes appear when we do not restrict ourselves to $\lambda = 1$, $\zeta = 1$.
This complex picture, as well as the fact that $\cR^*[n,d,1,1] \ll \cR^*_{est}[n,d,1,1]$ for $n \ll d$ - and also in many other configurations of $\lambda,\zeta$ - highlights the fact that the difficulty of estimating $\pi^*$ is not governed by the size of the space of permuted bi-isotonic matrices. As a consequence, even if we leave computational aspects aside, it is not clear that the least-square estimator of~\cite{mao2020towards} achieves optimal risk for estimating the permutation $\pi^*$ and, in any case, entropy-based arguments would lead to suboptimal bounds, at least if we use the same arguments as in \cite{mao2020towards}.

As a byproduct of our results, we also establish the minimax risk - and prove that it is achievable in polynomial time - for another loss function termed $l_{\infty}(\hat{\pi},\pi^*)$ (see~\eqref{eq:loss}) put forward in~\cite{chatterjee2019estimation,shah2019feeling,mao2020towards} - and we also disprove a conjecture regarding a computational-statistical gap for this loss. See Subsection~\ref{ss:maxloss}.

\medskip 

As our minimax results remain valid in the noiseless case $(\zeta=0)$ where one has access to partial observation of the matrix $M$ itself, we are able to tightly decipher the approximation error which is due to the partial sampling of the matrix $M$ from the stochastic error stemming from noisy observations. In some way, this complements the works of Pananjady et al.~\cite{pananjady2020worst} on the effect of the design in the specific case where the sampled entries are sampled uniformly.

\subsection{Proof techniques and further comparison with the literature}\label{ss:litmoi}

In order to build a polynomial-time procedure nearly achieving the minimax permutation risk in the partial observation setting~\eqref{eq:model_partial}, we first consider the so-called full observation setting where one has access to poly-logarithmic number $\Upsilon$ of samples $Y^{(0)},Y^{(1)},\ldots, Y^{(\Upsilon)}$ of the complete $n\times d$ matrix. This setting is akin to that of Liu and Moitra~\cite{liu2020better} when they handled the specific square case $n=d$ with noise level $\zeta=1$.

\begin{figure}
	\includegraphics[width=0.75\linewidth]{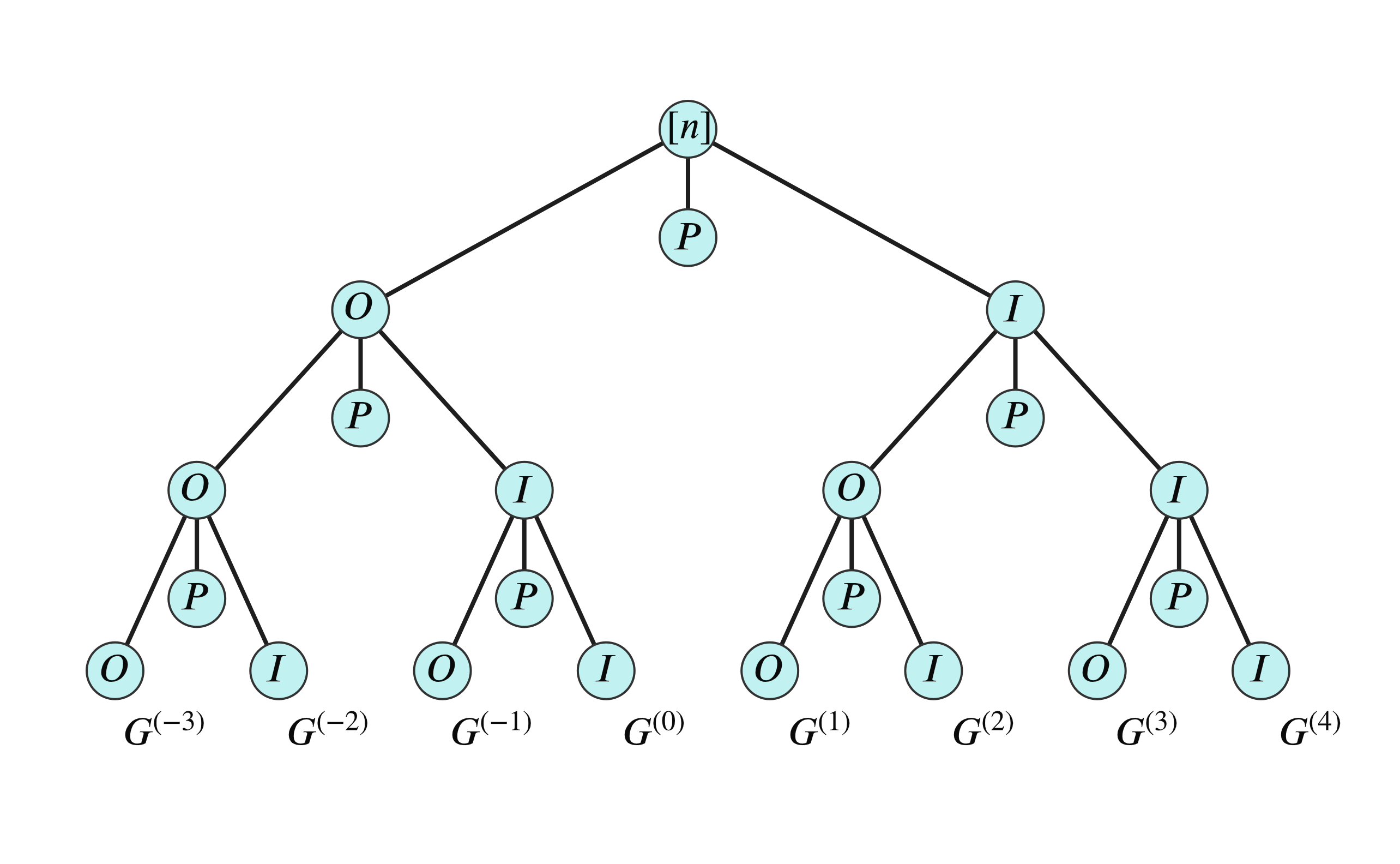}
	\caption{Example of a hierarchical sorting tree.}
	\label{fig:tree}
\end{figure}

For this reason, our estimators $\hat{\pi}_{HT}$ and $\hat{\pi}_{WM}$ introduced in Section~\ref{sec:algo_sketch} share some features with the procedure of~\cite{liu2020better}. From a broad perspective, our procedure and theirs build a hierarchical sorting tree using a top-down approach as depicted in~\Cref{fig:tree}. We start from the complete set $[n]$ of all experts and build a trisection $(O,P,I)$ of $[n]$, where $O$ (resp. $I$) contains experts that provably are below (resp. above) the median expert, whereas $P$ contains all the experts for which we cannot certify with high confidence whether they are above or below the median. Then, we recursively trisect the sets $O$ and $I$ as depicted in~\Cref{fig:tree}. At the end of the algorithm, we obtain a partial ordering on all the experts which can be used to estimate the oracle permutation $\pi^*$.

Then, the problem of building a suitable estimator boils down to introducing a suitable trisection procedure. We could naively do this by comparing the row-sums of the observed matrix $Y$ which amounts to comparing the mean ability of each expert, but this is well known to lead to suboptimal performances by a factor $\sqrt{d}$ --see e.g.~\cite{mao2020towards}. To improve over this rate, we need to compare the experts according to convex combinations of suitable questions. As in~\cite{liu2020better}, we start by selecting suitable blocks of questions by detecting the high-variation regions of the mean empirical expert and combine them with spectral algorithms to select suitable convex combinations of questions. Still, we have to refine 
significantly their spectral procedure to handle the rectangular case $n\ll d$. Equipped with these refinements, which are  involved technically, but are built on the ideas developed in~\cite{liu2020better}, we arrive at the estimator $\hat{\pi}_{HT}$ (see Section~\ref{sec:algo_sketch}) that turns out to be minimax optimal in some regimes of $(n,d,\zeta)$.

Unfortunately, this method turns out to be sub-optimal in many regimes, for instance for mild values of $n\in [d^{1/3}, d]$. Informally, this is due to the fact that our first estimator $\hat{\pi}_{HT}$ as well as that of Liu and Moitra~\cite{liu2020better} build an oblivious hierarchical sorting tree. This means that the trisection method decomposes a group $G^{(0)}$ of experts in the hierarchical sorting tree in $(O,P,I)$ only using the experts in $G^{(0)}$ of the matrix $Y$. In the related problem of hierarchical clustering, most top-down procedures also share this feature. It turns out that the observations of other experts can help improving the trisection of $G^{(0)}$. In particular, sets of experts that are close in the ordering --such as $G^{(1)}$ and $G^{(-1)}$ in Figure~\ref{fig:tree}-- are sometimes valuable to improve the selection of a suitable convex combination of questions. We emphasize this phenomenon and provide more intuition on it in Section~\ref{sec:algo_sketch}, when we introduce a new estimator $\hat{\pi}_{WM}$ that builds upon the memory of the sorting tree. This new procedure  $\hat{\pi}_{WM}$ turns out to be near minimax optimal
for all values of $(n,d,\zeta)$.

\medskip

Coming back to the partial observation setting~\eqref{eq:model_partial}, we introduce in Section~\ref{sec:partial} a reduction scheme which boils down to reducing the number of questions in order to come back to  a full observation model  for a sub-matrix of size $n\times d_-$ where $d_-$ is possibly much smaller than $d$. Then, relying on the full observation setting described above, we estimate the permutation $\pi^*$ based on the corresponding reduced matrix. In comparison to the full observation model, we can suffer from an additional bias 
terms which arises in the reduction process. To handle this, we develop a slight variant $\hat{\pi}_{WM-SR}$ of $\hat{\pi}_{WM}$ --see Appendix~\ref{sec:semi_random} for details. The resulting procedure turns out to nearly achieve minimax permutation recovery risk for all values $(n,d,\zeta)$ and all values of $\lambda$. Plugging our procedure to estimating the matrix $M$, we close all the computational gaps pointed out in Mao et al.~\cite{mao2020towards} for the problem of matrix estimation with a single unknown permutation - see Subsection~\ref{ss:recM}.

\subsection{Notation  and organization of the manuscript}

In the following, $c$, $c_1$, $\ldots$ stand for numerical positive constants that may change from line to line. Given a vector $u$ and $p\in [1,\infty]$, we write $\|u\|_{p}$ for its $l_p$ norm. For a matrix $A$, $\|A\|_F$ and $\|A\|_{\mathrm{op}}$ stand for its Frobenius and its operator norm. We write $\lfloor x\rfloor$ (resp. $\lceil x\rceil$) for the largest (resp. smallest) integer smaller than (resp. larger than) or equal to $x$.

Although $M$ stands for an $n\times d$ matrix, we extend it sometimes in an infinite matrix by setting $M_{i,k} = 0$ when either $i \leq 0$ or $k\leq 0$ and $M_{i,k} = 1$ when either $i \geq n+1$ and $k>0$ or $k \geq d+1$ and $i>0$. The corresponding infinite matrix $M_{\pi^{*-1}}$ which is obtained by permuting the $n$ original rows is still bi-isotonic and takes values in $[0,1]$. We shall often work with sub-matrices of $M$ that are restricted to a subset $P\subset [n]$ and $Q\subset [d]$ of rows and columns, in which case we write that the corresponding matrix $M'$ belongs to $\mathbb{R}^{P\times Q}$. More precisely, $M'$ is such that, $M'_{i,j}=M_{i,j}$ for any $i\in P$ and any $j\in Q$.

In the following, we write that two sequences or functions $u$ and $v$ satisfy $u\lesssim v$, if there exists a universal constant such that $u\leq c v$.

\medskip

In Section~\ref{sec:full}, we first consider the complete observation problem, where one has access to a poly-logarithmic number of independent samples of the complete noisy matrix $Y$. We characterize the minimax risk for permutation recovery and prove that it is achieved by a polynomial-time procedure. In section~\ref{sec:algo_sketch}, we describe the corresponding polynomial-time procedure. In Section~\ref{sec:partial}, we deal with the problem of partially observed matrix in the model~\eqref{eq:model_partial}. All the proofs are deferred to the appendix. 

\section{Analysis of the full observation problem}\label{sec:full}

As explained in the introduction, and following~\cite{liu2020better}, we first consider a slightly different problem where we fully observe a  $\Upsilon$-sample  $\cY = (Y^{(0)}, \dots, Y^{(\Upsilon-1)})$ of the noisy matrix according to the model $Y=M+E$ in~\eqref{eq:model_0}.
Here, $\Upsilon$ should be considered as a polylogarithms in $n$ and $d$. 
This is of course not very realistic in applications, but it is simpler to first present our algorithmic procedure in this setting, and it also enables more direct comparison to~\cite{liu2020better}. We will explain later in Section~\ref{sec:partial}, how one can transform data in the more realistic partial observation scheme from~\eqref{eq:model_partial} to this full observation scheme. We will then prove that the algorithm applied to the transformed data is near minimax optimal.

We recall that $M$ is a bi-isotonic matrix, up to an unknown permutation $\pi^*$ of its rows. Besides, the noise matrix $E$ is made of independent mean zero subGaussian entries, with Orlicz norm less than or equal to $\zeta$.

\subsection{Minimax lower bounds}\label{subsec:lower_bound}

Before considering ranking procedures, we characterize the minimax risk for the problem of ranking with full observations. For the purpose of the minimax lower bound, we assume that the noise matrix $E$ in~\eqref{eq:model_0} is made of independent  normal random variables with variance $\zeta^2$. For a permutation $\pi^*$ and a matrix $M$ such that $M_{\pi^*}\in 
\mathbb{C}_{\text{BISO}}$, we respectively denote $\P_{(\pi^*, M)}$ and $\E_{(\pi^*, M)}$ the corresponding probability and expectations with respect to the $\Upsilon$ independent observations of $Y$. 
Define 
\beq\label{eq:definition_RF}
	\cR_{F}(n,d,\zeta)= \zeta^2 \left[ \frac{nd^{1/6}}{\zeta^{1/3}} \wedge \frac{n^{3/4}d^{1/4}}{\zeta^{1/2}}\wedge n\sqrt{d}\wedge \frac{n^{2/3}\sqrt{d}}{\zeta^{1/3}}+  n \right]	\enspace . 
\eeq

The following minimax lower bound is stated in a setting where one has access to a polylogarithmic number $\Upsilon$  of full samples to be in line with the analysis of the next subsection. Still, we can forget about the dependency in $\Upsilon$  at first reading. 
\begin{theorem}\label{thm:lower_bound_full_observation}
There exists a universal constant $c$ such that the following holds for any $n\geq 2$, $d\geq 1$, $\zeta>0$, and $\kappa>2$. Provided that the sample size $\Upsilon$ is less than or equal to $\log^{\kappa}(2nd/\zeta)$, we have    
	 \begin{align}\label{eq:lower_bound_full_observation}
	\inf_{\hat \pi}  \sup_{\pi^*\in \Pi_n}\sup_{M:\,  M_{\pi^{*-1}}\in \mathbb{C}_{\text{BISO}}} \mathbb E_{(\pi^*, M)} \|M_{\hat \pi^{-1}} - M_{\pi^{*-1}}\|_F^2
	&\geq c\left[\log^{-\kappa}(nd/\zeta) \cR_F[n,d,\zeta]\bigwedge nd\right]\ . 
	\end{align}
	\end{theorem}
In fact, this theorem turns out to be a consequence of the minimax lower bound in the partial observation scheme --see Section~\ref{sec:partial}. Together with the risk upper bounds of the next section,~\eqref{eq:lower_bound_full_observation} characterizes, up to polylogarithmic terms, the minimax risk for estimating $\pi^*$.  The term $nd$ in~\eqref{eq:lower_bound_full_observation} is related to  the fact that the loss $\|M_{\hat \pi^{-1}} - M_{\pi^{*-1}}\|_F^2$ cannot be larger than $nd$ because the entries of $M$ are in $[0,1]$. 

\medskip 

\noindent 
{\bf Mild noise level}. 
The risk bound $\cR_{F}(n,d,\zeta)$ involves five different terms, some of them being significant only when $\zeta$ is small in comparison to $n$ and $d$. As these regimes with very small $\zeta$ are arguably quite specific, and to simplify the discussion, we will now detail the minimax lower bound in the specific case when $\zeta=1$. \beq\label{eq:risk_norm_1}
	\overline{\cR}_{F}(n,d)	= 	\cR_{F}(n,d,1)= (nd^{1/6})\bigwedge  (n^{3/4}d^{1/4}) + n\ . 
\eeq
In particular, we recognize three main regimes in~\eqref{eq:risk_norm_1} that depend on $n$ and $d$. 
When the number of experts is relatively small ($n\leq d^{1/3}$), the risk is proportional to $nd^{1/6}$. 
Specifying the result to $n=2$, one  checks that a square distance $d^{1/6}$ is necessary to distinguish two experts. As a consequence, a suitable estimator $\hat{\pi}$ should be able to coherently rank experts that are distant by more than $d^{1/6}$ in squared Frobenius norm, and then to achieve a risk smaller than $nd^{1/6}$. For larger $n\geq d^{1/3}$, it is in fact possible to build upon the large number of experts to improve the comparisons between experts using in particular spectral methods. For this reason, the optimal risk is proportional to $n$ for $n\geq d$. For an intermediary number of experts $n\in [d^{1/3},d]$, the risk is of the order of $n^{3/4}d^{1/4}$. Our main contribution is the construction of a polynomial-time procedure that achieves these risk bounds, see below.

\medskip 

\noindent 
{\bf Low noise level}. For mild values of $\zeta$, the minimax risk $\cR_{F}(n,d,\zeta)$ has the same form as $\overline{\cR}_{F}(n,d)$, up to some factors that depend on $\zeta$. However, for very small $\zeta$, the risk becomes qualitatively different. For example, we have $\cR_{F}(n,d,\zeta)\asymp \zeta^2 n\sqrt{d}$ when $\zeta\in (0, \tfrac{1}{n\vee d}]$. In fact, this rate is quite easy to achieve by a polynomial time algorithm in this extreme case. It is proven in various works -- see e.g.~\cite{shah2016stochastically} that ranking the experts according to the row sum of the matrix correctly compares two experts as long as their square distance is at least $\zeta^2\sqrt{d}$ (up to logarithmic  terms). As a consequence, this simple procedure leads to an error $\zeta^2 n\sqrt{d}$. While $\zeta^2n\sqrt{d}$ is highly suboptimal in most realistic regimes, it turns out to be tight for extremely low level of noise. 
Finally, the intermediary rate $\cR_{F}(n,d,\zeta)\asymp \zeta^{5/3}n^{2/3}\sqrt{d}$ is achieved for slightly larger values of $\zeta$, but it is less clear how to interpret it.

\subsection{Minimax upper bounds}\label{sec:minimax_upper}

In the following, we fix a parameter  $\delta\in (0,1)$ that will correspond to a small probability. We write $\zeta_-=\zeta\wedge 1$, where $\zeta$ is the noise level. In this section, we analyze two estimators $\hat{\pi}_{HT}$ and $\hat{\pi}_{WM}$ of $\pi^*$ that are described in Section~\ref{sec:algo_sketch} and more formally defined in Appendix~\ref{sec:algorithms}. The first estimator $\hat{\pi}_{HT}$ is based on the construction of an \emph{oblivious} hierarchical sorting tree. We will later explain all the ingredients of this procedure. In contrast, the second estimator $\hat{\pi}_{WM}$ relies on the construction of a hierarchical sorting tree \emph{with memory}. Both procedures have a computational complexity of the order of $\log^{c}(\frac{nd}{\zeta_-\delta})(n^3 + nd^2)$, for some $c>0$, which makes them polynomial time - unlike the least square procedure e.g.~from~\cite{mao2020towards}.

\begin{theorem}\label{th:first_estimator}
	There exist three numerical constants $c$, $c'$, and $c_0$ such that the following holds. Fix $\delta>0$ and assume that $\Upsilon \geq c_0 \log^{8}\left(nd/\delta\right)$. For any permutation $\pi^*\in \Pi_n$ and any matrix $M$ such that $M_{\pi^{*-1}}\in \mathbb{C}_{\text{BISO}}$, 
	 the oblivious hierarchical sorting tree estimator $\hat{\pi}_{HT}$ defined in the next section satisfies
	$$ \|M_{\hat \pi^{-1}_{HT}} - M_{\pi^{*-1}} \|_F^2 \leq c  \zeta^2 \log^{10.5}\left(\frac{2nd}{\delta\zeta_-}\right) \left[\frac{n^{2/3}d^{1/3}}{\zeta^{2/3}}\land \frac{nd^{1/6}}{\zeta^{1/3}}\wedge n\sqrt{d} + n \right] \enspace ,$$
	with probability at least $1- c'n\log^{9}(\tfrac{nd}{\delta\zeta_-}) \delta$.
\end{theorem}
If we take $\delta=\zeta^2(nd)^{-1}$ in the above expression, we easily deduce - reminding that the entries of $M$ are in $[0,1]$ - the following risk bound 
\[
\E\left[\|M_{\hat{\pi}_{HT}^{-1}} - M_{\pi^{*-1}} \|_F^2\right] \leq c   \zeta^{2}\log^{10.5}\left(\frac{2nd}{\zeta_-} \right) \left[\frac{n^{2/3}d^{1/3}}{\zeta^{2/3}}\land \frac{nd^{1/6}}{\zeta^{1/3}}\wedge n\sqrt{d} + n\right] \enspace .	
\]
Comparing this bound with~\eqref{eq:risk_norm_1} in the specific case where $\zeta=1$, we observe that $\hat{\pi}_{HT}$ achieves the optimal  risk $nd^{1/6}$ for small $n\leq d^{1/3}$ and the optimal risk $n$ for large $n\geq d$. Unfortunately,  for mild $n\in [d^{1/3}, d]$, the risk bound is of the order of $n^{2/3}d^{1/3}$, which is significantly higher than the minimax lower bound $n^{3/4}d^{1/4}$. To close this gap, we turn to the more refined estimator $\hat{\pi}_{WM}$.

\begin{theorem}\label{th:second_estimator_WM}
	There exist three numerical constants $c$, $c'$, and $c_0$ such that the following holds. Fix $\delta>0$ and assume that $\Upsilon \geq c_0 \log^{8}\left(nd/(\delta\zeta_-)\right)$.  For any permutation $\pi^*\in \Pi_n$ and any matrix $M$ such that $M_{\pi^{*-1}}\in \mathbb{C}_{\text{BISO}}$, the hierarchical sorting tree estimator with memory $\hat{\pi}_{WM}$ satisfies 
	\beq\label{eq:risk_upper_bound_WM}
	\|M_{\hat \pi_{\WM}^{-1}} - M_{\pi^{*-1}} \|_F^2 \leq \left[c
	\log^{11}\left(\frac{2nd}{\delta \zeta_-}\right)\cR_F[n,d,\zeta]\right] \bigwedge nd\enspace ,
	\eeq
    with probability at least $1- c'n\log^{9}(\tfrac{nd}{\delta\zeta_-}) \delta$. 
\end{theorem}
As for the previous theorem, this high probability result can be turned into a risk bound by taking $\delta=\zeta^2/(nd)$. In particular, this risk bound matches, up to polylogarithmic terms, the minimax lower bound~\eqref{eq:lower_bound_full_observation} for all possible values of $n$, $d$, and $\zeta$. As a consequence, the estimator $\hat{\pi}_{WM}$ is nearly minimax and this ranking problem does not exhibit any computational gap.

 In~\cite{liu2020better}, the polynomial-time estimator $\hat{\pi}_{LM}$ of Liu and Moitra achieves the minimax risk in the specific square where $n=d$ and $\zeta=1$.  In all the other regimes, no polynomial-time procedure was previously proved to achieve the minimax risk.  In fact, even if we do not restrict our attention to polynomial-time procedures, least-square type procedures studied e.g. in~\cite{mao2020towards} provably achieve the minimax risk only in the regime when $n\geq d$. As alluded in the introduction - see Equations~\eqref{eq:minimax_intro} and ~\eqref{eq:minimax_intro_reconstruction}, the minimax risks for estimating $\pi^*$ and $M$ differ when $n\leq d$, so that achieving the optimal risk for $\pi^*$ is not possible using the classical entropy arguments as in~\cite{shah2016stochastically,mao2020towards}. This highlights the fact that  estimating the permutation $\pi^*$ is significantly more challenging  in the regime  $n\leq d$ - both from a statistical and computational perspective - than in the regime $n \geq d$ handled in~\cite{liu2020better,mao2020towards}.

\medskip

\noindent
{\bf Consequences for the estimation of the matrix $M$.} Provided that we have estimated $\pi^*$ with $\Upsilon-1$ independent samples, we could use the last sample $Y^{(\Upsilon)}$ to estimate the matrix $M$ by minimizing the least-square criterion  $\widehat{B} = \argmin_{B\in \mathbb{C}_{\text{BISO}}}\|Y_{\hat{\pi}_{WM}}^{(\Upsilon)} - B\|_F^2 $ and setting $\widehat{M}= \widehat{B}_{\hat{\pi}_{WM}^{-1}}$.  Since the set of bi-isotonic matrices is convex, this estimator is computable efficiently~\cite{kyng2015fast}.
As argued in Proposition 3.3 of~\cite{mao2020towards} and often used in the ranking literature~\cite{shah2019feeling,chatterjee2019estimation,pananjady2022isotonic}, it turns out that, with high probability, the reconstruction error $\|M-\widehat{M}\|_F^2$ is (up to polylogarithmic terms) the sum of the expected permutation loss $\E\left[\|M_{\hat{\pi}_{WM}}-M_{\pi^*}\|_F^2\right]$ and the minimax reconstruction risk of a bi-isotonic matrix $\inf_{\widehat{B}}\sup_{B\in \mathbb{C}_{\text{BISO}}}\E[\|\widehat B-B\|_F^2]$ where $Y=B+E'$ and $E'$ is made of independent subGaussian random variables. Hence, based on $\hat{\pi}_{WM}$ and Theorem~\ref{th:second_estimator_WM}, it is easy to construct a polynomial-time estimator of $M$ that is also near minimax-optimal in the sense of Equation~\eqref{eq:minimax_reconstruction}. We will further build upon this remark in Section~\ref{sec:partial} when we come back to the problem of partial observations of the matrix.

\section{Description of the hierarchical sorting estimators}\label{sec:algo_sketch}

Let us now describe the construction of the estimators $\hat{\pi}_{HT}$ and $\hat{\pi}_{WM}$ of $\pi^*$. The construction is quite long and involves several subroutines. For this reason and to ease the understanding of proof details, we also provide a more formal and longer definition in~\Cref{sec:algorithms}. Afterwards, we comment on the different steps of the procedure and on their connection to the literature in Subsection~\ref{subsec:comments}.

\medskip 

Define $\tau_{\infty}=\lceil 4\cdot 10^7\log^7(\tfrac{nd}{\delta(\zeta_-)^2})\rceil$ and $t_{\infty}= \lceil \log(n)/\log(2)\rceil$. We define $\Upsilon^*= 6\tau_{\infty}t_{\infty}$ for the total number of independent samples required for the computation of these two estimators. 

Hence, we are given independent samples $\cY= (Y^{(0)}, \dots, Y^{(\Upsilon^*-1)})$. From a broad perspective, both procedures are based on the construction of the recursive sorting tree as illustrated in~\Cref{fig:tree}. Starting from the root of the tree which corresponds to the set $[n]$ of all experts, we build a partition $O$, $P$, $I$, of $[n]$ in such a way that, with high probability, all the experts in $O$ are below the median expert of $[n]$, all the experts in $I$ are above the median expert of $[n]$, while the remaining experts in $P$ are those for which we are not able to decipher whether they are below or above the median expert of $[n]$.

Having trisected $[n]$, we recursively trisect the subsets $O$ and $I$- see \Cref{fig:tree}. Each time, the size of the groups $O$ and $I$ is divided by at least $2$. Hence, at depth $t_{\infty}$, all the groups of $O$ and $I$ have size at most $1$. For each depth $t=0,\ldots, t_{\infty}-1$, we use $6 \tau_{\infty}$ new samples. The construction of the tree is described in $\algoTreeSorting$ --see Algorithm~\ref{alg:procedure_sketch} and is based on the routine 
$\algoBlockSorting$ which performs the trisection of a group into $(O,P,I)$.

\medskip 

Let us now explain how to deduce an estimator $\hat{\pi}$ from the final hierarchical sorting tree $\cT$. Indeed, the hierarchical sorting tree $\cT$ induces an order on its leaves as follows. For any groups ($O,P,I$) sharing the same parent, we say that any descendent of $O$ in the tree $\cT$ is below $P$, which, in turn, is below any descendent of $I$ in $\cT$. This endows a complete ordering on the leaves of the tree $\cT$. Denote $\cG = (G_1, \dots, G_{\alpha})$ the sequence of leaves of the final tree ranked according to this complete order. For any $a\in [\alpha]$, we define the lower bound $\pi^-_{\cG}(G_a)$ and the upper bound $\pi^+_{\cG}(G_a)$ of the ranks of experts in $G_a$ by $\pi^-_{\cG}(G_a) := \sum_{a' < a} |G_{a'}|$ and $\pi^+_{\cG}(G_a) := \sum_{a' \leq  a} |G_{a'}|$. Finally, we sample $\hat{\pi}$ arbitrarily in such a way that 
\begin{equation}\label{eq:pi_plus_sketch}
\hat{\pi}(G_a) = [\pi^-_{\cG}(a) + 1, \pi^+_{\cG}(a)]\enspace . 
\end{equation}
In other words, the estimator $\hat{\pi}$ ranks the groups $G_a$ according to the ordering of the groups endowed by $\cT$ and, given that, ranks the experts $G_a$ uniformly at random. 
See \Cref{sec:algorithms} for a more formal definition of the ordering.

\subsection{Description of the trisection of a leaf $G$ into $(O, P, I)$ with $\algoBlockSorting$}

The purpose of $\algoBlockSorting$ is to build a trisection of a group $G$ of experts into $(O,P,I)$ where $O$ is made of experts that are, with high probability, below the median expert in $G$ and $I$ is made of experts which are, with high probability, above this median expert. It turns out that this construction is based on $\tau_{\infty}$ iterations of a procedure called $\algoDoubleTrisection$ which is the backbone of our procedure. Intuitively, we shall iteratively detect subgroups of experts that are below (resp. above) the median expert of $G$ which, after $\tau_{\infty}$ iterations, will allow us to obtain $O$ and $I$.

For technical reasons, our definition is slightly more intricate. 
We shall simultaneously build two collections  $(O_{\tau},I_{\tau})$ and $(\overline{O}_{\tau},\overline{I}_{\tau})$ of groups, the second one being more conservative. We start with empty sets for $(O_0, I_0, \overline O_0, \overline I_0) = \emptyset$. 
Then, at each step $\tau$, we will consider the remaining set of experts $G\setminus(\overline O_{\tau} \cup \overline I_{\tau})$. Define $\gamma= \floor{|G|/2} - |\overline O_{\tau}|$ for the presumed rank of the median expert of $G$ inside $G\setminus(\overline O_{\tau} \cup \overline I_{\tau})$. Then, using $6$ independent samples, we apply $\algoDoubleTrisection(\cY,\cT, G\setminus(\overline O_{\tau} \cup \overline I_{\tau}), \gamma)$ to compute four subsets $(L_{\tau}, U_{\tau})$ and  $(\overline L_{\tau},\overline U_{\tau})$. With high probability, it turns out that $\overline{L}_{\tau}\subset L_{\tau}$ is made of experts below the median expert of $G$ and $\overline{U}_{\tau}\subset U_{\tau}$ is made of experts above the median expert of $G$. This allows us to update as follows
\begin{equation}
	O_{\tau + 1} = O_{\tau}\cup L_{\tau}, ~~I_{\tau + 1} = I_{\tau}\cup U_{\tau}, ~~\overline O_{\tau + 1} = \overline O_{\tau}\cup \overline L_{\tau}, ~~\overline I_{\tau + 1} = \overline I_{\tau}\cup \overline U_{\tau} \enspace .
\end{equation}
The procedure is summarized in Algorithm~\ref{alg:algoBlockSorting_sketch} below.

\begin{minipage}[c]{0.53\textwidth}
\begin{algorithm}[H]
	\caption{\label{alg:procedure_sketch} $\algoTreeSorting(\cY)$}
	\begin{algorithmic}[1]
		\Require $6\tau_{\infty}t_{\infty}$ samples $\cY = (Y^{(0)}, \dots, Y^{(6\tau_{\infty}t_{\infty}-1)})$ 
		\Ensure A tree $\cT$ and an estimator $\hat{\pi}$
		\State Initialize $\cT$ as the tree with only the root $[n]$
		\For{$t = 0, \dots, t_{\infty}-1$}
		\State Take $6\tau_{\infty}$ samples $\cY_t = (Y^{(6t\tau_{\infty})}, \dots, Y^{(6(t+1)\tau_{\infty}-1)})$ 
		\State Initialize $\cT'= \cT$
		\For{All the leaves $G$ at depth $t$ corresponding to $O$ or $I$ as in \Cref{fig:tree}}
		\State Set $(O_G,P_G,I_G) = \algoBlockSorting(\cY_t, \cT,G)$ 
		\State Add  $(O_G,P_G,I_G)$  to the tree $\cT'$
		\EndFor
		\State Update $\cT= \cT'$
		\EndFor
		\State Set $\hat{\pi}:=\hat{\pi}(\cT)$  as in \Cref{eq:pi_plus_sketch} 
		\State \Return{$\cT$ and $\hat{\pi}$}
	\end{algorithmic}
\end{algorithm}
\end{minipage}
\ ~ \
\begin{minipage}[c]{0.45\textwidth}
\begin{algorithm}[H]
	\caption{\label{alg:algoBlockSorting_sketch} $\algoBlockSorting(\cY, \cT, G)$ }
	\begin{algorithmic}[1]
		\Require $6\tau_{\infty}$ samples $\cY = (Y^{(0)}, \dots, Y^{(6\tau_{\infty}-1)})$, the tree $\cT$, a leaf $G$ in $\cT$
		\Ensure A partition of $G$ into  $(O,P, I)$
		\Statex
		\State Set $\gamma = \floor{|G|/2}$ and $O_0$, $I_0$, $\overline O_0$,  $\overline I_0$ $= \emptyset$
		\For{$\tau = 0, \dots, \tau_{\infty}-1$}
		\State Take $6$ samples $ \cY_\tau = (Y^{(6\tau)}, \dots, Y^{(6\tau + 5)})$ 
		\State set $\gamma = \floor{|G|/2} - |\overline O_{\tau}|$
		\State $(L_{\tau}, U_{\tau}), ~ (\overline L_{\tau},\overline U_{\tau}) = \algoDoubleTrisection(\cY_\tau,\cT, G\setminus(\overline O_{\tau} \cup \overline I_{\tau}), \gamma)$ as in Algorithm~\ref{alg:double_trisection_sketch}
		\State Update
		$O_{\tau+1} = O_{\tau}\cup L_{\tau}, ~~I_{\tau+1} = I_{\tau}\cup U_{\tau}, ~~\overline O_{\tau+1} = \overline O_{\tau}\cup \overline L_{\tau}, ~~\overline I_{\tau+1} = \overline I_{\tau}\cup \overline U_{\tau}$
		\EndFor
		\State \Return{$(O_{\tau_{\infty}}, G\setminus(O_{\tau_{\infty}} \cup I_{\tau_{\infty}}), I_{\tau_{\infty}})$}
	\end{algorithmic}
\end{algorithm}
\end{minipage}

\subsection{Description of the double trisection procedure}\label{subsec:description_procedure_trisection}

We now describe the trisection procedure $\algoDoubleTrisection$. For this purpose, we first provide a few definitions.

\subsubsection{Definitions}
In this subsection, we write $Y$ for one data set sampled according to \Cref{eq:model_0}. For the sake of simplicity, we often omit the dependence of $Y$ in the definitions. We write $\cD$ for the set of all dyadic numbers: $\cD = \{ 2^k ~:~ k \in \bbZ \}$ and we define the sets $\cR = \cD \cap [1,d]$ and $\cH = \cD \cap \left[\frac{\zeta^2}{nd},1\right]$. The collection $\cR$  corresponds to the possible scales, that is the number of questions under consideration, whereas the collection $\cH$ corresponds to the possible heights of variations. 

For all $r \in \cR$, we write $\cQ_r = \{1, r+1, 2r + 1,\dots ,\floor{\frac{d}{r}}r+1\}$ for the regular grid of questions with spacing $r$.
If $\overline P \subset [n]$ is a set of experts, we denote $\overline y(\overline P)$ as the mean of the vectors $Y_{i, \cdot}$ for $i \in \overline P$, that is, for all $k \in [d]$, we have $\overline y_k(\overline P) = \tfrac{1}{|\overline P|}\sum_{i \in \overline P} Y_{i, k}$. For any $r \in \cR$, we define $Z(Y,\overline P, r)$ as the aggregation of the matrix $Y$ on blocks of questions of size $r$ and with lines restricted to $\overline P$. More formally, for any $i \in \overline P$ and $l \in \cQ_r$, we have 
\begin{equation}\label{eq:aggregation}
	Z_{i, l}(Y,\overline P, r) = \frac{1}{\sqrt{r}} \sum_{k = l}^{l+r-1} Y_{i,k} \spaceAnd \overline Z_{i, l}(Y,\overline P, r) = \frac{1}{\sqrt{r}} \sum_{k = l}^{l+r-1} \overline y_{k}(\overline P)\enspace . 
\end{equation}
Both matrices are of size $|\overline P| \times |\cQ_r|$. Note that, in the above definition, $Z_{i, l}(Y,\overline P, r)$ and $\overline{Z}_{i, l}(Y,\overline P, r)$ are rescaled by $\sqrt{r}$ so that the subGaussian norm remains at most $\zeta$.  For any subset $Q \subset \cQ_r$, we also write $Z(Y,\overline P, Q, r)$ for the sub-matrix of $Z(Y,\overline P, r)$ restricted to columns in $Q$.

Given a matrix $Z \in \bbR^{\overline P \times Q}$, a vector $w \in \bbR_+^{Q}$ with non-negative components and $i$, $j$ in $\overline P$, we say $i$ is $(Z, w)$-above $j$ (or equivalently that $j$ is $(Z, w)$-below $i$) if the projection of $Z_{i, \cdot}$ on the direction $w$ is larger than the projection of $Z_{j, \cdot}$ on $w$, that is $\proscal<Z_{i, \cdot} - Z_{j, \cdot}, w> > 0$, where $\proscal<.,.>$ stands for the standard inner product between vectors.
Now, for $\gamma \in \{1, \dots, |\overline P| \}$, we can consider the $\gamma$-th expert $i_{\gamma} \in \overline P$ such that there are exactly $\gamma - 1$ experts which are $(Z,w)$-below $i_{\gamma}$. Given a tuning parameter $\beta> 0$ to be fixed below, we then define the $(Z, w, \gamma, \beta)$-trisection of $\overline P$ on direction $w$ with respect to pivot index $\gamma$ and matrix $Z$ as the sets:
\begin{equation}\label{eq:pivot_trisection}
	\begin{cases}
		U_{w} := U(Z, w, \gamma, \beta) & = \left\{ i \in \overline P ~:~ \proscal<Z_{i, \cdot} - Z_{i_{\gamma}, \cdot}, \frac{w}{\|w\|_2}> \geq \beta\sqrt{ \Log{2\frac{|\overline P|}{\delta}}}\right\}             \\
		L_{w} := L(Z, w, \gamma, \beta) & = \left\{ i \in \overline P ~:~ \proscal<Z_{i, \cdot} - Z_{i_{\gamma}, \cdot}, \frac{w}{\|w\|_2}> \leq -\beta\sqrt{ \Log{2\frac{|\overline P|}{\delta}}}\right\} \enspace .
	\end{cases}
\end{equation}
Hence a $(Z, w, \gamma, \beta)$-trisection on direction $w$ and pivot $\gamma$ consists of two possibly empty disjoint subsets $U$ and $L$ which are respectively taken among the $\gamma - 1$ experts (resp. the $|\overline P| - \gamma$) which are $(Z, w)$-above (resp. $(Z, w)$-below) the expert $i_{\gamma}$, with a margin of the order of $\sqrt{\log(|\overline P|/\delta)}$. Remark that if $\beta < \overline \beta$ then $U(Z, w, \gamma, \overline \beta) \subset U(Z, w, \gamma, \beta)$, which means that the trisection of $\overline P$ on direction $w$ becomes more conservative as $\beta$ increases.

In fact, \eqref{eq:pivot_trisection} turns out to be the cornerstone or our procedure. 
Since the coordinates of $w$ are non-negative, the corresponding row-wise weighted sums of the aggregation $\E[Z(Y,\overline P, r)]$ of the signal matrix $M$ are also ordered according to the oracle permutation. In other words, the $k$-th expert in $\overline P$ has the $k$-th highest value of the expectation of this  weighted sum.

For $r \in \cR$ and $Q \subset \cQ_r$, choosing $w = \1_Q$ in \Cref{eq:pivot_trisection} amounts to trisecting $\overline P$ according to the average of the observations over all questions in $\bigcup_{l \in Q} [l, l+r)$. In that case, we write for simplicity $(L_{Q}, U_{Q}) = (L_{\1_Q}, U_{\1_Q})$. When $Q = \cQ_r$ and $w = \1_{Q}$, then~\eqref{eq:pivot_trisection} simply amounts to ranking experts according to their average over all the questions.
As explained in the introduction, the global average does not lead to optimal performances.  This is why 
most following steps in the algorithm amount to selecting suitable blocks $Q$ of questions and directions $w$.

In the following, the tuning parameters $\beta$ are set as follows. 
\beq\label{eq:definition_beta}
\beta_{\tris} = 4\sqrt{2}\zeta\, ,\quad \quad \overline{\beta}_{\tris} = 8\sqrt{2}\zeta \enspace . 
\eeq

\subsubsection{Description of the double trisection procedure}

Recall that the purpose of $\algoDoubleTrisection$ is to select subsets $(L, U)$ and $(\overline L, \overline U)$ of a group $\overline P$ of experts in such a way that $\overline{L}\subset L$, $\overline{U}\subset U$, and experts in $L$ (resp. in $U$) are with high probability below (resp. above) the $\gamma$-th  expert of $\overline P$. 

For that purpose, we have $6$ independent samples $(Y^{(s)})_{s =1, \dots, 6}$ sampled from \Cref{eq:model_0} at our disposal.  Fix any height $h \in \cH$ and any scale $r \in \cR$. $\algoDoubleTrisection$ relies on the following steps also described in Algorithm~\ref{alg:double_trisection_sketch}.
\begin{enumerate}
    \item {\bf Selection of a suitable subset of questions.} Using the first sample $Y^{(1)}$, we first select
    a subset $\widehat{Q}\subset \cQ_r$. We postpone the definition of the selection procedure to the next subsection. We will introduce  two approaches for this $\widehat Q := \widehat Q_{\cp}(\overline P, h, r)$ as in \Cref{eq:significant_blocks_estim} or $\widehat Q := \widehat Q_{WM}(\cT, \overline P, h, r)$ as in \cref{eq:significant_blocks_estim_wm}.  These two definitions respectively correspond to the \emph{oblivious} estimator $\hat{\pi}_{HT}$ and to the estimator \emph{with memory} $\hat{\pi}_{\WM}$.
    \item {\bf Average-based trisection.} Using the second sample $Y^{(2)}$, we consider the corresponding aggregated matrix  $Z^{(2)}:=Z(Y^{(2)},\overline P, \widehat Q, r)$ as defined in \Cref{eq:aggregation} which focuses on the selected blocks of questions $\widehat{Q}$ . Then, we consider experts whose corresponding row sums on $Z^{(2)}$ is unusually large or small. More formally, we 
	compute the $(Z^{(2)}, \1_{\widehat Q}, \gamma, \beta_{\tris})$-trisection and the $(Z^{(2)}, \1_{\widehat Q}, \gamma, \overline \beta_{\tris})$-trisection of $\overline P$ as defined in \Cref{eq:pivot_trisection} and where the tuning parameters $\beta_{\tris}$ and $\overline{\beta}_{\tris}$ are defined in~\eqref{eq:definition_beta}. This allows us to obtain $(L_{\widehat Q}, U_{\widehat Q})$ 
	and $(\overline L_{\widehat Q}, \overline U_{\widehat Q})$.
\item {\bf PCA-based trisection.} Then, we focus on the conservative subset of remaining experts  $\widetilde P= \overline P \setminus \overline L_{\widehat Q} \cup \overline U_{\widehat Q}$. Relying on the samples $Y^{(3)}$, $Y^{(4)}$, $Y^{(5)}$, $Y^{(6)}$, we build the corresponding aggregated matrices $Z^{(s)} := Z(Y^{(s)},\widetilde P, \widehat Q, r)$ restricted to the subset $\widetilde P$ for $s=3,4,5$. In principle, we would like to aim at the right singular value of $\mathbb{E}[Z^{(3)}-\overline{Z}^{(3)}]$ as this would give us a nice direction $w$ on which we could apply~\eqref{eq:pivot_trisection}. For technical reasons to be explained later, we take a roundabout way, by first computing  a vector $\hat v$ indexed by $\widetilde P$ which, in principle, is not too far from the left singular value of $\mathbb{E}[Z^{(3)}-\overline{Z}^{(3)}]$. More precisely, we compute $\hat{v}$ as follows 
\begin{equation}\label{eq:ACP}
	\hat v := \hat v(\widetilde P, \widehat Q, r) = \argmax_{\| v \|_2 \leq 1} \Big[ \|v^T(Z^{(3)} - \overline{Z}^{(3)})\|_2^2 - \frac{1}{2}\| v^T(Z^{(3)} - \overline{Z}^{(3)} - Z^{(4)} + \overline{Z}^{(4)})\|_2^2\Big] \enspace .
\end{equation}
The right-hand side term in~\eqref{eq:ACP} allows us to deal with the fact that the entries of the noise matrix $E$ in~\eqref{eq:model_0} are possibly heteroskedastic. Although there exist more elegant workarounds for heteroskedastic noise (e.g. PCA~\cite{zhang2022heteroskedastic}), the analysis in those works  does not apply in our non-parametric setting.
Moreover,  $\hat v$ in~\eqref{eq:ACP} corresponds to the leading eigenvector of a square symmetric matrix and can therefore be computed efficiently. Then, we consider the image  $\hat z = \hat v^T (Z^{(5)}-\overline Z^{(5)}) \in \bbR^{\widehat Q}$  of $\hat{v}$. After this, we threshold $\hat{z}$ and take the absolute values of the components. Thus, we get $\hat w^+\in \mathbb{R}^Q$ defined by $(\hat w^+)_l= |\hat{z}_l|\1_{|\hat z_l|\geq 2\zeta\sqrt{2\log(2|\widehat Q|/\delta)}}$ for any $l  \in \widehat Q$.
Finally, we consider the last aggregated sample $Z^{(6)} := Z(Y^{(6)},\overline P, \widehat Q, r)$ on the set $\overline{P}\supset \widetilde{P}$ of experts. We apply these weights $\hat w^+$ to compute the row-wise weighted sums of $Z^{(6)}$ and discard experts whose corresponding weighted sums is unusually small or large. More formally, we apply $(Z^{(6)}, \hat w^+, \gamma, \beta)$-trisection and $(Z^{(6)}, \hat w^+, \gamma, \overline \beta)$-trisection of $\overline P$ as defined in \Cref{eq:pivot_trisection}. Doing so we obtain $(L_{\hat w^+}, U_{\hat w^+})$ and $(\overline L_{\hat w^+}, \overline U_{\hat w^+})$ respectively.

In the definition of $Z^{(6)}$ we consider the whole set of experts $\overline P$ instead of the remaining of experts $\widetilde P$ that have not been discarded because otherwise we should have needed to update the value of $\gamma$ when applying~\Cref{eq:pivot_trisection}.
\end{enumerate}

Finally, we define the trisections $(L, U)$ (resp. $(\overline L, \overline U)$) as the union of the corresponding discarded subsets of experts based on $\1_{\widehat{Q}}$ and $\hat w^+$, this for all possible height $h\in \cH$ and scale $r\in \cR$. We recall that the definition of $\widehat{Q}$ was depending on $h$ and $r$. 
\begin{equation}
	\begin{cases}
		(L, U) &= \left(\bigcup_{(h,r) \in \cH \times \cR}L_{ \widehat Q}(h,r)\cup L_{\hat w^+}(h,r), \bigcup_{(h,r) \in \cH \times \cR}U_{\widehat Q}(h,r)\cup U_{\hat w^+}(h,r)\right)\\	
		(\overline L, \overline U) &= \left(\bigcup_{(h,r) \in \cH \times \cR}\overline L_{ \widehat Q}(h,r)\cup \overline L_{\hat w^+}(h,r), \bigcup_{(h,r) \in \cH \times \cR}\overline U_{\widehat Q}(h,r)\cup \overline U_{\hat w^+}(h,r)\right).
	\end{cases}	
\end{equation}
This whole routine for computing $(L, U)$, $(\overline L, \overline U)$ is referred to as $\algoDoubleTrisection$ and is summarized in \Cref{alg:double_trisection_sketch}. We underline that $\overline L \subset L \subset \overline P$ and $\overline U \subset U \subset \overline P$ as we took $\beta_{\tris} < \overline \beta_{\tris}$.

\begin{algorithm}[H]
	\caption{$\algoDoubleTrisection((Y^{(s)})_{s= 1, \dots, 6}, \cT, \overline P, \gamma)$ \label{alg:double_trisection_sketch}}
	\begin{algorithmic}[1]
		\Require $6$ samples $(Y^{(s)})_{s=1,\dots,6}$, a set $\overline{P}$, a tree $\cT$, 	a pivot index $\gamma  \in [1, \dots, |\overline P|]$
		\Ensure Two trisections $(L, U)$ and $(\overline L, \overline U)$ of $\overline P$
		\Statex
		\State Start from $L, U, \overline L, \overline U = \emptyset$
		\For{$h \in \cH$, $r \in \cR$}
		\State \label{line:choice_Q}Compute $\widehat Q := \widehat Q_{\cp}(\overline P, h, r)$ as in \Cref{eq:significant_blocks_estim} or $\widehat Q := \widehat Q_{WM}(\cT, \overline P, h, r)$ as in \cref{eq:significant_blocks_estim_wm} using sample $Y^{(1)}$
		\State Let $Z^{(s)}:=Z(Y^{(s)},\overline P, \widehat Q, r)$, for $s\in \{2,6\}$ be the aggregated matrices of samples defined
		as in \Cref{eq:aggregation} 
		\State \label{line:first_stat_L1}Let $(L_{\widehat Q}, U_{\widehat Q})$, $(\overline L_{\widehat Q}, \overline U_{\widehat Q})$ be resp.~the $(Z^{(2)}, \1_{\widehat Q}, \gamma, \beta)$ and the $(Z^{(2)}, \1_{\widehat Q}, \gamma, \overline \beta)$-trisections of $\overline P$ as in \Cref{eq:pivot_trisection}
		\State Define $\widetilde P = \overline P \setminus (\overline L_{\widehat Q} \cup \overline U_{\widehat Q})$ and the aggregated samples $Z^{(s)}:= Z^{(s)}(Y^{(s)},\widetilde P, \widehat Q, r)$ for $s \in \{3,4,5\}$ 
		\State \label{line:ACP}Compute the PCA-like direction $\hat v := \hat v(\widetilde P, \widehat Q, r)$ as in~\eqref{eq:ACP}
		\State \label{line:direction}Compute $\hat z = \hat v^T (Z^{(5)}-\overline Z^{(5)})$ and define $\hat w^+$ by $(\hat w^+)_l= |\hat{z}_l|\1_{|\hat z_l|\geq 2\zeta\sqrt{2\log(2|\widehat Q|/\delta)}}$ for any $l  \in \widehat Q$
		\State \label{line:second_stat_L1}  Let $(L_{\hat w^+}, U_{\hat w^+})$, $(\overline L_{\hat w^+}, \overline U_{\hat w^+})$ be resp.~the $(Z^{(6)}, \hat w^+, \gamma, \beta)$ and the $(Z^{(6)}, \hat w^+, \gamma, \overline \beta)$-trisections of $\overline P$ as in \Cref{eq:pivot_trisection}
		\State Update $L = L\cup L_{\hat w^+}\cup L_{\widehat Q}$,~~ $U = U\cup U_{\hat w^+} \cup U_{\widehat Q}$,~~ $\overline L = \overline L\cup \overline L_{\hat w^+} \cup \overline L_{\widehat Q}$,~~ $\overline U = \overline U\cup \overline U_{\hat w^+} \cup \overline U_{\widehat Q}$
		\EndFor
		\State \Return{$( L, U)$, $(\overline L, \overline U)$}
	\end{algorithmic}
\end{algorithm}

To finish the definition of the estimator, it remains to describe the selection procedures for the suitable blocks of questions that are used in  Line~\ref{line:choice_Q} of Algorithm~\ref{alg:double_trisection_sketch}. 
As explained above, we consider two procedures $\widehat Q := \widehat Q_{\cp}(\overline P, h, r)$ as in \Cref{eq:significant_blocks_estim} or $\widehat Q := \widehat Q_{WM}(\cT, \overline P, h, r)$ as in \cref{eq:significant_blocks_estim_wm} - which respectively apply to the \emph{oblivious} estimator $\hat{\pi}_{HT}$ and to the estimator $\hat{\pi}_{WM}$ \emph{with memory}.

\subsubsection{Definition of $\widehat Q_{\cp}$}

We start with $\widehat Q_{\cp}(\overline P, h, r)$. The corresponding estimator $\hat{\pi}_{HT}$ is called an oblivious hierarchical sorting tree estimator because $\widehat Q_{\cp}$ only depends on the restriction of the data to $\overline P$. As a consequence, the corresponding $\algoBlockSorting$ procedure (see~Algorithm~\ref{alg:algoBlockSorting_sketch}) which builds a trisection of a group $G$ of experts into three subgroups $(O,P,I)$ only depends on the observations on this set $G$ of experts. In other words, the recursive construction of the hierarchical sorting tree estimator is completely oblivious of the rest of the tree. Up to our knowledge, this feature is shared by most hierarchical clustering algorithms.

Fix some height $h\in \cH$ and $r\in \cR$. Intuitively, $\widehat Q_{\cp}(\overline P, h, r)$ amounts to focusing on the subset of questions around which the  empirical mean expert $\overline y(\overline P)$ has a high-variation. We provide some intuition on the rationale of this approach in the next subsection. More precisely, we define
the CUSUM statistic:
\begin{equation}\label{eq:def_cusum_sketch}\tilde r = 8\left[\left(\frac{32\zeta^2 }{|\overline P|h^2}\log(\tfrac{2d}{\delta})\right)\lor r\right] \spaceAnd \widehat \bC_{k, \tilde r}(\overline P) = \frac{1}{\tilde r}\left(\sum_{k' = k}^{k+ \tilde r - 1}\overline y_{k'}(\overline P) -  \sum_{k' = k- \tilde r}^{k-1}\overline y_{k'}(\overline P)\right) \enspace .\end{equation}
In a nustshell, $ \widehat \bC_{k, \tilde r}(\overline P)$ is the empirical variation of $\overline y(\overline P)$ at question $k$ and at scale $\tilde r \geq r$.  Then, we define $\widehat D_{\cp}$ as the set of questions where the CUSUM statistic is larger than $h/4$, and $\widehat Q_{\cp}\subset \cQ_r$ for the corresponding subset of blocks  or questions of size $r$.
\begin{equation}\label{eq:significant_blocks_estim}
	\widehat D_{\cp} = \big\{k \in [d] ~:~ \widehat \bC_{k, \tilde r}(\overline y(\overline P)) \geq \frac{h}{4} \big\} \spaceAnd \widehat Q_{\cp} = \big\{l \in \cQ_r ~:~ \widehat D_{\cp} \cap [l, l+r) \neq \emptyset \big\} \enspace .
\end{equation}
In~\eqref{eq:def_cusum_sketch}, the choice of $\tilde{r}\geq r$ is due to the fact that we need to compute an empirical mean $\bC_{k, \tilde r}(\overline y(\overline P))$ on enough questions so that its standard deviation is small compared to $h$.

\subsubsection{Definition of $\widehat Q_{WM}$}

Finally, we describe $\widehat Q_{WM}(\cT, \overline P, h, r)$ which corresponds to the estimator $\hat \pi_{WM}$.
The set $\overline{P}$ is a subset of a leaf $G$ of the tree $\cT$ and we write $t$ for its depth. 
As illustrated in~\Cref{fig:tree}, there is a natural order on the nodes of $\cT$ at depth $t$ that have been either obtained as subsets of type $O$ or $I$ in $\algoBlockSorting$ (\Cref{alg:algoBlockSorting_sketch}). We can index these nodes according to the ordering by setting $G^{(0)}=G$ and then $G^{(1)}$, $G^{(2)}$,\dots as the following groups. Similarly, $G^{(-1)}$, $G^{(-2)}$,\dots stand for the groups preceding $G^{(0)}$. See~\Cref{fig:tree} for an illustration. In fact,  with high probability, for any $a$, all the experts in $G^{(a+1)}$ are above the expert in $G^{(a)}$. As a consequence, the observations in $G^{(1)}$ and $G^{(-1)}$ can bring some informations on the behaviour of the experts in $\overline{P}\subset G^{(0)}$.

Fix $r\in \cR$ and $h\in \cH$. Define $\tilde{r}\in \cR$ as $\tilde r = 4(\lceil 2^9\log(4d|\cR|/\delta)\frac{\zeta^2}{|\overline P|h^2} \rceil^{dya} \lor r)$, where $\lceil x\rceil^{dya}= 2^{\lceil \log_2(x)\rceil}$. As before, $\tilde{r}\geq r$ stands for the scale which is required if we want to estimate the variation of $\overline{y}(\overline{P})$ with a standard error small compared to $h$.

Now, we consider any scale $r_{\cp}\in [4r,2\tilde{r}]\cap \cR$. The rationale is that, if $r_{\cp}< \tilde{r}$, we can reduce the standard deviations of the empirical means by considering an average over experts in neighboring groups. Define the upper neighborhood $\cV^{+}_{r_{\cp}}$ and lower neighborhood $\cV^{-}_{r_{\cp}}$  as the set of groups above $G$ and below $G$ that are necessary to have enough experts at scale $r_{\cp}$.
\begin{align}\label{eq:nb_v_plus}
	a^{+}_{WM} = \min \{ a ~:~ |G^{(1)}| + \dots + |G^{(a)}| \geq \tfrac{2^{11}\zeta^2 \log(4d|\cR|/\delta)}{r_{\cp}h^2} \}
	 & \spaceAnd  \cV^{+}_{r_{\cp}} = \bigcup_{a =  1}^{a^{+}_{WM}}G^{(a)}  \enspace ;  \\
	 \label{eq:nb_v_minus}
	a^{-}_{WM} = \min \{ a ~:~ |G^{(-1)}| + \dots + |G^{(-a)}| \geq \tfrac{2^{11}\zeta^2 \log(4d|\cR|/\delta)}{r_{\cp}h^2} \}
	 & \spaceAnd \cV^{-}_{r_{\cp}} = \bigcup_{a \in - a^{-}_{WM}}^{-1}G^{(a)}  \enspace .
\end{align}

For a given subset $\overline P \subset G$, we define the corresponding CUSUM statistic $\widehat \bC^{(\ext)}_{k, r_{\cp}}$ computed on the questions $[k - r_{\cp}, k+r_{\cp})$ and using the empirical mean observations in $\cV^{+}_{r_{\cp}}\cup\cV^{-}_{r_{\cp}}$ if $r < 2\tilde r$ and in $\overline P$ if $r_{\cp} = 2\tilde r$:
\begin{equation}\label{eq:def_cusum_sketch_wm}
	\widehat \bC^{(\ext)}_{k, r_{\cp}} = \frac{1}{r_{\cp}}
	\begin{cases}
		\sum_{k' = k}^{k + r_{\cp} - 1} \overline y_{k'}(\cV^{+}_{r_{\cp}} \cup \cV^{-}_{r_{\cp}}) - \sum_{k' = k - r_{\cp}}^{k - 1}\overline y_{k'}(\cV^{+}_{r_{\cp}} \cup \cV^{-}_{r_{\cp}}) & \text{ if } r_{\cp} \in [8r, 2\tilde r) \\
		\sum_{k' = k}^{k + r_{\cp} - 1} \overline y_{k'}(\overline P) - \sum_{k' = k - r_{\cp}}^{k - 1}\overline y_{k'}(\overline P)                                                                               & \text{ if } r_{\cp} = 2\tilde r
	\end{cases}
\end{equation}
If  $r_{\cp} = 2\tilde r$, this new definition of the CUSUM with memory matches the  definition \Cref{eq:def_cusum_sketch} in the previous paragraph. For $r_{\cp} < 2\tilde r$, we are not able to average on enough expert in $\overline{P}$. To deal with this issue, we average on a suitable number of neighboring experts.

Beside considering questions around which the variations of $\overline{y}(\overline P)$ are large enough, we also check whether, on the corresponding regions, the width of $\overline P$, that is the difference between the best expert and the worst expert in $\overline P$ is high enough. Given a question 
 $k \in [d]$, we define $\widehat \bDelta^{(\ext)}_{k, r_{\cp}}$ as the difference between the locals average on $[k - r_{\cp}, k+r_{\cp})$ of the neighbourhoods of $G$ that is
\begin{equation}\label{eq:stat_env_wm}
	\widehat \bDelta^{(\ext)}_{k, r_{\cp}} = \frac{1}{2r_{\cp}}\sum_{k' = k - r_{\cp}}^{k + r_{\cp} - 1} \overline y_{k'}(\cV^+_{r_{\cp}}) - \overline y_{k'}(\cV^-_{r_{\cp}}) \enspace .
\end{equation}
Since the groups $G^{(1)}$,\ldots, $G^{(2)}$ are above the best expert in $\overline{P}$ and since the groups $G^{(-1)}$, $G^{(-2)}$,\dots are below the worst expert in $\overline{P}$, this statistic $\widehat \bDelta^{(\ext)}_{k, r_{\cp}}$ overestimates the width of $\overline{P}$. In the next subsection,  we will explain why it is relevant to consider the width of $\overline{P}$.

We are now equipped to define the subsets $\widehat D_{WM}(\cT, \overline P, h, r)$ of suitable questions and the corresponding $\widehat Q_{WM}(\cT, \overline P, h, r)$ of corresponding blocks. 
\begin{align}
	\widehat D_{WM} & = \big\{k \in [d] ~:~ \exists r_{\cp} \in [4r, \tilde r]\cap \cR \mbox{ s.t. } \widehat \bC^{(\ext)}_{k, 2r_{\cp}} \geq \frac{h}{16} \mbox{ and } \widehat \bDelta^{(\ext)}_{k, r_{\cp}} \geq \frac{h}{16}  \big\} \enspace ;  \\
	\widehat Q_{WM} & = \{l \in \cQ_r ~:~ \widehat D_{WM} \cap [l, l+r) \neq \emptyset \}\enspace .\label{eq:significant_blocks_estim_wm}
\end{align}
In other words, $\widehat D_{WM}$ is made of questions for which there exists a scale $r_{\cp}$ such that simultaneously the empirical variations $\widehat \bC^{(\ext)}_{k, 2r_{\cp}}$ at scale $2r_{\cp}$ is at least of order $h$ and the empirical width at scale $r_{\cp}$ is at least of order $h$.

\subsection{Comments on the procedure and relation to the literature}
\label{subsec:comments}

These twin procedure are quite involved and combine several ingredients, some of them being already used by Liu and Moitra~\cite{liu2019minimax}. In particular, they introduced the key ideas of localization of the suitable blocks of questions through change-point detection on the mean expert and of a spectral clustering scheme for dividing blocks of experts. Still, we need to add several key elements in order to deal with the arguably more involved setting where $n\leq d$. We describe below how our procedure compares to~\cite{liu2019minimax}  and highlight also the main differences and new ideas. Also, despite the fact that our procedure is very involved, it remains computationally efficient. Overall, the full procedure requires  $O[\log^{c}(\frac{nd}{\zeta_- \delta})(nd^2+ n^3)]$ operations for some $c>0$. Indeed, each of the main steps of the algorithm correspond to matrix multiplications and computations of the largest eigenvector of a square symmetric matrix.

\medskip 

In this subsection, we  discuss three key steps of the algorithm: (i) the selection of blocks of questions corresponding to the high-variation regions of the average expert in the group as in the definition of $\widehat{Q}_{\cp}$, (ii) construction of the weights vector $w^+$ by a spectral procedure, (iii) the use of neighboring groups in $\widehat{Q}_{WM}$.

\subsubsection{Detecting high-variation regions of the average expert}

Recall that, for a fixed $r\in \cR$ and $h \in \cH$, $\widehat Q_{\cp}$ selects blocks of questions in which the variations~\eqref{eq:def_cusum_sketch} $\widehat \bC_{k, \tilde r}(\overline P)$ of the average $\overline{y}(\overline P)$ at question $k$ and at scale $\tilde{r}\geq r$ is higher than $h/4$. 

\begin{figure}[h]
	\includegraphics[scale=0.6]{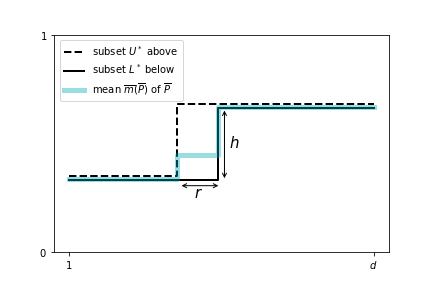}
	\caption{In this toy example, the group $\overline{P}$ is only made of two types of experts, those in $U^*$ and those in $L^*$. The high-variation region of $\overline{m}(\overline P)$  corresponds to the questions on which $U^*$ and $L^*$ differ.}
	\label{fig:simple_cp}
\end{figure}

To explain the rationale behind this choice, let us first consider a toy example depicted in Figure~\ref{fig:simple_cp}. Assume that the group $\overline{P}$ is made of two subgroups of experts $U^*$ and $L^*$ and that all the experts in $U^*$ and all the experts in $L^*$ are identical. Also, assume that the corresponding rows only differ on $r$ consecutive questions by $h$ and are otherwise identical. As illustrated in \Cref{fig:simple_cp}, it turns out that the expected average expert $\overline{m}(\overline{P})= \mathbb{E}[\overline y(\overline P)]$ needs to vary by $h$ at scale $r$ near the block of questions on which the two groups of experts are differing. This is due to the fact that both the rows corresponding to $U^*$ and $L^*$ are isotonic and that the row of $U^*$ is always larger or equal to that of $L^*$. As a consequence, by restricting our attention to the blocks of questions corresponding to high-variation regions of $\overline m(\overline P)$ (or in practice $\overline y(\overline P)$), we are able to much reduce the dimension of the problem and thereby to improve our ability to distinguish different experts.

Beyond this toy example, we show in~\Cref{lem:restriction_to_CP} that there exists a suitable scale $r\in \cR$ and a suitable height $h\in \cH$ such that, by restricting our attention to blocks of questions of size $r$ such that the expected average expert $\overline m(\overline P)$ varies by at least $h/2$, we are able to retain a significant proportion of the differences between experts in $\overline P$. In other words, focusing on regions of high-variation of $\overline{y}(\overline P)$ in the blocks $\widehat{Q}_{\cp}$ is, at least for some scale and some height, a suitable dimension reduction technique. This phenomenon was already observed in~\cite{liu2020better} and their procedure also uses such dimension detection techniques. In our paper, we also build upon this idea, which has also important consequences, in a related yet different manner, in the rectangular  case where $n\leq d$.

If we do not apply the spectral clustering sorting steps in $\hat{\pi}_{HT}$, that is, if we do not compute $\hat{v}$ and $\widehat{w}^+$ in  $\algoDoubleTrisection$, then  we would get a risk bound for $\hat{\pi}_{HT}$ of the order of $\zeta^{5/3}nd^{1/6}$ instead of that of Theorem~\ref{th:first_estimator}. In other words, the dimension reduction in $\widehat{Q}_{\cp}$ is alone sufficient to recover the optimal risk in the case where $d$ is quite large and $\zeta$ is mild - namely $n \leq \zeta^{2/3}d^{1/3}$ and $\zeta\in [1/d, \sqrt{d}]$.

\subsubsection{On the spectral estimation of the weights}

\begin{figure}[h]
\begin{minipage}[c]{0.45\textwidth}
	\includegraphics[scale=0.5]{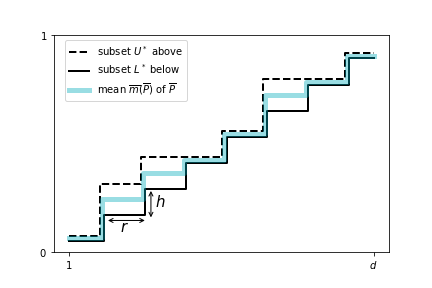}	
\end{minipage}
\ ~ \
\begin{minipage}[c]{0.5\textwidth}

	\[
		\mathbb{E}[Z-\overline{Z}]= \frac{1}{2}\sqrt{r}h  \begin{pmatrix}
				0 & -1 & -1 & 0 & 0 & -1 & 0 & 0\\
				0 & 1 & 1 & 0 & 0 & 1 & 0 & 0\\
				0 & 1 & 1 & 0 & 0 & 1 & 0 & 0\\
				0 & -1 & -1 & 0 & 0 & -1 & 0 & 0\\
				0 & -1 & -1 & 0 & 0 & -1 & 0 & 0\\
				0 & 1 & 1 & 0 & 0 & 1 & 0 & 0
			\end{pmatrix}
		\]
\end{minipage}
	\caption{In this toy example, the group $\overline{P}$ is only made of two types of experts $U^*$ and $L^*$ with $|L^*|=|U^*|=3$.}
	\label{fig:lower_bound}
\end{figure}

In this subsection, we explain how the computation of $\hat z$ in~\eqref{eq:ACP} and the corresponding weights $\hat w^+$ allow to improve over the $\zeta^{5/3}nd^{1/6}$ rate. Again, we start with a motivating toy example depicted in~\Cref{fig:lower_bound}. As previously, we consider a situation where $\overline P$ can be decomposed into two subgroups $U^*$ and $L^*$ of the same size. The corresponding rows $L^*$ are block-constant with blocks of questions of size $r$	 and increased by $h$ at the end of each block of questions. On the other hand, the corresponding lines of $U^*$ are, in each block of questions, either equal to the rows of $L^*$, or are exactly at a distance $h$ above. These last blocks of questions are the only ones which are informative when it comes to distinguishing the best experts in the group - namely $U^*$ - from the worst experts in the group - namely $L^*$. Some of the blocks corresponding to high-variation regions of the expected average row $\overline{m}(\overline{P})$ do not convey any information on the difference between $U^*$ and $L^*$ -- see~\Cref{fig:lower_bound}. In this example, at scale $r$, all the blocks of size $r$ are to be detected by the high variation dimension reduction step, that is $\widehat{Q}_{\cp}=\cQ_r$.
At the second step, we consider the corresponding aggregated matrix $Z- \overline{Z}$ at scale $r$ as defined in~\eqref{eq:aggregation}. To be more specific, let us assume that $|L^*|=|U^*|=3$. Then, $Z-\overline{Z}$ is a $6\times 8$ matrix whose expectation is of the form of the right panel in \Cref{fig:lower_bound}. 

In this specific example, the rank of this expected matrix is one and some of its columns are completely useless to decipher experts in $U^*$ from experts in $L^*$. In contrast, taking $w$ as a right singular vector associated to the largest singular value of this matrix would allow us to select the  significant blocks of questions while discarding the irrelevant ones. While this example is very specific, this still sheds some light on why spectral clustering procedure can be of interest for this problem and how it can help recover blocks of questions that are the most informative for dividing the experts.

Let us come back to a general matrix $M$ and to the spectral step of $\algoDoubleTrisection$ as described in the previous section. Up to a permutation of its rows, the expectation $\Theta^{(3)} -\overline{\Theta}^{(3)}$ of $Z^{(3)}-\overline{Z}^{(3)}$ is isotonic in each column. It turns out that the entries of any left singular vector associated to the largest singular value of $\Theta^{(3)} -\overline{\Theta}^{(3)}$ is, up to the permutation, either non-increasing or non-decreasing. As a consequence, the left-singular value of $\Theta^{(3)} -\overline{\Theta}^{(3)}$ can bring information on the underlying ranking. This property is at the heart of spectral ranking algorithms~\cite{vigna2016spectral}. 
 Unfortunately, contrary to the analysis of spectral methods in the Bradley-Luce-Terry model~\cite{chen2019spectral,chen2021optimal}, we cannot control the entry-wise deviations of the left singular eigenvector of $Z^{(3)}-\overline{Z}^{(3)}$ because the matrix $\Theta^{(3)}-\overline{\Theta}^{(3)}$ is non-parametric and does not necessarily exhibit any spectral gap. To handle this, Liu and Moitra~\cite{liu2020better} suggest to compute a right singular vector of $Z^{(3)}-\overline{Z}^{(3)}$ and, using another independent sample, to compare the experts based on the corresponding weighted average of the experts. Unfortunately, while their analysis provides near optimal results for $n = d$, this would not work for $n \leq d$. In $\algoDoubleTrisection$, we apply a more involved workaround (i) to handle possible heteroskedastic noise and (ii) to improve the convergence rates in comparison to Liu and Moitra~\cite{liu2020better}. Indeed, we first compute in~\eqref{eq:ACP} a debiased version $\hat{v}$ of the left-singular vector of $Z^{(3)}-\overline{Z}^{(3)}$. Then, we compute the image $[Z^{(5)}-\overline{Z}^{(5)}]^T \hat{v}$, threshold it, and take its absolute value to obtain our estimated weights $\hat w^+$. In principle, $\hat{w}^+$ aims at being close to the right first singular vector of $\Theta^{(3)} -\overline{\Theta}^{(3)}$. In comparison to Liu and Moitra \cite{liu2020better}, $\hat w^+$ better handles the situation where the matrix $\Theta^{(3)} -\overline{\Theta}^{(3)}$ is highly rectangular (with many columns) and where its corresponding right singular vector is nearly sparse.

\subsubsection{On the tree information and the definition of $\widehat{Q}_{WM}$}

The oblivious estimator $\hat{\pi}_{HT}$ based on $\widehat{Q}_{\cp}$ is only proved to achieve the suboptimal error of Theorem~\ref{th:first_estimator}. In this section, we explain how $\widehat{Q}_{WM}$ improves the performances of the procedure by relying on the neighboring experts to fix one possible weakness of $\widehat{Q}_{\cp}$ and so, improve the dimension reduction step.

\begin{figure}[h]
	\begin{minipage}[c]{0.48\textwidth}
		\centering
		\includegraphics[scale=0.55]{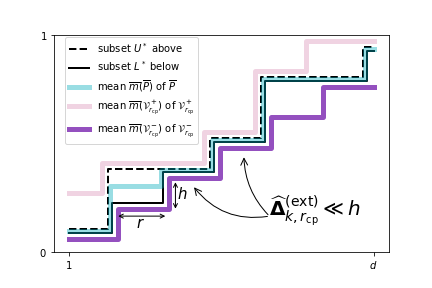}
	\end{minipage}
	\ ~ \
	\begin{minipage}[c]{0.48\textwidth}
		\centering
		\includegraphics[scale=0.55]{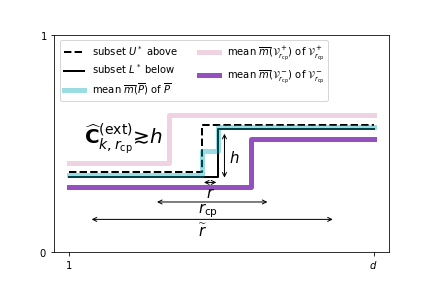}
		\end{minipage}

	\caption{In these two panels, the group $\overline{P}$ is only made of two types of experts, those in $U^*$ and those in $L^*$. The curves $\overline{m}(\cV_{r_{\cp}}^+)$ and $\overline{m}(\cV_{r_{\cp}}^-)$ respectively correspond to the expected average experts in the neighboring groups $\cV_{r_{\cp}}^+$ and $\cV_{r_{\cp}}^-$ defined in~\eqref{eq:nb_v_plus}.  In the left panel, the third and fourth blocks are not selected because the corresponding statistic $\widehat \bDelta^{(\ext)}_{k,r_{\cp}}$ is small.  In the right panel,  both the statistic $\widehat \bC^{(\ext)}_{k, r_{\cp}}$ and $\widehat \bDelta^{(\ext)}_{k,r_{\cp}}$ are large compared to $h$. If $\overline{P}$ is small, this block is selected using a scale $r_{\cp}\lesssim \tilde r$. } \label{fig:wm_cp}

\end{figure}

Indeed, $\widehat{Q}_{\cp}$ selects \emph{spurious} blocks of questions. In the previous toy example (\Cref{fig:lower_bound}), some of the blocks corresponding to high-variation values of the expected mean expert $\overline{m}(\overline{P})$ do not bring any suitable information for ordering the experts in $\overline{P}$ because, in these blocks, all the experts are close to each other. In other words, the width of $\overline{P}$, that is the difference between the best and worst experts in $\overline{P}$, is small. It is not possible to easily estimate this width from the observations in $\overline{P}$ since this would require to have sorted the experts in $\overline{P}$ in the first place. Still, we can estimate this width by comparing the average of experts that are above $\overline{P}$ with average of  experts that are below $\overline{P}$. A first idea would therefore be to consider a large enough number of experts above and below $\overline{P}$ in order to estimate the width with a small variance and to exclude regions such that the estimated width on a window of size $r$ is small compared to $h$. This is exactly the purpose of the 
statistic $\widehat \bDelta^{(\ext)}_{k, r}$. The selected blocks $\widehat Q_{WM}$ only contain regions such that the estimated width 
$\widehat \bDelta^{(\ext)}_{k, r}$ is large enough compared to $h$ --see the left panel in~\Cref{fig:wm_cp}. Unfortunately, the statistic $\widehat \bDelta^{(\ext)}_{k, r}$ may suffer from a large positive bias if the experts above or below $\overline{P}$ are away from $\overline{P}$. Moreover, considering only the scale $r$ is not sufficient because we are forced to average over many experts above and below $\overline{P}$ in order to have a small variance at this small scale, leading to a large bias. For this reason,  we consider all possible scales $r_{\cp}$ (in a dyadic grid) between $r$ and $\tilde r$.

Another important idea in the dimension reduction scheme is the following: If there is a region of questions in which, not only the mean experts of the group $\overline{P}$ but also the mean experts in neighboring groups of $\overline{P}$ have a high variation, it is interesting to detect this high-variation region by relying on all these neighboring groups in order to decrease the variance of the CUSUM statistic. With this idea, we are able to consider the CUSUM statistic at a smaller scale $r_{\cp} \leq \tilde r$ --see the right panel in~\Cref{fig:wm_cp}. This is exactly the purpose of the statistic $\widehat \bC^{(\ext)}_{k, r_{\cp}}$.

In our procedure, we build
a collection $\widehat{D}_{WM}$ that selects a question $k$ if there exists a scale $r_{\cp}$ in $[4r,\tilde{r}]$ such that both the CUSUM statistic $\widehat \bC^{(\ext)}_{k, 2r_{\cp}}$ at scale $2r_{\cp}$ is large and the empirical width $\widehat \bDelta^{(\ext)}_{k, r_{\cp}}$ is large. This combines the two ideas described in the previous paragraphs which, in turn, allows us to  further reduce the dimension in comparison to $\widehat{D}_{\cp}$ while ensuring that the selected questions in  $\widehat{D}_{WM}$ contains all the relevant regions to trisect $\overline{P}$, namely regions of size $r$, on which $\overline{P}$ has a variation at least of the order of  $h$ and the width of $\overline{P}$ is at least of the order of $h$.

Interestingly, in the square case where $n=d$ considered in~\cite{liu2020better} or more generally when $n\geq d$, this dimension reduction variant is not necessary to achieve the minimax risk as the oblivious estimator $\hat \pi_{HT}$ is already optimal. The dimension reduction scheme $\widehat{Q}_{WM}$ allows us to improve the risk bound from that Theorem~\ref{th:first_estimator} to that of Theorem~\ref{th:second_estimator_WM}. In the specific case where the noise level $\zeta$ is equal to one, the term  $n^{2/3}d^{1/3}$ in the risk bound is improved to the optimal one $n^{3/4}d^{1/4}$. Hence, building upon the neighboring groups in $\widehat{Q}_{WM}$ turns out to be the key ingredient to recover the minimax risk in the large $d$ regime where $n\in [d^{1/3},d]$.

\section{Partial observations}\label{sec:partial}

We now come back to the partial observation setting. Given $\lambda>0$, we are given $Poi(\lambda nd)$ independent observations $(x_t,y_{t})$ where $x_t$ is sampled uniformly in $[n]\times [d]$ and, conditionally to $x_t$, $y_{t}=M_{x_t}+ E_{x_t}$
is an observation of the full model~\eqref{eq:model_0} at position $x_t$. As noted above, $\lambda$ stands for the sampling effort and the larger $\lambda$, the more samples on average.

\subsection{Minimax Lower bound}

As in Section~\ref{subsec:lower_bound}, we first state a minimax lower bound in the case where the noise matrix $E$ is made of independent Gaussian random variables with variance $\zeta^2$. Note that the following minimax lower bound also handle the noise case where $\zeta=0$, i.e.~the noiseless case. 
\begin{theorem}\label{th:lower_bound_poisson}
	There exist universal constants $c$ and $c'$ such that the following holds for any $n\geq 2$, any $d\geq 1$, $\lambda>0$, and $\zeta\geq 0$:
	 \begin{align}
	\inf_{\hat \pi}  \sup_{\stackrel{\pi^*\in \Pi_n}{M:\,  M_{\pi^{*-1}}\in \mathbb{C}_{\text{BISO}}}} \mathbb E_{(\pi^*, M)} [\|M_{\hat \pi^{-1}} - M_{\pi^{*-1}}\|_F^2]
	\geq  c\left[ \left(\cR_F[n,d,\zeta/\sqrt{\lambda}] + \frac{n}{\lambda}e^{-2\lambda}\right)\bigwedge nd\right]\ . \label{eq:lower_bound_partia_observation}
	\end{align}
\end{theorem}
As in the previous minimax lower bound, the quantity $nd$ simply appears because the entries of $M$ lie in $[0,1]$. 
In~\eqref{eq:lower_bound_partia_observation}, we recognize two terms. First, $\cR_F[n,d,\zeta/\sqrt{\lambda}]$ corresponds to the minimax risk for recovering $\pi^*$ in a full observation model with noise $\zeta/\sqrt{\lambda}$. The second term $\frac{n}{\lambda}e^{-2\lambda}$ does not depend on $\zeta$ and is also present in the noiseless setting. It simply quantifies the fact that, for $\lambda<1$, observations are lacking so that it is impossible to correctly rank experts if there are no observations on the questions on which they are distinct. As the minimax lower bound in~\eqref{eq:lower_bound_partia_observation} turns out to be nearly tight in light of Theorem~\ref{th:second_estimator_WM_poisson_intro}, we refer to~\eqref{eq:lower_bound_partia_observation} as the minimax risk in the following. For the purpose of the discussion, we will first focus on the case where $\zeta=1$ and $\lambda <1$, which corresponds to the case where we really have partial observations on the matrix. We will then turn to $\zeta=1$ and $\lambda >1$, which corresponds to the case where we observe several times each entry of the matrix. Finally, we discuss the noiseless case where $\zeta=0$.

\medskip 

\noindent 
{\bf Low-sample size}. 
We first focus on the case where $\zeta=1$ and $\lambda <1$, which corresponds to the case where we really have partial observations on the matrix. If $\lambda \leq 1/d$, then the minimax risk is of the order of $nd$ and it is impossible to perform significantly better  than a random guess. This is not surprising as there are, in expectation, less than one observation on each row. For $\lambda\in [1/d, 1]$, the minimax risk is of the order of 
\[
\frac{n d^{1/6}}{\lambda^{5/6} } 	\bigwedge    \frac{n^{3/4}d^{1/4}}{\lambda^{3/4}} +  \frac{n }{\lambda }	
\]
In the rectangular case where $n\geq d$, the minimax risk is then of the order of $n/\lambda$ for $\lambda\in [1/d,1]$. When $n\in [d^{1/3},d]$, the minimax risk is of the order of $n/\lambda$ for $\lambda \in [1/d, n/d]$, and of the order of $\frac{n^{3/4}d^{1/4}}{\lambda^{3/4}}$ for $\lambda\in [n/d,1]$. For even smaller $n\leq d^{1/3}$, there is one more regime since 
\[
	\cR_F[n,d,1/\sqrt{\lambda}]\asymp \left\{
		\begin{array}{cc}
			\frac{n}{\lambda} & \text{ if }\lambda\in [\frac{1}{d},\frac{n}{d}] \ ; \\
			\frac{n^{3/4}d^{1/4}}{\lambda^{3/4}} & \text{ if }\lambda\in [\frac{n}{d},\frac{n^3}{d}] \ ; \\
			\frac{nd^{1/6}}{\lambda^{5/6}} &\text{ if }\lambda\in [\frac{n^3}{d},1]\ . 
		\end{array}
	\right. 
\]

\medskip 

\noindent 
{\bf Large-sample size}. In the setting where $\lambda>1$ and $\zeta=1$, there are several observations per entries. In this case, there are many regimes in~\eqref{eq:lower_bound_partia_observation} that depend on $n$, $d$, $\zeta$, and $\lambda$. To simplify the discussion, we focus here on the case $n=d$ and $\zeta=1$. Then, the minimax risk is of the order of  $\frac{n^{3/4}d^{1/4}}{\lambda^{3/4}}$ for $\lambda\in [1,n^2]$ and is of the order of $n\sqrt{d}/\lambda$ for $\lambda\geq n^2$. This 'easy rate' $n\sqrt{d}/\lambda$ is achieved by the simple procedure that ranks the experts according to the row sums~\cite{shah2016stochastically,mao2020towards}. This simple method turns out to be optimal in the regime where there are more than $n^2$ observations per entry. 

\medskip 

\noindent
{\bf Noiseless case}. In the extreme case where $\zeta=0$ and $\lambda\geq 1/d$, the minimax risk is of the order of $(n/\lambda) \times e^{-2\lambda}$, which, for some small $\lambda$ is of the order of $n/\lambda$. This minimax lower bound is quite simple to prove. Without loss of generality, suppose that $1/\lambda$ is an integer. Consider a matrix $M$ such that  all its columns, except its $1/\lambda$ first ones are constant and equal to one, so that it boils down to considering a ranking problem of size $n\times (1/\lambda)$. In this reduced model, there are two types of experts: (a) experts that are constant and equal to zero and (b) experts that are constant and equal to one. Obviously, if one is given at least one noiseless observation on a row, then it is possible to assign it to a group. However, on each row there is a probability $e^{-1}$ of having no observations. Hence, on expectations there are $n/e$ experts that are impossible to classify. For this reason, any estimator must suffer from a risk at least of the order $n/\lambda$.

\subsection{Reduction to the full observation model}

We now describe a scheme to adapt the estimators $\hat{\pi}_{HT}$ and  $\hat{\pi}_{WM}$ that we developed in the full observation setting of Section~\ref{sec:full}, to this more general Poissonian setting~\eqref{eq:model_partial}, which encompasses the partial observation setting as well as the over-complete observation setting where each entry is sampled several times. Roughly, if $\lambda$ is small, we simply decrease the number of columns of the matrix $M$ in order to obtain a reduced matrix with full observations. Conversely, if $\lambda$ is really large, which corresponds to the case of multiple observations per entry, we simply average the multiple observations per entry to reduce the noise levels.

As in Section~\ref{sec:full}, we fix $\delta\in (0,1)$ that will correspond to a small probability. Given this $\delta>0$, we denote $\Upsilon^*= \Upsilon^*(n,d,\zeta/\sqrt{\lambda\vee 1})$ the number of independent samples required in Section~\ref{sec:full} for the estimation through $\hat{\pi}_{HT}$ or $\hat{\pi}_{WM}$ of the $n\times d$ matrix $M$ with  a noise level equal to $\zeta/\sqrt{\lambda\vee 1}$. Recall that $\Upsilon^*$ is of  the order of $\log^{8}(nd(\lambda\vee 1)/(\delta\zeta_-))$. 

Define $\lambda_{-}= \lambda/[4 \Upsilon^*]$. For any $i\in [n]$ and any $S\subset [d]$, we write $n_{i,S}$ the number of observations in the sample falling in $\{i\}\times S$, that is $n_{i,S}= |\{t:\   x_t\in \{i\}\times S  \}|$. The following lemma is a simple consequence of Chernoff inequality for Poisson random variables. 

\begin{lemma}\label{lem:subsampling_1}
Assume that $\lambda_{-}\in [2/d,1]$, we fix $l(\lambda)= \lfloor 1/\lambda_{-}\rfloor$. With probability higher than 
$1 - \delta$, we have 
\[
\min_{i\in [n]}\min_{j\in [\lfloor d/l(\lambda)\rfloor] }n_{i,[(j-1)l(\lambda)+1, jl(\lambda)]} \geq \Upsilon^*\ . 	
\]
Now assume that $\lambda_{-}> 1$. With probability higher than $1-\delta$, we have 
\[
\min_{i\in [n]}\min_{j\in [d] }n_{i,\{j\}} \geq \lambda_{-}\Upsilon^*\ . 	
\]
\end{lemma}

Henceforth, we work under the event introduced in the previous lemma. If this event does not hold, we choose $\hat{\pi}_{WMP}$ arbitrarily. To build $\hat{\pi}_{WMP}$, we consider three subcases that depend on the value of $\lambda$:
\begin{enumerate}
	\item {\bf Very small sample size}. If $\lambda_{-}\leq 2/d$, then we simply choose $\hat{\pi}_{WMP}$ uniformly at random over the set of all possible permutations. While this choice does not depend on the data and could therefore seem sub-optimal, it is not the case, as the minimax lower bound states that it is impossible to perform better than random guess in this setting.
 \item {\bf Small sample size}. If $\lambda_{-}\in [2/d, 1]$, then we build $\Upsilon^*$ matrices $\cY^{\downarrow}= (Y^{\downarrow(0)}, Y^{\downarrow(1)}, Y^{\downarrow(\Upsilon^*-1)})$ of size $n\times \lfloor d/l(\lambda)\rfloor$ in the following way. For any $i\in [n]$, $j\in [\lfloor d/l(\lambda)\rfloor]$ and $s\in [0,\Upsilon^*-1]$, $Y^{\downarrow(s)}_{i,j}=y_t$ where $t$ is the $(s+1)$-th observation such that $x_{t}\in \{i\}\times [l(\lambda)j +1, l(\lambda)(j+1)]$. On the event of Lemma~\ref{lem:subsampling_1}, this definition is valid as we observe enough samples for any $i,j$. Then, we compute $\hat{\pi}_{WMP}$  as the variant $\hat{\pi}_{WM-SR}$, introduced in Section~\ref{sec:semi_random}, applied to this sample of reduced matrices.
 \item {\bf Large sample size}. If $\lambda_{-}\geq 1$, then we build $\Upsilon^*$ matrices $\cY^{\downarrow}= (Y^{\downarrow(0)}, Y^{\downarrow(1)}, Y^{\downarrow(\Upsilon^*-1)})$ of size $n\times d$ in the following way. For any $i\in [n]$, $j\in [d]$, $l\in [\lfloor \lambda_{-}\rfloor]$, and $s\in [0,\Upsilon^*-1]$, define $Y^{\downarrow(s)}_{i,j}=\frac{1}{\lfloor\lambda^-\rfloor}\sum_{t}y_t$ where the $y_t$'s are the $z$-th observations such that $x_{t}=(i,j)$ with $z\in [1+(s-1)\lfloor\lambda_-\rfloor ,s\lfloor\lambda_-\rfloor]$. In other words, we build the samples 
 $\cY^{\downarrow}$ be averaging $\lfloor\lambda_-\rfloor$ observations on each entries. Again, on the event of Lemma~\ref{lem:subsampling_1}, this definition is valid as we observed enough samples for any $i,j$. Then, we define $\hat{\pi}_{WMP}$ as  $\hat{\pi}_{WM}$ applied to this sample of averaged matrices. By averaging the independent observations, we reduce the noise level of each entry from $\zeta$ to $\zeta/\sqrt{\lambda_-}$. 
 
\end{enumerate}

For $\lambda_{-}\leq 2/d$, there are very few observations on each row so that it is very difficult to compare the experts. For $\lambda_{-}\in [2/d,1]$, we have access to less than $\Upsilon^{*}$ noisy observations of the matrix $M$. The rationale of our procedure is to group together $l(\lambda)$ consecutive questions together in such a way that there are enough observations on each of these groups. The resulting matrices of observations $Y^{\downarrow (s)}$ have around $\lambda d/\Upsilon^*$ columns. We could have applied the procedure $\hat{\pi}_{WM}$ defined in the previous section to $\cY^{\downarrow}$, but the corresponding subGaussian norm of the noise would be $1+\zeta$ (instead of $\zeta$) because there is additional variability coming from the fact that any entry in the reduced matrices has been sampled uniformly among $l(\lambda)$ entries in the original matrices. This would lead us to a procedure achieving the minimax rate with respect to $n$, $d$, and $\lambda$ but with a suboptimal dependency with respect to $\zeta$ since $\zeta$ would be replaced by $\zeta+1$. This is the reason why, for $\lambda_-\in [2/d,1]$, we rely on  a slight variant  $\hat{\pi}_{WM-SR}$ (see Section~\ref{sec:semi_random}) of  $\hat{\pi}_{WM}$ that builds upon the fact that the variations that are due to the aggregation of $M$ are very specific.

\begin{theorem}\label{th:second_estimator_WM_poisson_intro}
	There exist four numerical constants $c_1$--$c_4$ such that the following holds. Fix $\delta= \zeta^2_-[(\lambda\vee 1) nd]^{-2} $. For any permutation $\pi^*\in \Pi_n$ and any matrix $M$ such that $M_{\pi^{*-1}}\in \mathbb{C}_{\text{BISO}}$, 
	 the  sorting tree estimator $\hat{\pi}_{WMP}$ defined above satisfies
\beq\label{eq:upper_bound_partial_observation}
	 \E\left[\|M_{\hat \pi_{WMP}^{-1}} - M_{\pi^{*-1}} \|_F^2\right] \leq c_1\log^{c_2}\left(\tfrac{nd(\lambda\vee 1)}{\zeta_-}\right) \left[\cR_F(n, d,\zeta\lambda^{-1/2})  + \frac{n}{\lambda}e^{-c_3 \lambda\log^{-c_4}(\tfrac{nd(\lambda\vee 1)}{\zeta_- })}  \right]\enspace . 
\eeq 
\end{theorem}

 Up to logarithmic terms and up to the logarithmic term inside the exponential term~\eqref{eq:upper_bound_partial_observation}, both the minimax  upper bound~\eqref{eq:upper_bound_partial_observation} and lower bound~\eqref{eq:lower_bound_partia_observation} match for all values of $n$, $d$, $\lambda$ and $\zeta$. As a consequence, this problem of estimating a single permutation $\pi^*$ does not exhibit any significant computational gap.

\medskip

Let us further discuss and compare the exponential term $n\lambda^{-1}	e^{-c_3 \lambda\log^{-c_4}(\tfrac{nd(\lambda\vee 1)}{\zeta_- })}$ in~\eqref{eq:lower_bound_partia_observation} and $n\lambda^{-1}e^{-2\lambda} $ in~\eqref{eq:upper_bound_partial_observation}. First, observe that these two terms are larger than $R_F(n,d,\zeta \lambda^{-1/2})$ only when the noise level $\zeta$ is small, so that it is relevant to discuss them only when $\zeta\ll 1$. Second, note that there is a significant mismatch between these exponential terms only when $\lambda$ is close to one, up to a polylogarithmic factor, since otherwise, either the exponential is close to one  (for $\lambda\leq 1$) or the exponential is so small that it becomes negligible in comparison to $R_F(n,d,\zeta \lambda^{-1/2})$. One may object that the logarithmic term $\log^{-c_4}(\tfrac{nd(\lambda\vee 1)}{\zeta_- })$ may be large in case $\zeta$ is really small --think e.g. of $\zeta=e^{-nd}$. 
Let us consider this extremely low noise setting where, say $\zeta\leq 1/(nd)^2$. If one applies the procedure $\hat{\pi}_{WMP}$ with $\zeta_0=1/(nd)^2\geq \zeta$, then the logarithmic terms become bounded inside the exponential. Since $\cR_F(n, d,\zeta_0\lambda^{-1/2})$ is always smaller than $nd\bigwedge \frac{n}{\lambda}e^{-c_3 \lambda}$ provided that $\lambda\leq 1$, this  estimator achieves the risk bound $\frac{n}{\lambda}\wedge nd$, which is optimal for all $\zeta\in [0,1/(nd)^2]$ and all $\lambda \leq 1$.  To sum up, there is gap between our minimax lower and upper bounds only either (i) in the low-noise level with large but mild sampling effort, that is $\zeta=o(\log^{-c}(nd))$, $\zeta\geq (nd)^{-2}$, and $\lambda\in 	 [\log^{c}\log(nd),\log^{c'}(nd)]$ for some $c$ and $c'>0$ or (ii) in the extremely low noise level with large sampling effort, that is $\zeta\leq (nd)^{-2}$ and $\lambda\geq 1$.

\medskip

In $\hat{\pi}_{WMP}$, we have plugged in the hierarchical sorting tree estimator with memory $\hat{\pi}_{WM}$. If we had plugged in the oblivious hierarchical sorting tree estimator $\hat{\pi}_{HT}$, then the resulting estimator would satisfy a similar rate similar to \eqref{eq:upper_bound_partial_observation} except that the term $n^{3/4}d^{1/4}/\lambda^{3/4}$ would be replaced by the slower rate $n^{2/3}d^{1/3}/\lambda^{2/3}$.

\subsection{Reconstruction of the matrix $M$}\label{ss:recM}

In this subsection, we assume again that the noise level $\zeta=1$ to simplify the exposition. 
As alluded in Section~\ref{sec:full}, it is quite straightforward to estimate the matrix $M$ and control the 
corresponding loss $\|\widehat{M}- M\|_F^2$ by a simple subsampling step explained e.g. in~\cite{mao2020towards} that we recall here. First, we split the sample into two part by assigning independently each observation to the first subsample with probability $1/2$ and the second subsample with probability $1/2$. Then, we use the first subsample  to estimate the permutation $\hat{\pi}$ of the experts. As for the second subsample $(x^{(2)}_t,y_t^{(2)})$, we define the empirical observed matrix $Y^{(2)}$ by 
\[
Y^{(2)}_{i,j}= \frac{1}{\lambda}\sum_{t}y_t^{(2)}\1_{x_t^{(2)}=(i,j)}. 
\]
Then, we compute the least-square estimator $\widehat{M}_{\hat{\pi}}$ of $M_{\hat{\pi}}$ in the class of bi-isotonic matrix $\widehat{M}_{\tilde{\pi}}= \arg\min_{B\in \mathbb{C}_{\text{BISO}}}\|B - Y^{(2)}_{\tilde{\pi}}\|_F^2$. This estimator can be computed in near linear-time~\cite{kyng2015fast}. Then, Proposition~3.3 in \cite{mao2020towards} states, that with high probability, the loss $\|\widehat{M}-M\|_F^2$ is, up to logarithmic terms, smaller than the sum of the minimax risk for estimating a bi-isotonic matrix $B$ and the loss $\|M_{\tilde{\pi}^{-1}}-M_{\pi^{*-1}}\|_F^2$. Plugging this proposition with our estimator $\hat{\pi}_{WMP}$ with $\delta=(\lambda\vee 1)/(np)$, we readily arrive to the following risk bound for the corresponding estimator $\widehat{M}_{WMP}$. 

Define $\cR_1(n,d,\lambda)= \sqrt{\frac{nd}{\lambda}}\bigwedge \frac{nd}{\lambda^{2/3}(n\vee d)^{2/3}} \bigwedge  \frac{nd}{\lambda}$. Mao et al.~\cite{mao2020towards} have proved that, up to polylogarithmic factor and up to a possible additive term $(n\wedge d)/\lambda$,  the minimax risk in square Frobenius norm for estimating a bi-isotonic matrix with partial observations  is $\cR_1(n,d,\lambda)$. 
\begin{corollary}\label{cor:reconstruction}
	There exist two numerical constants $c$ and $c'$ such that the following holds. For any permutation $\pi^*\in \Pi_n$ and any matrix $M$ such that $M_{\pi^*}\in \mathbb{C}_{\text{BISO}}$, we have 
\begin{eqnarray}\nonumber
	 \E\left[\|\widehat{M}_{WMP} - M\|_F^2\right] &\leq& (nd) \bigwedge \left[ c   \log^{c'}((\lambda\vee 1)nd)\left( \cR_1(n,d,\lambda)+  \cR_F(n,d,\lambda^{-1/2})\right) \right]\\
	 &\leq & (nd) \bigwedge \left[ c   \log^{c'}((\lambda\vee 1)nd)\left( \cR_1(n,d,\lambda)+  \frac{n}{\lambda}\right) \right]\ . 
	 \label{eq:upper_bound_partial_observation_reconstruction}
\end{eqnarray}
\end{corollary}
The proof is a straightforward consequence of Proposition 3.3 in~\cite{mao2020towards} and Theorem~\ref{th:second_estimator_WM_poisson_intro} and is therefore omitted. It turns out that $\cR_F(n,d,\lambda^{-1/2})$ is always smaller than $\cR_1(n,d,\lambda)+\frac{n}{\lambda}$, so that the cost of reconstruction for not knowing $\pi^*$ is $n/\lambda$.

This risk bound~\eqref{eq:upper_bound_partial_observation_reconstruction} is minimax optimal, up to polylogarithms, and this for all possible values of $n\geq 2$, $d$, and $\lambda>0$. Indeed, in their Theorem 3.1, Mao et al.~\cite{mao2020towards} provide a matching minimax lower bound in $\cR_1(n,d,\lambda)$ in the specific case where $n\geq d$, but their proof easily extends to the case where $n\leq d$.  Besides, our proof of the minimax lower bound $\frac{n}{\lambda}$ in Theorem~\ref{th:lower_bound_poisson} for the problem of estimating $\pi^*$ straightforwardly extends to the problem of matrix estimation (recall that we consider $\zeta=1$ here).

The least-square estimator $\hat{\pi}_{LS}$ of Mao et al. has also been proved to achieve the minimax risk for $n\geq d$ --see their theorem 3.1 in~\cite{mao2020towards}. However, no efficient algorithm is known for computing this estimator in  $\hat{\pi}_{LS}$, so that our estimator $\widehat{M}_{WMP}$ is, to the best of our knowledge, the first efficient minimax-optimal estimator for estimating $M$ in this context, for any values of $n,d,\lambda$.

\subsection{Bounds for the max loss of Mao et al.~\cite{mao2020towards}}\label{ss:maxloss}

In~\cite{mao2020towards}, Mao et al. control, for an estimator $\hat \pi$ of the permutation, a different loss from ours. Up to normalization factors, they indeed focus on the maximum $l_2$ norm of the rows of $(M_{\hat \pi^{-1}})_{i,.} - (M_{\pi^{*-1}})_{i,.}$, that is
\beq\label{eq:loss}
l_{\infty}(\hat {\pi},\pi^*)= \sup_{i\in [n]}\|(M_{\hat \pi^{-1}})_{i,.} - (M_{\pi^{*-1}})_{i,.}\|_2^2\ .
\eeq
This loss also considered in~\cite{shah2019feeling,chatterjee2019estimation} corresponds to some maximum error of the estimated permutation so that $l_{\infty}(\hat {\pi},\pi^*)\geq \|M_{\hat \pi^{-1}} - M_{\pi^{*-1}}\|_F^2/n$. Alternatively, we can define the loss $l_{err}$
\[
l_{err}(\hat {\pi},\pi^*)=	\max_{\stackrel{i,j\in [n]\ :}{ \hat {\pi}(i)< \hat {\pi}(j)\text{ and }\pi^*(i)> \pi^*(j)  } }\|M_{i,.}-M_{j,.}\|_2^2 \ ,
\]
which quantifies the maximum distance between two experts that have not been ranked in a consistent manner. The loss $l_{\infty}$ and $l_{err}$ turn out to be equivalent as stated in the following lemma. 
\begin{lemma}\label{lem:l_infty}
For any permutation $\hat{\pi}$, we have 
\beq\label{eq:upper_l_infty}
l_{\infty}(\hat{\pi},\pi^*)\leq l_{err}(\hat{\pi},\pi^*)\leq 4l_{\infty}(\hat {\pi},\pi^*)
\eeq
\end{lemma}

To simplify the discussion in this section, we assume  again that the noise level $\zeta$ equals one. 
Mao et al.~\cite{mao2020towards} provide a simple polynomial time $\hat{\pi}_{\text{ref}}$ achieving  
\beq\label{eq:risk_mao_sup}
	\E[l_{\infty}(\hat{\pi}_{\text{ref}},\pi^*)]\lesssim d \bigwedge  \frac{d^{1/4}}{\lambda^{3/4}}\log^{3/4}( n)\ . 
\eeq
Conversely, they prove in their Theorem 3.7 that any estimator $\hat {\pi}$ that only ranks the experts $i$ and $j$ according to the differences of the observations on the rows $i$ and $j$ must incur this risk bound-- see~\cite{mao2020towards} for further details. Besides, they conjecture that the risk bound~\eqref{eq:risk_mao_sup} cannot be improved. In~\cite{liu2020better}, Liu and Moitra already pointed out that the max loss $l_{\infty}(\hat{\pi},\pi^*)$ is less suited than the loss $\|M_{\hat  \pi} - M_{\pi^*}\|_F^2$ for the purpose of estimating the matrix $M$ --see the discussion in the previous subsection. Still, controlling the max loss $l_{\infty}(\hat {\pi},\pi^*)$ may be an objective per se, and the study of its minimax value and of the existence of related minimax estimators is relevant. In the following proposition, which is mainly a consequence of our results and proof techniques, we disprove Mao et al.'s conjecture by introducing an estimator $\hat{\pi}_{PC}$ achieving a faster rate than~\eqref{eq:risk_mao_sup}. Besides, this rate turns out to be minimax-optimal.

\begin{proposition}\label{prp:risk_minimax_sup}
There exist numerical constants $c$, $c'$, and $c''$ such that the following result holds. There exists a polynomial-time estimator $\hat{\pi}_{PC}$ that performs pair-wise comparisons between the experts and that achieves the risk bound 
\beq\label{eq:upper_bound_l_infty}
\E[l_{\infty}(\hat{\pi}_{PC},\pi^*)] \leq c   \log^{c'}\left(nd(\lambda\vee 1)\right)\left[\frac{ d^{1/6}}{\lambda^{5/6} }   \bigwedge \frac{\sqrt{d}}{\lambda} \right]\bigwedge d  \enspace  .
\eeq
Conversely, for any $n\geq 2$, any $d\geq 1$, and $\lambda>0$, we have 
\beq\label{eq:lower_bound_l_infty}
\inf_{\hat \pi}  \sup_{\pi^*\in \Pi_n}\sup_{M:\,  M_{\pi^{*-1}}\in \mathbb{C}_{\text{BISO}}} \mathbb E_{(\pi^*, M)} [l_{\infty}(\hat{\pi},\pi^*)]
\geq c''  \left[\frac{ d^{1/6}}{\lambda^{5/6} } \bigwedge \frac{\sqrt{d}}{\lambda} \bigwedge d\right] \ . 
\eeq
\end{proposition}
For $\lambda\leq 1/d$, it is not possible to perform significantly better than random guess. Then, in the interesting regime $\lambda\in [1/d,d^2]$, the risk is of the order of $\frac{ d^{1/6}}{\lambda^{5/6}}$. 
It turns out that this rate corresponds, up to polylogarithmic terms, to the minimal distance between two experts so that one is able to consistently compare them. 
For very large sample size $\lambda\geq d^2$, we arrive at the easy regime which is of the order of $ \frac{\sqrt{d}}{\lambda}$. 

\medskip 

The estimator $\hat{\pi}_{PC}$ is based on pairwise comparisons. For any two experts $i$ and $j$, we apply the procedure $\hat{\pi}_{WMP}$ to $i$ and $j$ with $\delta= [(\lambda\vee 1) (n^2d)]^{-2}$ . If the trisection $(O, P, I)$ is of the form $(\emptyset, \{i\},\{j\})$, we return $i \prec j$.   If the trisection $(O, P, I)$ is of the form $(\emptyset, \{j\},\{i\})$, we return $j \prec i$. Otherwise, we return nothing. Applying this comparison algorithm to all $(i,j)$, we recover a set of pairwise comparisons  $\cP\cC= \{(i,j): \ i\prec j\}$. With high probability --see the proof for more details--, it turns that $\cP\cC$ satisfies two properties:
\begin{enumerate}
	\item[(i)] $\cP\cC$ is consistent. For any $(i,j)\in \cP\cC$, we have $\pi^*(i)< \pi^*(j)$. 
	\item[(ii)] $\cP\cC$ contains all $2$-tuple of experts that are far apart. More precisely, $\cP\cC$ contains all $(i,j)$ such that $\pi^*(i)< \pi^*(j)$, and 
	\beq
	\|M_{i,.}-M_{j,.}\|^2_2\geq 	c   \log^{c'}\left(nd(\lambda\vee 1)\right)\left[\frac{ d^{1/6}}{\lambda^{5/6} }   \bigwedge \frac{\sqrt{d}}{\lambda} \right]\bigwedge d  \enspace  ,
	\eeq
	for suitable constants $c$ and $c'$. 
\end{enumerate}
Then, define the function $\phi:[n]\mapsto \mathbb{N}$ by $\phi(j)=|\{(i,j):\ (i,j)\in \cP\cC\}|$ which simply counts the number of experts $i$ that are detected to be lower than $j$. Finally, we build $\hat{\pi}_{PC}$ as any permutation that ranks the experts consistently with $\phi$.

\medskip 
In fact, the procedure for computing  $\hat{\pi}_{PC}$ could be greatly simplified. Indeed, as we only perform pairwise comparisons, some parts of $\algoBlockSorting$ turn out to be irrelevant. For instance, the PCA steps are not required. Besides, the sample splits could be avoided and it could even  be possible to work with a single observation. As the problem of optimal permutation recovery with respect to the $l_{\infty}$ loss is not the main scope of this paper, we do not provide a simplified and dedicated algorithm. Besides, we conjecture that our original estimator $\hat{\pi}_{WMP}$ also achieves the minimax risk~\eqref{eq:upper_bound_l_infty} with respect to the $l_{\infty}$ loss.

\paragraph*{Acknowledgements.} The work of A. Carpentier is partially supported by the Deutsche Forschungsgemeinschaft (DFG) Emmy Noether grant MuSyAD (CA 1488/1-1), by the DFG - 314838170, GRK 2297 MathCoRe, by the FG DFG, by the DFG CRC 1294 'Data Assimilation', Project A03, by the Forschungsgruppe FOR 5381 "Mathematical Statistics in the Information Age - Statistical Efficiency and Computational Tractability", Project TP 02, by the Agence Nationale de la Recherche (ANR) and the DFG on the French-German PRCI ANR ASCAI CA 1488/4-1 "Aktive und Batch-Segmentierung, Clustering und Seriation: Grundlagen der KI" and by the UFA-DFH through the French-German Doktorandenkolleg CDFA 01-18 and by the SFI Sachsen-Anhalt for the project RE-BCI. The work of E. Pilliat and N. Verzelen has been partially supported by ANR-21-CE23-0035 (ASCAI).

\appendix

\section{Full description of the procedures}\label{sec:algorithms}

In this section, we provide a fuller description of the estimators $\hat{\pi}_{HT}$ and $\hat{\pi}_{WM}$ as a collection of algorithms. We will rely on this description in the analysis of these estimators. To ease its understanding, we make this section completely self-contained. As a consequence, the material presented here is partly redundant with Section~\ref{sec:algo_sketch}.

\subsection{Sorting a group of experts}

Some of the notation have already been introduced in Section~\ref{sec:algo_sketch}. Still we define them again here for the sake of completeness. 
We write $\cD$ for the set of all dyadic numbers, that is $\cD = \{ 2^k ~:~ k \in \bbZ \}$. Equipped with $\cD$, let 
\begin{equation}\label{eq:definition_R_D}
	\cR = \cD \cap [1,d]  \spaceAnd \cH = \cD \cap \left[\frac{\zeta^2}{nd},1\right] \enspace ,
\end{equation}
respectively  denote the dyadic collection of numbers beween $1$ and $d$ and the dyadic  collection of numbers between $1/nd$ and $1$. 

Besides for an integer $r\in \cR$, we write $\cQ_r$ for the regular grid of $[d]$ of width $r$: 
$$ \cQ_r = \left\{1, r+1, 2r + 1,\dots \floor{\frac{d}{r}}r+1 \right\} \enspace .$$

In contrast to Section~\ref{sec:algo_sketch}, we start by describing the simple comparison routine before moving to the dimension reduction techniques and to the general architecture of the procedures.

Given a collection $\overline P$ of experts, some data $Z\in \mathbb{R}^{\overline P\times Q}$ and a direction $w\in (\mathbb{R}^{+})^{Q}$ and a pivot $\gamma\in [1:|\overline P|]$, the following pivoting algorithm sorts the experts in $\overline P$ according to the projection of the data onto the vector $w$.  More precisely, it returns four subsets $\overline{L}\subset L$ and $\overline{U}\subset U$ of experts such that the $\gamma$-th best expert according to the $(Z,w)$-order - as defined above Equation \Cref{eq:pivot_trisection} - is significantly above all experts in $L$ and below all experts in $U$. The subsets $L$ and $\overline{L}$ (resp. $U$ and $\overline{U}$) differ in the level of significance we require. 
We define the tuning parameters $\beta_{\tris}$ and $\overline{\beta}_{\tris}$ for $\algoPivot$
\beq\label{eq:definition_betris}
\beta_{\tris}= 4\sqrt{2}\zeta\ , \quad \quad \overline{\beta}_{\tris}= 8\sqrt{2}\zeta.
\eeq

\begin{algorithm}[H]
	\caption{$\algoPivot(Z,w,\gamma)$ \label{alg:pivot}}
	\begin{algorithmic}[1]
		\Require A matrix $Z \in \bbR^{\overline P \times Q}$ with $\overline P \subset [n]$ a set of experts and $Q \subset [d]$ a set of blocks, a direction $\nw \in \bbR_+^{Q}$, $w \neq 0$ and a pivot index $\gamma$
		\Ensure Two couples of subsets $(L, U)$ and $(\overline L, \overline U)$ of $\overline P$
		\Statex
		\For{$i \in \overline P$}
		\State Compute the statistic $ \psi (i, w) = \proscal<Z_{i,\cdot}, \frac{w}{\|w\|_2}>$
		\EndFor
		\State Sort the statistics $\psi(i,w)$ : $\psi(i_1,w) \leq \dots \leq \psi(i_{|\overline P|}, w)$
		\State $\nU = \{ i \in \overline P ~:~ \psi (i,w) > \psi (i_{\gamma}, w) + \beta_{\tris}\sqrt{ \Log{2\frac{|\overline P|}{\delta}}}\}$
		\State $\overline \nU = \{ i \in \overline P ~:~ \psi (i,w) > \psi (i_{\gamma}, w) + \overline \beta_{\tris}\sqrt{ \Log{2\frac{|\overline P|}{\delta}}}\}$
		\State $\nL= \{ i \in \overline P ~:~ \psi (i,w) < \psi (i_{\gamma}, w) - \beta_{\tris}\sqrt{ \Log{\frac{2|\overline P|}{\delta}}}\}$
		\State $\overline \nL= \{ i \in \overline P ~:~ \psi (i,w) < \psi (i_{\gamma}, w) - \overline \beta_{\tris}\sqrt{ \Log{\frac{2|\overline P|}{\delta}}}\}$
		\State \Return{$(\nL,\nU)$, $(\overline \nL,\overline \nU)$}
	\end{algorithmic}
\end{algorithm}
When the vector $w$ is equal to $\1_{Q}$, we simply write $\algoPivot(Z,\gamma)$ instead of $\algoPivot(Z,\1_{Q},\gamma)$ for the sake of simplicity. 
\begin{equation}\label{eq:pivot_overloading}
	\algoPivot(Z,\gamma) = \algoPivot(Z, w=\1_{Q}, \gamma)\enspace .
\end{equation}
In fact, $\algoPivot(Z,\gamma)$ simply amounts to comparing the row sums of $Z$ for each of the experts in $P$.

\medskip

In the next two pages, we redefine in more detail the Double Trisection algorithm of Section~\ref{sec:algo_sketch}. First, $\algoDoubleTrisectionPCA$ relies on a PCA-type argument to find a suitable direction $\hat w^+$ and then provides two trisections of the subset $\overline P$ of experts using the $\algoPivot$ sub-routine.

\begin{algorithm}[H]
	\caption{$\algoDoubleTrisectionPCA(\cZ,\gamma)$\label{alg:double_trisection_pca}}
	\begin{algorithmic}[1]
		\Require $4$ reduced samples $\cZ = (Z^{(1)}, Z^{(2)}, Z^{(3)}, Z^{(4)})$ where $Z^{(1)}$, $Z^{(2)}$, $Z^{(3)}$ $\in \bbR^{\widetilde P \times Q}$ and $Z^{(4)} \in \bbR^{\overline P\times Q}$ with some $\widetilde P \subset \overline P$, and a pivot index $\gamma$ 
		\Ensure Four subsets $(\nL_{\PCA}, \nU_{\PCA})$ and $(\overline \nL_{\PCA}, \overline \nU_{\PCA})$	 of $\overline P$
		\Statex
		\State Compute the following vector with coefficients in $\widetilde P$: $$ \hat v = \argmax_{\| v \| \leq 1} \Big[ \|v^T(Z^{(1)} - \overline{Z}^{(1)})\|_2^2 - \frac{1}{2}\| v^T(Z^{(1)} - \overline{Z}^{(1)} - Z^{(2)} + \overline{Z}^{(2)})\|_2^2\Big] \in \bbR^{\widetilde P}\enspace $$
		\State $\hat z = v^T Z^{(3)} \in \bbR^{Q}$
		\State Define $\hat w^+$ by $(\hat w^+)_l= |\hat{z}_l|\1_{|\hat z_l|\geq 2\zeta\sqrt{2\log(2|Q|/\delta)}}$
		\State $(L_{\PCA}, U_{\PCA}), (\overline L_{\PCA}, \overline U_{\PCA})= \algoPivot(Z^{(4)}, \hat w^+,\gamma)$
		\State \Return{$( L_{\PCA}, U_{\PCA})$, $(\overline L_{\PCA}, \overline U_{\PCA})$}
	\end{algorithmic}
\end{algorithm}

Next, $\algoDoubleTrisectionLocal(\cZ,\gamma)$ builds two trisections of $\overline P$ based on the reduced samples. First, it builds these trisections by simply using the row sums on the data and then it improves them thanks to $\algoDoubleTrisectionPCA$.

\begin{algorithm}[H]
	\caption{$\algoDoubleTrisectionLocal(\cZ,\gamma)$\label{alg:double_trisection_local}}
	\begin{algorithmic}[1]
		\Require $5$ reduced samples $\cZ = (Z^{(1)}, Z^{(2)}, Z^{(3)}, Z^{(4)}, Z^{(5)})$ in $\bbR^{\overline P \times Q}$, a pivot index $\gamma$ and a threshold $\beta_{\tris}$
		\Ensure Two couples of subsets $(L, U)$ and $(\overline L, \overline U)$ of $\overline P$
		\Statex
		\State $(L_{\cp}, U_{\cp}), (\overline L_{\cp}, \overline U_{\cp}) = \algoPivot(Z^{(1)},\1_Q,\gamma)$
		\State Set  $\widetilde P = \overline P \setminus (\overline L_{\cp}\cup \overline U_{\cp})$
		\State Set $\cZ' = (Z^{(2)}(\widetilde P), Z^{(3)}(\widetilde P), Z^{(4)}(\widetilde P), Z^{(5)}(\overline P))$ be the sequence of reduced samples where the three first samples are restricted to $\widetilde P$.
		\State $(L_{\PCA}, U_{\PCA}), (\overline L_{\PCA}, \overline U_{\PCA})= \algoDoubleTrisectionPCA(\cZ', \gamma)$ \State Set $ L =  L_{\cp}\cup L_{\PCA} \spaceAnd \overline L =  \overline L_{\cp}\cup \overline L_{\PCA} \spaceAnd U =  U_{\cp}\cup U_{\PCA} \spaceAnd \overline U =  \overline U_{\cp}\cup \overline U_{\PCA} $
		\State \Return{$(L, U), (\overline L, \overline U)$}
	\end{algorithmic}
\end{algorithm}

To finish defining $\algoDoubleTrisection$, we simply need to plug a dimension reduction procedure to select a subset of questions $Q\subset[d]$ and then to sum the data on these questions. 

The two following algorithms are mainly definitions. For some data $Y \in \bbR^{[n] \times [d]}$, a set of experts $\overline P$ and a set of blocks $Q \subset \cQ_r$, and a scale $r$, the $|\overline P|\times |Q|$ matrix $\algoEncodeMatrix(Y,\overline P,Q,r)$ is simply a reduced matrix where we consider the normalized row sums of $Y$ around the questions of $Q$ at scale $r$. 
\begin{algorithm}[H]
	\caption{$\algoEncodeMatrix(Y,\overline P,Q,r)$}
	\begin{algorithmic}[1]
		\Require A matrix $Y \in \bbR^{[n] \times [d]}$, a set of experts $\overline P$ and a set of blocks $Q \subset \cQ_r$, a scale $r$
		\Ensure A reduced matrix $Z \in \bbR^{\overline P \times Q}$
		\Statex
		\For{$i \in \overline P$ and $l \in Q$}
		\State Define $Z_{i, l} = \frac{1}{\sqrt{r}}\sum_{k\in [l, l + r)}Y_{i, k}$ \Comment{$Y_{i,k} = 1$ for $k \geq d+1$}
		\EndFor
		\State \Return{$Z \in \bbR^{\overline P \times Q}$} \Comment{the restriction of $Z$ to $\overline P$ and $Q$}
	\end{algorithmic}
\end{algorithm}
Second, $\algoEncodeSet(D,r)$ transforms a subset $[d]$ of questions into a subset $Q\subset \cQ_r$ of blocks of questions at scale $r$. 

\begin{algorithm}[H]
	\caption{$\algoEncodeSet(D,r)$}
	\begin{algorithmic}[1]
		\Require A set of questions $D\subset[d]$, a scale $r \in \cR$
		\Ensure A set of blocks $Q \subset  \cQ_r$
		\Statex
		\Return{$Q = \{l \in \cQ_r ~:~ [l, l+r) \cap D \neq \emptyset\}$}
	\end{algorithmic}
\end{algorithm}

Then, we are in position to redefine this  version of \Cref{alg:double_trisection_sketch}. As in the original definition in Section~\ref{sec:algo_sketch}, there are two variations of this procedure depending on whether we are building the estimator $\hat{\pi}_{HT}$ or the estimator $\hat{\pi}_{WM}$ that uses the memory of the tree. Algorithm $\algoDoubleTrisection(\cY, \cT, \overline P, \gamma)$ takes some original data and then reduces the dimension of the problem to build two trisections of the set $\overline P$ of experts. 

\begin{algorithm}[H]
	\caption{$\algoDoubleTrisection(\cY, \cT, \overline P, \gamma)$}\label{alg:double_trisection}
	\begin{algorithmic}[1]
		\Require $6$ samples $\cY = (Y^{(1)}, \dots, Y^{(6)})$, a tree $\cT$, a set of expert $\overline P$ included in a leaf $G$ of $\cT$ at maximal depth and a pivot index $\gamma$
		\Ensure Two couples of subsets $(L, U)$, $(\overline L, \overline U)$ of $\overline P$
		\Statex
		\State Initialize $L, U, \overline L, \overline U = \emptyset$
		\For{$h \in \cH,r \in \cR$}
		\If{Not using the memory of the tree} 
		\State Set $\widehat Q := \widehat Q_{\cp}(h, r) = \algoDimensionReduction(Y^{(1)}, \overline P,h,r)$ - see \Cref{algo:dimension_reduction} or \Cref{eq:significant_blocks_estim}
		\ElsIf{Using the memory of the tree}
		\State Set $\widehat Q := \widehat Q_{WM}(h, r) = \algoDimensionReductionWM(Y^{(1)}, \cT,\overline P,h,r)$ - See \Cref{algo:dimension_reduction_wm} or \Cref{eq:significant_blocks_estim_wm}
		\EndIf 
		
		\State Consider the five samples $\cY' = (Y^{(2)},Y^{(3)}, Y^{(4)}, Y^{(5)}, Y^{(6)})$
		\State Consider the five reduced samples $\cZ = \algoEncodeMatrix(\cY', \overline P, \widehat Q_{WM}, r)$
		\State Compute $(L_{\loc}, U_{\loc}), (\overline L_{\loc}, \overline U_{\loc}) = \algoDoubleTrisectionLocal(\cZ, \gamma)$
		\State Update $L = L \cup L_{\loc} \spaceAnd \overline L = \overline L \cup \overline L_{\loc} \spaceAnd U = U \cup U_{\loc} \spaceAnd \overline U = \overline U \cup \overline U_{\loc}$
		\EndFor
		\State \textbf{return} $(L,U)$, $(\overline L,\overline U)$
	\end{algorithmic}
\end{algorithm}

Finally, we reproduce $\algoBlockSorting$ here that was originally defined in \Cref{alg:algoBlockSorting_sketch}. We recall that $\algoBlockSorting$ iteratively applies a logarithmic number of times the procedure $\algoDoubleTrisection$ to build two suitable trisections of a set $G$ of experts. Although implicit in this description, there are two different versions of the corresponding procedure whether we use the memory of the tree - estimator $\hat{\pi}_{WM}$ - or not - estimator $\hat{\pi}_{HT}$ in $\algoDoubleTrisection$. In the following,  $\tau_{\infty}=\lceil 4\cdot 10^7\log^7(\tfrac{nd}{\delta(\zeta_-)^2})\rceil$  stands for the number of iterations in $\algoBlockSorting$.

\begin{algorithm}[H]
	\caption{\label{alg:algoBlockSorting} $\algoBlockSorting(\cY, \cT, G)$}
	\begin{algorithmic}[1]
		\Require $6\tau_{\infty}$ samples $\cY = (Y^{(0)}, \dots, Y^{(6\tau_{\infty}-1)})$, the tree $\cT$, a leaf $G \subset [n]$ in $\cT$ at maximal depth
		\Ensure A partition of $G$ into three groups $(O,P, I)$
		\Statex
		\State Set $\gamma = \floor{|G|/2}$ and $O_0$, $I_0$, $\overline O_0$,  $\overline I_0$ $= \emptyset$
		\For{$\tau = 0, \dots, \tau_{\infty}-1$}
		\State Consider $6$ fresh samples $ \cY_\tau = (Y^{(6\tau)}, \dots, Y^{(6\tau + 5)})$ 
		\State set $\gamma = \floor{|G|/2} - |\overline O_{\tau}|$
		\State $(L_{\tau}, U_{\tau}), ~ (\overline L_{\tau},\overline U_{\tau}) = \algoDoubleTrisection(\cY_\tau,\cT, G\setminus(\overline O_{\tau} \cup \overline I_{\tau}), \gamma)$ as in Algorithm~\ref{alg:double_trisection}
		\State Update
		$O_{\tau+1} = O_{\tau}\cup L_{\tau}, ~~I_{\tau+1} = I_{\tau}\cup U_{\tau}, ~~\overline O_{\tau+1} = \overline O_{\tau}\cup \overline L_{\tau}, ~~\overline I_{\tau+1} = \overline I_{\tau}\cup \overline U_{\tau}$
		\EndFor
		\If{$O_{\tau_{\infty}}\cap I_{\tau_{\infty}}\neq \emptyset $} 
		\State Set $O_{\tau_{\infty}} := O_{\tau_{\infty}}\setminus I_{\tau_{\infty}}$ and $I_{\tau_{\infty}} := I_{\tau_{\infty}}\setminus O_{\tau_{\infty}}$
		\EndIf
		\State \Return{$(O_{\tau_{\infty}}, G\setminus(O_{\tau_{\infty}} \cup I_{\tau_{\infty}}), I_{\tau_{\infty}})$}
	\end{algorithmic}
\end{algorithm}
Under an event of high probability (to be later discussed), we have $O_{\tau_{\infty}}\cap I_{\tau_{\infty}}= \emptyset$. The correction at the end of the algorithm simply forces the algorithm to return a partition of $G$.

\subsection{Hierarchical sorting Trees and $\algoTreeSorting$ Algorithm}

In this subsection, we formally describe how we build and navigate into a hierarchical tree. In the following, a node $G \in \Nodes$ is a labelled subset of $[n]$. Its  label belongs to $\{\0, \bp, \1 \}$. For a node $G$, we write $\algoType(G)$ for the label (also called type) of $G$.

\begin{definition}(Hierarchical sorting Trees)\label{eq:def:HST0}
	A hierarchical sorting tree $\cT$ is a rooted tree that satisfies the three following properties: 
	\begin{itemize}
		\item The root $G$ of $\cT$ corresponds to the set $[n]$ and its label is $\0$. 
        \item Any node $G$ of type $\bp$ is a leaf.
        \item Any  node $G$ of type in $\{\0,\1\}$ is either a  leaf or has three children $(\nO,P,\nI)$ with type $\0$, $\bp$, $\1$ respectively. Besides, $(\nO,P,\nI)$  correspond to a partition of $G$. 
	\end{itemize}
\end{definition}
We write $\cT_0$ for the tree of depth $0$. The procedure $\algoTreeSorting$  iteratively builds a hierarchical sorting tree. Hence, we need to define the operation of adding children to a leaf in a tree $\cT$. For a specific leaf $G$ of type $\0$ or $\1$, we consider three labelled subsets $\nO$, $P$, $\nI$ of type $\0$, $\bp$, $\1$, respectively. Besides,  those subsets satisfy the third condition in Definition~\ref{eq:def:HST0}. Then,   $\cT' = \algoAddChild(\cT,G, (\nO,P,\nI))$ is the supertree of $\cT$ where we have added the nodes $(\nO, P, \nI)$ as children of $G$.
Finally, we observe that for any $t>0$, all the nodes at depth $t$ of a hierarchical sorting tree $\cT$  are disjoint.

In fact, we shall prove in \Cref{prop:algo_tree_sorting} and in \Cref{cor:tree_sorting} that, with high probability, the final tree $\cT_{t_{\infty}}$ turns out to be a valid hierarchical sorting tree as defined below.

\begin{definition}(Valid hierarchical sorting Tree)\label{eq:def:HST}
	A hierarchical sorting tree $\cT$ is valid if non-terminal nodes $G$ of $\cT$ satisfy  the two following additional properties: if we denote $(\nO,P,\nI)$ their children of type $\0$, $\bp$, $\1$ respectively, then   
	\begin{itemize}
		\item  All the experts in $\nO$ are below those of $\nI$. In other words, for any $i\in \nO$ and any $j\in \nI$, we have $\pi^*(i)< \pi^*(j)$. 
 		\item $|\nO|< |G|$ and $|\nI|< |G|$.
	\end{itemize}
\end{definition}
The second property ($|\nO|< |G|$ and $|\nI|< |G|$) forces the tree to be finite.

For a node $G$ in a such valid hierarchical sorting tree $\cT$, $\algoDepth(\cT, G)$ stands for the depth of $G$ in $\cT$. 
In light of this definition of valid hierarchical sorting trees, a labelled subset $G$ cannot appear twice in a tree $\cT$, so that $\algoDepth(\cT, G)$ is well-defined.

 We are now equipped to provide a more formal definition of $\algoTreeSorting$, although the procedure is in fact the same as the one described in Algorithm~\ref{alg:procedure_sketch}. Let $t_{\infty}= \lceil \log(n)/\log(2)\rceil$.

\begin{algorithm}[H]
	\caption{\label{alg:algoTreeSorting}$\algoTreeSorting(\cY)$}
	\begin{algorithmic}[1]
		\Require $6\tau_{\infty}t_{\infty}$ samples $\cY = (Y^{(0)}, \dots, Y^{(6\tau_{\infty}t_{\infty}-1)})$
		\Ensure A final tree $\cT$
		\State $\cT = \cT_0$ \Comment{The root is at depth $0$}
		\For{$t = 0, \dots, t_{\infty}-1$}
		\State Consider $6\tau_{\infty}$ fresh samples $\cY = (Y^{(6t\tau_{\infty})}, \dots, Y^{(6(t+1)\tau_{\infty}-1)})$
		\For{$G \in \cL^{(\0, \1)}(\cT)$} \Comment{See \Cref{eq:leaves_0_1} for the definition of $\cL^{(\0, \1)}$}
		\State $(\nO_G, P_G, \nI_G) = \algoBlockSorting(\cY, \cT, G)$
		\State Set $\algoType(\nO_G) = \0 \spaceAnd \algoType(P_G) = \bp \spaceAnd \algoType(\nI_G) = \1$
		\EndFor
		\For{$G \in \cL^{(\0, \1)}(\cT)$}
		\State $\algoAddChild(\cT, G, (\nO_G, P_G, \nI_G))$
		\EndFor
		\EndFor
		\State \Return{$\cT$}
	\end{algorithmic}
\end{algorithm}
As explained in Section~\ref{sec:algo_sketch}, the final estimators $\hat{\pi}_{HT}$ or $\hat{\pi}_{WM}$ are computed from their corresponding hierarchical sorting tree $\cT$.

In order to define the $\algoDimensionReductionWM$ algorithm in the next subsection, we need to introduce a few more notation. First, we define 
\begin{equation}\label{eq:leaves_0_1}
	\cL^{(\0, \1)}(\cT) = \{ G \in \Leaves(\cT) ~:~ \algoType(G) \in \{ \0, \1 \} \}\ ,
\end{equation}
as the collection of leaves of $\cT$ that are either of type $\0$ or of type $\1$. In the algorithm $\algoTreeSorting$, these are the leaves to be partitionned. In particular at step $t$ of $\algoTreeSorting$, $\cL^{(\0, \1)}(\cT)$ is only made of leaves at depth $t$.

For a subset $P\subset [n]$,   $\algoLeaf(\cT, P)$ is defined as the leaf $G \in \Leaves(\cT)$ containing $P$ (if it exists). Finally, the groups $G\in \cL^{(\0, \1)}(\cT)$ inherit from a natural order provided that $\cT$ is a valid hierarchical sorting tree.
We can enumerate the  groups  $G_1$, $G_2$,\ldots, $G_{|\cL^{(\0,\1)}(\cT)|}$ in such a way that all the experts in $G_{s}$ are below those of $G_{s'}$ for $s< s'$. To ease the presentation, we also introduce, for any positive integer $s$ the groups $G_{|\cL^{(\0,\1)}(\cT)|+s}=\{n+s\}$. The corresponding data and signal for the $n+s$-th expert satisfies  $Y_{n+s,j}=1=M_{n+s,j}$ almost-surely for any $j\in [d]$. Also, for any positive integer $s$ we introduce the groups $G_{1-s}=\{1-s\}$. The corresponding data and signal for this synthetic expert satisfy $Y_{1-s,j}=0=M_{1-s,j}=0$ almost-surely for any $j\in [d]$. 

Then, for a specific leaf $G_s\in \cL^{(\0, \1)}(\cT)$, $\algoOrder(\cT, G)$ stands for the collection $(G^{(a)})$, $a\in \mathbb{Z}$ of leaves  where $G^{(a)}= G_{a+s}$. In other words,  we have $G^{(0)}=G_s$ and $G^{(1)}$ is the following group, and so on.

\subsection{Dimension Reduction Algorithms}\label{sec:algorithms_WM}

To finish the description of the two procedures, we fully describe the two dimension reduction algorithms both for the oblivious estimator $\hat{\pi}_{HT}$ and for the estimator $\hat{\pi}_{WM}$ with memory. These procedures were already introduced in Section~\ref{sec:algo_sketch}. First,  $\algoDimensionReduction(Y, \overline P, h, r)$ considers the columns-wise mean of the restriction of $Y$ to the group $\overline P$ and detects high-variation regions of this vector. 
\begin{algorithm}[H]
	\caption{$\algoDimensionReduction(Y, \overline P, h, r)$ \label{algo:dimension_reduction}}
	\begin{algorithmic}[1]
		\Require A sample $Y \in \bbR^{\overline P \times [d]}$, a set of experts $\overline P$, $h \in \cH$ and $r$ in $\cR$
		\Ensure An encoded set including the high-variation regions $\widehat Q_{\cp}:= \widehat Q_{\cp} (Y, \overline P, h, r) \subset \cQ_r$
		\Statex
		\State $\overline y(\overline P) = \frac{1}{|\overline P|}\sum_{i \in \overline P} Y_{i, \cdot}$
		\State $\tilde r = 8\left[\left\lceil \left(\frac{32\zeta^2 }{|\overline P|h^2}\log(\tfrac{2d}{\delta})\right)\right\rceil \lor r\right]$
		\State Initialize $\widehat D_{\cp} = \emptyset$
		\For{$k \in [d]$}
		\State Compute \begin{equation}\label{eq:def_cusum}\widehat \bC_k(\overline y(\overline P)) = \frac{1}{\tilde r}\left(\sum_{k' = k}^{k+ \tilde r - 1}\overline y_{k'}(\overline P) -  \sum_{k' = k- \tilde r}^{k-1}\overline y_{k'}(\overline P)\right) \enspace ;\end{equation}

		\EndFor
	\State $\widehat D_{\cp} = \{k \in [d] ~:~ \widehat \bC_k(\overline y(\overline P))\geq h/4\}$
	\State $\widehat Q_{\cp} =\algoEncodeSet(\widehat D_{\cp}, r)$
	\State \Return{$\widehat Q_{\cp}$}
	\end{algorithmic}
\end{algorithm}

For the more involved dimension reduction procedure with memory $\algoDimensionReductionWM$, we compute the CUSUM statistic in larger groups $\cV\supset \overline P $ to reduce its variance and we also require that the estimated "width" of the group of experts is high enough. More precisely, given three sets of expert $\cV$, $\cV^+$ and $\cV^-$ and a sample $Y$, we consider the two following statistics, for any $k = 1, \dots, d$ and $r' \in \cR$:
\begin{align}
	\widehat \bDelta^{(\ext)}_{k, r'}(\cV^+, \cV^-) & = \frac{1}{2r'}\sum_{k' = k - r'}^{k + r' - 1} \overline y_{k'}(\cV^+) - \overline y_{k'}(\cV^-)   \ ; \quad \quad \widehat \bC^{(\ext)}_{k, r'}(\cV)               = \frac{1}{r'}\left(\sum_{k' = k}^{k + r' - 1} \overline y_{k'}(\cV) - \sum_{k' = k - r'}^{k - 1}\overline y_{k'}(\cV)\right)\label{eq:stat_cusum_wm}\enspace . 
\end{align}
Here, $\widehat \bDelta^{(\ext)}_{k, r'}(\cV^+, \cV^-)$ computes the width - i.e.~the difference - between the mean of experts in $\cV^+$ and the mean of experts in $\cV^{-}$. Since $\cV^+$ and $\cV^-$ are built in the algorithm below in such a way that experts in $\overline P$ are below those of $\cV^+$ and above those of $\cV^{-}$, $\widehat \bDelta^{(\ext)}_{k, r'}(\cV^+, \cV^-)$ provides an upper bound of the width between the best expert in $\overline P$ and the worst expert in $\overline P$.

The algorithm $\algoDimensionReductionWM$ described below builds a collection of sets $\cV^+$, $\cV^-$, and $\cV$ and detects questions such that both the CUSUM $\bC^{(\ext)}_{k, r'}(\cV)$ and the width $\widehat \bDelta^{(\ext)}_{k, r'}(\cV^+, \cV^-)$ are large enough. Further explanations are postponed to the analysis of the algorithm in Section~\ref{sec:proof_trisection_wei}.	
Below, we write $\lceil x\rceil^{dya}$ for $2^{\lceil \log_2(x)\rceil}$. 

\begin{algorithm}[H]
	\caption{$\algoDimensionReductionWM(Y,\cT,\overline P,h, r)$\label{algo:dimension_reduction_wm}}
	\begin{algorithmic}[1]
		\Require A sample $Y \in \bbR^{n \times d}$, a tree $\cT$, a set $\overline P$ included in a leaf $G$ of $\cT$ of type $\0$ or $\1$, $h \in \cH$ and $r \in \cR$
		\Ensure A set of blocks $\widehat Q_{WM} := \widehat Q_{WM}(Y, \cT, \overline P, h, r) \subset \cQ_r$
		\Statex
		\State $r_0 = 2^9\log(4d|\cR|/\delta)\frac{\zeta^2}{|\overline P|h^2} \spaceAnd \tilde r = 4(\lceil r_0\rceil^{dya} \lor r)$
		\State $(G^{(a)})_{a \in \bbZ} = \algoOrder(\cT, G)$
		\For{$r_{\cp} \in [4r, 2\tilde r]\cap \cR$}
		\State \label{line:awm}Set $ a^{+}_{WM} = \min \{ a ~:~ |G^{(1)}| + \dots + |G^{(a)}| \geq 2^{11}\log(4d|\cR|/\delta)\tfrac{\zeta^2}{r_{\cp}h^2} \} $
		\State Set  $a^{-}_{WM} = \min \{ a ~:~ |G^{(-1)}| + \dots + |G^{(-a)}| \geq 2^{11}\log(4d|\cR|/\delta)\tfrac{\zeta^2}{r_{\cp}h^2} \} $
		\State \label{line:Gres}Set
		\begin{equation}
			\cV^{+}_{r_{\cp}} := \cV^{+}_{r_{\cp}}(\cT, G, h) = \bigcup_{a =  1}^{a^{+}_{WM}}G^{(a)}
			\spaceAnd \cV^{-}_{r_{\cp}} := \cV^{-}_{r_{\cp}}(\cT, G, h) = \bigcup_{a \in - a^{-}_{WM}}^{-1}G^{(a)}
		\end{equation}
		\If{$r_{\cp} > \tilde r$}
		\State Set $\cV_{r_{\cp}} := \cV_{r_{\cp}}(\cT, \overline P, h) = \overline P$
		\ElsIf{$r_{\cp} \leq \tilde r$}
		\State Set $\cV_{r_{\cp}} := \cV_{r_{\cp}}(\cT, \overline P, h) = \cV^{-}_{r_{\cp}} \cup \cV^{+}_{r_{\cp}}$
		\EndIf
		\EndFor
		\State $\widehat Q_{WM} = \emptyset$
		\For{$r_{\cp} \in [4r, \tilde r] \cap \cR$}
		\State $\widehat D_{WM} = \emptyset$
		\For{$k = 1, \dots, d$}
		\State Compute $\widehat \bDelta^{(\ext)}_{k, r_{\cp}} := \widehat \bDelta^{(\ext)}_{k, r_{\cp}}(\cV^+_{ r_{\cp}}, \cV^-_{ r_{\cp}})$
		\State Compute $\widehat \bC^{(\ext)}_{k, 2r_{\cp}} := \widehat \bC^{(\ext)}_{k,2r_{\cp}} (\cV_{2r_{\cp}})$
		\EndFor
		\State Update $\widehat D_{WM} = \{k \in [d] ~:~ \widehat \bDelta^{(\ext)}_{k, r_{\cp}} \geq h/16$ and $\widehat \bC^{(\ext)}_{k, 2r_{\cp}} \geq h/16 \}$
		\State Update $\widehat Q_{WM} = \widehat Q_{WM} \cup \algoEncodeSet(\widehat D_{WM}, r)$
		\EndFor
		\State \Return{$\widehat Q_{WM}$}\Comment{The same set is defined in \Cref{eq:significant_blocks_estim_wm}}
	\end{algorithmic}
\end{algorithm}

\section{Overview and organization of the proofs of Theorems~\ref{th:first_estimator} and~\ref{th:second_estimator_WM}}\label{sec:deterministic_analysis}

In this section, we divide the analysis of the procedures into several properties that will be proved to hold with high probability in the next sections.

\subsection{Definitions}

Since we build our estimator using a hierarchical tree, we need to quantify the error that we suffer at each depth of the tree. For $i \in [n]$, we write $M_i = M_{i, \cdot}$ for the expert $i$. By definition of $\pi^*$, we recall that 
\[
	M_{\pi^{*-1}(1)} \leq M_{\pi^{*-1}(2)}\leq \ldots\leq M_{\pi^{*-1}(n)}\enspace .
\]
For a given group of experts $G$, we write $\pi^*_{\{G\}}$ for the oracle ordering in $[1, |G|]$ of the group $G$ according to $\pi^*$, that is for all $i,j \in G$, $\pi^*_{\{G\}}(i)$ and $\pi^*_{\{G\}}(j)$ belong to  $[1,|G|]$ and
$$\pi^*_{\{G\}} (i) < \pi^*_{\{G\}} (j) \quad \text{iff} \quad \pi^*(i) < \pi^*(j) \enspace .$$

We say that a sequence of sets $\cG = (G_1, \dots, G_{\alpha})$ is an \emph{ordered partition} of a set $S$ if $\{G_1, \dots, G_{\alpha} \}$ is a partition of $S$. 
For a given ordered partition $\{G_1, \dots, G_{\alpha} \}$ and $a \in [1, \alpha]$ and any  $i \in G_a$ we write
\begin{equation}\label{eq:pi_plus}
	\pi^-_{\cG}(G_a) = \pi^-_{\cG}(i) :=  \sum_{a' < a} |G_{a'}| \quad \text{ and } \quad \pi^+_{\cG}(G_a) = \pi^+_{\cG}(i) := \sum_{a' \leq a} |G_{a'}| \enspace . 
\end{equation}
If we are to build a permutation $\pi$ which is consistent with this ordered partition, then one easily checks that $\pi(i)\in [\pi^-_{\cG}(i)+1, \pi^+_{\cG}(i)]$.  For simplicity, we write $G(i)$ for the group $G_a$ such that $i \in G_a$.
For a given ordered partition $\cG = (G_1, \dots, G_{\alpha})$, we define the oracle permutation associated to $\cG$ by
\begin{equation} \label{eq:pi_partition}
	\pi^*_{\cG}(i) = \pi^-_{\cG}(i) + \pi^*_{\{G(i)\}}(i) \enspace .
\end{equation}
For example, $\pi^*_{\{[n]\}} = \pi^*$ is simply the true permutation. By definition, we have $\pi^*_{\cG}(G(i)) = [\pi^-_{\cG}(i)+1, \pi^+_{\cG}(i)]$. Given an ordered partition, $\pi^*_{\cG}$ is the best permutation we could hope for after any statistical treatment.

Given an ordered partition $\cG = (G_1, \dots, G_{\alpha})$, we define the random estimation of $\pi^*$ given $\cG$ as $\hat \pi_{\cG} (i)$ which is uniformly distributed in $[\pi^-_{\cG}(i) +1 , \pi^+_{\cG}(i)]$:
$$\hat \pi_{\cG}(i) \in [\pi^-_{\cG}(i) +1 , \pi^+_{\cG}(i)]\enspace .$$
Note that $\hat \pi_{\cG}$ is not necessarily bijective.

\subsection{Deterministic Analysis}\label{subsec:deterministic_analysis}

In this subsection, we analyze $\algoTreeSorting$ (\Cref{alg:algoTreeSorting}) and we characterize the loss of the estimator $\hat{\pi}$ in terms of that of the trisections that are computed inside the subroutine $\algoBlockSorting(\cY, \cT, G)$. This algorithm takes a subset  $G$ of experts and computes two trisections of $G$. The first one 
\[
(O,P,I)= (O_{\tau_{\infty}}, G\setminus(O_{\tau_{\infty}} \cup I_{\tau_{\infty}}), I_{\tau_{\infty}})
\]
is returned by the algorithm.  The second one 
\[
(\overline \nO, \overline P, \overline \nI)= (\overline{O}_{\tau_{\infty}}, G\setminus(\overline{O}_{\tau_{\infty}} \cup \overline{I}_{\tau_{\infty}}), \overline{I}_{\tau_{\infty}})
\]
is important for our analysis. From the definitions of the different procedures, one readily checks that $\overline\nO\subset O$ and $\overline\nI\subset I$. In fact, we shall prove later that, with high probability, the subsets $(\nO,P,\nI)$ and $(\overline \nO, \overline P, \overline \nI)$ satisfy the following stronger property.

\begin{property}\label{property:block-sorting} \
	\begin{enumerate}
		\item $\{ \nO, P, \nI\}$ and $\{ \overline \nO, \overline P, \overline \nI\}$ are partitions of the leaf $G$ with $\overline \nO \subset \nO$, $\overline \nI \subset \nI$, and $P \subset \overline P$, \label{item:point1_assum_block_sorting} \label{item:point2_assum_block_sorting}
		\item For $\omega = \pi^{*-1}_{\{\nO, P, \nI\}}\pi^{*}_{\{G\}}$, we have  $\omega(i) = i$ for any $i \in \overline \nO \cup \overline \nI$.\label{item:point3_assum_block_sorting}
		\item For any $i \in \nO$ and $j \in \nI$, we have $\pi^*(i) < \pi^*(j)$. \label{item:point4_assum_block_sorting}
  		\item We have $|\nO|\leq |G|/2$ and $|\nI|\leq |G|/2$.   \label{item:point5_assum_block_sorting}
	\end{enumerate}
\end{property}
The last claim states that all experts in $\nO$ are below all experts of $\nI$. The second claim can be understood as the fact that, if an expert $i$ belongs to $\overline \nO$, then all experts below $i$ belong to $\nO$.

\medskip 

Let $\cY = (Y^{(0)}, \dots, Y^{(6\tau_{\infty}-1)})$ be a sequence of $6\tau_{\infty}$ matrices in $\bbR^{n \times d}$. We say that $\algoBlockSorting$ satisfies \Cref{property:block-sorting} on $(\cY,\cT,G)$ if the two partitions $(\nO,P,\nI)$ and $(\overline \nO, \overline P, \overline \nI)$ worked out in \Cref{alg:algoTreeSorting} satisfy \Cref{property:block-sorting}. We recall that by definition $\overline \nO \subset \nO$,  $\overline \nI \subset \nI$, $P \subset \overline P$ so that $\overline P$ corresponds to the collection of experts that are either not sorted by $\algoTreeSorting$ or are sorted with a small confidence.

For each $t = 0, \dots, t_{\infty}$, we write $\cT_t$ for the hierarchical sorting tree at the beginning of step $t$ of $\algoTreeSorting$. Besides, we write $\cG_{t}$ for the corresponding ordered partition obtained by taking the leaves of the tree $\cT_{t}$ in increasing order in the ternary base $\{\0, \bp, \1\}$.
We define the tree $\overline \cT_{t_{\infty}}$ as the tree $\cT_{t_{\infty}}$ where we replaced all the leaves $P$ of type $\bp$ - at any depth - by $\overline P$, where $(\overline \nO,\overline P, \overline \nI)$ has been worked out by $\algoTreeSorting$ at the same time as $(O,P,I)$.

We also define
\begin{eqnarray}\nonumber
	\cL_{t}(\cT_{t_{\infty}})& =& \{ P \in \cG_{t_{\infty}} ~:~ P \text{ is a nonempty leaf at depth } t \text{ of } \cT_{t_{\infty}} \} \enspace ;\\
	\label{eq:leaf_of_tree_bar}
	\overline \cL_{t}(\overline \cT_{t_{\infty}}) &= &\{ \overline P \in \cG_{t_{\infty}} ~:~ \overline P \text{ is a nonempty leaf at depth } t \text{ of } \overline \cT_{t_{\infty}} \} \enspace .
\end{eqnarray}
For simplicity, we sometimes write $\cL_t$ for $ \cL_{t}(\cT_{t_{\infty}})$ and $\overline \cL_{t}$ for $\overline \cL_{t}(\overline \cT_{t_{\infty}})$. $\cL_t$ stands for  the collection of experts that have not been sorted at the $t$-th iteration $\algoTreeSorting$. The sets in the collection $\overline \cL_{t}$ are strictly larger and correspond to the collections of experts in $\overline{P}$ that are either not sorted by $\algoTreeSorting$ or are sorted with less  confidence.
Let $M(\overline P)$ be defined as the restriction of $M$ to the experts in $\overline P$, and $\overline{M}(\overline P)$ the $|\overline P|\times d$ matrix with constant columns which correspond to the mean row of $M(\overline P)$. 
	The following proposition characterizes the loss of the final estimator estimator $\hat{\pi}_{\cG_{t_{\infty}}}$ which is obtained from a hierarchical sorting tree in terms of the variance of the experts $M$ within the groups $\overline{P}$ in $\overline \cL_{t}(\overline \cT_{t_{\infty}})$.

\begin{proposition}[Deterministic Analysis of $\algoTreeSorting$]\label{prop:algo_tree_sorting}
	Assume that at each step of $\algoTreeSorting$, the routine $\algoBlockSorting$ applied to the data satisfies \Cref{property:block-sorting}. Then, the error of $\hat{\pi}= \hat \pi^{-1}_{\cG_{t_{\infty}}}$ is controlled as follows
	\beq\label{eq:loss_general}
	\| M_{\hat \pi^{-1}} - M_{\pi^{*-1}} \|_F^2 \leq 10t_{\infty}\sum_{t=1}^{t_{\infty}}\sum_{\overline P \in \overline \cL_{t}}\|M(\overline P) - \overline M(\overline P)\|_F^2 \enspace .
	\eeq
Besides, the hierarchical tree $\cT_{t_{\infty}}$ is valid (as in \Cref{eq:def:HST}) and all its non-empty leaves are of type $\bp$. 
\end{proposition}
Up to a normalization, $\|M(\overline P) - \overline M(\overline P)\|_F^2$ corresponds to the variance of $M$ within the group $\overline P$. The bound~\eqref{eq:loss_general} expresses that the loss of a hierarchical sorting tree is controlled by the variance of the set $\overline{P}$ that are not sorted with confidence at each step of the algorithm. Also, we recall that $t_{\infty}= \lceil \log(n)/\log(2)\rceil$.
This proposition only relies on \Cref{property:block-sorting} and on the construction of the tree. Hence, it applies both to the estimators $\hat{\pi}_{HT}$ and $\hat{\pi}_{WM}$.

\medskip

The sets $(O,\overline{O},I,\overline{I})$ built in $\algoBlockSorting$ arise as unions of set $(L,U)$ and $(\overline{L},\overline{U})$ that are computed by $\algoDoubleTrisection$ for a set $\overline P$ and a pivot $\gamma\in [1, |\overline P|]$. For this reason, we now state a desired property of the result of the algorithm that will enforce \Cref{property:block-sorting}. 
\begin{property}[Property on $(L,U)$ and $(\overline{L},\overline{U})$]\label{property:trisection}
	For
	$\overline P' = \overline P \setminus (\overline \nL \cup \overline \nU)$ and $P' = \overline P \setminus (\nL \cup \nU$), we have
	\begin{enumerate}
		\item  $\overline \nL \subset \nL$, and  $\overline \nU \subset \nU$,\label{item:point2_assumDoubleTrisection}
		\item if $\omega = \pi^{*-1}_{\{\nL, P', \nU\}}\pi^{*}_{\{\overline P\}}$ then for any $i \in \overline \nL \cup \overline \nU$ it holds that $\omega(i) = i$, \label{item:point3_assumDoubleTrisection}
		\item For any $i \in \nL$ and $j \in \nU$ we have $\pi^*_{\{\overline P\}}(i) < \gamma < \pi^*_{\{\overline P\}}(j)$. \label{item:point4_assumDoubleTrisection}
	\end{enumerate}
\end{property}
We say that $\algoDoubleTrisection$ with $(\cY, \cT, \overline P, \gamma)$ satisfies~\Cref{property:trisection} if the corresponding subsets $(L,U)$ and $(\overline{L},\overline{U})$ satisfy~\Cref{property:trisection}.

\begin{proposition}[Deterministic Analysis of  $\algoBlockSorting$]\label{prop:algo_block_sorting}
	$\algoBlockSorting$ satisfies \Cref{property:block-sorting} on $(\cY,\cT,G)$ if, at each step of \Cref{alg:algoBlockSorting}, each call  of $\algoDoubleTrisection$  satisfies \Cref{property:trisection}.
\end{proposition}

In light of Propositions~\ref{prop:algo_tree_sorting} and \ref{prop:algo_block_sorting}, it suffices to show that, with high probability, all applications of $\algoDoubleTrisection$ in the construction of the hierarchical sorting tree satisfy~\Cref{property:trisection}, and then to control the sum of within-group variances in~\eqref{eq:loss_general}.

\subsection{High probability Control of Property~\ref{property:trisection}}

We write in this part of the proof (this subsection), for simplicity,  $\cY = (Y^{(1)}, \dots, Y^{(6)})$ for $6$ independent matrices that are identically distributed as $Y=M+E$ in~\eqref{eq:model_0}, where we recall that  the entries of $E$ are centered, independent and $\zeta$-subgaussian.

Fix a hierarchical sorting tree $\cT$ (recall~\Cref{eq:def:HST0}), a leaf $G$ of $\cT$, a set $\overline P \subset G$, a pivot $\gamma \in \{1, \dots, |\overline{P}| \}$. Let $\cP_2:= \cP_2(\cT, \overline P, \gamma, \beta_{\tris}, \overline \beta_{\tris})$ be the event holding true if $\algoDoubleTrisection$ satisfies \Cref{property:trisection} on $\cY$ for $(\cT, \overline P, \gamma, \beta_{\tris}, \overline \beta_{\tris})$. 
The following proposition states that
$\cP_2$ holds with uniformly high probability.
\begin{proposition}\label{cor:trisection}
	For any $\cT$, any leaf $G$, any $\overline P\subset G$, any pivot $\gamma \in [|\overline P|]$, 
	 we have $\Prob( \cP_2 ) \geq 1 - 3|\cH||\cR|\delta$.
\end{proposition}
This result is valid for both versions of $\algoDoubleTrisection$ where we use the memory of the tree (estimator $\hat{\pi}_{WM}$) or not (estimator $\hat{\pi}_{HT}$). Recall that in $\algoBlockSorting$ there are at most $\tau_{\infty}$ calls of $\algoDoubleTrisection$.  Since the construction of the hierarchical tree requires at most $2^{t_{\infty}+1}$ applications of $\algoBlockSorting$, we arrive at the following straightforward corollary of Propositions~\ref{prop:algo_block_sorting}, \ref{prop:algo_tree_sorting} and \ref{cor:trisection}.

\begin{corollary}\label{cor:tree_sorting}
	There exists an event $\xi$ of probability higher than  $ 1 - 2^{t_{\infty} + 1}3\tau_{\infty}|\cH||\cR|\delta$ such that all results of $\algoBlockSorting$ within $\algoTreeSorting$ satisfy Property~\ref{property:block-sorting}. In particular, the tree $\cT_{t_{\infty}}$ is a valid hierarchical sorting tree (as in \Cref{eq:def:HST}) whose non-empty leaves are all of type $\bp$.
	Besides, on this event we also have 
	\begin{equation}\label{eq:majoration_intermediary_groups}
		\| M_{\hat \pi^{-1}_{\cG_{t_{\infty}}}} - M_{\pi^{*-1}} \|_F^2 \leq 10t_{\infty}\sum_{t=1}^{t_{\infty}}\sum_{\overline P \in \overline \cL_{t}}\|M(\overline P) - \overline M(\overline P)\|_F^2 \enspace .
	\end{equation}
\end{corollary}
Again, this results applies to both variants of our procedure - with or without memory.

\subsection{Control of the Loss function}

In contrast to the previous subsection, we now need to specify the dimension reduction scheme $\algoDimensionReduction$ (which corresponds to $\hat{\pi}_{HT}$) or $\algoDimensionReductionWM$ (which corresponds to $\hat{\pi}_{WM}$) inside $\algoDoubleTrisection$ as the convergence rates depend on these quantities.

First we state the results for the method without memory: $\hat{\pi}_{HT}$.
\begin{proposition}\label{prop:rate_of_trisection}
	Consider the oblivious hierarchical sorting tree estimator $\hat{\pi}_{HT}$. On the intersection of event $\xi$ (defined in \Cref{cor:tree_sorting}) and an event of probability higher than {$1-5\cdot2^{t}\tau_{\infty}\delta$}, it holds  that
	$$\sum_{\overline P \in \overline \cL_{t}}\|M(\overline P) - \overline M(\overline P)\|^2 \, \lesssim \, \zeta^2\log^{8.5}\left(\frac{6nd}{\delta\zeta_-}\right) \left[\frac{n^{2/3}d^{1/3}}{\zeta^{2/3}}\wedge \frac{nd^{1/6}}{\zeta^{1/3}}\wedge n\sqrt{d} + n \right] \enspace .$$
\end{proposition}

Then we state the results for the method with memory: $\hat{\pi}_{WM}$.
\begin{proposition}\label{prop:rate_of_trisection_WEI}
	Consider the hierarchical sorting tree estimator $\hat{\pi}_{WM}$. 
	On the intersection of event $\xi$ (defined in \Cref{cor:tree_sorting}) and an event of probability higher {$1-5\cdot2^{t}\tau_{\infty}\delta$}, it holds  that
	$$\sum_{\overline P \in \overline \cL_{t}}\|M(\overline P) - \overline M(\overline P)\|^2 \lesssim  \zeta^2 \log^{9}\left(\frac{6nd}{\delta \zeta_-}\right)\left[\left(\frac{n^{3/4}d^{1/4}}{\zeta^{1/2}}\wedge \frac{nd^{1/6}}{\zeta^{1/3}}\wedge  n\sqrt{d}\wedge \frac{n^{2/3}\sqrt{d}}{\zeta^{1/3}} \right)+ n\right] \enspace .$$
\end{proposition}

Now, we are in position to easily conclude the proof of Theorems~\ref{th:first_estimator} and~\ref{th:second_estimator_WM}. 

\begin{proof}[Proof of \Cref{th:first_estimator}]
Let $\hat{\pi}_{HT} := \hat \pi_{\cG_{t_{\infty}}}$ denote the oblivious hierarchical sorting tree estimator. Combining \Cref{cor:tree_sorting} with \Cref{prop:rate_of_trisection} and a union bound over all $t = 0, \dots, t_{\infty} - 1$, it holds with probability higher than $1 - 8\cdot 2^{t_{\infty} + 1}\tau_{\infty}|\cH||\cR|\delta$ that
	\begin{eqnarray*}
		\| M_{\hat{\pi}_{HT}^{-1}} - M_{\pi^{*-1}} \|_F^2 &\lesssim& t^2_{\infty} \, \zeta^2\log^{8.5}\left(\frac{2nd}{\delta\zeta_-}\right) \left[\frac{n^{2/3}d^{1/3}}{\zeta^{2/3}}\wedge \frac{nd^{1/6}}{\zeta^{1/3}}\wedge n\sqrt{d} + n \right]\\
		&    \lesssim&  \, \zeta^2 \log^{10.5}\left(\frac{2nd}{\delta\zeta_-}\right) \left[\frac{n^{2/3}d^{1/3}}{\zeta^{2/3}}\land \frac{nd^{1/6}}{\zeta^{1/3}}\wedge n\sqrt{d} + n \right]\enspace .
	\end{eqnarray*}
\end{proof}

\begin{proof}[Proof of \Cref{th:second_estimator_WM}]
	Let $\hat \pi_{\WM} := \hat \pi_{\cG_{t_{\infty}}}$ denote the hierarchical sorting tree where we use the memory to reduce the dimension (Algorithm $\algoDimensionReductionWM$). Combining \Cref{cor:tree_sorting} with \Cref{prop:rate_of_trisection_WEI} and a union bound on $t = 0, \dots, t_{\infty} - 1$, it holds with probability higher than $1 - 8\cdot 2^{t_{\infty} + 1}\tau_{\infty}|\cH||\cR|\delta$ that
	$$ \| M_{\hat{\pi}_{WM}^{-1}} - M_{\pi^{*-1}} \|_F^2 \lesssim \zeta^2 \log^{11}\left(\frac{6nd}{\delta \zeta_-}\right)\left[\left(\frac{n^{3/4}d^{1/4}}{\zeta^{1/2}}\wedge \frac{nd^{1/6}}{\zeta^{1/3}} \wedge  n\sqrt{d}\wedge \frac{n^{2/3}\sqrt{d}}{\zeta^{1/3}}  \right)+ n\right] \enspace .$$
\end{proof}

In the next four sections, we prove the intermediary results. Propositions~\ref{prop:algo_tree_sorting}--\ref{cor:trisection} are relatively simple. The main difficulty and the key arguments lie in the proofs of Proposition~\ref{prop:rate_of_trisection} and~\ref{prop:rate_of_trisection_WEI} which are respectively in Sections~\ref{sec:proof_rate_trisection} and~\ref{sec:proof_trisection_wei}.

\section{Proofs of Propositions~\ref{prop:algo_tree_sorting}, \ref{prop:algo_block_sorting}, and \ref{cor:trisection}}
\begin{proof}[Proof of \Cref{prop:algo_tree_sorting}]
First, we prove by induction that $\cT_{t_{\infty}}$ is a valid hierarchical sorting tree. Besides, the last part of \Cref{property:block-sorting} enforces that the cardinality of any non-terminal node $G$ of $\cT_{t_{\infty}}$ of depth $t$ is at most $n/2^{t}$. As a consequence, the cardinality of any non-terminal node at depth $t_{\infty-1}$ is at most $1$ and its children $O$ and $I$ are therefore empty.

	We control the error using a telescopic sum. Recall that, by convention, $\pi^*_{\cG_0}=\pi^*$. We start with the following inequality:
	\begin{equation}\label{eq:telescopic_sum}
		\| M_{\hat \pi^{-1}_{\cG_{t_{\infty}}}} - M_{\pi^{*-1}} \|_F^2 \leq 2\| M_{\hat \pi_{\cG_{t_{\infty}}}} - M_{\pi^{*-1}_{\cG_{t_{\infty}}}} \|_F^2 + 2t_{\infty}\sum_{t = 1}^{t_{\infty}}\| M_{\pi^{*-1}_{\cG_t}} - M_{\pi^{*-1}_{\cG_{t-1}}} \|_F^2 \enspace .
	\end{equation}
	Since, for any group $P$ in $\cG_{t_{\infty}}$, $\hat \pi_{\cG_{t_{\infty}}}$ sorts the elements of $P$ uniformly at random and $\pi^{*-1}_{\cG_{t_{\infty}}}$ acts as another permutation of $P$, we deduce from the triangular inequality that
	\begin{align} \nonumber
		\| M_{\hat \pi^{-1}_{\cG_{t_{\infty}}}} - M_{\pi^{*-1}_{\cG_{t_{\infty}}}} \|_F^2 & \leq  \sum_{P \in \cG_{t_{\infty}}}2\|M(P) - \overline M(P)\|_F^2 = 2\sum_{t=1}^{t_{\infty}}\sum_{P \in \cL_t}\|M(P) - \overline M(P)\|_F^2              \\
	& \leq 2\sum_{t=1}^{t_{\infty}}\sum_{\overline P \in \overline \cL_{t}}\|M(\overline P) - \overline M(\overline P)\|_F^2 \enspace , \label{eq:telescopic2}
	\end{align}
	where we used in the last line that $P\subset \overline P$.
	For the second term in~\eqref{eq:telescopic_sum}, remark that $\pi^*_{\cG_{t-1}}(P) = \pi^*_{\cG_{t}}(P)$ for any $P \in \cL_{t-1}$ so that the error at step $t$ in the telescopic sum can be restricted to the groups $G$ that are trisected at step $t-1$:
	\begin{align*}
		\sum_{t = 1}^{t_{\infty}} \| M_{\pi^{*-1}_{\cG_t}} - M_{\pi^{*-1}_{\cG_{t-1}}} \|_F^2 & = \sum_{t=1}^{t_{\infty}}\sum_{G \in \cG_{t-1} \setminus (\cup_{s\geq 1}\cL_{t-s})}\sum_{i \in G}\|M_{\pi^{*-1}_{\cG_t}(\pi^{*}_{\cG_{t-1}}(i))} - M_{i}\|_2^2 \enspace .
	\end{align*}
	Let $(\nO, P, \nI)$ be the trisection obtained at the $t$-th iteration when we apply  $\algoBlockSorting$ to a group $G\in \cG_{t-1}\setminus (\cup_{s\geq 1}\cL_{t-s})$. We also write  $(\overline \nO, \overline P, \overline \nI)$ for the more conservative  trisection obtained at the end of $\algoBlockSorting$. For short, we write $\omega = \pi^{*-1}_{\cG_{t}}\pi^{*}_{\cG_{t-1}}$. We decompose the sum over $i\in G$:
	\begin{align*}
		\sum_{i \in G}\|M_{\omega(i)} - M_{i}\|^2
		= \sum_{i \in \overline \nO} \|M_{\omega(i)} - M_i\|^2
		+ \sum_{i \in \overline \nI} \|M_{\omega(i)} - M_i\|^2
		+ \sum_{i \in \overline P}\|M_{\omega(i)} - M_i\|^2 \enspace ,
	\end{align*}
	By \Cref{property:block-sorting}, all the experts in $\overline \nO$ and in $\overline \nI$ are perfectly sorted within $G$ by $\pi^{*-1}_{\cG(t)}$. As a consequence, the two first sums in the right-hand side term of the above equality are equal to zero. To handle the last term, we introduce the row vector $m(\overline P)$ as the mean of the experts of $M$ over $\overline P$:
	\begin{align*}
		\sum_{i \in G}\|M_{\omega(i)} - M_{i}\|_2^2 & = \sum_{i \in \overline P} \| M_{\omega (i)} - M_i  \|_2^2
		\leq 2\sum_{i \in \overline P} (\| M_i - m(\overline P) \|_2^2 +\|m(\overline P) - M_{\omega (i)} \|_2^2)        \\
		                                            & = 4 \| M(\overline P) - \overline M(\overline P) \|_F^2 \enspace ,
	\end{align*}
	where we used in the last line that $\omega$ acts as a permutation of $\overline{P}$.  Since $\overline P \in \overline \cL_t$, we obtain
	$$\sum_{t = 1}^{t_\infty} \| M_{\pi^{*-1}_{\cG_t}} - M_{\pi^{*-1}_{\cG_{t-1}}} \|^2 \leq
		4\sum_{t=1}^{t_{\infty}}\sum_{\overline P \in \overline \cL_{t}}\|M(\overline P) - \overline M(\overline P)\|_F^2 \enspace .$$
	Together with~\eqref{eq:telescopic_sum} and \eqref{eq:telescopic2}, this concludes the proof since $t_{\infty}\geq 1$.
\end{proof}

\begin{proof}[Proof of \Cref{prop:algo_block_sorting}]
	Consider any data $\cY$, any tree $\cT$ and any leaf $G$ of $\cT$. Let $(O,P,I)$ and $(\overline{O},\overline{P},\overline{I})$ denote the trisections built in   $\algoBlockSorting(\cY, \cT, G)$. For any $\tau< \tau_{\infty}$, let $(\nL_{\tau}, \nU_{\tau})$, $(\overline \nL_{\tau}, \overline \nU_{\tau})$, $(\nO_{\tau}, \nI_{\tau})$ and $(\overline \nO_{\tau}, \overline \nI_{\tau})$ be defined as in \Cref{alg:algoBlockSorting}. We also write $P_{\tau} = G \setminus (O_{\tau}\cup I_{\tau})$ and $\overline P_{\tau} = G \setminus (\overline O_{\tau}\cup \overline I_{\tau})$.
	We only need to prove that, for all $\tau$,  $(\nO_{\tau}, P_{\tau}, \nI_{\tau})$, and  $(\overline \nO_{\tau}, \overline P_{\tau}, \overline \nI_{\tau})$ satisfy \Cref{property:block-sorting}.
	Since
	\begin{equation*}\label{eq:union_trisection}
		\nO_{\tau} = \bigcup_{\tau'< \tau}\nL_{\tau'} \spaceAnd \nI_{\tau} = \bigcup_{\tau'< \tau}\nU_{\tau'} \spaceAnd \overline \nO_{\tau} = \bigcup_{\tau'< \tau}\overline \nL_{\tau'} \spaceAnd \overline \nI_{\tau} = \bigcup_{\tau'< \tau}\overline \nU_{\tau'}\enspace ,
	\end{equation*}
	we easily deduce from \Cref{property:trisection} for $(\nL_{\tau}, \nU_{\tau})$ and $(\overline \nL_{\tau}, \overline \nU_{\tau})$ that the first part  of \Cref{property:block-sorting} is satisfied for $(\nO_{\tau}, P_{\tau}, \nI_{\tau})$, and  $(\overline \nO_{\tau}, \overline P_{\tau}, \overline \nI_{\tau})$.

	\medskip

	Let us turn to the third and fourth parts of \Cref{property:block-sorting}. Let us call $i_{m}$ the expert such that $\pi^*_{\{G\}}(i_{m})= \lfloor |G|/2\rfloor$. In fact, we only need to  prove that $\max_{i\in \nO_{\tau}}\pi^*_{\{G\}}(i)\leq |G|/2$ and $\min_{i\in \nI_{\tau}}\pi^*_{\{G\}}(i)\geq |G|/2$. For this purpose, we prove by induction on $\tau$ that the pivot always satisfies  $\pi^{*-1}_{\{\overline P_{\tau}\}}(\gamma)=i_{m}$ and that all the experts of $O_{\tau}$ (resp. $I_{\tau}$) are below (resp. above) $i_m$, where $\gamma$ depends on $\tau$ and is defined in \Cref{alg:algoBlockSorting}. Assume that this property holds at step $\tau$. Since $\overline{\nO}_{\tau}$ only contains experts that are below the  median expert and since $\gamma = \floor{|G|/2} - |\overline O_{\tau}|$, it follows that $\pi^{*-1}_{\{\overline P\}}(\gamma)=i_{m}$. Consider any $i\in \nO_{\tau+1}$. If $i\in \nO_{\tau}$, then $\pi^*_{\{G\}}(i)\leq |G|/2$ by induction. If $i\in \nL_{\tau}$, then it follows from \Cref{property:trisection} that $i$ is below $i_{m}$, which in turn implies that  $\pi^*_{\{G\}}(i)\leq |G|/2$. By symmetry, the property also holds for $I_{\tau}$. We have proved  the third and the fourth parts of \Cref{property:block-sorting}.

	\medskip

	Finally, we consider the second part of \Cref{property:block-sorting}. Assume that the property holds at step $\tau$. This implies that, for any $i\in \overline{\nO}_{\tau}$, all experts below $i$ belong to $\nO_{\tau}$. Consider any expert $i \in \overline \nO_{\tau + 1}$. If $i \in \overline \nO_{\tau}$, then, by induction, we have  $\pi^*_{\{G\}}(i) = \pi^*_{\{\nO_{\tau}\}}(i) = \pi^*_{\{\nO_{\tau + 1}\}}(i)$. Then, we turn to the case where $i$ belongs to $\overline \nL_{\tau}\subset \overline P_{\tau}$. Consider any $j\in G$ such that  $\pi^{*}_{\{G\}}(j)\leq \pi^{*}_{\{G\}}(i)$.  If $j\in \overline{\nO}_{\tau}$, then we obviously have $j\in \overline{\nO}_{\tau+1}$. If  $j\in \overline P_{\tau}$, then the second part of \cref{property:trisection}
	enforces that $j\in \nL_{\tau}$ and therefore $j\in \nO_{\tau+1}$. Finally, it is not possible that $j\in \overline{\nI}_{\tau}$ since this enforces that $\pi^*_{\{G\}}(j)>|G|/2 >\pi^*_{\{G\}}(i)$ and contradicts the hypothesis. We prove similarly that, for any expert $i \in \overline \nI_{\tau + 1}$, all experts $j$ above $i$ belong to $\nI_{\tau+1}$. 

\end{proof}

\begin{proof}[Proof of \Cref{cor:trisection}]

	As $\algoDoubleTrisection$ is based on multiple applications of the $\algoPivot$ algorithm, we start by considering the latter procedure.

	Consider two sets $|\overline P|\subset[n]$ and $Q\subset[d]$ and  a matrix $\Theta \in \bbR^{|\overline P| \times |Q|}$ which, up to the permutation $\pi^*_{\{\overline{P}\}}$, is bi-isotonic. Let  $Z$ be a noiy observation of $\Theta$,
	\beq\label{eq:model_aggregated}
	Z = \Theta + \NOISEC \ ,
	\eeq
	where the noise matrix $N$ is made of independent, centered, $\zeta$-subGaussian random variables.
	Let $w \in \bbR_+^Q$ be a non-zero vector with nonegative coordinates.  we write
	$(\nL, \nU)$ and $(\overline \nL, \overline \nU)$ for the result of $\algoPivot(Z, w, \gamma)$. 
	We define the event $\cP_3 := \cP_3(\overline P, Q, w, \gamma)$ as the event on $Z$ such that $(L,U)$ and $(\overline{L},\overline{U} )$ satisfy \Cref{property:trisection}.

We remind  that $P'= \overline{P}\setminus (\nL\cup \nU)$, and $\overline{P'}= \overline{P}\setminus (\overline{\nL}\cup \overline{\nU})$.
	\begin{lemma}\label{cor:meta_analysis_of_pivot}
		For any non-zero vector $w \in \bbR_+^Q$, any pivot $\gamma \in \{1, \dots, |\overline P| \}$
		, we have $\P[\cP_3]\geq 1-\delta$. Besides, on the same event of probability at least $1-\delta$, we have
		\begin{equation}\label{eq:condition_espacement_pivot}
			\abs{\proscal<\Theta_{i, \cdot} - \Theta_{i_{\gamma}, \cdot},\frac{w}{\|w\|_2}>} \leq (2\zeta\sqrt{2}+ \overline{\beta}_{\tris})\sqrt{\Log{\frac{2|\overline P|}{\delta}}} \quad  \text{ if } i\in P'\enspace .
		\end{equation}

	\end{lemma}

	Before proving the lemma, let us explain why Proposition~\ref{cor:trisection} is easily deduced from it. The procedure $\algoDoubleTrisection$ calls at most $3|\cH||\cR|$ times $\algoPivot$. Note that, each time, we rely on an independent sample to choose the direction $w$ and to apply  $\algoPivot$. Then, applying the Lemma, we derive that, with probability higher than $1-3|\cH||\cR|\delta$, each of these $3|\cH||\cR|$ sets  $(\nL, \nU)$ and $(\overline \nL, \overline \nU)$ satisfy \Cref{property:trisection}. Hence, we only need to check that \Cref{property:trisection} is stable by union. If, both $(\nL^{(1)}, \nU^{(1)} )$ and $(\overline \nL^{(1)}, \overline \nU^{(1)})$  and $(\nL^{(2)}, \nU^{(2)} )$ and $(\overline \nL^{(2)}, \overline \nU^{(2)})$ satisfy \Cref{property:trisection}, then one easily checks that the first and third part of \Cref{property:trisection} are also true for $(\nL,\nU)=(\nL^{(1)}\cup \nL^{(2)}, \nU^{(1)}\cup \nU^{^{(2)}})$  and $(\overline{\nL},\overline{\nU})=(\overline \nL^{(1)}\cup \overline \nL^{(2)}, \overline \nU^{(1)}\cup \overline \nU^{(2)})$. Consider any expert $i$ in $\overline{L}$. Without loss of generality, we may assume that  $i\in \overline{\nL}^{(1)}$ so that all experts below $i$ in $\overline{P}$ belong to $\nL^{(1)}$ by the second part of \Cref{property:trisection}. As a consequence, all these experts below $i$ belong to $\nL$ and we deduce that the second part of~\Cref{property:trisection} holds. Similarly, we deal with experts $i \in \overline{\nU}$.  This concludes the proof of Proposition~\ref{cor:trisection}.

\end{proof}

\begin{proof}[Proof of \Cref{cor:meta_analysis_of_pivot}]
	Since the noise matrix in~\eqref{eq:model_aggregated} is made of independent $\zeta$-subGaussian random variables, it follows from a union bound, that  with probability higher than $1-\delta$, we have
	\[
		\left|  \langle Z_{i, \cdot},\frac{w}{\|w\|_2}\rangle   -   \langle \Theta_{i, \cdot},\frac{w}{\|w\|_2}\rangle  \right|\leq \zeta \sqrt{2\log\left(\frac{2|\overline{P}|}{\delta}\right)} \ .
	\]
	simultaneously for all $i$ in $\overline{P}$. Since the entries of $w$ are non-negative, the quantities $\langle \Theta_{i, \cdot},\frac{w}{\|w\|_2}\rangle$ are ordered according the permutation $\pi^*_{\{\overline P \}}$. Denote $i_{\gamma}= \pi^{*-1}_{\{\overline{P}\}}(\gamma)$ and $\widehat{i}_{\gamma}$ the index of  $\gamma$-th value of  $\langle Z_{i, \cdot},\frac{w}{\|w\|_2}\rangle$ for $i\in \overline{P}$.
	Since at least $\gamma$ experts satisfy $\langle \Theta_{i, \cdot},\frac{w}{\|w\|_2}\rangle\leq \langle \Theta_{i_{\gamma}, \cdot},\frac{w}{\|w\|_2}\rangle$, we deduce from the above uniform deviation inequality that 
	\[
		\langle Z_{\widehat{i}_{\gamma}, \cdot},\frac{w}{\|w\|_2}\rangle\leq \langle \Theta_{i_{\gamma}, \cdot},\frac{w}{\|w\|_2}\rangle + \zeta \sqrt{2\log\left(\frac{2|\overline{P}|}{\delta}\right)} \ .
	\]
	By symmetry, we  deduce that 
	\[
		\big|\langle Z_{\widehat{i}_{\gamma}, \cdot},\frac{w}{\|w\|_2}\rangle -\langle \Theta_{i_{\gamma}, \cdot},\frac{w}{\|w\|_2}\rangle \big|\leq \zeta \sqrt{2\log\left(\frac{2|\overline{P}|}{\delta}\right)} \ .
	\]
	As a consequence, we have
	\begin{eqnarray}\nonumber
		\proscal<\Theta_{i, \cdot},\frac{w}{\|w\|_2}> &< & \proscal<\Theta_{i_{\gamma}, \cdot},\frac{w}{\|w\|_2}> -(\beta_{\tris}-2\zeta\sqrt{2})\sqrt{\Log{\frac{2|\overline P|}{\delta}}} \quad \text{ if } i\in \nL \, ; \\ \nonumber
		\proscal<\Theta_{i, \cdot},\frac{w}{\|w\|_2}> &> & \proscal<\Theta_{i_{\gamma}, \cdot},\frac{w}{\|w\|_2}>+(\beta_{\tris}-2\zeta\sqrt{2})\sqrt{\Log{\frac{2|\overline P|}{\delta}}}  \quad  \text{ if } i\in \nU \enspace ; \\
		\abs{\proscal<\Theta_{i, \cdot} - \Theta_{i_{\gamma}, \cdot},\frac{w}{\|w\|_2}>} &\leq& (2\zeta\sqrt{2}+ \beta_{\tris})\sqrt{\Log{\frac{2|\overline P|}{\delta}}} \quad  \text{ if } i\in P'\enspace .\label{eq:aa1}
	\end{eqnarray}
	The same inequalities hold for $\overline{\nL}$, $\overline{\nU}$, and $\overline{P}'$ provided that we replace $\beta_{\tris}$ by $\overline{\beta}_{\tris}$. It remains to show that $(\nL,\nU)$ and $(\overline{\nL},\overline{\nU})$ satisfy \Cref{property:trisection}. The first part of the property is obvious. Since $\beta_{\tris}\geq 2\sqrt{2}\zeta$, one observes that $\pi^*_{\{\overline{P}\}}(i)<  \pi^*_{\{\overline{P}\}}(i_{\gamma})=\gamma$ if $i\in \nL$. Similarly, $\pi^*_{\{\overline{P}\}}(i)>\gamma$ if $i\in \nU$ and the third part of \Cref{property:trisection} follows. Turning to the second part of the property, we consider without loss of generality some $i\in \overline{\nL}$ and we need to show that all $j$ satisfying $\pi^*_{\{\overline{P}\}}(j)\leq \pi^*_{\{\overline{P}\}}(i)$ belong to $\nL$. First, such a $j$ does not belong to $\nU$ since  $ \pi^*_{\{\overline{P}\}}(i) < \gamma$. Since $i\in \overline{\nL}$, we deduce that
	\[
		\proscal<\Theta_{i, \cdot},\frac{w}{\|w\|_2}> <  \proscal<\Theta_{i_{\gamma}, \cdot},\frac{w}{\|w\|_2}> -(\overline{\beta}_{\tris}-2\sqrt{2}\zeta)\sqrt{\Log{\frac{2|\overline P|}{\delta}}}\ ,
	\]
	which implies that
	\[
		\big|\proscal<\Theta_{i_{\gamma}, \cdot}-\Theta_{j, \cdot} ,\frac{w}{\|w\|_2}> \big|> (\overline{\beta}_{\tris}-2\sqrt{2}\zeta)\sqrt{\Log{\frac{2|\overline P|}{\delta}}}\geq (2\zeta\sqrt{2}+ \beta_{\tris})\sqrt{\Log{\frac{2|\overline P|}{\delta}}}\ .
	\]
	which in light of~\eqref{eq:aa1} implies that $j\notin P'$. We have proved that $j$ belongs to $\nL$. Hence, \Cref{property:trisection} holds, which concludes the proof.
\end{proof}

\section{Proof of Proposition~\ref{prop:rate_of_trisection}}\label{sec:proof_rate_trisection}
In this section, we prove Proposition~\ref{prop:rate_of_trisection} and thereby control the loss of the estimator $\hat{\pi}$ with simple dimension reduction. For this purpose, we analyze each step of the algorithm. In~\Cref{subsec:algofindcp}, we first prove that, by detecting the high-variation regions of $M$, we are able to aggregate $M$ at some scale $r$ without decreasing much the variation of $M$. This allows us to drastically reduce the dimension of the problem. Then, in Sections~\ref{sec:analysis_pivot} and~\ref{sec:analysis_pivot}, we show that, unless this aggregated matrix $\Theta$ has small  variations,  $\algoDoubleTrisectionPCA$ and $\algoPivot$ will remove some experts so that the corresponding new aggregated matrix $\Theta'$ exhibit significantly smaller variations.  As a consequence, after a polylogarithmic number of iterations of the procedure, the variations of the matrix $M$ restricted to the remaining experts of $\overline P$ is small enough.

\subsection{Notation}
As the arguments rely on considering aggregation of the matrix at different scales, we recall some notation. Let $Y=M+\Noise$ denote a sample of the original matrix. For a set $P\subset[n]$ of experts and a set $\nQ$ of blocks of questions and a scale $r\in R$, we respectively denote
\begin{align*}
	Z(P, \nQ,r)        & = \algoEncodeMatrix(Y,P, \nQ, r) \in \bbR^{P \times \nQ} \\
	\Theta(P, \nQ,r)   & = \algoEncodeMatrix(M,P, \nQ, r)                         \\
	\NOISEC (P, \nQ,r) & = \algoEncodeMatrix(\Noise,P, \nQ, r)\enspace\enspace ,
\end{align*}
the aggregations of $Y$, $M$, and $\Noise$ at scale $r$ so that
$$Z(P, \nQ,r) = \Theta(P, \nQ,r) + \NOISEC(P, \nQ,r) \enspace .$$
By definition of $\algoEncodeMatrix$, all the entries of $\NOISEC$ are independent and $\zeta$-subGaussian.
For any $p\times q$ matrix $A$, we define $\overline a$ as the row vector corresponding to the column-wise mean of $A$, that is  $\overline a_j = \frac{1}{q}\sum_{i=1}^{q} A_{i,j}$. Besides,  we write $\overline{A}$ for $p\times q$ matrix whose experts are all equal to $\overline{a}$.

\subsection{Analysis of $\algoDimensionReduction$}\label{subsec:algofindcp}
In this subsection, we mainly state, that for any $h\in \cH$, $r\in \cR$, the set $\widehat Q_{\cp}$ (which depends on $h$, $r$) detects the high-variation regions of $M$ with high probability. Then, we show in \Cref{lem:constraints_on_CP_blocks}, that for, for some $(h,r)\in \cH\times \cR$, the aggregation of $M$ at scale $r$ and at these high-variation regions contains most of the variance of $M$. This motivates us to work with this aggregated matrix henceforth. Consider any set $\overline P$  of experts  and a sample
$$Y(\overline P)  = M(\overline P) + \Noise(\overline P) \enspace .$$
Fix any scale $r \in \cR$ and any height $h \in \cH$. Recall the two quantities  $r_0$ and $\tilde r$ defined in $\algoDimensionReduction$ by
\begin{equation}\label{eq:def_tilde_r}
	r_0 = 32\zeta^2\log\left(\frac{2d}{\delta}\right)\frac{1}{|\overline P|h^2} \spaceAnd \tilde r = 8(\lceil r_0\rceil  \lor r) \enspace .
\end{equation}
In a nustshell, $r_0$ stands for the minimal scale at which a variation of order $h$ in the mean $\overline{m}(\overline P)=\E[\overline y(\overline P)]$ can be statistically detected. This is why we consider empirical variations of $\overline{y}(\overline P)$ at scale $\tilde{r}\geq (r_0\vee r)$ in Algorithm $\algoDimensionReduction$ to possibly detect variations at scale $r$.

The purpose of this subsection is to prove, that with high probability, the collections $\widehat Q_{\mathrm{\cp}}(h,r) = \algoDimensionReduction(Y,\overline P,h,r)$ of selected blocks of length $r$  is not too large and that there exists at least one $(h,r)\in \cH\times \cR$ such that the aggregation of $M(\overline P)$ at scale $r$ restricted to the blocks $\widehat Q_{\mathrm{\cp}}(h,r)$ captures most of the variance of $M(\overline P)$.

For this purpose, we recall the CUSUM statistics introduced in $\algoDimensionReduction$ and we introduce its population counterpart.
Given positive integers $k\in [d]$ and $r>0$, consider 
\[
	\widehat \bC_{k,r} = \frac{1}{r}\left(\sum_{k' = k}^{k+ r - 1}\overline y_{k'}(\overline P) -  \sum_{k' = k- r}^{k-1}\overline y_{k'}(\overline P)\right)\quad \text{ and }\quad
	\bC_{k,r}^* = \frac{1}{r}\left(\sum_{k' = k}^{k+ r-1}\overline m_{k'}(\overline P) -  \sum_{k' = k- r}^{k-1}\overline m_{k'}(\overline P)\right) \enspace .
\]
Equipped with this notation, we define $\widehat D_{\cp}$ as  in the algorithm as the set of positions $d$ such that the association CUSUM statistic is above the threshold,
and $D_{\cp}^*$ and $\overline D^*_{\cp}$ as some population versions of $\widehat D_{\cp}$, but with different tuning parameters:
\begin{align}
	\widehat D_{\cp}(h,r) & = \left\{k \in [d] ~:~  \widehat \bC_{k,\tilde r} \geq  \frac{1}{4}h \right\} \enspace , \\
	D^*_{\cp}(h,r)        & = \left\{k \in [d] ~:~  \bC^*_{k, 8r} \geq \frac{1}{2}h \right\}\ ; \quad
	\overline D^*_{\cp}(h,r) = \left\{k \in [d] ~:~  \bC^*_{k,\tilde r} \geq \frac{1}{8} h\right\}\enspace .\label{eq:Dcp*}
\end{align}

Then, we consider the collection of blocks $\widehat Q_{\mathrm{\cp}}(h,r)$ , $Q^*_{\mathrm{\cp}}(h,r)$, and 
$\overline Q^*_{\mathrm{\cp}}(h,r)$ of size $r$ that are associated with these positions.  In terms of our algorithms, this means that
$Q^*_{\mathrm{\cp}}(h,r) = \algoEncodeSet(D^*_{\cp},r)$,
$\overline Q^*_{\mathrm{\cp}}(h,r) = \algoEncodeSet(\overline D^*_{\cp},r)$, and $\widehat Q_{\mathrm{\cp}}(h,r) = \algoEncodeSet(\widehat D_{\cp},r)$.
The first proposition states that, with high probability, $\widehat Q_{\mathrm{\cp}}(h,r)$ is sandwidched between $Q^*_{\mathrm{\cp}}(h,r)$ and $\overline Q^*_{\mathrm{\cp}}(h,r)$, so that, on the corresponding event, it is sufficient to study these two quantities.

\begin{lemma}\label{lem:analysis_of_findCP}
	For all $h,r$, the event
	$\xi_{\cp} := \xi_{\cp}(\overline P,h,r)$ defined by
	\begin{equation}\label{eq:blocks_CP}
		Q^*_{\mathrm{\cp}} \subset \widehat Q_{\mathrm{\cp}} \subset \overline Q^*_{\mathrm{\cp}} \enspace ,
	\end{equation}
	holds true with probability at least $1- \delta$.
\end{lemma}

Then, we show that there are not too many significant blocks in $\overline Q^*_{\mathrm{\cp}}$. The proof is based on the fact that the row vector $\overline{m}(\overline P)$ is isotonic and lies in $[0,1]$. As a consequence, there cannot exist two many regions where the variations of $\overline{m}(\overline P)$ is large.

\begin{lemma}\label{lem:constraints_on_CP_blocks} For all $h\in \cH$ and all $r\in \cR$, we have
	\begin{equation}\label{eq:constraints_on_CP_blocks}
		|\overline Q^*_{\cp}| \leq \frac{64\tilde r}{rh} \enspace .
	\end{equation}
\end{lemma}

The next lemma states that, at least for a height $h\in \cH$ and a scale $r\in \cR$, the aggregation of $M$ at scale $r$ and restricted to the regions $Q^*_{\cp}(h,r)$ of significant variations contains almost all the variance of the signal.

For any number $\theta$ and any $\eta>0$,  we define $[\theta]_{\eta} = (-1)^{\mathrm{sgn}(\theta)}\eta\1_{|\theta|\geq \eta}$. For any matrix $\Theta$, we write
$[\Theta]_\eta$ for the thresholded matrix with coefficients $[\Theta_{i,j}]_{\eta}$.

\begin{lemma}\label{lem:restriction_to_CP}
	For any set $\overline P\subset [n]$ and any bi-isotonic matrix $M\in [0,1]^{n\times d}$, there exist
	$r \in \cR$ and  $h \in \cH$ such that
	\begin{equation}\label{eq:restriction_to_CP2}
		\| M(\overline P) - \overline M(\overline P) \|^2_F \leq 16\zeta^2 + 96|\cR||\cH|\left\| \left[\Theta(\overline P, Q^*_{\cp}) - \overline \Theta(\overline P,Q^*_{\cp})\right]_{\sqrt{r}h}\right\|^2_F \enspace ,
	\end{equation}
\end{lemma}

The proof of the above lemmas is postponed to~\Cref{sec:proof:algofindcp}.

\subsection{Analysis of $\algoPivot$ based on the row sums}\label{sec:analysis_pivot}

We consider a specific subset $\overline P$ of experts, a subset $Q$ of blocks of questions, the corresponding aggregated model
\begin{equation}\label{eq:definition_aggregated}
	Z(\overline  P,Q) = \Theta(\overline P,Q) + \NOISEC(\overline P,Q) \in \bbR^{\overline P \times Q} \enspace ,
\end{equation}
and a pivot $\gamma \in [1, |\overline P|]$. Let $(\overline L,\overline U) = \algoPivot(Z,\1_{Q}, \gamma)$ be the conservative result of $\algoPivot$ based on the row sums of $Z(\overline P,Q)$  and let  $\overline P' = \overline P \setminus (\overline L \cup \overline U)$ be the subgroup of experts which have not been classified by $\algoPivot$. The following proposition states that, provided that for some $\eta$ the norm $\|\left[\Theta(\overline P,Q) - \overline \Theta(\overline P,Q)\right]_{\eta}\|^2_F$ is large enough compared  to $|\overline P|\sqrt{|Q|}$, the resulting matrix $\Theta(\overline P',Q) - \overline \Theta(\overline P',Q)$ after $\algoPivot$ has a significantly smaller norm.
We shall often use the following quantity.
\begin{equation}\label{eq:definition_phi_l_1}
	\phi_{l_1}= 2( 2\zeta\sqrt{2}+ \overline{\beta}_{\tris})\leq 29\zeta   \enspace .
\end{equation}

\begin{proposition}\label{prop:rate_analysis_of_pivot}
	Consider any $\overline P\subset [n]$, any $r\in \cR$, and any  subset $Q\subset Q_r$. Also, fix  any $\eta>0$ and any  $\phi>0$. If
	$$\|\left[\Theta(\overline P,Q) - \overline \Theta(\overline P,Q)\right]_{\eta}\|^2_F \geq \frac{1}{\phi}\|\Theta(\overline P,Q) - \overline \Theta(\overline P,Q)\|_F^2
		\geq 8\phi_{l_1}\eta\sqrt{\log(\tfrac{2|\overline P|}{\delta})}|\overline P|\sqrt{|Q|} \enspace,$$
	then, with probability higher than $1-\delta$, we have
	$$\|\Theta(\overline P',Q) - \overline \Theta(\overline P',Q)\|^2_F \leq \left(1 - \frac{1}{16 \phi }\right) \|\Theta(\overline  P,Q) - \overline \Theta(\overline  P,Q)\|^2_F \enspace .$$
\end{proposition}

\subsection{Analysis of $\algoDoubleTrisectionPCA$}\label{sec:analysis_pca}

In this subsection, we state the main result regarding the trisection of a set $\overline P$ based on the first singular vector of a suitable matrix. 
We start from a subset $\overline{P}$ of experts. In $\algoDoubleTrisectionLocal$, we start applying $\algoPivot$ and define $\widetilde{P}=\overline{P}\setminus (\overline{L}_{\cp}\cup \overline{U}_{\cp})$ as the set of experts that have not been classified by $\algoPivot$. 
We are given four independent samples   $\cZ = (Z^{(1)}, Z^{(2)}, Z^{(3)},Z^{(4)})$ according to the aggregated model~\eqref{eq:definition_aggregated}. 
The first three samples are restricted to $\widetilde{P}$, whereas the last one concerns $\overline P$. 
Fix $\gamma \in [1, |\overline  P|]$. We consider $(\overline L,\overline U) = \algoDoubleTrisectionPCA(\cZ,\gamma)$ and $\overline P' = \widetilde{P} \setminus (\overline L_{\pca} \cup \overline  U_{\pca})$ the set of experts that have not been classified by $\algoDoubleTrisectionPCA$.

Recall the definition~\eqref{eq:definition_phi_l_1} of $\phi_{l_1}$. Henceforth, the matrix  $\Theta(\widetilde{P}, Q)$ is said to be undistinguishable in $l_1$-norm if it satisfies
\begin{equation}\label{eq:undistinguishable_in_L1_norm}
	\max_{i,j\in \overline P}  \|\Theta_{i, \cdot}(\widetilde{P},Q) - \Theta_{j, \cdot}(\widetilde P,Q)\|_1 \leq \phi_{l_1} \sqrt{|Q|\log\left(\tfrac{2|\overline P|}{\delta}\right)} \enspace .
\end{equation}
Since, up to permutation of its experts, the matrix $\Theta(\widetilde P, Q)$ is bi-isotonic, the $l_1$ norm $\|\Theta_{i, \cdot}(\widetilde  P,Q) - \Theta_{j, \cdot}(\widetilde  P,Q)\|_1$ is simply the difference of the row sums of $\Theta(\widetilde  P,Q)$. Since $\widetilde  P$ has been deduced from $\overline{P}$ by applying  $\algoPivot(Z,\gamma)$, we can safely assume that
$\Theta(\widetilde P, Q)$  is undistinguishable in $l_1$-norm with high probability -- see the next subsection for a proper justification.

The next result states that, if $\Theta(\widetilde P,Q)$ is undistinguishable in $l_1$-norm and if the Frobenius norm of $\Theta(\widetilde P,Q)-\overline{\Theta}(\widetilde P,Q)$ is large enough, then the corresponding matrix $\Theta(\overline P',Q)$ obtained after trisection has a significantly smaller Frobenius norm.

\begin{proposition}\label{prop:rate_analysis_of_trisectionACP}
	Let $\overline P\subset[n]$ and $Q\subset[d]$.  If $\Theta(\widetilde{P},Q)$ is undistinguishable in $l_1$-norm and if
	\begin{align}\label{eq:condition_signal_trisection_pca}
		\|\Theta(\widetilde P,Q) - \overline \Theta(\widetilde P,Q)\|_F^2\geq  2\cdot 10^5 \zeta^2 \log^3\left(\frac{6nd}{\delta \zeta_-}\right) \left(\sqrt{|\widetilde P||Q|} + |\widetilde P|\right) \enspace ,
	\end{align}
	then, with probability higher than $1-3\delta$, we have
	$$\|\Theta(\overline P',Q) - \overline \Theta(\overline P',Q)\|^2_F \leq \left(1 - \frac{1}{200\log^2(nd/\zeta_-)}\right) \|\Theta(\overline P,Q) - \overline \Theta(\overline P,Q)\|^2_F \enspace .$$
\end{proposition}

Then, we gather the two previous results to analyze the routine $\algoDoubleTrisectionLocal$.
Fix any $\overline P\subset[n]$ and $Q\subset[p]$.
Let $\cZ = (Z^{(1)}(\overline P,Q), Z^{(2)}(\overline P,Q), Z^{(3)}(\overline P,Q), Z^{(4)}(\overline P,Q),Z^{(5)}(\overline P,Q))$ be five independent samples of the model \Cref{eq:model_aggregated}. Fix any $\gamma\in [1,|\overline P|]$.
Let $(\overline L,\overline U)$ be the conservative result of $\algoDoubleTrisectionLocal(\cZ,\gamma)$ and $\overline P' = \overline P \setminus (\overline L \cup \overline U)$.
In the following, we write $\Theta(\overline P,\cQ_{r})$ for the aggregation of $M(\overline P)$ at all blocks of size $r$.

\begin{corollary}\label{cor:rate_analysis_of_trisection_local}
	Fix any $r\in \cR$.   If, for some $\overline P\subset [n]$, $Q\subset \cQ_r$, and   $\eta > 0$, $\Theta(\overline P,Q)$ satisfies
	\begin{equation}\label{eq:energy_trisection_local}
		\left\{\begin{array}{c}\|\left[\Theta(\overline P,Q) - \overline \Theta(\overline P,Q)\right]_{\eta}\|^2_F \geq \frac{1}{120\log^{2}\left(\tfrac{nd}{\delta\zeta_-}\right)}\|\Theta(\overline P,\cQ_r) - \overline \Theta(\overline P,\cQ_r)\|_F^2 \\
			\|\Theta(\overline P,Q) - \overline \Theta(\overline P,Q)\|_F^2  \geq 4\cdot 10^5 \zeta^2 \log^3\left(\frac{6nd}{\delta \zeta_-}\right)\left[\frac{\eta}{\zeta} |\overline{P}|\sqrt{|Q|}\wedge  \left(\sqrt{|\overline P||Q|} + |\overline P|\right)\right] \enspace,
		\end{array}\right.
	\end{equation}
	then, with probability higher than $1-4\delta$, we have
	$$\|\Theta(\overline P',\cQ_r) - \overline \Theta(\overline P',\cQ_r)\|^2_F \leq \left(1 - \frac{1}{3\cdot 10^5 \log^4\left(\frac{nd}{\delta \zeta_-}\right)}\right) \|\Theta(\overline P,\cQ_r) - \overline \Theta(\overline P,\cQ_r)\|^2_F \enspace .$$
\end{corollary}

\subsection{Analysis of $\algoBlockSorting$}\label{sec:analysis_algoBlockSorting}

Next, we combine the results of the previous sections to control the error of $\algoBlockSorting$. We are given a collection $\cY$ of $6\tau_{\infty}$ samples of the model~\eqref{eq:model_0} and a valid hierarchical sorting tree $\cT$ of depth $t$ that we consider as fixed. Then, we take a leaf $G$ of $\cT$ with maximal depth and we consider $(O,P,I)= \algoBlockSorting(\cY, \cT, G)$ the trisection of $G$, as well as $\overline{P}\supseteq  P$ the more conservative intermediary set. In this section, we  provide a high-probability control of $\overline{P}$.

For any height $h\in \cH$ and scale $r\in \cR$, recall that $Q^*_{\cp}$ is the subset(defined in 
\Cref{subsec:algofindcp}) of block of questions at scale $r$ such that the mean $\overline{m}(\overline{P})$ increases by at least $h/2$. Also recall the superset $\overline{Q}^*_{\cp}\supset Q^*_{\cp}$.

At a high level, the next proposition states that, after $\tau_{\infty}$ iterations of the $\algoDoubleTrisectionLocal$ routines at all scales $r\in \cR$ and all heights $h\in \cH$, the size
$\|M( \overline{P})-\overline{M}(\overline P)\|^2_F$ is quite small. This is mainly due to the fact that, by Lemma~\ref{lem:restriction_to_CP}, at each step $\tau$, there exists some $(r,h)\in \cR\times\cH$ such that the norm of   the thresholded aggregated matrix
$\left[\Theta(\overline P_{\tau}, Q^{*}_{\cp}) - \overline \Theta(\overline P_{\tau}, Q^*_{\cp})\right]_{\sqrt{r}h}$ is of the same order as $\|M( \overline{P})-\overline{M}(\overline P)\|^2_F$. By Lemma~\ref{lem:analysis_of_findCP}, the estimated blocks $\widehat{Q}^*_{\cp}$ contain $Q^*_{\cp}$ with high probability. Hence, unless the norm of the thresholded aggregated matrix is small, we derive from corollary~\ref{cor:rate_analysis_of_trisection_local} that the norm of this aggregated matrix has contracted at step $\tau+1$. Hence, after $\tau_{\infty}$ steps, one could expect that the norm of $\|M( \overline{P})-\overline{M}(\overline P)\|^2_F$ is small. In fact, both the statement and the proof of this proposition are slightly more involved because we need to keep track of the scales and heights of interest.
Define the function $\Psi(p,r,h,q)$ by
\begin{equation}\label{eq:definition_psi}
	\Psi(p,r,h,q)= \frac{hp\sqrt{rq}}{\zeta}\land \left(\sqrt{pq}\right)+p\enspace .
\end{equation}

\begin{proposition}\label{prop:rate_analysis_of_block_sorting}
	With probability higher than $1- 5\tau_{\infty}\delta$, there exists a subset $\overline P^{\dagger}$ such that  $\overline{P}\subseteq \overline P^{\dagger}\subseteq G$ and the following property  holds. For some $r^{\dagger} \in \cR$ and some $h^{\dagger}\in \cH$, upon writing
	$Q^{\dagger}_{\cp} = Q^*_{\cp}(\overline P^{\dagger}, h^{\dagger}, r^{\dagger})$ and $\overline Q^{\dagger}_{\cp} = \overline Q^*_{\cp}(\overline P^{\dagger}, h^{\dagger}, r^{\dagger})$,  we have simultaneously
	\begin{align}\label{eq:block_sorting_1}
		\|\left[\Theta(\overline P^{\dagger}, Q^{\dagger}_{\cp}) - \overline \Theta(\overline P^{\dagger}, Q^{\dagger}_{\cp})\right]_{\sqrt{r^{\dagger}}h^{\dagger}}\|^2_F
		& \leq 4\cdot 10^5 \zeta^2 \log^{3}\left(\frac{6nd}{\delta \zeta_-}\right) \Psi(|\overline P^{\dagger}|,r^{\dagger},h^{\dagger},|\overline Q^{\dagger}_{\cp}|) \enspace  ;                        \\
		\label{eq:block_sorting_2}
		\|M(\overline P^{\dagger})-\overline{M}(\overline P^{\dagger})\|^2_F & \leq 16\zeta^2 + 96|\cR||\cH|\|[\Theta(\overline P^{\dagger}, Q^{\dagger}_{\cp}) - \overline \Theta(\overline P^{\dagger}, Q^{\dagger}_{\cp})]_{\sqrt{r^{\dagger}}h^{\dagger}}\|^2_F \enspace .
	\end{align}
\end{proposition}
In other words, there exists a superset $\overline P^{\dagger}$ of $\overline P$ such that, for a suitable height and scale, at the high-variation regions, both the original matrix $M(\overline P^{\dagger})$ and the thresholded aggregated matrix are controlled at the level $\Psi(|\overline P^{\dagger}|,r^{\dagger},h^{\dagger},|\overline Q^{\dagger}_{\cp}|)$. The virtue of the above result is that it easily adapts to the block sorting variant with memory. Unfortunately, the rate $\Psi(|\overline P^{\dagger}|,r^{\dagger},h^{\dagger},|\overline Q^{\dagger}_{\cp}|)$ is a bit difficult to handle. In the next corollary, we replace it by a simpler but cruder bound that only depends on $|G|$, $h^{\dagger}$ and $d$.

\begin{corollary}\label{cor:rate_analysis_of_block_sorting}
	Under the same event of probability higher than $1-5\tau_{\infty}\delta$ as in the previous proposition, the set $\overline P^{\dagger}$, the scale $r^{\dagger}$, and the height $h^{\dagger}$ also satisfy
	\begin{equation}\label{eq:block_sorting_3}
		\|\left[\Theta(\overline P^{\dagger}, Q^{\dagger}_{\cp}) - \overline \Theta(\overline P^{\dagger}, Q^{\dagger}_{\cp})\right]_{\sqrt{r^{\dagger}}h^{\dagger}}\|^2_F
		\lesssim \zeta^2\log^{3.5}\left(\frac{6nd}{\delta \zeta_-}\right) \left[\frac{ h^{\dagger}|G|\sqrt{d}}{\zeta}\land \sqrt{|G|d} \land \sqrt{\frac{|G|}{h^{\dagger}}} + |G|\right]\ ,
	\end{equation}
	where we recall that $G$ is the initial group.
\end{corollary}

\subsection{Analysis of the complete procedure $\algoTreeSorting$}\label{sec:analysis_treesorting}
We are now equipped to prove \Cref{prop:rate_of_trisection}.

\begin{proof}[Proof of \Cref{prop:rate_of_trisection}]
	Let us fix an integer $t\leq t_{\infty}$ and let us consider the collection  $\overline{\cL}_t$   of the $2^{t}$ groups $\overline{P}$ that are not sorted with confidence. Let us apply \Cref{prop:rate_analysis_of_block_sorting} to each of these sets $\overline{P}$.
	In view of this proposition, we define $(\overline P^{\dagger}, h^{\dagger}, r^{\dagger})$ as well as $ Q_{\cp}^{\dagger}$. We also define
	$$ s^{\dagger} = \left|\left\{l ~:~ \exists i,j \in \overline P^{\dagger} \text{ s.t. } \left[\Theta_{i,l}(\overline P^{\dagger}, Q^{\dagger}_{\cp}) - \Theta_{j,l}(\overline P^{\dagger}, Q^{\dagger}_{\cp})\right]_{\sqrt{r^{\dagger}}{h^{\dagger}}} \neq 0\right\}\right|\enspace . $$
	In a nustshell, $s^{\dagger}$ is the number of columns of the thresholded aggregated matrix which are not equal to zero.

	\begin{definition}
		Define the dyadic collection $\cS= \{1,2,4,\ldots, 2^{\lceil \log_2(d)\rceil}\}$. For any $s\in \cS$, $r\in \cR$, and $h\in \cH$, we consider the collection $ \cP^*(h,r,s)\subset \overline \cL_t$ satisfying $r^{\dagger}=r$, $h^{\dagger}=h$, and $s^{\dagger}\in [s,2s)$.
	\end{definition}

	The following lemma controls the cardinality of $\cP^*(h,r,s)$. This bound mainly relies on the facts that the matrix $M$ is, up to a row permutation, bi-isotonic and that its entries lie in $[0,1]$.
	\begin{lemma}\label{lem:constraints_on_peeling}
		Assume that there exists an ordering $\sigma$ of $\overline{\cL}_t$ that orders all groups $\overline{P}$'s. In other words, for any $r\leq s$, any expert $i\in \overline{P}_{\sigma(r)}$ is below any expert $j\in \overline{P}_{\sigma(s)}$. Then, upon this assumption,
		\[
			|\cP^*(h,r,s)|\leq \frac{2d}{hrs}\land 2^{t}\leq \sqrt{\frac{d2^{t+1}}{hrs} }\enspace ,
		\]
		for any  $h\in \cH$, $r\in \cR$, and $s\in \cS$.
	\end{lemma}
	In fact, all the collections $\overline{\cL}_t$ with $t=0,\ldots, t_{\infty}$ satisfy the assumption in the above under the event $\xi$ defined in \Cref{cor:tree_sorting} --see the proof of~\Cref{cor:trisection}.

	Putting everything together and summing over the groups $\cP^*(h,r,s)$, we derive from~\Cref{prop:rate_analysis_of_block_sorting} and \Cref{cor:rate_analysis_of_block_sorting} that, with probability higher than $1-5\cdot2^{t}\tau_{\infty}\delta$, we have
	\begin{align*}
		\sum_{\overline{P} \in \overline{\cL}_{t}}\|M(\overline{P}) - \overline M(\overline{P})\|_F^2 & \leq 16\zeta^2|\overline{\cL}_t| + 96 |\cR||\cH|\sum_{h,r,s}\sum_{\overline P \in \cP^*(h,r,s)}\|\left[\Theta(\overline P, Q^{\dagger}_{\cp}) - \overline \Theta(\overline P, Q^{\dagger}_{\cp})\right]_{\sqrt{r}h}\|_F^2                                 \\
		 & \stackrel{(a)}{\lesssim}  \zeta^2\log^{5.5}\left(\frac{6nd}{\delta\zeta_-}\right)\sum_{h,r,s}\sum_{\overline P \in \cP^*(h,r,s)}\left[\big(\frac{n}{2^{t}\zeta^2}s^{\dagger}rh^2\big)\land \big(\sqrt{\frac{n}{2^{t}h}}\big) + \frac{n}{2^{t}}\right]\\
		  & \stackrel{(b)}{\lesssim} \zeta^2 \log^{5.5}\left(\frac{6nd}{\delta\zeta_-}\right)  \sum_{h,r,s}\left[ \frac{n srh^2}{\zeta^2} \land \sqrt{\frac{dn}{srh^2}} + n \right]                                                                                       \\
		  & \stackrel{(c)}{\lesssim} \zeta^2 \log^{8.5}\left(\frac{6nd}{\delta\zeta_-}\right) \left[\frac{n^{2/3}d^{1/3}}{\zeta^{2/3}} + n \right]\ ,
	\end{align*}
	where in (a), we combined \Cref{cor:rate_analysis_of_block_sorting} with the fact that the size of each group is at most $n/2^{t}$, the crude bound $\|[A]_{\eta}\|_F^2\leq \eta^2 d_1d_2$ for any $d_1\times d_2$ matrix and that $n/2^{t}\geq 1$.
	In (b), we relied on Lemma~\ref{lem:constraints_on_peeling}, whereas in (c) we used that $x\land y\leq x^{1/3}y^{2/3}$.

	We have proved the desired $n^{2/3}d^{1/3}/\zeta^{2/3} + n$ upper bound. The rate $n d^{1/6}/\zeta^{1/3}$ is proved using the same scheme except that we apply  \Cref{cor:rate_analysis_of_block_sorting} differently in (a).
	More precisely, we have
	\begin{align*}
		\sum_{\overline{P} \in \overline{\cL}_{t}}\|M(\overline{P}) - \overline M(\overline{P})\|_F^2 & \leq 16\zeta^2 |\overline{\cL}_t| + 96 |\cR||\cH|\sum_{h,r,s}\sum_{\overline P \in \cP^*(h,r,s)}\|\left[\Theta(\overline P, Q^{\dagger}_{\cp}) - \overline \Theta(\overline P, Q^{\dagger}_{\cp})\right]_{\sqrt{r}h}\|_F^2 \\
		& \leq     \zeta^2\log^{5.5}\left(\frac{6nd}{\delta\zeta_-}\right)\sum_{h,r,s}\sum_{\overline P \in \cP^*(h,r,s)}\left[\big(\frac{n}{2^{t}\zeta}h\sqrt{d}\big)\land \big(\sqrt{\frac{n}{2^{t}h}}\big)\land \sqrt{\frac{n}{2^t}d} + \frac{n}{2^{t}}\right]                        \\
		& \lesssim    \zeta^2\log^{5.5}\left(\frac{6nd}{\delta\zeta_-}\right)\sum_{h} \sum_{r,s}|\cP^*(h,r,s)|\left[\big(\frac{n}{2^t\zeta }h\sqrt{d}\big) \land \big(\sqrt{\frac{n}{2^{t}h}}\big)\land \sqrt{\frac{n}{2^t}d} + \frac{n}{2^{t}} \right]                                        \\
		& \stackrel{(a')}{\lesssim}   \zeta^2\log^{5.5}\left(\frac{6nd}{\delta\zeta_-}\right) \sum_{h}\left[ \frac{nh\sqrt{d}}{\zeta} \land \sqrt{\frac{n^2}{h}}\land \sqrt{n2^td} + n \right]                                                                                        \\
		& \stackrel{(b')}{\lesssim}\zeta^2\log^{6.5}\left(\frac{6nd}{\delta\zeta_-}\right) \left[\frac{nd^{1/6}}{\zeta^{1/3}}\land n\sqrt{d} + n \right]\ ,
	\end{align*}
	where in $(a')$, we used that $\sum_{r,s}|\cP^*(h,r,s)|\leq 2^{t}\leq 2n$ and in (b') that $xy\leq x^{1/3}y^{2/3}$.
\end{proof}

\begin{proof}[Proof of \Cref{lem:constraints_on_peeling}]
	To ease the notation, we write $\cP^*=\cP^{*}(h,r,s)$ in this proof. Since $\cP^*\subset \overline{\cL}_t$, we straightforwardly derive that $|\cP^*|\leq |\overline{\cL}_t|\leq 2^{t}$.
	Let us introduce the width of a matrix $\Theta \in \bbR^{[n] \times \cQ_r}$ on the set $P \subset [n]$ and $Q \subset \cQ_r$:
	$$\Capa (\Theta, P, Q) = \max_{i,j \in P} \sum_{l \in Q}|\Theta_{i, l} - \Theta_{j, l}| \enspace .$$
	Consider any set $\overline{P}$ and the corresponding quantities $\overline P^{\dagger}$, $s^{\dagger}$, $r^{\dagger}$, and $h^{\dagger}$. By definition of $s^{\dagger}$, we have
	$ \Capa (\Theta, \overline P^{\dagger}, Q^{\dagger}_{\cp}) \geq s^{\dagger}\sqrt{r^{\dagger}}h^{\dagger}$.
	Recall that the matrix $\Theta$ is, up to a permutation of its rows, bi-isotonic. Besides, all the groups $\overline{P}$ in $\overline{\cL}_t$ are perfectly ordered by assumption. As a consequence, the width of $\Theta$ on $[n]$ is larger or equal to the sum of the width on each set $\overline{P}$.
	Since $\cP^*$ is an ordered sub-partition, it holds that
	\begin{equation}\label{eq:constraints_on_peeling_2}
		\Capa (\Theta, [n], \cQ_{r}) \geq \sum_{\overline P \in \cP^*}\Capa (\Theta, \overline P^{\dagger}, \cQ_{r}) \geq \sum_{\overline P \in \cP^*} \Capa (\Theta, \overline P^{\dagger}, Q^{\dagger}_{\cp})\geq |\cP^*|s^{\dagger}\sqrt{r}h\geq |\cP^*|s\sqrt{r}h \enspace .
	\end{equation}
	By definition of $\cQ_{r}$, we have $|\cQ_r|\leq 2d/r$. Since the values of $\Theta$ lie in $[0,\sqrt{r}]$, we deduce that $\Capa (\Theta, [n], \cQ_{r}) \leq 2d/\sqrt{r}$. Together with \Cref{eq:constraints_on_peeling_2}, this yields
	\begin{equation*}
		|\cP^*|\leq \frac{2d}{rsh} \enspace ,
	\end{equation*}
	which concludes the proof.

\end{proof}

\section{Remaining proofs for Proposition~\ref{prop:rate_of_trisection}}

\subsection{Proofs of the results on $\algoDimensionReduction$ (\Cref{subsec:algofindcp})}\label{sec:proof:algofindcp}

\begin{proof}[Proof of \cref{lem:analysis_of_findCP}]
	It is sufficient to prove that $ D^*_{\cp}(h,r) \subset \widehat D_{\cp}(h,r)\subset \overline{D}^*_{\cp}(h,r)$. Recall that we use the convention that $\overline{y}_i=\overline{m}_i=0$ if $i\leq 0$ and $\overline{y}_i=\overline{m}_i=1$ if $i>d$.
	Since the CUSUM statistic is linear, we have the decomposition
	\[
		\widehat{\bC}_{k,\tilde{r}} =     \bC^*_{k,\tilde{r}} + \frac{1}{\tilde r}\left(\sum_{k' = k}^{k+ \tilde{r} - 1}\overline e_{k'}(\overline P) -  \sum_{k' = k- \tilde{r}}^{k-1}\overline e_{k'}(\overline P)\right)\ ,
	\]
	where the latter random variable is centered and $\zeta(|\overline P|\tilde{r}/2)^{-1/2}$-subGaussian. By a union bound, we derive that, with probability higher than $1-\delta$, we have
	\[
		\max_{k\in [d]}\big|\widehat{\bC}_{k,\tilde{r}}- \bC^*_{k,\tilde{r}}\big|\leq \zeta\sqrt{\frac{2\cdot 2}{|\overline P|\tilde{r}}\log\left(\frac{2d}{\delta}\right)}\ .
	\]
	Since $\tilde{r}$ is defined in such a way that
	\[
		\zeta\sqrt{\frac{4}{|\overline P|}\log\left(\frac{2d}{\delta}\right)}\leq \frac{1}{8}\sqrt{\tilde{r}} h \ ,
	\]
	we deduce that $\widehat D_{\cp}(h,r)\subset \overline{D}^*_{\cp}(h,r)$. Conversely, if $k$ belongs to $D^*_{\cp}(h,r)$, we have
	$\bC^*_{k,8r}\geq h/2$. Since $\overline{m}(\overline P)$ is an isotonic vector and $\tilde{r}\geq 8r$, it follows that $\bC^*_{k,\tilde{r}}\geq \bC^*_{k,8r}\geq h/2$. We deduce that
	\[
		\widehat{\bC}_{k,\tilde{r}}\geq h[\tfrac{1}{2}-\tfrac{1}{8}]\geq \frac{h}{4}\ ,
	\]
	which implies that $D^*_{\cp}(h,r)\subset \widehat D_{\cp}(h,r)$.
\end{proof}

\begin{proof}[Proof of \Cref{lem:constraints_on_CP_blocks}]
	If an index $k$ belongs to $\overline{D}^*(h,r)$, this implies that $\overline{m}_{k+\tilde{r}}-\overline{m}_{k-\tilde{r}}\geq h/8$, since the vector $\overline{m}$ is isotonic. Define $\kappa= 1+ \lceil \tilde{r}/r\rceil$.
	Since $\overline{m}$ is an isotonic vector, for $l\in \overline{Q}^*_{\cp}(h, r)$, we deduce that $\overline{m}_{l+\kappa r }- \overline{m}_{l- \kappa r}\geq h/8$.  Consider the regular grid $\cQ_{\kappa r}$  of width $\kappa r$ and define $\overline{Q}^*(h, r, \kappa)=\{l \in \cQ_{\kappa r}: \overline{Q}_{\cp}^*(h,r)\cap [l,l+\kappa r)\neq \emptyset\}$. Since, for $l\in  \overline{Q}^*(h, r, \kappa)$, we have $\overline{m}_{l+2\kappa r }- \overline{m}_{l- 2\kappa r}\geq h/8$ and since the total variation of $\overline{m}$ is at most one, this implies
	\[
		\frac{h}{8}|\overline{Q}^*(h,r,\kappa)|\leq  \sum_{l\in D(\kappa,r,h) }\overline{m}_{l+2\kappa r }- \overline{m}_{l- 2\kappa r}\leq \sum_{l\in \cQ_{\kappa r} }\overline{m}_{l+2\kappa r }- \overline{m}_{l- 2\kappa r}  \leq 4 \enspace .
	\]
	Since $|\overline{Q}^*(h,r)|\leq \kappa |\overline{Q}^*(h,r,\kappa)|$ and since $\kappa\leq 2\tilde{r}/r$, we obtain the desired result.

\end{proof}

\begin{proof}[Proof of \Cref{lem:restriction_to_CP}]
	For any height $h \in \cH$ --recall the definition of the dyadic class  $\cH$ in~\eqref{eq:definition_R_D}-- and any expert $i\in \overline P$, we consider the $h$-level set $M_{i,.}-\overline{m}$, that is
	\begin{align}\label{eq:level_set}
		\nF(i,h) & = \{k \in [d] ~:~ M_{i,k} - \overline m_{k} \geq h \}\ ;\quad \quad \nF(i,-h) = \{k \in [d] ~:~ M_{i,k} - \overline m_{k} \leq - h \} \enspace .
	\end{align}
	Since $\nF(i,h)$ and $\nF(i,-h)$ are subsets of $[d]$, we can decompose them into unions of disjoint intervals. For any positive integer $r\in \cR$, we write $\nF(i,h,r)$ as the union of intervals of $\nF(i,h)$ whose size belongs $[2r-1,4r-1)$.
	Finally, we consider the subset $\nF(i,h,r;2h)\subset \nF(i,h,r)$ of all intervals of $\nF(i,h,r)$ that intersect $\nF(i,2h)$. In other words, any maximal interval $I$ in $\nF(i,h,r;2h)$ is a $h$-level set whose size belongs to $[2r-1,4r-1)$ and such that  $M_{i,.}-\overline{m}$ crosses the level $2h$ in $I$. We define similarly $\nF(i,h,r)$ and $\nF(i,h,r;2h)$ when $h$ is negative and $-h\in \cH$.
	It follows from these definitions that, for any $h$ such that either $h\in \cH$ or $-h\in \cH$, we have
	\begin{equation}\label{eq:level_set_jump_set}
		\nF(i, 2h) \subset \bigcup_{r \in \cR} F(i,h,r; 2h) \enspace .
	\end{equation}
	We define  $F^*(h,r,2h)$ as the union of those intervals for $i\in \overline P$.
	$$ F^*(h,r,2h) = \bigcup_{i \in \overline P} F(i,h,r;2h) \enspace .$$
	First, we claim that this collection of intervals  $F^*(h,r,2h)$ is contained in the significant regions of variation of $\overline{m}$. This result heavily relies on the monotonicity assumptions.
	\begin{lemma}\label{lem:jump_set_detected_as_CP}
		For any $h\in \cH$ and any $r\in \cR$.
		$$\left[F^*(h,r,2h)\bigcup F^*(-h,r,-2h)\right]\subset D^*_{\cp}(h,r) \enspace .$$
	\end{lemma}
	Next, we quantify $\| M(\overline P) - \overline M(\overline P)\|_F^2$ using regions of large variation of $M_{i,.}- \overline{m}$.
	\begin{lemma}\label{lem:peeling_in_r_h}
		For any $\overline P$, it holds that
		$$\| M(\overline P) - \overline M(\overline P)\|_F^2 \leq 16\left[ \zeta^2 +  \sum_{i \in \overline P}\sum_{r \in \cR, h \in \cH} h^2\left(|F(i,h,r;2h)|+ |F(i,-h,r;-2h)|\right) \right] \enspace .$$
	\end{lemma}
	The last lemma connects these sets $|F(i,h,r;2h)|$ to the norm of the thresholded aggregated matrix.

	\begin{lemma}\label{lem:majoration_theta}
		For any $r \in \cR$ and $h \in \cH$, we consider $\Theta(\overline P, Q^*_{\cp}(h,r))$ the aggregation of $M$  at scale $r$ and at $Q^*_{\cp}(h,r)$. We have
		$$ h^2\sum_{i \in \overline P} [|F(i,h,r;2h)|+|F(i,-h,r;-2h)|]  \leq 3\|\left[\Theta(\overline P, Q^*_{\cp}(h,r)) - \overline \Theta(\overline P,Q^*_{\cp}(h,r))\right]_{\sqrt{r}h}\|_F^2 \enspace .$$
	\end{lemma}
	Combining Lemmas \ref{lem:peeling_in_r_h} and \ref{lem:majoration_theta}, we conclude that
	\begin{align*}
		\|M(\overline P) - \overline M(\overline P)\|_F^2 & \leq 16\left[ \zeta^2 +\sum_{i \in \overline P}\sum_{r \in \cR,\  h \in \cH} h^2(|F(i,h,r;2h)|+ |F(i,-h,r;-2h)|) \right]                                                                              \\
		                              & \leq 16\left[ \zeta^2 + 3\sum_{r \in \cR,\  h \in \cH}\|\left[\Theta(\overline P,Q^*_{\cp}(h,r)) - \overline \Theta(\overline P,Q^*_{\cp}(h,r))\right]_{\sqrt{r}h}\|_F^2 \right]                              \\
		                              & \leq 16\left[ \zeta^2  +6|\cR||\cH|\max_{r \in \cR,\  h\in \cH}\left\|\left[\Theta(\overline P,Q^*_{\cp}(h,r)) - \overline \Theta(\overline P,Q^*_{\cp}(h,r))\right]_{\sqrt{r}h}\right\|_F^2\right] \enspace ,
	\end{align*}
	which concludes the proof of \Cref{lem:restriction_to_CP}.

\end{proof}

\begin{proof}[Proof of \Cref{lem:jump_set_detected_as_CP}]
	Consider any  $i \in \overline P$, any height $h \in \cH$, and any scale  $r \in \cR$. Without loss of generality, we only focus on $F(i,h,r;2h)$; the case of $F(i,-h,r;-2h)$ being analogous.
	Let $I$ be an interval of $F(i,h,r;2h)$.
	Fix any question $k\in I$ such that $|M_{i,k}-\overline{m}_k|\geq 2h$. Since $k\in F(i,h,r)$, it follows that there exists $l<4r$ such that $M_{i,k+l}-\overline{m}_{k+l}\leq h$. Since both the vectors $M_{i,\cdot}$ and $\overline{m}$ are isotonic, it follows that
	$\overline{m}_{k+l}- \overline{m}_{k}\geq h$. Now consider
	any $k_0\in I$. Using again the monotonicity of $\overline{m}$, we deduce that,
	\[
		\bC^*_{k_0, 8r}\geq \frac{4rh}{8r}\geq \frac{1}{2} h  \ ,
	\]
	and $k_0$ therefore belongs to $D^*_{\cp}(h,r)$. We have proved the desired result.
\end{proof}

\begin{proof}[Proof of \Cref{lem:peeling_in_r_h}]
	Consider any expert  $i \in \overline P$. We decompose the norm of $[M_{i, \cdot} - \overline m(\overline P)]$ using the level sets of this vector. We recall that $\cH$ is  of the form $\{h_{\min},2h_{\min}, 4h_{\min}, \ldots \}$ where $h_{\min}\in [\zeta^2/nd, 2\zeta^2/nd]$.
	\begin{align*}
		\| M_{i, \cdot} - \overline m(\overline P)\|_2^2 & \leq \sum_{h \in \cH}\sum_{k = 1}^d (M_{i, k} - \overline m_k(\overline P))^2 \1\{ 2h \leq |M_{i, k} - \overline m_k(\overline P)| < 4h\} + 4dh_{\min}^2 \\
		                                       & \leq  16\sum_{h \in \cH}\sum_{k = 1}^d h^2 \1\{ 2h \leq |M_{i, k} - \overline m_k(\overline P)| < 4h\} + \frac{16 \zeta^2}{n^2d}                       \\
		                                       & \leq  16\sum_{h \in \cH} h^2 |F(i,2h)|+ \frac{16\zeta^2}{n^2d}                                                                              \\
		                                       & \leq 16\sum_{h \in \cH}\sum_{r \in \cR} h^2|F(i,h,r;2h)| + \frac{16 \zeta^2}{n^2d} \ ,
	\end{align*}
	where in the last line, we used~\eqref{eq:level_set_jump_set}.
	Then, we sum  over $i \in \overline P$ to conclude.
\end{proof}

\begin{proof}[Proof of \Cref{lem:majoration_theta}]
	Consider any  $i \in \overline P$, any height $h \in \cH$, and any scale  $r \in \cR$. Without loss of generality, we only consider $F(i,h,r;2h)$ the case of $F(i,-h,r;-2h)$ being analogous. Let $I$ be a maximal interval of $F(i,h,r;2h)$.  We deduce from Lemma~\ref{lem:jump_set_detected_as_CP} that $I$ is included in $D^*_{\cp}(h,r)$. Let $I_0$ be the largest sub-interval of $I$ of the form $[qr, q'r)$ where $q$ and $q'\in \cQ_r$. Since $|I|\geq 2r-1$, it follows that $|I|\leq 3|I_0|$.  We write $L_0$ the subset of columns of the aggregated matrix $\Theta(\overline P,Q^*_{\cp}(h,r))$ corresponding to $I_0$ so that $|L_0|= |I_0|/r$. On each column $l$ of $L_0$, we have $\Theta_{i, l}(\overline P,Q^*_{\cp}(h,r)) - \overline \theta_{l}(\overline P,Q^*_{\cp}(h,r))\geq \sqrt{r}h$. Putting everything together, we get
	\begin{eqnarray*}
		h^2 |I|&\leq& 3h^2 |I_0| = 3 h^2r |L_0| \leq  3\sum_{l\in L_0}\left(\left[\Theta(\overline P,Q^*_{\cp}(h,r))_{i, l} - \overline \theta(\overline P,Q^*_{\cp}(h,r))_{l}\right]_{\sqrt{r}h}\right)^2\ .
	\end{eqnarray*}
	Summing over all intervals $I$ and over all experts $i\in \overline P$ and also accounting for the $F[i,-h,r;-2h]$ concludes the proof.
\end{proof}

\subsection{Proof of \Cref{prop:rate_analysis_of_pivot}}

To simplify the notation, we  define $\Phi_{l_1}= 2(2\zeta\sqrt{2}+\overline \beta_{\tris}) \sqrt{\log\left(\frac{2|\overline P|}{\delta}\right)}= \phi_{l_1}\sqrt{\log\left(\frac{2|\overline P|}{\delta}\right)}$.

For simplicity,  we respectivly write $\Theta(\overline P) = \Theta(\overline P,Q)$ and $\Theta(\overline P') = \Theta(\overline P',Q)$  in this proof. Recall that $\overline{\theta}(\overline P)$ stands the mean row of $\Theta(\overline P)$ whereas $\overline{\theta}(\overline P')$ stands for the mean row of $\Theta(\overline P')$.

Invoking~\cref{cor:meta_analysis_of_pivot} with $w=\1_{Q}$ and since the matrix $\Theta$ is isotonic, we deduce that outside an event of probability smaller than $\delta$, we have
\begin{equation}\label{eq:condition_l_1_pivot}
	\max_{i,j\in \overline P'} \|\Theta(\overline P')_{i, \cdot} - \Theta(\overline P')_{j, \cdot}\|_1 \leq  \Phi_{l_1} \sqrt{|Q|} \ .
\end{equation}
since the matrix $\Theta$ is isotonic. We shall deduce from this inequality the desired bound. We consider two cases depending on the difference between $ \overline \theta(\overline P)$ and  $\overline \theta(\overline P')$ the mean rows in $\overline P$ and $\overline P'$.

\medskip

{\noindent \bf Case 1}: $|\overline P'|\cdot\| \overline \theta(\overline P)- \overline \theta(\overline P')\|_2^2 > \frac{1}{16\phi}\|\Theta(\overline P) - \overline \Theta(\overline P)\|_F^2 $.
Since $\overline P'\subset \overline P$, we deduce that
\begin{eqnarray*}
	\|\Theta(\overline P)- \overline \Theta(\overline P)\|^2_F - \|\Theta(\overline P')- \overline \Theta(\overline P')\|^2_F &\geq& \sum_{i \in \overline P'}\|\Theta(\overline P)_{i, \cdot}- \overline \theta(\overline P)\|_2^2 - \|\Theta(\overline P)_{i, \cdot}- \overline \theta(\overline P')\|_2^2\\
	&= &|\overline P'|\cdot \|\overline \theta(\overline P) - \overline \theta(\overline P')\|_2^2 \\
	&\geq & \frac{1}{16\phi}\|\Theta(\overline P) - \overline \Theta(\overline P)\|_F^2  \enspace ,
\end{eqnarray*}
where we used the condition in the last line. We have proved the desired result.

	{\noindent \bf Case 2}: $|\overline P'|\cdot  \| \overline \theta(\overline P)- \overline \theta(\overline P')\|_2^2 \leq \frac{1}{16\phi}\|\Theta(\overline P) - \overline \Theta(\overline P)\|_F^2$.
We start with the decomposition
\begin{equation}\label{eq:step_0}
	\|\Theta(\overline P') - \overline \Theta(\overline P')\|_F^2  \leq \|\Theta(\overline P') - \overline \Theta(\overline P)\|_F^2 =  \|\Theta(\overline P) - \overline \Theta(\overline P)\|_F^2 - \|\Theta(\overline P \setminus \overline P') - \overline \Theta(\overline P)\|_F^2\ ,
\end{equation}
so that we only have to control $\|\Theta(\overline P \setminus \overline P') - \overline \Theta(\overline P)\|_F^2$ from below. By definition of the operator $[\cdot]_{\eta}$, we have
\[
	\|\Theta(\overline P \setminus \overline P') - \overline \Theta(\overline P) \|_F^2\geq  \|[\Theta(\overline P \setminus \overline P') - \overline \Theta(\overline P)]_{\eta}\|_F^2 = \|[\Theta(\overline P) - \overline \Theta(\overline P)]_{\eta}\|_F^2 - \|[\Theta(\overline P') - \overline \Theta(\overline P)]_{\eta}\|_F^2.
\]
By assumption, we have $\|[\Theta(\overline P) - \overline \Theta(\overline P)]_{\eta}\|_F^2\geq \phi^{-1}\|\Theta(\overline P) - \overline \Theta(\overline P)\|_F^2$. Hence, as long as we prove that
\begin{equation}\label{eq:main_step_pivot}
	\|[\Theta(\overline P') - \overline \Theta(\overline P)]_{\eta}\|_F^2\leq (2\phi)^{-1}\|[\Theta(\overline P) - \overline \Theta(\overline P)\|_F^2\ ,
\end{equation}
we can safely conclude from~\eqref{eq:step_0} that
\[
	\|\Theta(\overline P') - \overline \Theta(\overline P')\|_F^2\leq (1-(2\phi)^{-1}) \|\Theta(\overline P) - \overline \Theta(\overline P)\|_F^2\ .
\]
Thus, we only have to prove~\eqref{eq:main_step_pivot}. Again, by definition of the thresholding operator, we have
\begin{eqnarray}
	\|[\Theta(\overline P') - \overline \Theta(\overline P)]_{\eta}\|_F^2&=& \eta^2\sum_{i\in \overline P',\ l\in Q}\1_{|\Theta_{i,l}-\overline \theta(\overline P)_l|\geq \eta } \nonumber \\
	&\leq &   \eta^2\sum_{i\in \overline P',\ l\in Q}\1_{|\Theta_{i,l}-\overline \theta(\overline P')_l|\geq \eta/2 }+ \1_{|\theta(\overline P')_l-\overline \theta(\overline P)_l|\geq \eta/2 }\label{eq:pivot_step_3}\ .
\end{eqnarray}
By Markov inequality, the condition that defines Case 2 above implies that
\begin{equation}\label{eq:proof_pivot_majoration_tchebychev}
	\sum_{l\in Q}\1\{ |\overline \theta(\overline P')_l - \overline \theta(\overline P)_l| \geq \eta/2\} \leq \frac{1}{4\phi}\frac{\|\Theta(\overline P) - \overline \Theta(\overline P)\|_F^2}{|\overline P'| \eta^2} \enspace .
\end{equation}

From~\eqref{eq:condition_l_1_pivot} and a convexity argument, we deduce that, for any $i\in \overline P'$,
$\|\Theta(\overline P')_{i,\cdot} - \overline \theta(\overline P')\|_1 \leq  \Phi_{l_1} \sqrt{|Q|}$.
Then, applying again Markov inequality, we deduce that, for any expert $i$ in $\overline P'$ and any $\eta>0$, we have
$$\sum_{l\in Q}\1\{|\Theta(\overline P')_{i,l} - \overline \theta(\overline P')_l| \geq \eta/2\} \leq 2\Phi_{l_1}\frac{\sqrt{|Q|}}{\eta} \enspace .$$
Since we assume that $\|\Theta(\overline P) - \overline \Theta(\overline P)\|^2_F \geq 8\phi \Phi_{l_1}\eta|\overline P|\sqrt{|Q|} \geq 8\phi \Phi_{l_1} \eta|\overline P'|\sqrt{|Q|}$,
we deduce that
\begin{equation}\label{eq:proof_pivot_majoration_markov}
	\sum_{l\in Q}\1\{|\Theta(\overline P')_{i,l} - \overline \theta(\overline P')_l| \geq \eta/2\} \leq \frac{1}{4\phi} \frac{\|\Theta(\overline P) - \overline \Theta(\overline P)\|_F^2}{|\overline P'| \eta^2} \enspace .
\end{equation}
So that, combining \eqref{eq:pivot_step_3}, \Cref{eq:proof_pivot_majoration_tchebychev} and \Cref{eq:proof_pivot_majoration_markov}, we arrive at
\[
	\|[\Theta(\overline P') - \overline \Theta(\overline P)]_{\eta}\|_F^2\leq \frac{1}{2\phi} \|\Theta(\overline P) - \overline \Theta(\overline P)\|_F^2\ .
\]
We have proved \eqref{eq:main_step_pivot}.

\subsection{Proof of \Cref{prop:rate_analysis_of_trisectionACP}}
\label{subsec:rate_analysis_of_trisectionACP}

For simplicity, we write in this proof $\Theta := \Theta(\widetilde P,Q)$ and $\Theta(\overline P') := \Theta(\overline P',Q)$. Without loss of generality, we assume that the rows of $\Theta$ are already ordered according to the oracle order so that $\Theta$ is bi-isotonic.

First, the following lemma states that, the first singular value of $(\Theta - \overline{\Theta})$ is, up to polylogarithmic terms, of the same order as its Frobenius norm. This is mainly due to the fact that the entries of $\Theta$ lies in $[0,\sqrt{r}]$ and that $\Theta$ is a bi-isotonic matrix. 

\begin{lemma}\label{lem:structure_isotonic_matrices}
	Assume that $\|\Theta - \overline{\Theta} \|_F \geq 2\zeta$. For any sets $\widetilde P$ and  $Q$, we have
	\[
		\|\Theta - \overline{\Theta} \|^2_{\mathrm{op}} \geq \frac{1}{16\log^2(nd/\zeta_-)} \| \Theta - \overline{\Theta} \|^2_F \enspace .
	\]
\end{lemma}

Now, write
$\hat v = \argmax_{\|v\|_2 \leq 1} \Big[ \|v^T(Z^{(1)} - \overline{Z}^{(1)})\|_2^2 - \frac{1}{2}\| v^T(Z^{(1)} - \overline{Z}^{(1)} - Z^{(2)} + \overline{Z}^{(2)})\|_2^2\Big]  \enspace .$
\begin{lemma}\label{lem:concentration_pca}
	Fix any $\delta\in (0,1)$.
	If
	\begin{equation}\label{eq:condition_pca}
		\|\Theta - \overline{\Theta}\|_{\mathrm{op}}^2 \geq 1600 \zeta^2\left[\sqrt{|Q|(5|\widetilde P|+\log(6/\delta))} +7|\widetilde  P|+2\log(6/\delta)\right] \enspace ,
	\end{equation}
	then, with probability higher than $1-\delta$, we have  
	$$ \|\hat v^T \left(\Theta - \overline \Theta\right) \|_2^2 \geq \frac{1}{2}\| \Theta - \overline \Theta \|_{\mathrm{op}}^2 \enspace .$$
\end{lemma}
In light of Lemma~\ref{lem:structure_isotonic_matrices} and Condition~\eqref{eq:condition_signal_trisection_pca}, the Condition~\eqref{eq:condition_pca} in Lemma~\ref{lem:concentration_pca} is valid. Consequently, there exists an event of probability higher than $1-\delta$ such that
\begin{equation}\label{eq:pca_signal_with_frobenius}
	\|\hat v^T \left(\Theta - \overline \Theta\right) \|^2_2 \geq \frac{1}{32\log^2(nd/\zeta_-)}\|\Theta - \overline \Theta\|^2_F \enspace .
\end{equation}

Next, we show that a thresholded version of $\hat{z}= ( Z^{(3)}- \overline Z^{(3)} )^T \hat{v}$ is almost aligned  with $z^*= (\Theta-\overline{\Theta})^T \hat{v}$. We define the sets $S^*\subset Q$ and $\hat{S}\subset Q$ of blocks of questions by
\beqn
S^* &=& \left\{ l \in Q ~:~ |z^*_l| \geq 3\zeta\sqrt{2\log(2|Q|/\delta)}\right\} \ ; \quad \hat S = \left\{ l \in Q ~:~ |\hat{z}_l| \geq 2\zeta\sqrt{2\log(2|Q|/\delta)}\right\} \enspace .
\eeqn
$S^*$ stands for the collection of blocks of questions $l$ such that $z^*_l$ is large whereas $\hat{S}$ is the collection of blocks $l$ with large $\hat{z}_l$.  Finally, we consider the vectors  $w^*$ and $\hat{w}$ defined as theresholded versions of $z^*$ and $\hat{z}$ respectively, that is $w^*_i= z^*_i1_{i\in S^*}$ and $\hat{w}_i= \hat{z}_i\1_{i\in \hat{S}}$. Note that, up to the sign, $\hat{w}$ stands for the active coordinates computed in $\algoDoubleTrisectionPCA$.

\medskip
We write $v$ for any unit vector in $\mathbb{R}^{|\widetilde P|}$. Since the noise matrix $N^{(2)}$ is made of independent $\zeta$-subGaussian random variables, it follows that $(v^T (\NOISEC^{(3)} - \overline \NOISEC^{(3)}))_l$ is a $\zeta$-subGaussian random variables. Hence, we deduce that, for any fixed matrix $\Theta$, subsets $\overline P$ and $Q$, and any unit vector $v$, we have
\[\P\left[\max_{l\in Q} \left|(v^T (\NOISEC^{(3)} - \overline \NOISEC^{(3)}))_l\right| \leq \zeta\sqrt{2\log(2|Q|/\delta)}\right]\geq 1-\delta\enspace .
\]
Observe that $\hat{z}=z^* + \hat{v}^T (\NOISEC^{(3)} - \overline \NOISEC^{(3)})$.
Conditioning on $\hat{v}$, we deduce that, on an event of probability higher than $1-\delta$, we have
\begin{equation}\label{eq:upper_noisec2}
	\|\hat{z}-z^* \|_{\infty}\leq \zeta\sqrt{2\log(2|Q|/\delta)}\ .
\end{equation}
Under this event, we have $S^*\subset \hat{S}$ and for $l\in \hat{S}$, we have $z^*_l/\hat{z}_l\in [1/2,2]$. Next, we shall prove that, under this event,  $\hat v^T (\Theta - \overline \Theta) \hat w/\|\hat w\|_2$ is large (in absolute value):
\[
	\abs{\hat v^T (\Theta- \overline \Theta) \hat w}
	= \abs{(z^*)^T \hat{w}}= \sum_{l\in \hat{S}} z^*_l \hat{z}_l\geq \frac{2}{5}\sum_{l\in \hat{S}} (z^*_l)^2 + (\hat{z}_l)^2\geq \frac{2}{5} [\|w^*\|_2^2 + \|\hat{w}\|_2^2]\geq \frac{4}{5}\| \hat w \|_2 \| w^* \|_2\ ,
\]
where we used in the first inequality that $z^*_l/\hat{z}_l\in [1/2,2]$ and in the second inequality that $S^*\subset \hat{S}$. Thus, it holds that
\begin{equation}\label{eq:lb_estw_truew}
	\abs{\hat v^T (\Theta - \overline \Theta) \frac{\hat w}{\|\hat w\|_2}}^2 \geq \frac{16}{25}\|w^* \|_2^2 \enspace .
\end{equation}
It remains to prove that $\|w^* \|_2$ is large enough. Writing $S^{*c}$ for the complementary of $S^*$ in $Q$, it holds that
\begin{equation}\label{eq:wstar_on_S}
	\|w^* \|_2^2 = \|z^*\|_2^2 - \sum_{l\in S^{*c}}(z^*_l)^2\ ,
\end{equation}
so that we need to upper bound the latter quantity. Write $z^*_{S^{*c}}= z^*-w^*$. Coming back to the definition of $z^*$,
\beqn
\left[\sum_{l\in S^{*c}}(z^*_l)^2\right]^2&=& \left[\sum_{l\in S^{*c}}[\hat v^T(\Theta - \overline \Theta )]_lz^*_l\right]^2\\ &\leq&  \|\left(\Theta - \overline \Theta \right)z^*_{S^{*c}}\|^2_2=
\sum_{i \in \widetilde P}\left(\sum_{l \in S^{*c}}(\Theta_{i,l} - \overline \theta_{l})z^*_{l} \right)^2 \\
&\leq &  \frac{18\zeta^2}{|\widetilde  P|^2} \log\left(\frac{2|Q|}{\delta}\right) \sum_{i \in \widetilde  P}\left(\sum_{l \in S^{*c}}\sum_{j\in \widetilde  P}|\Theta_{i,l} - \Theta_{j,l}|\right)^2\\
&\leq & \frac{18\zeta^2}{|\widetilde P|^2} \log\left(\frac{2|Q|}{\delta}\right) \sum_{i \in \widetilde P}\left(\sum_{j\in \widetilde P}\|\Theta_{i,\cdot} - \Theta_{j,\cdot}\|_1\right)^2\\
&\leq & 18\zeta^2 \phi_{l_1}^{2} \log\left(\frac{2|Q|}{\delta}\right)\log\left(\frac{2|\widetilde P|}{\delta}\right) |\widetilde P||Q|\\
&\leq & \left[145 \zeta^2 \log\left(\frac{2|Q||\widetilde P|}{\delta}\right)\sqrt{|\widetilde{P}||Q|}\right]^2 \enspace  ,
\eeqn
where we used the definition of $S^*$ in the third line as well as the Condition~\eqref{eq:condition_l_1_pivot} in the fifth line. We recall that $\phi_{l_1}=2(2\zeta\sqrt{2}+ \overline{\beta}_{\tris})\leq 29\zeta$ is defined in~\eqref{eq:definition_phi_l_1}. Recall that  $z^*= \hat{v}^T(\Theta-\overline{\Theta})$. Combining \Cref{eq:pca_signal_with_frobenius}, \Cref{eq:wstar_on_S}, and Condition~\eqref{eq:condition_signal_trisection_pca}, we deduce that 
\[
	\|w^* \|_2^2 \geq \frac{1}{64\log^2(nd/\zeta_-)}\|\Theta - \overline \Theta\|_F^2\ ,
\]
which, together with \Cref{eq:lb_estw_truew}, yields
\[
	\left\|(\Theta - \overline \Theta) \frac{\hat w}{\|\hat w\|_2} \right\|_2^2 \geq \abs{\hat v^T (\Theta - \overline \Theta) \frac{\hat w}{\|\hat w\|_2}}^2 \geq \frac{1}{100\log^2(nd/\zeta_-)}\|\Theta - \overline \Theta\|_F^2 \enspace .
\]
Write $\hat{w}^{(1)}$ and $\hat{w}^{(2)}$ the positive and negative parts of $\hat{w}$
respectively so that $\hat{w}= \hat{w}^{(1)}- \hat{w}^{(2)}$ and $\hat{w}^{+}= \hat{w}^{(1)}+ \hat{w}^{(2)}$. We obviously have $\|\hat{w}\|_2= \|\hat{w}^{+}\|_2$. Besides, if the rows of $\Theta$ are ordered according to the oracle permutation, then
$(\Theta - \overline \Theta)\hat{w}^{(1)}$ and $(\Theta - \overline \Theta)\hat{w}^{(2)}$ are increasing vectors with mean zero. It then follows from Harris' inequality that these two vectors have a nonegative inner product. We have proved that
\begin{equation}\label{eq:lower_former_signal}
	\left\|(\Theta - \overline \Theta) \frac{\hat w^+}{\|\hat w^+\|_2} \right\|_2^2\geq \left\|(\Theta - \overline \Theta) \frac{\hat w}{\|\hat w\|_2} \right\|_2^2  \geq  \frac{1}{100\log^2(nd/\zeta_-)}\|\Theta - \overline \Theta\|_F^2 \enspace .
\end{equation}
Equipped with this bound, we are now in position to show that the set $\overline P'$ of experts obtained from $\widetilde{ P}$ when applying the pivoting algorithm with $\hat w^+/\|\hat w^+\|_2$ has a much smaller variance.

By  \Cref{cor:meta_analysis_of_pivot}, there exists an event of probability higher than $1-\delta$ such that
$$\max_{i,j\in \overline P'}\abs{\proscal<\Theta(\overline P')_{i,\cdot} - \Theta(\overline P')_{j,\cdot},\frac{\hat w^+}{\|\hat w^+\|_2}>}  \leq \phi_{l_1}\sqrt{\log\left(\tfrac{2|\overline P|}{\delta}\right)}\ , $$
where we recall that $\phi_{l_1}= 2(2\zeta\sqrt{2}+ \overline{\beta}_{\tris})$. By convexity, it follows that
\[
	\left\|(\Theta(\overline P') - \overline \Theta(\overline P'))\tfrac{\hat w^+}{\|\hat w^+\|_2} \right\|_2^2 \leq \phi_{l_1}^{2}\log\left(\tfrac{2|\overline P|}{\delta}\right)|\overline P'| \leq \phi_{l_1}^{2}\log\left(\tfrac{2|\overline P|}{\delta}\right)|\widetilde P|\ .
\]
In light of Condition~\eqref{eq:condition_signal_trisection_pca}, this quantity is small compared to $\| \Theta - \overline{\Theta}\|_F^2$:
\begin{equation}\label{eq:upper_new_signal}
	\| (\Theta(\overline P') - \overline \Theta(\overline P'))\tfrac{\hat w^+}{\|\hat w^+\|_2} \|_2^2 \leq\frac{1}{200\log^2(nd/\zeta_-)} \| \Theta - \overline{\Theta}\|_F^2\ ,
\end{equation}
which together with~\eqref{eq:lower_former_signal} leads to
\begin{equation}\label{eq:upper_new_signal_2}
	\|(\Theta - \overline \Theta )\tfrac{\hat w^+}{\|\hat w^+\|_2} \|_2^2 - \|(\Theta(\overline P') - \overline \Theta(\overline P') )\tfrac{\hat w^+}{\|\hat w^+\|_2} \|_2^2\geq \frac{1}{200\log^2(nd/\zeta_-)}  \| \Theta - \overline{\Theta}\|_F^2\enspace .
\end{equation}
Since $\overline P'\subset \widetilde P$, we deduce that, for any vector $w'\in \mathbb{R}^q$, we have
$\|(\Theta - \overline \Theta)w' \|_2^2 \geq \|(\Theta(\overline P') - \overline \Theta(\overline P') )w' \|^2$. It then follows from the Pythagorean theorem that
\[
	\|\Theta - \overline \Theta\|_F^2 - \|\Theta(\overline P') - \overline \Theta(\overline P') \|_F^2 \geq \|(\Theta - \overline \Theta )\tfrac{\hat w^+}{\|\hat w^+\|_2} \|_2^2 - \|(\Theta(\overline P') - \overline \Theta(\overline P') )\tfrac{\hat w^+}{\|\hat w^+\|_2} \|_2^2 \enspace .
\]
Then, together with  \eqref{eq:upper_new_signal_2}, we arrive at
$$\|\Theta(\overline P') - \overline \Theta(\overline P')\|_F^2 \leq \left(1 - \frac{1}{200\log^2(nd/\zeta)}\right)\|\Theta - \overline \Theta\|_F^2 \enspace .$$

\begin{proof}[Proof of \Cref{lem:structure_isotonic_matrices}]

	The proof mainly relies on a discretisation argument. Given any $a\in \mathbb{R}$ and any matrix $U$, we define the matrix $[U]_{a}^\peel$ by
	$([U]_a^\peel)_{i,j} = U_{i,j}\mathbf 1\{U_{i,j} \in (a, 2a]\}$.  If $a$ is negative, then the interval should be understood as $[2a,a)$.
	Recall that  all the  entries of $\Theta - \overline \Theta$ lie in
	$[-\sqrt{r},\sqrt{r}]$. This allows us to decompose this matrix as follows
	\begin{align*}
		r^{-1/2}(\Theta - \overline \Theta) & = \sum_{i \in \mathbb N^*} [r^{-1/2}(\Theta - \overline \Theta)]_{2^{-i}}^\peel + [r^{-1/2}(\Theta - \overline \Theta)]_{-2^{-i}}^\peel\enspace .
	\end{align*}
	All the matrices in this decomposition have disjoint support.
	For all $i > \log_2(nd/\zeta)$, all the entries of the discretised matrices in the decomposition are smaller than $\zeta (nd)^{-1}$. Since $|\widetilde{P}|=p\leq n$ and  $|Q|=q\leq d/r$, this implies that
	$$\Big\|\sum_{i \in \mathbb N^*, i> \log_2(nd/\zeta)} [r^{-1/2}(\Theta - \overline \Theta)]_{2^{-i}}^\peel +[r^{-1/2}(\Theta - \overline \Theta)]_{-2^{-i}}^\peel\Big\|_F^2 \leq \zeta^2 \frac{nd}{(nd)^2r} \leq\frac{\zeta^2}{r} \enspace .$$
	Coming back to the previous bound, we arrive at
	\begin{align*}
		\|r^{-1/2}(\Theta - \overline \Theta)\|^2_F & \leq  \sum_{i \in \mathbb N^*, i\leq \log_2(nd/\zeta)} \|[r^{-1/2}(\Theta - \overline \Theta)]_{2^{-i}}^\peel\|^2_F +  \|[r^{-1/2}(\Theta - \overline \Theta)]_{-2^{-i}}^\peel\|^2_F + \frac{\zeta^2 }{r}\enspace .
	\end{align*}
	As we assume that $\|\Theta - \overline \Theta\|_F \geq 2\zeta \geq 2\zeta /\sqrt{r}$,
	\begin{align*}
		\|r^{-1/2}(\Theta - \overline \Theta)\|^2_F & \leq  \frac{4}{3}\sum_{i \in \mathbb N^*, i\leq \log_2(nd/\zeta)} \|[r^{-1/2}(\Theta - \overline \Theta)]_{2^{-i}}^\peel\|^2_F  + \|[r^{-1/2}(\Theta - \overline \Theta)]_{-2^{-i}}^\peel\|^2_F\enspace .
	\end{align*}
	Hence, there exists an integer $i_0\in [1,\log_2(nd/\zeta)]$ such that
	\begin{align*}
		\frac{3\|\Theta - \overline \Theta\|^2_F}{8 r \log_2(nd/\zeta)} & \leq   \|[r^{-1/2}(\Theta - \overline \Theta)]_{2^{-i_0}}^\peel\|^2_F \lor \|[r^{-1/2}(\Theta - \overline \Theta)]_{-2^{-i_0}}^\peel\|^2_F\enspace .
	\end{align*}
	Assume w.l.o.g.~that, for this $i_0 \leq \log_2(nd/\zeta)$, we have
	\begin{align*}
		\frac{3\|\Theta - \overline \Theta\|^2_F}{8 r \log_2(nd/\zeta)} & \leq   \|[r^{-1/2}(\Theta - \overline \Theta)]_{2^{-i_0}}^\peel\|^2_F\enspace .
	\end{align*}
	Now, we define a different discretised version. For a matrix $U$ and some  $a \in \mathbb R^+$, let $[U]_a$ be defined by
	$([U]_a)_{ij} = (a \mathbf 1\{U_{i,j} \geq a\})_{i,j}$. We readily deduce that
	\begin{align}\label{eq:peel1}
		\|\Theta - \overline \Theta\|^2_F & \leq   \frac{32}{3} r \log_2(nd/\zeta) \|[r^{-1/2}(\Theta - \overline \Theta)]_{2^{-i_0}}\|^2_F.
	\end{align}
	The entries of the matrix $[r^{-1/2}(\Theta - \overline \Theta)]_{2^{-i_0}}$ lie in $\{0, 2^{-i_0}\}$. Up to a permutation of the rows of $\Theta$, we can assume that each column of $[r^{-1/2}(\Theta - \overline \Theta)]_{2^{-i_0}}$ is isotonic. One can easily check that a matrix that only takes two values and such that each column is isotonic can be transformed into a bi-isotonic matrix by applying a suitable permutation $\pi_0$ to its columns.
	We denote $B$ the corresponding permuted matrix. Recall that we denote $p$ and $q$ the dimensions of $B$. Then, define the function $\phi:[p]\rightarrow \{0,\ldots, q\}$ such that $\phi(i)$ is the number of non-zero entries in the $(p-i+1)$-th row of $B$. Since $B$ is bi-isotonic, the function $\phi$ is non-increasing.  Besides, we have
	\[
		\sum_{i=1}^{ p} \phi(i) = 2^{i_0} \sum_{i,j} B_{i,j} = 2^{2i_0} \sum_{i,j} B_{i,j}^2 = 2^{2i_0} \|B\|_F^2\enspace .
	\]
	\begin{lemma}\label{lem:peelmonfunc}
		Let $d_1$ and $d_2$ be two positive integers and consider a non-increasing function $f: [d_1]\rightarrow \mathbb{R}_+$. Then, there exists $m\in [d_1]$ such that $\sum_{i=1}^{d_1} f(i)\leq \log(ed_1)m f(m)$.
	\end{lemma}
	Applying this lemma to $\phi$, we deduce that, for some $m\in [p]$, we have
	\begin{align}\label{eq:applilempeel}
		2^{2i_0} \|B\|_F^2\leq   \log(ep)m \phi(m) \ .
	\end{align}
	Since $\phi(m)$ is the number of non-zero entries on the $p+1-m$-th row of $B$, since $B$ is bi-isotonic and since $B$ only takes two values, this implies that  $B$ contains in the lower right a rectangle of size $m\times \phi(m)$ with value $2^{-i_0}$. Define the vector $u\in \mathbb{R}^p$ such that $u_i=m^{-1/2}$ if $i\geq p-m+1$ and $u_i=0$, otherwise. Define also the vector $v\in \mathbb R^{q}$ such $v_j = 1/\sqrt{\phi(m)}$ if $j \geq q - \phi(m)+1$, and $v_j = 0$ otherwise. It follows from these definitions that $u^TBv =  2^{-i_0} \sqrt{m \phi(m)}$.
	Recall that $[r^{-1/2}(\Theta - \overline \Theta)]_{2^{-i_0}}$ corresponds to a row and column permutation of $B$. Hence, there exist two permutations $\pi_1$ and $\pi_2$ such that
	\[
		u^T_{\pi_1}[r^{-1/2}(\Theta - \overline \Theta)]_{2^{-i_0}}v_{\pi_2} = u^TBv\ .
	\]
	By construction, the entries of $\Theta - \overline \Theta$ are higher than $2^{-i_0}$ for all entries such that $(u_{\pi_1})_i\neq 0$ and $(v_{\pi_1})_j\neq 0$. We deduce that
	\[
		\left\|\Theta - \overline \Theta\right\|_{\mathrm{op}}\geq \sqrt{r} u^T_{\pi_1}r^{-1/2}(\Theta - \overline \Theta) v = \sqrt{r}2^{-i_0}\sqrt{m\phi(m)}\geq \frac{\sqrt{r}}{\sqrt{\log(ep)}}\|B\|_F\ .
	\]
	Finally, we come back to~\eqref{eq:peel1} to conclude that
	$\|\Theta - \overline \Theta\|_{\mathrm{op}} \geq [32\log(e p)\log_2(nd/\zeta)/3]^{-1/2} \|\Theta - \overline \Theta\|_F$, where we recall that $\zeta_- < \zeta$.
\end{proof}

\begin{proof}[Proof of Lemma~\ref{lem:peelmonfunc}]
	Define $a= \sup_{m=1}^{d_1}mf(m)$. As a consequence, we have $f(m)\leq a/m$. This implies that
	\begin{align*}
		\sum_{i=1}^{d_1} f(i) \leq \sum_{i=1}^{d_1} \frac{a}{i}\leq a\log(e d_1) \enspace .
	\end{align*}
	We have proved that $\log(ed_1) \sup_{m=1}^{d_1}mf(m)\geq   \sum_{i=1}^{d_1} f(i)$.    \end{proof}

\begin{proof}[Proof of \Cref{lem:concentration_pca}]

	We start with the two following lemmas. For short, we write $p=|\overline P|$ and $q=|Q|$ in this proof. 

	\begin{lemma}\label{lem:conc:scal}
		Let $\NOISEC'$ denote a random $d_1\times d_2$ matrix whose entries follow independent, centered and $\zeta$-subGaussian distributions. Let $\Omega\subset \mathbb{R}^{d_2}$ be a subspace of dimension $d_2'$. With probability larger than $1-\delta$, one has
		$$\sup_{u\in \mathbb R^{d_1},\  v\in \Omega:\  \|u\|_2\leq 1,\ \|v\|_2\leq 1} | u^T (\NOISEC' - \overline{\NOISEC'}) v | \leq 10 \zeta \sqrt{d_1+d_2'+ \log(2/\delta)}\enspace ,$$
		where $\overline{\NOISEC'}= d_1^{-1}\1_{d_1} \1_{d_1}^{T}\NOISEC'$ is made of the mean row of $\NOISEC'$.
	\end{lemma}

	\begin{lemma}\label{lem:conc:eigen}
		Let $\NOISEC'$ be a random $d_1\times d_2$ matrix whose entries follow independent, centered and $\zeta$-subGaussian distributions. It holds with probability larger than $1-\delta$ that
		\[\sup_{u \in \mathbb R^{d_1}: \|u\|_2\leq 1} |\|u^T(\NOISEC' - \overline{\NOISEC'})\|_2^2 - \mathbb{E} \|u^T(\NOISEC' - \overline{\NOISEC'})\|_2^2| \leq 64 \zeta^2 \left[\sqrt{d_2(5d_1+\log(2/\delta)}+ (5d_1+  \log(2/\delta))\right]\enspace .
		\]
	\end{lemma}

	We have
	\[
		Z^{(1)} - \overline{Z}^{(1)} = \Theta - \overline{\Theta} + \NOISEC^{(1)} - \overline{\NOISEC}^{(1)}\enspace ,
	\]
	so that, for any $v\in\mathbb{R}^p$,
	\begin{align*}
		\|v^T(Z^{(1)} - \overline{Z}^{(1)})\|_2^2 = \|v^T(\Theta - \overline{\Theta})\|_2^2 +  \|v^T\NOISEC^{(1)} - v^T\overline{\NOISEC}^{(1)}\|_2^2 + 2\langle v^T\NOISEC^{(1)} - v^T\overline{\NOISEC}^{(1)}, v^T(\Theta - \overline{\Theta}) \rangle\enspace ,
	\end{align*}
	which, in turn, implies that
	\begin{eqnarray}  \label{eq:upper_1}
		\lefteqn{    \Big|\|v^T(Z^{(1)} - \overline{Z}^{(1)})\|_2^2 - \|v^T(\Theta - \overline{\Theta})\|_2^2 - \E\big[\|v^T\NOISEC^{(1)} - v^T\overline{\NOISEC}^{(1)}\|_2^2\big] \Big|\leq } & &\\ &&  \Big|\|v^T\NOISEC^{(1)} - v^T\overline{\NOISEC}^{(1)}\|_2^2 - \E \big[\|v^T\NOISEC^{(1)} - v^T\overline{\NOISEC}^{(1)}\|_2^2\big] \Big|
		+ 2|\langle v^T\NOISEC^{(1)} - v^T\overline{\NOISEC}^{(1)}, v^T(\Theta - \overline{\Theta}) \rangle|\enspace .\nonumber
	\end{eqnarray}

	Write $W\subset \mathbb{R}^q$ for the image of $(\Theta-\overline{\Theta})^T$. Then, we apply~\Cref{lem:conc:scal} to derive that
	\begin{eqnarray}\nonumber
		\sup_{v \in \mathbb R^{p}:\ \|v\|_2\leq 1}|\langle v^T(\NOISEC^{(1)} - \overline{\NOISEC}^{(1)}), v^T(\Theta - \overline{\Theta}) \rangle &\leq& \|\Theta - \overline{\Theta}\|_{\mathrm{op}}\sup_{v \in \mathbb R^{p}:\ \|v\|_2\leq 1,\  u \in W:\ \|u\|_2\leq 1} | v^T(\NOISEC^{(1)} - \overline{\NOISEC}^{(1)}) u|\\
		&\leq& 10 \zeta \|\Theta - \overline{\Theta}\|_{\mathrm{op}} \sqrt{2p+\log(6/\delta)}\enspace , \label{eq:scalprodfin}
	\end{eqnarray}
	with probability higher than $1-\delta/3$ since the dimension of $W$ is no larger than $p$.
	We deduce from~\Cref{lem:conc:eigen} that, with probability higher than $1-\delta/3$, we have
	$$\sup_{v\in \mathbb R^{p}:\ \|v\|_2\leq 1}\left|\Big|\|v^T(\NOISEC^{(1)} - \overline{\NOISEC}^{(1)})\|_2^2 - \mathbb E \|v^T(\NOISEC^{(1)} - \overline{\NOISEC}^{(1)})\|_2^2\right|\leq 64 \zeta^2 \left[\sqrt{q(5p+\log(6/\delta))}+ 5p+  \log(6/\delta)\right] \enspace .$$
	Together with~\eqref{eq:upper_1} and~\eqref{eq:scalprodfin}, we have that with probability larger than $1-2\delta/3$,
	\begin{align*}
		\sup_{v\in \mathbb R^{p}:\ \|v\|_2\leq 1}
		 & \Big|\|v^T(Z^{(1)} - \overline{Z}^{(1)})\|_2^2 - \|v^T(\Theta - \overline{\Theta})\|_2^2 - \mathbb E \|v^T\NOISEC^{(1)} - v^T\overline{\NOISEC}^{(1)}\|_2^2 \Big| \\ &\leq 10\zeta\|\Theta - \overline{\Theta}\|_{\mathrm{op}} \sqrt{2p+\log(6/\delta)}+ 64 \zeta^2\left[\sqrt{q(5p+\log(6/\delta))}+ (5p+  \log(6/\delta))\right] \enspace .
	\end{align*}
	In the same way, we have that, with probability larger than $1-\delta/3$, 
	\begin{align*}
		\sup_{v\in \mathbb R^{p}:\ \|v\|_2\leq 1} & \Big|\frac{1}{2}\|v^T(Z^{(1)} - \overline{Z}^{(1)} - Z^{(2)} + \overline{Z}^{(2)})\|_2^2  - \mathbb E \|v^T(N^{(1)} - \overline{N}^{(1)})\|_2^2 \Big|                                                                \\
	& = \frac{1}{2}\sup_{v\in \mathbb R^{p}:\ \|v\|_2\leq 1}\Big|\|v^T(Z^{(1)} - \overline{Z}^{(1)} - Z^{(2)} + \overline{Z^{(2)}})\|_2^2  - \mathbb E \|v^T(Z^{(1)} - \overline{Z}^{(1)} - Z^{(2)} + \overline{Z}^{(2)})\|_2^2 \Big| \\ &\leq  128 \zeta^2\left[\sqrt{q(5p+\log(6/\delta)}+ (5p+  \log(6/\delta))\right] \enspace.
	\end{align*}
	Putting everything together we conclude that, on an event  of probability higher than $1-3\delta$, we have simultaneously for all	$v \in \mathbb R^{p}$ with $\|v\|_2\leq 1$ that
	\begin{align*}
		 & \Big|\|v^T(Z^{(1)} - \overline{Z}^{(1)})\|_2^2 - \|v^T(\Theta - \overline{\Theta})\|_2^2 - \frac{1}{2}\|v^T(Z^{(1)} - \overline{Z}^{(1)} - Z^{(2)} + \overline{Z}^{(2)})\|_2^2 \Big| \\ &\leq 10\zeta\|\Theta - \overline{\Theta}\|_{\mathrm{op}} \sqrt{2p+\log(6/\delta)}+ 192 \zeta^2\left[\sqrt{q(3p+\log(6/\delta))}+ (3p+  \log(6/\delta))\right] \enspace .
	\end{align*}
	Since $\|\Theta - \overline{\Theta}\|_{\mathrm{op}}^2 \geq 1600 \zeta^2[\sqrt{q(5p+\log(6/\delta))} +7p+2\log(6/\delta)]$, we deduce that, on the same
	event, we have
	\[
		\sup_{v\in \mathbb{R}^p: \ \|v\|_2\leq 1} \Big|\|v^T(Z^{(1)} - \overline{Z}^{(1)})\|_2^2 - \|v^T(\Theta - \overline{\Theta})\|_2^2 - \frac{1}{2}\|v^T(Z^{(1)} - \overline{Z}^{(1)} - Z^{(2)} + \overline{Z}^{(2)})\|_2^2 \Big|\leq \frac{1}{4}\|\Theta - \overline{\Theta}\|_{\mathrm{op}}^2\enspace .\]
	Writing  $\psi(v)= \big|\|v^T(Z^{(1)} - \overline{Z}^{(1)})\|_2^2  - \frac{1}{2}\|v^T(Z^{(1)} - \overline{Z}^{(1)} - Z^{(2)} + \overline{Z}^{(2)})\|_2^2\big|$, we deduce that, for $v$ such that $\|v^T(\Theta - \overline{\Theta})\|_{2}^2 = \|\Theta - \overline{\Theta}\|_{\mathrm{op}}^2$, we have $\Psi(v)\geq \tfrac{3}{4}\|\Theta - \overline{\Theta}\|_{\mathrm{op}}^2$, whereas, for $v$ such that $\|v^T(\Theta - \overline{\Theta})\|_{2}^2< \frac{1}{2} \|\Theta - \overline{\Theta}\|_{\mathrm{op}}^2$, we have $\Psi(v)<\tfrac{3}{4}\|\Theta - \overline{\Theta}\|_{\mathrm{op}}^2$. We conclude that $\hat{v}$ satisfies $\|\hat{v}^T(\Theta - \overline{\Theta})\|_{2}^2> \frac{1}{2} \|\Theta - \overline{\Theta}\|_{\mathrm{op}}^2$.

\end{proof}

\begin{proof}[Proof of Lemma~\ref{lem:conc:scal}]
	We start with a classical result. Variants of it can be found in random matrix textbooks (see e.g~\cite{vershynin2018high}). Still, we provide a simple dedicated proof below for the sake of completeness.
	\begin{lemma}\label{lem:conc:scalpur}
		Let  $\NOISEC'$ be a $d_1\times d_2$ matrix whose entries follow independent, centered, and  $\zeta$-subGaussian distributions. Consider any vector subspace $\Omega\subset\mathbb R^{d_2}$  with dimension $d_2'$. With probability higher than $1-\delta$, one has
		$$\sup_{u\in \mathbb R^{d_1},\  v\in \Omega:\ \|u\|_2\leq 1,\  \|v\|_2\leq 1} |u^T \NOISEC' v | \leq 5 \zeta \sqrt{d_1+d_2'+ \log(1/\delta)}\enspace .$$
	\end{lemma}
	We have the following decomposition
	\begin{eqnarray*}
		\sup_{\stackrel{u\in \mathbb R^{d_1}, v\in \Omega}{\|u\|_2\leq 1, \|v\|_2\leq 1} } | u^T (\NOISEC' - \overline{\NOISEC'}) v |& \leq &\sup_{\stackrel{u\in \mathbb R^{d_1}, v\in \Omega}{\|u\|_2\leq 1, \|v\|_2\leq 1} } | u^T \NOISEC' v |+  \sup_{\stackrel{u\in \mathbb R^{d_1}, v\in \Omega}{\|u\|_2\leq 1, \|v\|_2\leq 1} } | u^T \overline{\NOISEC'} v |  .
	\end{eqnarray*}
	The first expression in the right-hand side is handled with~\Cref{lem:conc:scalpur}. Regarding the second one, we observe that  $\overline{\NOISEC'} v $ is a constant vector. As a consequence,
	\[
		\sup_{\stackrel{u\in \mathbb R^{d_1}, v\in \Omega}{\|u\|_2\leq 1, \|v\|_2\leq 1} } | u^T \overline{\NOISEC'} v |\leq \sqrt{d_1}\sup_{v\in \Omega}\|\overline{n}' v \|_2\ ,
	\]
	where $\overline{n}'$ is a $\zeta/\sqrt{d_1}$-subGaussian random vector. Then, we control this expression applying Lemma~\ref{lem:conc:scalpur} to a $1\times d'_2$ matrix. All in all, we have proved that, with probability higher than $1-\delta$, we have
	\[
		\sup_{\stackrel{u\in \mathbb R^{d_1}, v\in \Omega}{\|u\|_2\leq 1, \|v\|_2\leq 1} } | u^T (\NOISEC' - \overline{\NOISEC'}) v |\leq 10 \zeta\sqrt{(d_1+d'_2)+\log\left(\frac{2}{\delta}\right)}\enspace .
	\]
\end{proof}

\begin{proof}[Proof of Lemma~\ref{lem:conc:scalpur}]
	Let $\cU_{d}(\epsilon)$ denote the $\epsilon$-covering number of the $d$-dimensional unit ball and let  $\cU_{d}(\epsilon)$ denote a corresponding minimal covering set. For $\Omega$ a $d'_2$-dimensional subspace of $\mathbb{R}^{d'_2}$, we also write with a slight abuse of notation $\cU_{d'_2}(\epsilon)$ for a corresponding minimal covering set of its unit ball. Consider any $d_1\times d_2$ matrix $W$. Write
	$w^*= \sup_{u: \|u\|_2\leq 1}\sup_{v\in \Omega, \ \|v\|_2\leq 1} |u^T W v|$ and $w= \sup_{u\in \cU_{d_1}(1/4)}\sup_{v\in \cU_{d'_2}(1/4)}
		|u^{T}W v|$. Given $u\in \mathbb{R}^{d_1}$, let $\pi(u)$ denote any closest point of $u$ in $\cU_{d_1}(1/4)$. Similarly, for $v\in \Omega$, $\pi'(v)$ stands for a closest point of $v$ in $\cU_{d'_2}(1/4)$.
	By triangular inequality, we have
	\begin{eqnarray*}
		w^* &\leq & w  + \sup_{u: \|u\|_2\leq 1}\sup_{v\in \Omega, \ \|v\|_2\leq 1}
		|u^{T}W v|-|\pi(u) W \pi'(v)|\\
		&\leq & w+ \sup_{u: \|u\|_2\leq 1}\sup_{v\in \Omega, \ \|v\|_2\leq 1}  |(u^{T}-\pi(u)^{T})W v| +   |\pi(u)^{T}W (v-\pi'(v))|\\
		&\leq &  w +w^*/2\enspace .
	\end{eqnarray*}
	We have proven that
	\begin{equation}\label{eq:operator+covering}
		\sup_{u: \|u\|_2\leq 1}\,\sup_{v\in \Omega, \ \|v\|_2\leq 1} |u^T W v|\leq 2\sup_{u\in \cU_{d_1}(1/4)}\, \sup_{v\in \cU_{d'_2}(1/4)}
		|u^{T}W v|\enspace .
	\end{equation}

	Since $\log(\cU_{d}(\epsilon))\leq d\log(3/\epsilon)$ (see e.g.~\cite{wu2017lecture}), we deduce from triangular inequality that, with probability higher than $1-\delta$, we have
	\[
		\sup_{u\in \mathbb R^{d_1},\  v\in \Omega:\ \|u\|_2\leq 1,\  \|v\|_2\leq 1} |u^T \NOISEC' v | \leq 2\zeta\sqrt{2(d_1+d'_2)\log(12)+2\log(1/\delta)}\ .
	\]
\end{proof}

\begin{proof}[Proof of Lemma~\ref{lem:conc:eigen}]
	Relying on~\eqref{eq:operator+covering} with $W= (\NOISEC'-\overline{\NOISEC}')(\NOISEC'-\overline{\NOISEC}')^T- \mathbb{E}\left[(\NOISEC'-\overline{\NOISEC}')(\NOISEC'-\overline{\NOISEC}')^T\right]$, we derive that $\sup_{u: \|u\|_2\leq 1}    |\|u^T(\NOISEC'-\overline{\NOISEC}')\|_2^2 - \mathbb E \|u^T(\NOISEC'-\overline{\NOISEC}')\|_2^2|$ is less than or equal to
	\[
		2 \sup_{u\in \cU_{d_1}(1/4)}\, \sup_{v\in \cU_{d_1}(1/4)}
		u^T (\NOISEC'-\overline{\NOISEC}')(\NOISEC'-\overline{\NOISEC}')^T v - \E\left[u^T (\NOISEC'-\overline{\NOISEC}')(\NOISEC'-\overline{\NOISEC}')^T v\right]
	\]
	As a consequence, it amounts to simultaneously  control $|\cU_{d_1}(1/4)|^2$ quadratic forms of subGaussian random variables. For this purpose, we use the Hanson-Wright inequality~\cite{vershynin2018high}. Below we provide a version of this inequality with explicit numerical constants.
	\begin{lemma}\label{lem:HW}
		Let $x$ be $d$-dimensional $\zeta$-subGaussian centered random vector with independent components. For any $d\times d$ matrix $A$ and any $t>0$, we have
		\[
			\P\left[x^T A x - \E[x^TAx]\geq 32 \zeta^2 \left(\|A\|_F\sqrt{t}+\|A\|_{\mathrm{op}}t\right)\right]  \leq 2e^{-t}
		\]
	\end{lemma}

	For any fixed $u$ and $v$, we interpret $u^T (\NOISEC'-\overline{\NOISEC}')(\NOISEC'-\overline{\NOISEC}')^T v$ as a quadratic form of $d_1d_2$ independent random variables where the corresponding matrix $B$ of the quadratic form satisfies $\|B\|_{\mathrm{op}}\leq 1$ and $\|B\|_{F}\leq \sqrt{d_2}$. Putting everything together we deduce that, with probability higher than $1-\delta$, we have
	\begin{eqnarray*}
		\lefteqn{    \sup_{u: \|u\|_2\leq 1}    \left|\|u^T(\NOISEC'-\overline{\NOISEC}')\|_2^2 - \mathbb E \|u^T(\NOISEC'-\overline{\NOISEC}'))\|_2^2\right|}&&\\  &\leq &
		64 \zeta^2 \left[\sqrt{d_2(2d_1\log(12)+ \log(2/\delta)}+2d_1\log(12)+\log(2/\delta)\right]\\
		&\leq & 64 \zeta^2 \left[\sqrt{d_2(5d_1+\log(2/\delta)}+ (5d_1+  \log(2/\delta))\right]
		\enspace .
	\end{eqnarray*}

\end{proof}
\begin{proof}[Proof of Lemma~\ref{lem:HW}]
	We consider separately the diagonal terms of $A$ and the non-diagonal terms. Write $A^{-}$ for the matrix  such that $A^{-}_{ij}= A_{ij}\1_{i\neq j}$.  First, we use Section 2.8 in~\cite{pollard_book} to handle $x^T A^{-}x$. We know that
	\[
		\P\left[ x^{T}A^{-}x  \geq 8\zeta^2\left(\|A^{-}\|_F\sqrt{t}+\sqrt{2}\|A^{-}\|_{\mathrm{op}}t\right)\right]\leq e^{-t}\enspace ,
	\]
	for any $t>0$. Regarding the diagonal part, we know from Rudelson and Vershynin~\cite{rudelson2013hanson} (Step 1 of the main proof) that $\|x_i^2 - \E[x_i^2]\|_{\psi_1}\leq 4\zeta^2$ (see~\cite{vershynin2018high} for a definition of $\|.\|_{\psi_1}$). Then, we are in position to apply Bernstein's inequality~\cite{boucheron2013concentration} (Theorem 2.10) to $\sum_{i}a_{ii}x_i^2$ with $v= (16\zeta^2)^{2}\sum_i a^2_{ii}$ and $c= 16\zeta^2\max_{i}|a_{ii}|$. For any $t>0$, we have
	\[
		\P\left[\sum_{i=1}^d a_{ii} (x_i^2-\E[x_i^2]) \geq 16\zeta^2\left( \sqrt{2\sum_{i}a_{ii}^2 t} + \max_{i}|a_{ii}|t\right)\right]\leq e^{-t}\ ,
	\]
	which implies that
	\[
		\P\left[\sum_{i=1}^d a_{ii} (x_i^2-\E[x_i^2]) \geq 16\zeta^2 \left(\|A\|_F\sqrt{2 t} + \|A\|_{\mathrm{op}}t\right)\right]\leq e^{-t}\ .
	\]
	We combine the two deviation inequalities and use $\|A^{-}\|_{\mathrm{op}}\leq 2\|A\|_{\mathrm{op}}$ to conclude that
	\[
		\P\left[ x^{T}Ax -\E[x^TA x] \geq 32\zeta^2\left[\|A\|_F\sqrt{t}+\|A\|_{\mathrm{op}}t\right]\right]\leq 2e^{-t}\enspace .
	\]

\end{proof}

\subsection{Proof of \Cref{cor:rate_analysis_of_trisection_local}}
Let $(\overline L_{\cp}, \overline U_{\cp})$ denote the conservative result of $\algoPivot(Z^{(1)},\1_{Q},\gamma)$  and
$\widetilde P = \overline  P\setminus (\overline L_{\cp}\cup  \overline U_{\cp})$.~\\  Let $(\overline L_{\pca}, \overline U_{\pca}) = \algoDoubleTrisectionPCA(\cZ,\gamma)$ with $\cZ = (Z^{(2)}, Z^{(3)}, Z^{(4)}, Z^{(5)})$. Here,  $(Z^{(2)}, Z^{(3)}, Z^{(4)})$ restricted to the experts in $\widetilde P$, whereas $Z^{(5)}$ is restricted to experts in $P$. Finally, we write $\overline P' = \widetilde P \setminus (\overline{L}_{\pca} \cup \overline{U}_{\pca})$. We first prove the following intermediary result
\begin{equation}\label{eq:cor_intermediary_1}
	\|\Theta(\overline  P',Q) - \overline \Theta(\overline  P',Q)\|^2_F
	\leq \left(1 - \frac{1}{1920\log^2(\frac{nd}{\delta\zeta_-})}\right) \|\Theta(\overline  P,Q) - \overline \Theta(\overline  P,Q)\|^2_F       \enspace .
\end{equation}

We consider two cases. First, we assume that $\frac{\eta}{\zeta}|\overline  P|\sqrt{|Q|}\leq \sqrt{|\overline  P||Q|}+|\overline  P|$. Then, it follows from Equation~\eqref{eq:energy_trisection_local} that we are in position to apply~\Cref{prop:rate_analysis_of_pivot} with $\phi= 120\log^2(\tfrac{nd}{\delta\zeta_-})$. Since $\overline  P' \subset \widetilde P$, it follows that
$\|\Theta(\overline  P',Q) - \overline \Theta(\overline  P',Q)\|^2_F\leq \|\Theta(\widetilde P,Q) - \overline \Theta(\widetilde P,Q)\|^2_F$ and~\eqref{eq:cor_intermediary_1} follows from~\Cref{prop:rate_analysis_of_pivot}.

Now, we assume that  $ \sqrt{|\overline  P||Q|}+|\overline  P|\leq \frac{\eta}{\zeta}|\overline  P|\sqrt{|Q|}$. If $\|\Theta(\widetilde P,Q) - \overline \Theta(\widetilde P,Q)\|^2_F\leq 0.5\|\Theta(\overline  P,Q) - \overline \Theta(\overline  P,Q)\|^2_F$, then the result obviously holds. Otherwise, it follows from~\eqref{eq:energy_trisection_local} that
\[
	\|\Theta(\widetilde P,Q) - \overline \Theta(\widetilde P,Q)\|_F^2  \geq 2\cdot 10^5\log^3\left(\frac{6nd}{\delta \zeta_-}\right)\left[\eta|\overline  P|\sqrt{|Q|}\land \left(\sqrt{|\overline  P||Q|} + |\overline  P|\right)\right]\enspace,
\]
Besides, with probability higher than $1-\delta$, $\Theta(\widetilde P)$ is undistinguishable in $l_1$-norm by \eqref{eq:condition_l_1_pivot}. Hence, we are in position to apply~\cref{prop:rate_analysis_of_trisectionACP} and it follows that
\beqn
\|\Theta(\overline  P',Q) - \overline \Theta(\overline  P',Q)\|^2_F& \leq &\left(1 - \frac{1}{200\log^2(\tfrac{nd}{\delta \zeta_-})}\right) \|\Theta(\widetilde P,Q) - \overline \Theta(\widetilde P,Q)\|^2_F \\
&\leq& \left(1 - \frac{1}{1920\log^2(\tfrac{nd}{\delta\zeta_-})}\right) \|\Theta(\overline  P,Q) - \overline \Theta(\overline  P,Q)\|^2_F       \enspace ,
\eeqn
which is exactly Equation~\eqref{eq:cor_intermediary_1}.
\medskip

It remains to conclude from Equation~\eqref{eq:cor_intermediary_1}. We start from
\beqn
\|\Theta(\overline  P',\cQ_r) - \overline \Theta(\overline  P',\cQ_r)\|^2_F&=& \|\Theta(\overline  P',Q) - \overline \Theta(\overline  P',Q)\|^2_F+ \|\Theta(\overline  P',\cQ_r\setminus Q) - \overline \Theta(\overline  P',\cQ_r\setminus Q)\|^2_F\\
&\leq & \left(1 - \frac{1}{1920\log^2(\frac{nd}{\delta\zeta_-})}\right) \|\Theta(\overline  P,Q) - \overline \Theta(\overline  P,Q)\|^2_F \\ &&   +  \|\Theta(\overline  P,\cQ_r\setminus Q) - \overline \Theta(\overline  P,\cQ_r\setminus Q)\|^2_F\\
&\leq & \left(1 - \frac{1}{3\cdot 10^5 \log^4\left(\frac{nd}{\delta \zeta_-}\right)}\right) \|\Theta(\overline  P,\cQ_r) - \overline \Theta(\overline  P,\cQ_r)\|^2_F \ ,
\eeqn
where we used in the last line that $\|\Theta(\overline  P,Q) - \overline \Theta(\overline  P,Q)\|^2_F\geq 1/(120\log^2(\tfrac{nd}{\delta \zeta_-})) \|\Theta(\overline  P,\cQ_r) - \overline \Theta(\overline P,\cQ_r)\|^2_F$.

\subsection{Proof of \cref{prop:rate_analysis_of_block_sorting}}
For all $\tau = 0, \dots, \tau_{\infty}$, let $(\overline O_{\tau}, \overline P_{\tau}, \overline I_{\tau})$ be the sets defined in $\algoBlockSorting$. Let also
$(h^{\dagger \tau} , r^{\dagger \tau}) = \argmax_{h,r}\|\left[\Theta(\overline P_{\tau}, Q^*_{\cp}) - \overline \Theta(\overline P_{\tau}, Q^*_{\cp})\right]_{\sqrt{r}h}\|^2_F$, where we recall that $Q^*_{\cp}$ depends on $h$ and $r$. For simplicity, we write $Q^{\dagger \tau}_{\cp} = Q^*_{\cp}(h^{\dagger \tau}, r^{\dagger \tau})$. Equipped with this notation, we readily deduce from \Cref{lem:restriction_to_CP} that
\begin{equation}\label{eq:uper_deconcaten1}
	\| M(\overline P_{\tau}) - \overline M(\overline P_{\tau}) \|^2_F \leq 16\zeta^2 + 96|\cR||\cH|\left\| \left[\Theta(\overline P_{\tau}, Q^{\dagger \tau}_{\cp}) - \overline \Theta(P,Q^{\dagger \tau}_{\cp})\right]_{\sqrt{r^{\dagger \tau}}h^{\dagger \tau}}\right\|^2_F \enspace .
\end{equation}
If, for some $\tau<  \tau_{\infty}$, we have
\begin{equation}\label{eq:upper_p_dagg}
	\left\|\left[\Theta(\overline P_{\tau}, Q^{\dagger \tau}_{\cp}) - \overline \Theta(\overline P_{\tau}, Q^{\dagger \tau }_{\cp})\right]_{\sqrt{r^{\dagger\tau }}h^{\dagger\tau}}\right\|^2_F
	\leq 4\cdot 10^5 \zeta^2 \log^{3}\left(\frac{6nd}{\delta \zeta_-}\right)\Psi(|\overline P_{\tau}|,r^{\dagger\tau },h^{\dagger\tau },|\overline Q^{\dagger\tau }_{\cp}|)\enspace ,
\end{equation}
then we can fix, for any such $\tau$,  $\overline{P}^{\dagger}= \overline{P}_{  \tau}$, $h^{\dagger}= h^{\dagger \tau}$, $Q^{\dagger}= Q^{\dagger \tau}_{\cp}$,  and $r^{\dagger}= r^{\dagger \tau}$ so that both the properties~\eqref{eq:block_sorting_1} and \eqref{eq:block_sorting_2} hold.

Hence, we assume henceforth that, for all $\tau$, Equation~\eqref{eq:upper_p_dagg} does not hold and we shall arrive at a contradiction.   In particular, this implies that $\| [\Theta(\overline P_{\tau}, Q^{\dagger \tau}_{\cp}) - \overline \Theta(\overline P_{\tau},Q^{\dagger \tau}_{\cp})]_{\sqrt{r^{\dagger \tau}}h^{\dagger \tau}}\|^2_F \geq \zeta^2\geq \zeta^2(|\cR||\cH|)^{-1}$. We have  $112|\cR||\cH|\leq 120\log^2(\tfrac{nd}{\delta\zeta_-})$  provided that $n$ is a large enough constant. In light of~\eqref{eq:uper_deconcaten1}, this implies that, for all $\tau$,
\begin{align}
	\left\| \left[\Theta(\overline P_{\tau}, Q^{\dagger \tau}_{\cp}) - \overline \Theta(\overline P_{\tau},Q^{\dagger \tau}_{\cp})\right]_{\sqrt{r^{\dagger \tau}}h^{\dagger \tau}}\right\|^2_F
	 & \geq \frac{1}{120\log^2(\tfrac{nd}{\delta\zeta_-})}\| M(\overline P_{\tau}) - \overline M(\overline P_{\tau}) \|^2_F \label{eq:signal_Q_dagger}                                          \\ \label{eq:signal_Q_dagger2}
	 & \geq \frac{1}{120\log^2(\tfrac{nd}{\delta\zeta_-})}\| \Theta(\overline P_{\tau}, \cQ_{r^{\dagger \tau}}) - \overline \Theta(\overline P_{\tau}, \cQ_{r^{\dagger \tau}}) \|^2_F\enspace ,
\end{align}
where we recall that $\cQ_{r^{\dagger \tau}}$ is the collection of all blocks at scale $r^{\dagger \tau}$ and we use the Pythagorean equality in the second line. Applying Lemma~\ref{lem:analysis_of_findCP} at the scale $r^{\dagger \tau}\in \cR$, at the height $h^{\dagger \tau}\in \cH$, and at all steps $\tau=0,\ldots, \tau_{\infty}-1$, we deduce that  the event $\overline \xi_{\cp}:=\bigcap_{\tau = 0}^{\tau_{\infty} - 1} \xi_{\cp}(\overline P_{\tau}, h^{\dagger \tau}, r^{\dagger \tau})$ holds with probability at least $1 - \tau_{\infty} \delta$. Under this event, we write $\widehat Q^{\tau}_{\cp}$ for the estimated set defined at step $\tau$ and scales $(h^{\dagger \tau}, r^{\dagger \tau})$ in $\algoBlockSorting$. Then, it holds that $Q^{\dagger \tau}_{\cp} \subset \widehat Q^{\tau}_{\cp} \subset  \overline{Q}^{\dagger \tau}_{\cp}$ and we deduce from~\eqref{eq:signal_Q_dagger2} that
\begin{align*}
	\left\| \left[\Theta(\overline P_{\tau}, \widehat Q^{\tau}_{\cp}) - \overline \Theta(\overline P_{\tau},\widehat Q^{\tau}_{\cp})\right]_{\sqrt{r^{\dagger \tau}}h^{\dagger \tau}}\right\|^2_F
	\geq \frac{1}{120\log^2(\tfrac{nd}{\delta\tau^-})}\| \Theta(\overline P_{\tau}, \cQ_{r^{\dagger \tau}}) - \overline \Theta(\overline P_{\tau}, \cQ_{r^{\dagger \tau}}) \|^2_F\enspace .
\end{align*}
Since \eqref{eq:upper_p_dagg} is not satisfied, we also have
\begin{align*}
	\left\| \Theta(\overline P_{\tau}, \widehat Q^{\tau}_{\cp}) - \overline \Theta(\overline P_{\tau},\widehat Q^{\tau}_{\cp})\right\|^2_F
	 & \geq         \left\| \left[\Theta(\overline P_{\tau}, \widehat Q^{\tau}_{\cp}) - \overline \Theta(\overline P_{\tau},\widehat Q^{\tau}_{\cp})\right]_{\sqrt{r^{\dagger \tau}}h^{\dagger \tau}}\right\|^2_F \\
	 & \geq
	\left\| \left[\Theta(\overline P_{\tau}, Q^{\dagger \tau}_{\cp}) - \overline \Theta(\overline P_{\tau},Q^{\dagger \tau}_{\cp})\right]_{\sqrt{r^{\dagger \tau}}h^{\dagger \tau}}\right\|^2_F           \\
	 & \geq 4\cdot 10^5  \zeta^2 \log^{3}\left(\frac{6nd}{\delta \zeta_-}\right)\Psi(|\overline P_{\tau}|,r^{\dagger\tau },h^{\dagger\tau },|\widehat Q^{\tau }_{\cp}|) \ ,
\end{align*}
since $\widehat Q^{\tau}_{\cp} \subset  \overline{Q}^{\dagger \tau}_{\cp}$. Hence, we are in position to apply
\Cref{cor:rate_analysis_of_trisection_local} at all steps $\tau$ with $r^{\dagger \tau}$, $\overline P_{\tau}$, $\widehat Q^{\tau}_{\cp}$, and $\eta = \sqrt{r^{\dagger \tau}}h^{\dagger \tau}$. There exists an event of probability higher than  $1 - 4\tau_{\infty}\delta$ such that, at all steps $\tau$, we have
\begin{equation*}
	\|\Theta(\overline P_{\tau+1},\cQ_{r^{\dagger \tau}}) - \overline \Theta(\overline P_{\tau + 1},\cQ_{r^{\dagger \tau}})\|^2_F \leq \left(1 - \frac{1}{3\cdot 10^5 \log^4\left(\frac{nd}{\delta\zeta_-}\right)}\right) \|\Theta(\overline P_{\tau};\cQ_{r^{\dagger \tau}}) - \overline \Theta(\overline P_{\tau};\cQ_{r^{\dagger \tau}})\|^2_F \enspace .
\end{equation*}
Together with Equation~\eqref{eq:signal_Q_dagger}, we deduce that
\begin{align*}
	\| M(\overline P_{\tau}) & - \overline M(\overline P_{\tau}) \|^2_F - \| M(\overline P_{\tau+1}) - \overline M(\overline P_{\tau+1}) \|^2_F                                                                                                                                                 \\
	                         & \geq \|\Theta(\overline P_{\tau},\cQ_{r^{\dagger \tau}}) - \overline \Theta(\overline P_{\tau + 1},\cQ_{r^{\dagger \tau}})\|^2_F - \|\Theta(\overline P_{\tau+1},\cQ_{r^{\dagger \tau}}) - \overline \Theta(\overline P_{\tau },\cQ_{r^{\dagger \tau}})\|^2_F \\
	                         & \geq \frac{1}{4\cdot 10^7\log^6(\tfrac{nd}{\delta\zeta_-})}\| M(\overline P_{\tau}) - \overline M(\overline P_{\tau}) \|^2_F \enspace .
\end{align*}
Hence,
\begin{eqnarray*}
	\| M(\overline P_{\tau_{0}})\|_F^2 &\geq&   \| M(\overline P_{\tau_{0}}) - \overline M(\overline P_{\tau_{0}}) \|^2_F  \\  &\geq&
	\| M(\overline P_{\tau_{\infty}}) - \overline M(\overline P_{\tau_{\infty}}) \|^2_F \left(1-  \frac{1}{4\cdot 10^7\log^6(\frac{nd}{\delta\zeta_-})}\right)^{-\tau_{\infty}} \enspace .
\end{eqnarray*}
Since~\eqref{eq:upper_p_dagg} does not hold at $\tau=\tau_\infty$, this implies that the Frobenius norm in the right-hand side of the above inequality is larger than $2\zeta^2$ and, in light of the definition of $\tau_{\infty}=4\cdot 10^7\log^7(\tfrac{nd}{\delta(\zeta_-)^2}) $, the right-hand side is larger than $2nd$. This contradicts the fact that  $\| M(\overline P_{\tau_{0}})\|_F^2\leq nd$ since the entries of $M$ lie in $[0,1]$.

\begin{proof}[Proof of Corollary~\ref{cor:rate_analysis_of_block_sorting}]

To ease the notation in this proof, we simply write $P$ for $\overline{P}^{\dagger}$, $r$ for $r^{\dagger}$, $\overline Q^*_{\cp}$ for $\overline Q^{\dagger}_{\cp}$, and $h$ for $h^{\dagger}$.
	Since $\overline Q^*_{\cp}$ corresponds to a set of blocks of questions of size $r$, it follows that $|\overline Q^*_{\cp}| \leq d/r$.  This, in turn, implies that $\sqrt{r}h|P|\sqrt{|\overline{Q}^*_{\cp}|}\leq |P|h\sqrt{d}$ and $\sqrt{|P||\overline Q^*_{\cp}|}\leq \sqrt{|P|d}$. We have proven that
	\beq\label{eq:1_d8}
		\frac{\sqrt{r}h|P|\sqrt{|\overline Q^*_{\cp}|}}{\zeta}\land\left( \sqrt{|P||\overline Q^*_{\cp}|}+ |P|\right)\leq \frac{|G|h\sqrt{d}}{\zeta}\land (\sqrt{|G|d} + |G|)\enspace .
	\eeq
	Second, we know from~\Cref{lem:constraints_on_CP_blocks} that $|\overline Q^*_{\cp}| \leq 64 \frac{\tilde r }{rh}$ so that
	\beq\label{eq:upper_intermediaire}
	\frac{\sqrt{r}h|P|\sqrt{|\overline Q^*_{\cp}|}}{\zeta}\land\left( \sqrt{|P||\overline Q^*_{\cp}|}+ |P|\right)\lesssim  \frac{|P|\sqrt{\tilde r h}}{\zeta}\land \left(\sqrt{|P|\frac{\tilde r }{rh}}+ |P|\right)\ .
	\eeq
	If $r > 32\zeta^2 \log(\tfrac{2d}{\delta})\frac{1}{|P| h^2}$ then, it follows from the definition~\eqref{eq:def_tilde_r} of $\tilde{r}$ that $\tilde r= 8r$ so that the right-hand side of~\eqref{eq:upper_intermediaire} is at most of the order of $\sqrt{|P|/h}+|P|$. For a smaller $r$, we know from~\eqref{eq:def_tilde_r} that $\tilde{r}\lesssim \zeta^2 \log(\tfrac{2d}{\delta})/(|P| h^2)$, which in turn implies that
	\[
		\frac{|P|\sqrt{\tilde r h}}{\zeta}\lesssim  \sqrt{\log(\tfrac{2d}{\delta})} \sqrt{\frac{|P|}{h}}\enspace .
	\]
	Hence, we deduce from~\eqref{eq:upper_intermediaire} that 
	\[
		\frac{\sqrt{r}h|P|\sqrt{|\overline Q^*_{\cp}|}}{\zeta}\land\left( \sqrt{|P||\overline Q^*_{\cp}|}+ |P|\right)\lesssim \sqrt{\log(\tfrac{2d}{\delta})} \left(\sqrt{\frac{|P|}{h}} + |P|\right) \enspace .
	\]
	Together with~\eqref{eq:1_d8}, this leads us to
	\[
		\frac{\sqrt{r}h|P|\sqrt{|\overline Q^*_{\cp}|}}{\zeta}\land\left( \sqrt{|P||\overline Q^*_{\cp}|}+ |P|\right)\lesssim  \sqrt{\log(\tfrac{2d}{\delta})}\left[ \frac{|G|h\sqrt{d}}{\zeta}\land |G|\sqrt{d} \land \sqrt{\frac{|G|}{h}} + |G| \right] \enspace ,
	\]
	which, together with \eqref{eq:block_sorting_1} concludes the proof.

\end{proof}

\section{Proof of Proposition~\ref{prop:rate_of_trisection_WEI}}\label{sec:proof_trisection_wei}
In this section, we prove \Cref{prop:rate_of_trisection_WEI} which states a tighter bound than \Cref{prop:rate_of_trisection} on $\algoTreeSorting$ when we use the variant $\algoDimensionReductionWM$ to compute $\hat{\pi}_{WM}$. Recall that, for any $t=1, \dots, t_{\infty}$, $\cT_t$ stands for the hierarchical sorting tree built by $\algoTreeSorting$ at the beginning of step $t$. Thus,  $\cT_t$ has depth $t$. 

The main difference with the analysis of \Cref{prop:rate_of_trisection} lies in the analysis of the algorithm $\algoDimensionReductionWM$, which is the purpose of the next subsection. Then, we combine it with the general scheme of the proof of \Cref{prop:rate_of_trisection_WEI} to get the desired bound.

\subsection{Analysis  of $\algoDimensionReductionWM$}\label{subsec:dimension_reduction_wm}

The key idea of $\algoDimensionReductionWM$ is to examine the high-variation regions of the observations not only in a set of experts $\overline P$ but also in the neighboring sets of experts. For this reason, we remind the reader of the notation of $\algoDimensionReductionWM$. Through this subsection, we fix the step $t \geq 0$ of $\algoTreeSorting$. For simplicity, we write $\cT := \cT_t$. Recall that  $\cL^{(\0, \1)}(\cT)$ stands for the set of leaves of $\cT$ of type $\0$ or $\1$. By definition, those leaves are all at depth $t$. Let us focus on a specific leaf $G\in  \cL^{(\0, \1)}(\cT)$, and we consider a subset $\overline P$ of $G$. 

Finally, we recall that we consider an ordering of the leaves $\cL^{(\0, \1)}(\cT)$ at depth $t$ and centered on $G$ as:
$$( G^{(a)})_{a \in \bbZ} = \algoOrder(\cT, G) \enspace ,$$
where $G^{(0)}=G$.

Also, we fix any $h \in \cH$ and $r \in \cR$. As in $\algoDimensionReductionWM$, define
\begin{equation}\label{eq:def_r0_wm}
	r_0 = 2^{9}\log(4d|\cR|/\delta)\frac{\zeta^2}{|\overline P|h^2} \spaceAnd \tilde r =  4 (\lceil r_0\rceil^{dya}  \lor r)\enspace ,
\end{equation}
where $\lceil r_0\rceil^{dya} =2^{\lceil \log_2(r_0)\rceil}$ is the smaller power of $2$ which is larger than $r_0$. Up to numerical constants, $\tilde{r}$ is defined as for the original procedure $\algoDimensionReduction$. If $r\geq r_0$, then we can simply rely on CUSUM statistics at the scale $8r$ and on the set $\overline P$ to detect high variation regions in $\overline P$. If $h$ (or $|\overline P|$) is so small that $r_0> r$, we applied the CUSUM statistic at a larger scale $\tilde{r}$ in $\algoDimensionReduction$. In this version, we compute the CUSUM statistics at a scale smaller than $\tilde{r}$ to the price of considering more experts than those in $\overline P$.

If $r < \lceil r_0\rceil^{dya}$, let us consider any $r_{\cp} \in [8r, 2\tilde r]\cap \cR$. We respectively define
\begin{align*}
	a^{+}_{WM} := a^{+}_{WM}(\cT, G, h, r_{\cp}) & = \min \left\{ a ~:~ |G^{(1)}| + \dots + |G^{(a)}| \geq 2^{11}\log(4d|\cR|/\delta)\frac{\zeta^2}{r_{\cp}h^2} \right\}  \enspace ;            \\
	a^{-}_{WM} := a^{-}_{WM}(\cT, G, h, r_{\cp}) & = \min \left\{ a ~:~ |G^{(-1)}| + \dots + |G^{(-a)}| \geq 2^{11}\log(4d|\cR|/\delta)\frac{\zeta^2}{r_{\cp}h^2} \right\} \enspace ,
\end{align*}
as the minimum number of groups above and below $G$ in such a way that there are enough experts to detect a $h$-variation in the mean at the scale $r_{\cp}$. Then, $\cV^{+}_{r_{\cp}}$ and  $\cV^{-}_{r_{\cp}}$ stand for the collection of experts in the corresponding groups:
\begin{align}\label{eq:vplus}
	\cV^{+}_{r_{\cp}}:=\cV^{+}_{r_{\cp}}(\cT, \overline P, h)   = \bigcup_{a = 1}^{a^{+}_{WM}}G^{(a)} \spaceAnd  \cV^{-}_{r_{\cp}}:=\cV^{-}_{r_{\cp}}(\cT, \overline P, h)  = \bigcup_{a = - a^{-}_{WM}}^{ - 1}G^{(a)} \enspace ,
\end{align}
Finally, we define
\begin{align}\label{eq:v}
	\cV_{r_{\cp}}:=\cV_{r_{\cp}}(\cT,\overline P, h)   = \left\{ \begin{array}{ll}
		                                                   \cV^{+}_{r_{\cp}} \cup \cV^{-}_{r_{\cp}} & \mbox{ if } r_{\cp} \leq  \tilde r \\
		                                                   \overline P                                        & \mbox{ if } r_{\cp} = 2\tilde r
	                                                   \end{array}\right.
\end{align}
which  exactly corresponds to the definition at \Cref{line:awm} and \Cref{line:Gres} of $\algoDimensionReductionWM$. For any $k\in [d]$, we recall here the definition of the statistic $\widehat \bDelta^{(\ext)}_{k, r_{\cp}}$ its deterministic counterpart:
$$\widehat \bDelta^{(\ext)}_{k, r_{\cp}} = \frac{1}{2r_{\cp}}\sum_{k' = k - r_{\cp}}^{k + r_{\cp} - 1} \overline y_{k'}(\cV^+_{r_{\cp}}) - \overline y_{k'}(\cV^-_{r_{\cp}}) \spaceAnd \bDelta^{*(\ext)}_{k, r_{\cp}} = \frac{1}{2r_{\cp}}\sum_{k' = k - r_{\cp}}^{k + r_{\cp} - 1} \overline m_{k'}(\cV^+_{r_{\cp}}) - \overline m_{k'}(\cV^-_{r_{\cp}})\enspace .$$
In the notation of $\widehat \bDelta^{(\ext)}_{k, r_{\cp}}$, we remove the dependency on $\cV^+_{r_{\cp}}$ and $\cV^-_{r_{\cp}}$ to simplify the notation. Here, $\widehat \bDelta^{(\ext)}_{k, r_{\cp}}$ stands for the width between the empirical means of the  groups above $\overline P$ and below $\overline P$. Recall also the definition of the statistic $\widehat \bC^{(\ext)}_{k, 2r_{\cp}}$ and introduce its deterministic counterpart:
\begin{align*}
	\widehat \bC^{(\ext)}_{k, 2r_{\cp}} & = \frac{1}{2r_{\cp}}\left(\sum_{k' = k}^{k + 2r_{\cp} - 1} \overline y_{k'}(\cV_{2r_{\cp}}) - \sum_{k' = k - 2r_{\cp}}^{k - 1}\overline y_{k'}(\cV_{2r_{\cp}})\right)            \\
	\bC^{*(\ext)}_{k, 2r_{\cp}}         & = \frac{1}{2r_{\cp}}\left(\sum_{k' = k}^{k + 2r_{\cp} - 1} \overline m_{k'}(\cV_{2r_{\cp}}) - \sum_{k' = k - 2r_{\cp}}^{k - 1}\overline m_{k'}(\cV_{2r_{\cp}})\right) \enspace .
\end{align*}
Here, $\widehat \bC^{(\ext)}_{k, 2r_{\cp}}$ stands for the mean CUSUM statistic over the experts in $\cV_{2r_{\cp}}$. 
Consider any $r_{\cp}\in  [4r, \tilde r]\cap \cR$. Then, as in the algorithm $\algoDimensionReductionWM$, we define the collection of positions where both the width and the CUSUM statistic are large:
$$\widehat D_{WM} := \widehat D_{WM}(\cT, \overline P, h, r, r_{\cp}) = \left\{k ~:~ \widehat \bDelta^{(\ext)}_{k, r_{\cp}} \geq \frac{h}{16} ~ \mathrm{and} ~ \widehat \bC^{(\ext)}_{k, 2r_{\cp}} \geq \frac{h}{16} \right\} \enspace .$$
See Figure~\ref{fig:wm_cp} for illustrations. Then, we define $D^*_{WM}$ and  $\overline D^*_{WM}$ as the population counterparts of $\widehat D_{WM}$ with different constants
\begin{align*}
	D^*_{WM}(\cT, \overline P, h, r, r_{\cp})           & = \left\{ k \in 1, \dots, d ~:~ \bC^{*(\ext)}_{k, 2r_{\cp}}  \geq \frac{h}{8} ~ \mathrm{and} ~ \bDelta^{*(\ext)}_{k, r_{\cp}} \geq \frac{h}{8}  \right\} \enspace ;  \\
	\overline D^*_{WM}(\cT, \overline P, h,r, r_{\cp}) & = \left\{ k \in 1, \dots, d ~:~ \bC^{*(\ext)}_{k, 2r_{\cp}} \geq \frac{h}{32} ~ \mathrm{and} ~ \bDelta^{*(\ext)}_{k, r_{\cp}} \geq \frac{h}{32} \right\}\enspace . 
\end{align*}
Then, we consider the collections of blocks $\widehat Q_{WM}(\cT, \overline P, h, r, r_{\cp})$, $Q^*_{WM}(\cT, \overline P, h, r, r_{\cp})$, and $\widehat Q_{WM}(\cT, \overline P, h, r, r_{\cp})$ of size $r$. With our notation, this means that 
$Q^*_{WM}(\cT, \overline P, h, r, r_{\cp}) = \algoEncodeSet(D^*_{WM},r)$,
$\overline Q^*_{WM}(\cT, \overline P, h, r, r_{\cp}) = \algoEncodeSet(\overline D^*_{WM},r)$, and $\widehat Q_{WM}(\cT, \overline P, h, r, r_{\cp}) = \algoEncodeSet(\widehat D_{WM},r)$.
Finally, we consider the unions over all possible  $r_{\cp} \in \cR$ with  $4r \leq r_{\cp} \leq \tilde r$:
\begin{align*}
	\widehat Q_{WM} & := \widehat Q_{WM}(\cT, \overline P, h, r) = \bigcup_{r_{\cp} = 4r}^{\tilde r} \widehat Q_{WM}(\cT, \overline P, h, r, r_{\cp})\enspace ;  \\
	Q^*_{WM}        & := Q^*_{WM}(\cT, \overline P, h, r) = \bigcup_{r_{\cp} = 4r}^{\tilde r} Q^*_{WM}(\cT, \overline P, h, r, r_{\cp}) \enspace ; \\
	\overline Q^*_{WM} &:= \overline Q^*_{WM}(\cT, \overline P, h, r) = \bigcup_{r_{\cp} = 4r}^{\tilde r} \overline Q^*_{WM}(\cT, \overline P, h, r, r_{\cp})\enspace .
\end{align*}

The following lemma states that, with high probability, $\widehat Q_{WM}$ is sandwiched between $Q^*_{WM}$ and $\overline Q^*_{WM}$, so that, on the corresponding event, it is sufficient to study these two quantities.

\begin{lemma}\label{lem:analysis_of_residual_dim}
	Consider any valid hierarchical sorting tree $\cT$, any subset $\overline P$ of a leaf $G$ of $\cT$, any  $h \in \cH$, and any  $r \in \cR$. With probability at least $1- \delta$, it holds that 
	\begin{equation}\label{eq:blocks_residu}
		Q^*_{WM} \subset \widehat Q_{WM} \subset \overline Q^*_{WM} \enspace .
	\end{equation}
\end{lemma}

Next, we show that the aggregation of $M(\overline P)$ at $Q^*_{WM}$ captures most of the variance of $M(\overline P)$.
\begin{lemma}\label{lem:restriction_to_residual_dim} Assume that $\cT$ is a valid hierarchical sorting tree. 
	Then, there exist $h\in \cH$ and  $r\in \cR$ such that
	\begin{equation}\label{eq:restriction_to_residual_dim}
		\| M(\overline P) - \overline M(\overline P) \|^2_F \leq 16\zeta^2  + 96|\cR||\cH|\left\| \left[\Theta(\overline P, Q^*_{WM}) - \overline \Theta(\overline P,Q^*_{WM})\right]_{\sqrt{r}h}\right\|^2_F \enspace . 
	\end{equation}
\end{lemma}
Recall that $\cT_{t_{\infty}}$ (and in particular also $\cT=\cT_{t}$) is a valid hierarchical sorting tree under the event $\xi$ of high probability defined in \Cref{cor:tree_sorting}. 
This lemma is the counterpart of~\Cref{lem:restriction_to_CP} for the oblivious $\algoDimensionReduction$ algorithm.

\subsection{Analysis of the variant $\algoBlockSorting$ with $\algoDoubleTrisectionWEI$}
Recall the definition~\eqref{eq:definition_psi} of the function $\Psi$ by $\Psi(p,r,h,q)= \frac{hp\sqrt{rq}}{\zeta}\land \sqrt{pq}+p$.
In \Cref{prop:rate_analysis_of_block_sorting}, we stated a high probability control for the result of $\algoBlockSorting$ when fed with $\algoDimensionReduction$. In particular, this proposition only used the properties of  $\algoDimensionReduction$ stated in Lemmas~\ref{lem:analysis_of_findCP} and~\ref{lem:restriction_to_CP}. As we have proven in Lemmas~\ref{lem:analysis_of_residual_dim} and~\ref{lem:restriction_to_residual_dim} (their counterparts for $\algoDimensionReductionWM$), we readily obtain the following result whose proof is omitted.

\begin{proposition}\label{prop:rate_analysis_of_block_sorting_wm} Assume that $\cT_t$ is a valid hierarchical sorting tree. 
	Consider a leaf $G$ of $\cT$ of type $\0$ or $\1$ at depth $t$. 
	With probability higher than $1- 5\tau_{\infty}\delta$, there exists a subset $\overline P^{\dagger}$ such that  $\overline{P}\subseteq \overline P^{\dagger}\subseteq G$ and the following property  holds. For some $r_{\cp}^{\dagger} \geq r^{\dagger} \in \cR$ and some $h^{\dagger}\in \cH$, upon writing
	$Q^{\dagger}_{WM} = Q^*_{WM}$ and $\overline Q^{\dagger}_{WM} = \overline Q^*_{WM}$, we have simultaneously
	\begin{align}\label{eq:block_sorting_wm_1}
		\|\left[\Theta(\overline P^{\dagger}, Q^{\dagger}_{WM}) - \overline \Theta(\overline P^{\dagger}, Q^{\dagger}_{WM})\right]_{\sqrt{r^{\dagger}}h^{\dagger}}\|^2_F
		& \leq  4\cdot 10^5 \zeta^2\log^{3}\left(\frac{6nd}{\delta \zeta_-}\right) \Psi(|\overline P^{\dagger}|,r^{\dagger},h^{\dagger},|\overline Q^{\dagger}_{WM}|) \enspace  ;                     \\
		\label{eq:block_sorting_wm_2}
		\|M(\overline P^{\dagger})-\overline{M}(\overline P^{\dagger})\|^2_F & \leq 16\zeta^2 + 96|\cR||\cH|\|[\Theta(\overline P^{\dagger}, Q^{\dagger}_{WM}) - \overline \Theta(\overline P^{\dagger}, Q^{\dagger}_{WM})]_{\sqrt{r^{\dagger}}h^{\dagger}}\|^2_F \enspace .
	\end{align}
\end{proposition}

Since $\|M(\overline P)-\overline{M}(\overline P)\|^2_F \leq \|M(\overline P^{\dagger})-\overline{M}(\overline P^{\dagger})\|^2_F$, the above proposition controls $\|M(\overline P)-\overline{M}(\overline P)\|^2_F$ in terms of $\Psi(|\overline P^{\dagger}|,r^{\dagger},h^{\dagger},|\overline Q^{\dagger}_{WM}|)$.

\subsection{Analysis of the complete procedure $\algoTreeSorting$ with $\algoDoubleTrisectionWEI$}

In light of~\Cref{prop:rate_analysis_of_block_sorting_wm}, we need to control the cardinality of $|\overline Q^{\dagger}_{WM}|$. In comparison to the oblivious procedure analyzed in the previous section, the main improvement here is that the typical cardinalities $|\overline Q^{\dagger}_{WM}|$ are smaller than $|\overline Q^{\dagger}_{\cp}|$ thanks to the refined dimension reduction procedure $\algoDimensionReductionWM$.

Unfortunately, it is not possible to get a tight control of the cardinality of each $\overline Q^{\dagger}_{WM}$ individually. Still, we are able to show that among all groups $G\in \cL^{(\0,\1)}(\cT_t)$ that are refined in the $t$-th iteration of $\algoTreeSorting$, many of them will correspond to small $|\overline Q^{\dagger}_{WM}|$. To formalize this argument, we need to be careful about the dependencies of the quantities under consideration.

We start from the ordered collection $\cL^{(\0,\1)}(\cT_t)$ of $v\leq 2^{t}$ leaves of types  $\0$ or $\1$. 
We write $G_1,\ldots, G_v$ for these  groups and we are given a collection $\overline{P}_1$,\ldots, $\overline{P}_v$ of subgroups such that $\overline{P}_i\subset G_i$ for $i=1,\ldots, v$. Later, we will specify $\overline{P}_v= \overline{P}^{\dagger}_v$, but those sets can be considered arbitrarily. 

For a specific group $\overline{P}_v\subset G_v$, we write $\overline Q^{\dagger}_{WM}(\overline{P}_v,h,r)$ instead of $\overline Q^{\dagger}_{WM}$ to emphasize its dependency on $\overline{P}_v$, $r$ and $h$. Given a positive integer $p>0$, we define $\cP^*(p)=\{\overline{P}_v:\ |\overline{P}_v|\in [p,2p)\}$ the collection of groups $\overline{P}_v$ of size in $[p,2p)$.

\begin{lemma}\label{lem:constraints_on_residual_dim}
Assume that $\cT_t$ is a valid hierarchical sorting tree. 
For any $h\in \cH$, $r\in R$, any integer $p$, any sequence $\overline{P}_v$ of subsets of $G_v$, it holds that 
	\begin{equation}
		\sum_{\overline P \in \cP^*(p)} |\overline Q^*_{WM}(\overline P,h,r)| \lesssim 
		\log(d)	\left[  \frac{\sqrt{nd(r_0\vee r)}}{rh\sqrt{p}} \wedge\frac{d(r_0\vee r)}{r^2 h}\wedge  \frac{nd}{pr} \wedge \frac{n(r_0\vee r)}{p rh} \right]\ .
	\end{equation}
\end{lemma}
We are now equipped to prove \Cref{prop:rate_of_trisection_WEI}.

\begin{proof}[Proof of Proposition~\ref{prop:rate_of_trisection_WEI}]
We work under the event $\xi$ (\Cref{cor:tree_sorting}) ensuring $\cT_{t_{\infty}}$ and in particular $\cT_t$ is a valid hierarchical sorting tree. 
For each group $G_s\in \cL^{(\0,\1)}(\cT_t)$ we apply \Cref{prop:rate_analysis_of_block_sorting_wm} and define a corresponding subgroup $\overline{P}^{\dagger}$, with $r^{\dagger} \in \cR$, $h^{\dagger}\in \cH$ and a corresponding collection of blocks $\overline Q^{\dagger}_{WM}$. 
Define the collection $\cD_n= \{1,2,4,\ldots, 2^{\lceil \log_2(n)\rceil }\}$. 
For $p\in \cD_n$, we define $\cP^*(p,h,r)$ as the collection of groups $\overline{P}^{\dagger}$ satisfying $|\overline{P}^{\dagger}|\in [p,2p)$, $h^{\dagger}=h$, and $r^{\dagger}=r$.

Then, we derive from~\Cref{prop:rate_analysis_of_block_sorting_wm} that, on an additional event of probability higher 
than $1-5\cdot2^{t}\tau_{\infty}\delta$, we have 
\begin{align}\nonumber
\lefteqn{	\sum_{\overline{P} \in \overline{\cL}_{t}}\|M(\overline{P}) - \overline M(\overline{P})\|_F^2 } &\\ & \nonumber \leq 16\zeta^2|\overline{\cL}_t| + 96 |\cR||\cH|\sum_{p,h,r}\sum_{\overline P^{\dagger} \in \cP^*(p,h,r)}\|\left[\Theta(\overline P^{\dagger}, \overline Q^{\dagger}_{WM}) - \overline \Theta(\overline P^{\dagger}, \overline Q^{\dagger}_{WM})\right]_{\sqrt{r}h}\|_F^2 \\ \nonumber
& \stackrel{(a)}{\lesssim}   \zeta^2\log^{5}\left(\frac{6nd}{\delta  \zeta_-}\right)\sum_{p,h,r}\sum_{\overline P^{\dagger} \in \cP^*(p,h,r)}\left[ \sqrt{[\frac{h^2pr}{\zeta^2}\wedge 1] p|\overline Q^*_{WM}(\overline P^{\dagger},h,r)|} +p\right]                                                \\ \nonumber
& \stackrel{(b)}{\lesssim} \zeta^2 \log^{5.5}\left(\frac{6nd}{\delta  \zeta_-}\right)  \sum_{p,h,r}
	\left[ \sqrt{ \frac{r}{r\vee r_0} p|\cP^*(p,h,r)|\sum_{\overline P^{\dagger} \in \cP^*(p,h,r)}|\overline Q^*_{WM}(\overline{P}^{\dagger},h,r)\big|}+ 	p|\cP^*(p,h,r)|\right]\\ \nonumber
	& \stackrel{(c)}{\lesssim}  \zeta^2 \log^{5.5}\left(\frac{6nd}{\delta  \zeta_-}\right)  \sum_{p,h,r}
	\left[ \sqrt{ \frac{n r}{r\vee r_0}  \sum_{\overline P^{\dagger} \in \cP^*(p,h,r)}|\overline Q^*_{WM}(\overline{P}^{\dagger},h,r)\big|}+ 	n\right]\\  
 \nonumber 
	& \stackrel{(d)}{\lesssim} \zeta^2 \log^{6}\left(\frac{6nd}{\delta  \zeta_-}\right)\sum_{p,h,r}\left[\left(n^{3/4}d^{1/4}\left(\frac{1}{(r_0\vee r)ph^2}\right)^{1/4} \wedge  n\sqrt{\frac{d}{p(r_0\vee r)}} \wedge \frac{n}{\sqrt{ph}} \wedge\sqrt{\frac{nd}{rh}}\right) + n\right]\\	\nonumber
 & \lesssim  \zeta^2 \log^{7}\left(\frac{6nd}{\delta  \zeta_-}\right)\sum_{p,h}\left[\left(\frac{n^{3/4}d^{1/4}}{\zeta^{1/2}}\wedge  n\sqrt{d}\wedge \frac{nh}{\zeta}\sqrt{d} \wedge \frac{n}{\sqrt{h}}\wedge \frac{\sqrt{nd}}{\sqrt{h}} \right)+ n\right]\\ \nonumber
 & \stackrel{(e)}{\lesssim} \zeta^2 \log^{7}\left(\frac{6nd}{\delta  \zeta_-}\right)\sum_{p,h}\left[\left(\frac{n^{3/4}d^{1/4}}{\zeta^{1/2}}\wedge  n\sqrt{d}\wedge \frac{n^{2/3}\sqrt{d}}{\zeta^{1/3}} \wedge \frac{nd^{1/6}}{\zeta^{1/3}} \right)+ n\right]\\ \nonumber
 &\lesssim \zeta^2 \log^{9}\left(\frac{6nd}{\delta \zeta_-}\right)\left[\left(\frac{n^{3/4}d^{1/4}}{\zeta^{1/2}}\wedge  n\sqrt{d}\wedge \frac{n^{2/3}\sqrt{d}}{\zeta^{1/3}} \wedge \frac{nd^{1/6}}{\zeta^{1/3}} \right)+ n\right]\ ,
\end{align}
where we applied \Cref{prop:rate_analysis_of_block_sorting_wm} in (a), Jensen inequality and the definition of $r_0$ in (b),  as well as the bound $|\cP^*(p,h,r)|\leq n/p$ in (c), \Cref{lem:constraints_on_residual_dim} in (d), and $x\wedge y \leq x^{1/3}y^{2/3}$ in (e).

\end{proof}

\subsection{Remaining proofs}

\begin{proof}[Proof of \cref{lem:analysis_of_residual_dim}]
	It is sufficient to prove that with high probability, $ D^*_{WM}(\cT, P, h, r_{\cp}) \subset \widehat D_{WM}(\cT, P, h, r_{\cp})\subset \overline{D}^*_{WM}(\cT, P, h, r_{\cp})$ for all $r_{\cp} \in [4r, \tilde r] \cap\cR$. Recall that we use the convention that $\overline{y}_i=\overline{m}_i=0$ if $i\leq 0$ and $\overline{y}_i=\overline{m}_i=1$ if $i>d$.
	Since the CUSUM and the envelope statistics are linear, we have the decompositions
	\begin{align*}
		\widehat{\bC}^{(\ext)}_{k, 2r_{\cp}}(\cV_{2r_{\cp}}) & =     \bC^{*(\ext)}_{k, 2r_{\cp}}(\cV_{2r_{\cp}}) + \frac{1}{2r_{\cp}}\left(\sum_{k' = k}^{k+ 2r_{\cp} - 1}\overline e_{k'}(\cV_{2r_{\cp}}) -  \sum_{k' = k- 2r_{\cp}}^{k-1}\overline e_{k'}(\cV_{2r_{\cp}})\right)         \\
		\widehat{\bDelta}^{(\ext)}_{k,r_{\cp}}(\cV^{+}_{r_{\cp}}, \cV^{-}_{r_{\cp}})
		                                                     & = \bDelta^{*(\ext)}_{k,r_{\cp}}(\cV^{+}_{r_{\cp}}, \cV^{-}_{r_{\cp}}) + \frac{1}{2r_{\cp}}\sum_{k' = k - r_{\cp}}^{k+ r_{\cp} - 1}\left(\overline e_{k'}(\cV^{+}_{r_{\cp}}) - \overline e_{k'}(\cV^{-}_{r_{\cp}})\right)\ ,
	\end{align*}
	where the two latter random variables are centered and respectively $\zeta(r_{\cp}|\cV_{2r_{\cp}}|)^{-1/2}$-subGaussian and $\zeta[r_{\cp}(|\cV^{+}_{r_{\cp}}| \land |\cV^{-}_{r_{\cp}}|)]^{-1/2}$-subGaussian. By a union bound, we deduce that, with probability higher than $1-\delta$, we have simultaneously
	\begin{align}\label{eq:0_1}
		\max_{r_{\cp}\in [4r, \tilde r] \cap \cR}\max_{k\in [d]}\big|\widehat{\bC}^{*(\ext)}_{k,2r_{\cp}}- \bC^{*(\ext)}_{k,2r_{\cp}}\big|& \leq \zeta\sqrt{\frac{2}{r_{\cp}|\cV_{2r_{\cp}}|}\log\left(\frac{4d|\cR|}{\delta}\right)} \ ; \\
		\max_{r_{\cp}\in [4r, \tilde r] \cap \cR}\max_{k\in [d]}\big|\widehat{\bDelta}^{*(\ext)}_{k,r_{\cp}}- \bDelta^{*(\ext)}_{k,r_{\cp}}\big|&\leq \zeta\sqrt{\frac{2}{(|\cV^{+}_{r_{\cp}}| \land |\cV^{-}_{r_{\cp}}|)r_{\cp}}\log\left(\frac{4d|\cR|}{\delta}\right)} \ .	\label{eq:0_2}
	\end{align}
	To conclude, it suffices to check that $|\cV^{+}_{r_{\cp}}|$, $|\cV^{-}_{r_{\cp}}|$, and $|\cV_{r_{\cp}}|$ have been chosen large enough so that the right-hand side of the two above equations is at most $h/32$. 
	
	By definition of $\cV^{+}_{r_{\cp}}$ and $\cV^{-}_{r_{\cp}}$, we know that $|\cV^{+}_{r_{\cp}}| \land |\cV^{-}_{r_{\cp}}| \geq 2^{11}\log(4d|\cR|/\delta)\tfrac{\zeta^2}{r_{\cp}h^2}$ which implies that~\eqref{eq:0_2} is at most $h/32$. 
	
	If $r_{\cp}\leq \tilde{r}/2$, then $\cV_{2r_{r_{\cp}}}= |\cV^+_{2r_{r_{\cp}}}|+ |\cV^-_{2r_{r_{\cp}}}|\geq 2^{11}\log(4d|\cR|/\delta)\tfrac{\zeta^2}{r_{\cp}h^2}$, which implies that~\eqref{eq:0_1} is at most $h/32$. Finally, for $r_{\cp}=\tilde{r}$, we use 
	\[r_{\cp}\geq 4r_0\geq 2^{11}\log(4d|\cR|/\delta)\frac{\zeta^2}{|P|h^2}
	\]
and that $|\cV|=|P|$ to conclude that~\eqref{eq:0_1} is at most $h/32$. 
\end{proof}

\begin{proof}[Proof of \Cref{lem:restriction_to_residual_dim}]
	In the analysis of $\algoDimensionReduction$, we introduced in \Cref{eq:Dcp*} the sets $D^*_{\cp}(P, h, r)$ of questions such that the corresponding CUSUM of the mean expert in $P$ is above $h/2$ at scale $8r$. Recall the set $Q^*_{\cp} := Q^*_{\cp}(P, h, r) = \algoEncodeSet (D^*_{\cp}(P, h, r), r)$. In \Cref{lem:restriction_to_CP}, we stated that, for some $h \in \cH$ and $r \in \cR$, we have  
	\begin{equation}\label{eq:energy_cp_simple}
		\| M(P) - \overline M(P) \|^2_F \leq 16\zeta^2 + 96|\cR||\cH|\left\| \left[\Theta(P, Q^*_{\cp}) - \overline \Theta(P,Q^*_{\cp})\right]_{\sqrt{r}h}\right\|^2_F \enspace .
	\end{equation}
	Define $D^*_{\env} := D^*_{\env}(\cT, P, h, r) = \{ k\in [d]~:~ \bDelta^{*(\ext)}_{k, r}(\cV^+_r, \cV^-_r) \geq h/2\}$ for the questions where the population width between $\cV^{+}_r$ and $\cV^{-}_r$ at scale $r$ is at least $h/2$. Besides, we define  $Q^*_{\env} := Q^*_{\env}(\cT, P, h, r) = \algoEncodeSet(D^*_{\env}, r)$.
	If $l \in Q^*_{\cp} \setminus Q^*_{\env}$, then for any $i, j \in P$, we have 
	\[|\Theta_{i, l} - \Theta_{j, l}| \leq \frac{1}{\sqrt{r}}\sum_{k' = lr }^{(l+1)r - 1} \overline m_{k'}(\cV^+_{r}) - \overline m_{k'}(\cV^-_{r})\leq 2\sqrt{r}\bDelta^{(*\ext)}_{lr, r} < \sqrt{r}h\enspace .
	\]
	Hence, it follows that 
	\begin{equation}\label{eq:Qstarcpenv}
		\left\| \left[\Theta(P, Q^*_{\cp}) - \overline \Theta(P,Q^*_{\cp})\right]_{\sqrt{r}h}\right\|^2_F = \left\| \left[\Theta(P, Q^*_{\cp}\cap Q^*_{\env}) - \overline \Theta(P,Q^*_{\cp}\cap Q^*_{\env})\right]_{\sqrt{r}h}\right\|^2_F \enspace .
	\end{equation}
	In light of~\eqref{eq:energy_cp_simple} and~\eqref{eq:Qstarcpenv}, we only have to prove that, for any fixed $\cT$, $P$, $h$, and $r$, we have 
	\begin{equation}\label{eq:Dstarwm}
		D^*_{\cp}(P, h, r) \cap D^*_{\env}(\cT, P, h, r) \subset \bigcup_{\substack{r_{\cp} \in [4r, \tilde r]\cap \cR}}D^*_{WM}(\cT, P, h, r, r_{\cp}) \enspace .
	\end{equation}
Since the remainder of the proof heavily relies on the comparisons between CUSUM statistics for different subsets of experts,  we respectively write $\bC^{*(\ext)}_{k, r}(\cV_r)$ and  $\Delta^{*(\ext)}_{k, r}(\cV^+_r,\cV^-_r)$ instead of $\bC^{*(\ext)}_{k, r}$ and $\Delta^{*(\ext)}_{k, r}$ to better keep track of the dependencies. 
	Fix any question $k \in D^*_{\cp}(P, h, r) \cap D^*_{\env}(\cT, P, h, r) $ and define 
 $$r_{\min} = \max \{r' \in \cR ~:~ \bC^{*(\ext)}_{k, r'}(\cV_{r'}) < h/8  \} \enspace ,$$
with the convention that $\max (\emptyset)=1$. 
 $r_{\min}$ can be interpreted as the largest scale $r'$ in $\cR$ such that the population CUSUM at scale $r'$ applied to $\cV_{r'}$ is smaller than $h/8$. By definition, we have $\cV_{2\tilde{r}}=P$. As a consequence, for any $r'\geq 2\tilde{r}$, we have  $\bC^{*(\ext)}_{k, r'}(\cV_{r'}) = \bC^{*}_{k, r'}(P)\geq \bC^*_{k, 2\tilde{r}}(P)\geq \bC^*_{k, 8r}(P)\geq h/8$ since $k\in  D^*_{\cp}(P, h, r)$ and since $\tilde{r}\geq 8r$ (see~\eqref{eq:def_r0_wm}). This implies that  $r_{\min}\leq \tilde r$. 
 We consider two distinct cases.

 \medskip 

 \noindent 
 {\bf Case 1: $r_{\min} \leq 4r$.} Then, we simply choose $r_{\cp} = 4r$. By definition of $r_{\min}$, we have   $\bC^{*(\ext)}_{k, 2r_{\cp}}(\cV_{2r_{\cp}}) \geq h/8 $. Since $k \in D^*_{\env}(\cT, P, h, r)$,  we can lower bound the envelope statistic as
	$$\bDelta^{*(\ext)}_{k, r_{\cp}}(\cV^{+}_{r_{\cp}}, \cV^{-}_{r_{\cp}}) \geq \frac{1}{4}\bDelta^{*(\ext)}_{k, r}(\cV^{+}_r, \cV^{-}_r) \geq h/8 \enspace .$$
We have proved that $k \in D^*_{WM}(\cT, P, h,r,  r_{\cp})$.

	\medskip 

	\noindent 
	{\bf Case 2: $r_{\min}\in (4r, \tilde{r}]$.} In that case, we choose $r_{\cp} = r_{\min}\geq 8r$ (since $r_{\min}$ is a power of 2). By definition of $r_{\min}$, we have both $\bC^{*(\ext)}_{k, 2r_{\cp}} \geq h/8 $ and $\bC^{*(\ext)}_{k, r_{\cp}} < h/8 $. Since $k \in D^*_{\cp}(P, h, r)$ and $r_{\cp} \geq 8r$, we also deduce by monotonocity that the CUSUM of the mean expert in $P$ at scale $r_{\cp}$ is higher than $h/2$, this is $\bC^*_{k, r_{\cp}} \geq \bC^*_{k, 8r} \geq h/2$ since $k\in D^*_{\cp}(P,h,r)$ -- see \Cref{eq:Dcp*}.

	Remark that, since $r_{\cp} \leq \tilde r$, we have $\cV_{r_{\cp}} = \cV^{+}_{r_{\cp}}\cup\cV^{-}_{r_{\cp}}$. Without loss of generality, we can assume that $|\cV^{+}_{r_{\cp}}| \geq |\cV^{-}_{r_{\cp}}|$. This implies in particular that 
	\begin{align*}
		r_{\cp}\bC^{*(\ext)}_{k,r_{\cp}}(\cV_{r_{\cp}}) & \geq \frac{|\cV^{+}_{r_{\cp}}|}{|\cV_{r_{\cp}}|}\sum_{k' = k}^{k + r_{\cp} - 1}\overline m_k(\cV^ {+}) - \frac{|\cV^{+}_{r_{\cp}}|}{|\cV_{r_{\cp}}|}\sum_{k' = k - r_{\cp}}^{k - 1}\overline m_k(\cV^ {+}) \\
		                                       &\geq \frac{1}{2}\left( \sum_{k' = k}^{k + r_{\cp} - 1}\overline m_k(\cV^ {+}_{r_{\cp}}) - \sum_{k' = k - r_{\cp}}^{k - 1}\overline m_k(\cV^{+}_{r_{\cp}})\right)\enspace . 
	\end{align*}
	Since $\bC^*_{k,r_{\cp}}(P)\geq h/2$ and $\bC^{*(\ext)}_{k,r_{\cp}}(\cV_{r_{\cp}})\leq h/8$, this implies that 
	\begin{align*}
		h/4 \leq \bC^*_{k,r_{\cp}}(P) - 2\bC^{*(\ext)}_{k,r_{\cp}}(\cV_{r_{\cp}}) & \leq \frac{1}{r_{\cp}}\left(\sum_{k' = k}^{k + r_{\cp}-1}\overline m_{k'}(P) - \overline m_{k'}(\cV^{+}_{r_{\cp}})\right) + \frac{1}{r_{\cp}}\left(\sum_{k' = k - r_{\cp}}^{k-1}\overline m_{k'}(\cV^{+}_{r_{\cp}}) - \overline m_{k'}(P)\right) \\
		& \leq \frac{1}{r_{\cp}}\left(\sum_{k' = k - r_{\cp}}^{k-1}\overline m_{k'}(\cV^{+}_{r_{\cp}}) - \overline m_{k'}(P)\right)                                                                                                              \\
		& \leq \frac{1}{r_{\cp}}\left(\sum_{k' = k - r_{\cp}}^{k + r_{\cp} - 1}\overline m_{k'}(\cV^{+}_{r_{\cp}}) - \overline m_{k'}(\cV^{-}_{r_{\cp}})\right)                                                                                            \\
		& = 2\bDelta^{*(\ext)}_{k, r_{\cp}}(\cV^+_{r_{\cp}},\cV^-_{r_{\cp}})\enspace .
	\end{align*}
	Hence, we have proved that $\bDelta^{*(\ext)}_{k, r_{\cp}}(\cV^+_{r_{\cp}},\cV^-_{r_{\cp}})\geq h/8$ and $\bC^{*(\ext)}_{k, 2r_{\cp}} (\cV_{r_{\cp}})\geq h/8$. Thus, $k \in D_{WM}(\cT, P, h, r_{\cp})$. We have shown \Cref{eq:Dstarwm} and the proof is finished. 

\end{proof}

\begin{proof}[Proof of \Cref{lem:constraints_on_residual_dim}]
	We fix $h \in \cH$ and $r \in \cR$. Let us consider a subgroup $\overline{P}\subset G\in \cL^{(\0,\1)}(\cT)$ 
	Recall that the blocks $\overline Q^*_{WM}(\overline P, h, r)= \bigcup_{r_{\cp} = 4r}^{\tilde r} \overline Q^*_{WM}(\overline P, h, r,r_{\cp})$ -- see the definitions in Section~\ref{subsec:dimension_reduction_wm}. Again, we remove the dependency on $\cT$ in $Q^*_{WM}(\overline P, h, r,r_{\cp})$ for the ease of exposition. 
	First, we bound $\sum_{\overline P \in \cP^*(p)} |\overline Q^*_{WM}(\overline P,h,r,r_{\cp})|$ before summing over the range over all possible $r_{\cp}$. 

	\medskip 
	
	Let us consider some  $l \in \overline Q^*_{WM}(\overline P, h, r,r_{\cp})$. By definition, there exists at least one question $k(l) \in [lr, (l+1)r)$ such that we have simultaneously $\overline \bC^{*(\ext)}_{k(l), 2r_{\cp}} \geq h/32$ and $\overline \bDelta^{*(\ext)}_{k(l), r_{\cp}} \geq h/32$. For $l \in \cQ_r \setminus \overline Q^*_{WM}(r_{\cp})$, we simply define $k(l) = l$. 
	We deduce from this definition that
	\begin{equation}\label{eq:first_maj}
		|\overline Q^*_{WM}(\overline{P}, h, r, r_{\cp})|
		\leq \sum_{l \in \cQ_r}\1 \{ \overline \bC^{*(\ext)}_{k(l), 2r_{\cp}} \geq h/32\}\1\{\overline \bDelta^{*(\ext)}_{k(l), r_{\cp}} \geq h/32\} \enspace .
	\end{equation}
This implies that 
	\begin{equation}\label{eq:capacity_cp_first}
		|\overline Q^*_{WM}(\overline{P}, h, r, r_{\cp})|
		\leq \frac{32}{h} \sum_{l \in \cQ_r} \overline \bC^{*(\ext)}_{k(l), 2r_{\cp}}
		\leq 2^8\frac{r_{\cp}}{rh}
	\end{equation}
	where the last inequality comes from the fact that the total variation of $\overline m(\cV_{2r_{\cp}})$ is at most $1$ and that, for any $l \in \cQ_r$, the interval $[k(l) - 2r_{\cp}, k(l) + 2r_{\cp})$ intersects at most $8r_{\cp}/r$ intervals of the form $[k(l') - 2r_{\cp}, k(l') + 2r_{\cp})$ with $l' \in \cQ_r$.

\medskip 
	Let $p$ be an integer and assume that $|\overline{P}|\in  [p,2p)$.
	Let us introduce $\Gamma := \Gamma(p, h, r_{\cp}) = \frac{\tilde r}{r_{\cp}} \geq 1$, where we recall that $\tilde r\geq 4r_0$ is defined by $r_0 = 2^9\log(4d|\cR|/\delta)\frac{ \zeta^2}{ph^2}$ in \Cref{eq:def_r0_wm}. 
	Intuitively, $\Gamma$ would correspond to the number $a_{WM}^+$ and $a_{WM}^-$ of sets of experts above $\overline{P}$ or below $\overline{P}$ that would be considered if those sets were of size $p$. More generally,  $\cV^{+}_{r_{\cp}}(\cT, \overline{P}, h) = \cup_{a=1}^{a_{WM}^+} G^{(a)}$ contains at most $\Gamma$ groups of size at least $p$ among $G^{(1)}, \dots, G^{(a_{WM} - 1)}$ since the total size of the groups $G^{(a)}$ with $a \leq a_{WM} - 1$ must be less than $2^{11}\log(4d|\cR|/\delta)\frac{\zeta^2}{r_{\cp}h^2}$. Thus, we deduce that $\cV^{-}_{r_{\cp}}(\cT, \overline{P}, h)\cup \overline{P} \cup \cV^{+}_{r_{\cp}}(\cT, \overline{P}, h)$ contains at most $2\Gamma + 3$ groups of size at least $p$.

	The following lemma states that the neighbourhoods $\cV^{-}_{r_{\cp}}(\cT, \overline{P}, h)\cup \overline{P} \cup \cV^{+}_{r_{\cp}}(\cT, \overline{P}, h)$ of groups $\overline{P}$ in $\cP^*(p)$ only intersect on a few groups. 
	\begin{lemma}\label{lem:at_most_gamma_groups}
		Consider any group $\overline P \in \cP^*(p)$.  There exists at most $4\Gamma + 3$ groups $\overline P' \in \cP^*(p)$ such that 
		\begin{equation}\label{eq:gamma_intersection}
			\left(\cV^{-}_{r_{\cp}}(\cT, \overline P,h)\cup \overline P \cup \cV^{+}_{r_{\cp}}(\cT,\overline P,h)\right) \cap \left(\cV^{-}_{r_{\cp}}(\cT, \overline P',h)\cup \overline P' \cup \cV^{+}_{r_{\cp}}(\cT, \overline P',h) \right)  \neq \emptyset \enspace .
		\end{equation}
	\end{lemma}

	As in the proof of \Cref{lem:constraints_on_peeling}, we introduce the width of the matrix $M$ on a set $A$ of experts and an interval of questions  $[k_1, k_2]$ by 
	$$W_{\infty, 1}(M, A, [k_1, k_2]) := \max_{i,j \in A}\sum_{k = k_1}^{k_2}|M_{i,k} - M_{j, k}| \enspace .$$
	From \Cref{eq:first_maj} again, we deduce  that
	\begin{align*}
		\sum_{\overline P \in \cP^*(p)}|\overline Q^*_{WM}(\overline P, h, r, r_{\cp})|
		 & \leq \frac{32}{h} \sum_{\overline P \in \cP^*(p)}\sum_{l \in \cQ_r} \overline \bDelta^{*(\ext)}_{k(l), r_{\cp}}(\cV^{+}_{r_{\cp}}(\cT,\overline P,h), \cV^{-}_{r_{\cp}}(\cT,\overline P,h))                                       \\
		 & \leq \frac{32}{h} \sum_{l \in \cQ_r} \sum_{\overline P \in \cP^*(p)} \frac{1}{2r_{\cp}}W_{\infty, 1}(\cV^{+}_{r_{\cp}}(\cT,\overline P,h) \cup \overline P \cup \cV^{-}_{r_{\cp}}(\cT,\overline P,h), [k - r_{\cp}, k + r_{\cp}]) \\
		 & \leq \frac{32}{rh}(4\Gamma + 3)d \enspace ,
	\end{align*}
	where the last inequality comes \Cref{lem:at_most_gamma_groups} and the fact that
	the sum over disjoints sets $\cV^{-}_{r_{\cp}}(\overline P) \cup \overline P \cup \cV^{+}_{r_{\cp}}(\overline P)$ of $W_{\infty, 1}(\cV^{+}_{r_{\cp}}(\overline P) \cup \overline P \cup \cV^{-}_{r_{\cp}}(\overline P), [k - r_{\cp}, k + r_{\cp}))$ is upper bounded by $2r_{\cp}$ since the total variation of any column of $M$ is at most $1$.

	Combining \Cref{eq:capacity_cp_first} with the latter upper bound together with  $|\cP^*(p)| \leq n/p$ we 
	deduce that 
	\beq\label{eq:inegalite_Q_WM}
		\sum_{\overline P \in \cP^*(p)}|\overline Q^*_{WM}(\overline P, h, r, r_{\cp})|\lesssim \frac{n r_{\cp}}{prh} \land \frac{\Gamma d}{rh} \enspace .
	\eeq
	If $r_0 > r$, then we have $\Gamma \leq \frac{8r_0}{r_{\cp}}$. This implies that
	$$ \sum_{\overline P \in \cP^*(p)}|\overline Q^*_{WM}(\overline P, h, r, r_{\cp})| \lesssim \frac{n r_{\cp}}{prh} \land \frac{ r_0d }{r_{\cp}rh}  \lesssim \frac{\sqrt{ndr_0}}{rh\sqrt{p}} \enspace.$$
	Since $r_{\cp}\leq [4r, \tilde{r}]\cap \cR$, there are at most $c\log(d)$ possible values for $r_{\cp}$, we conclude that 
	$$\sum_{r_{\cp}} \sum_{\overline P \in \cP^*(p)}|\overline Q^*_{WM}(\overline P, h, r, r_{\cp})| \lesssim \frac{n r_0}{prh} \land \frac{ r_0d }{r^2h} \land \frac{\sqrt{ndr_0}}{rh\sqrt{p}} \enspace.$$

	Otherwise, if $r_0 \leq r$, then $\Gamma \leq 8$ and $ r_{\cp}\in [4r,8r]$.  We deduce from~\eqref{eq:inegalite_Q_WM} that 
	$$ \sum_{\overline P \in \cP^*(p)}|\overline Q^*_{WM}(\overline P, h, r, r_{\cp})| \lesssim  \frac{n}{ph} \land \frac{d}{rh}  \leq \frac{\sqrt{nd}}{\sqrt{r}h\sqrt{p}}\enspace .$$
	We have proved that, in any case, 
	\beq\label{eq:upper_1_Q_wm}
		\sum_{r_{\cp}} 	\sum_{\overline P \in \cP^*(p)}|\overline Q^*_{WM}(\overline P, h, r, r_{\cp})| \lesssim \log(d)	\left[  \frac{d(r_0\vee r)}{r^2 h}\wedge \frac{\sqrt{nd(r_0\vee r)}}{rh\sqrt{p}}\right]\ . 
	\eeq
	To establish the remaining bound for the sum of $|\overline Q^*_{WM}(\overline P, h, r, r_{\cp})|$, we control each  $|\overline Q^*_{WM}(\overline P, h, r, r_{\cp})|$ individually in a similar fashion to what we did for the analysis of the oblivious hierarchical sorting estimator $\hat{\pi}_{HT}$. First,  we have $\overline Q^*_{WM}(\overline P, h, r, r_{\cp})\subset \cQ_r$ so that  $|\overline Q^*_{WM}(\overline P, h, r, r_{\cp})|\leq d/r$. Besides, arguing as in the proof of Lemma~\ref{lem:constraints_on_CP_blocks}, $|\overline Q^*_{WM}(\overline P, h, r, r_{\cp})|\lesssim  r_{\cp}/(rh)\lesssim (r_0\vee r)/(rh)$. 
	\beq\label{eq:upper_2_Q_wm}
		\sum_{r_{\cp}} 	\sum_{\overline P \in \cP^*(p)}|\overline Q^*_{WM}(\overline P, h, r, r_{\cp})| \lesssim \log(d)	\left[  \frac{nd}{pr} \wedge  \frac{n(r_0\vee r)}{p r h} \right]\ . 
	\eeq
	Combining~\eqref{eq:upper_1_Q_wm} and~\eqref{eq:upper_2_Q_wm} concludes the proof.
\end{proof}

\begin{proof}[Proof of \Cref{lem:at_most_gamma_groups}]
	Consider two distinct groups  $\overline P$ and $\overline P'$ in $\cP^*(p)$. Let $(G^{(a)}(\overline P))_{a \in \bbZ} = \algoOrder(\cT, \overline P)$ be the ordering of $\cL^{(\0, \1)}(\cT)$ centered on $\overline P$ and $a' \in \bbZ$ the index of the leaf $G^{(a')}(\overline P)$ containing $\overline P'$. Obviously, $|G^{(a')}|\geq |\overline P'|\geq p$. 
	
	Without loss of generality, we  assume that $a' > 0$. 
	In that case, if \Cref{eq:gamma_intersection} is satisfied then necessarily $$(\cV^{+}_{r_{\cp}}(\cT,\overline P,h) \cup \overline P) \cap \cV^{-}_{r_{\cp}}(\cT,\overline P',h) \neq \emptyset \enspace .$$
	This can only happen if the number of leaves $G^{(a)}(\overline P)$ for $0<a< a'$ that are of size at least $p$ is less than or equal to $2\Gamma$. The same holds if $a' < 0$ and this proves the lemma.
\end{proof}

\section{Proof of \Cref{{lem:subsampling_1}} and \Cref{th:second_estimator_WM_poisson_intro}}

\subsection{Proof of \Cref{{lem:subsampling_1}}}

We start with the case $\lambda_-\in [2/d,1]$. The random variable $n_{i,[(j-1)l(\lambda)+1, jl(\lambda)]}$ is distributed as a Poisson random variable with parameter $\lambda l(\lambda)$. Let us apply Chernoff's inequality for Poisson random variable (e.g.~\cite{boucheron2013concentration}, section 2.2). We have 
\[
\bbP\left[n_{i,[(j-1)l(\lambda)+1, jl(\lambda)]} \leq \lambda l(\lambda)/2\right] \leq \exp\left[- \frac{3}{28}\lambda l(\lambda) \right]\leq \frac{\delta}{nd}\ , 
\]
provided that $\lambda l(\lambda)\geq \tfrac{28}{3}\log(nd/\delta)$. Since $\lambda_-\leq 1$, we have  $\lambda l(\lambda)/2\geq \Upsilon^*$. In view of the definition of $\Upsilon^*$, the condition $\lambda l(\lambda)\geq \tfrac{28}{3}\log(nd/\delta)$ is therefore valid and we conclude that 
\[
	\bbP\left[n_{i,[(j-1)l(\lambda)+1, jl(\lambda)]} \leq \Upsilon^*\right] \leq \frac{\delta}{nd}\ , 
\]
and the first result follows. Turning to the second result, we observe that $n_{i,\{j\}}$ is distributed as a Poisson random variable. We apply again Chernoff's inequality to derive that 
\[
	 \bbP\left[n_{i,\{j\}} \leq \frac{\lambda}{2} \right] \leq \exp\left[- \frac{3}{28}\lambda \right]\leq \frac{\delta}{nd}\ , 
\]
since $\lambda\geq \frac{28}{3}\log(nd/\delta)$. Since $\lambda \geq 2\lambda_-\Upsilon^*$, the result follows.

\subsection{Proof of \Cref{th:second_estimator_WM_poisson_intro}}
If $\lambda_-\leq 2/d$, we use the trivial bound $\|M_{\hat \pi_{WMP}^{-1}} - M_{\pi^{*-1}} \|_F^2\leq nd$, which ensures that 
\[
	\E\left[\|M_{\hat \pi_{WMP}^{-1}} - M_{\pi^{*-1}} \|_F^2\right]\leq \frac{n}{\lambda_-}\leq c\log^{c'}\left(\frac{nd\lambda^{1/2}}{\zeta_-}\right) \frac{n}{\lambda}\ . 
\]
If $\lambda_-\geq 1$, then Lemma~\ref{lem:subsampling_1} ensures that, with probability higher than $1-\delta$, we are able to build the $\Upsilon^*$ subsamples and we are in position to apply Theorem~\ref{th:second_estimator_WM} with subGaussian norm $\zeta/\lfloor \lambda_-\rfloor^{1/2}$. Hence, with probability higher than $1-c'\log^{9}(nd\lambda^{1/2}_-/(\delta \zeta_-))\delta$, we have  
\beqn 
\|M_{\hat \pi_{WMP}^{-1}} - M_{\pi^{*-1}} \|_F^2&\leq& c \log^{11}\left(\frac{nd\lfloor\lambda_-\rfloor^{1/2}}{\delta\zeta_-}\right) \cR_{F}(n,d,\zeta\lfloor\lambda_-\rfloor^{-1/2})\\
&\leq & c' \log^{c''}\left(\frac{nd\lambda^{1/2}}{\zeta_-}\right) \cR_{F}(n,d,\zeta\lfloor\lambda\rfloor^{-1/2})\ , 
\eeqn
where we use the definition of $\delta$ and $\lambda_-$ in the last line. On the complementary event, we simply use that $\|M_{\hat \pi_{WMP}^{-1}} - M_{\pi^{*-1}} \|_F^2\leq nd$. Since $\delta$ has been chosen small enough, we can conclude that 
\[
\E[\|M_{\hat \pi_{WMP}^{-1}} - M_{\pi^{*-1}} \|_F^2]\leq c' \log^{c''}\left(\frac{nd\lambda^{1/2}}{\zeta_-}\right) \cR_{F}(n,d,\zeta\lambda^{-1/2})\ . 
\]
\bigskip

It remains to consider the case where $\lambda_-\in [2/d,1]$. Working under the event of probability higher than $1-\delta$ ensured by Lemma~\ref{lem:subsampling_1}, we have $\Upsilon^*$ independent samples $Y^{{\downarrow(0)}},\ldots, Y^{{\downarrow(\Upsilon^*-1)}}$ of size $n\times \lfloor d/l(\lambda)\rfloor$. Define the matrix $M^{{\downarrow}}$ of size $n\times \lfloor d/l(\lambda)\rfloor$ by $M^{{\downarrow}}_{i,j}=M_{i,l(\lambda)(j-1)+1}$. Obviously, $M^{{\downarrow}}_{\pi^{*-1}}$ is a bi-isotonic matrix. Besides, for $s=0,\ldots, \Upsilon^*-1$, $(i,j)\in [n]\times \lfloor d/l(\lambda)\rfloor$, we have the decomposition 
\[
	Y_{ij}^{{\downarrow(s)}}= M^{{\downarrow}(s)}_{ij}+ E^{{\downarrow(s)}}_{ij}\ , 
\]
where $M^{{\downarrow}(s)}_{ij}$ belongs to $[M^{{\downarrow}}_{ij},M^{{\downarrow}}_{ij+1}]$ with the convention $M^{{\downarrow}(s)}_{i,\lfloor d/l(\lambda)\rfloor+1}=1$ and the $E^{{\downarrow(s)}}_{ij}$'s are independent and, for fixed $i$ and $j$, are i.i.d.~distributed and $\zeta$-subGaussian. In fact, the  $M^{{\downarrow}(s)}_{ij}$ are random since $M^{{\downarrow}(s)}_{ij}$ has been sampled uniformly in  $\{M_{i,l(\lambda)(j-1)+1},M_{i,l(\lambda)(j-1)+2},\ldots, M_{i,l(\lambda)(j-1)+l(\lambda)}\}$. Besides, those are correlated with the noise $E^{{\downarrow(s)}}_{ij}$. For the sake of the analysis, it is in fact easier to consider that $M^{{\downarrow}(s)}_{ij}$ has been set by an adversary. Hence, we fall into the semi-random model of Section~\ref{sec:semi_random} and we are in position to apply Theorem~\ref{th:second_estimator_WM-SR} to $\hat{\pi}_{WM-SR}$. With probability at least $1- c'n\log^{9}(\tfrac{nd}{\delta\zeta_-}) \delta$, we have 
\[
	\|M^{\downarrow}_{\hat \pi_{\WM-SR}^{-1}} - M^{\downarrow}_{\pi^{*-1}} \|_F^2 \leq c
	 \log^{11}\left(\frac{2nd}{\delta \zeta_-}\right)\left[\cR_F(n,\lfloor d/l(\lambda)\rfloor,\zeta) + n\right]\enspace ,
\]

Define the matrix $M^{\downarrow\uparrow}$ of size $n\times d$ such that each column is duplicated $l(\lambda)$ times, except the last one which has been duplicated $l(\lambda)-l(\lambda)\lfloor d/l(\lambda)\rfloor$. We readily deduce that 
\beq\label{eq:risk_upper_bound_WM_1}
	\|M^{\downarrow\uparrow}_{\hat \pi_{\WM-SR}^{-1}} - M^{\downarrow\uparrow}_{\pi^{*-1}} \|_F^2 \leq c'l(\lambda)	
	 \log^{11}\left(\frac{2nd}{\delta \zeta_-}\right)\left[\cR_F(n,\lfloor d/l(\lambda)\rfloor,\zeta) + n\right]\enspace ,
\eeq
By triangular inequality, we have
\[
	\|M_{\hat \pi_{\WM-SR}^{-1}} - M_{\pi^{*-1}} \|_F^2\leq 2\|M^{\downarrow\uparrow}_{\hat \pi_{\WM-SR}^{-1}} - M^{\downarrow\uparrow}_{\pi^{*-1}} \|_F^2+ 8 \|M-M^{\downarrow\uparrow}\|_F^2\ .
\]
Thus it remains to upper bound the square Euclidean norm of each row of $M-M^{\downarrow\uparrow}$:
\beqn 
\sum_{j=1}^{d}[M-M^{\downarrow\uparrow}]^2_{i,j}&=& \sum_{k=1}^{\lfloor d/l(\lambda)\rfloor} \sum_{r=1}^{l(\lambda)} [M_{i,(k-1)l(\lambda)+r}- M_{i,(k-1)l(\lambda)+1}]^2\\ & & + \sum_{r=1}^{d-l(\lambda)\lfloor d/l(\lambda)\rfloor } \left[M_{i,(\lfloor d/l(\lambda)\rfloor -1) l(\lambda) + r}  - M_{i,(\lfloor d/l(\lambda)\rfloor - 1) l(\lambda)+1}\right]^2\\
&\leq & \sum_{k=1}^{\lfloor d/l(\lambda)\rfloor} l(\lambda) [M_{i,kl(\lambda)}-M_{i,(k-1)l(\lambda)+1}]^2 + l(\lambda)[M_{i,d}-M_{i,(\lfloor d/l(\lambda)\rfloor - 1)l(\lambda)+1} ]^2\\
&\leq & 2l(\lambda)\ , 
\eeqn
since the total variation of the $i$-th row of $M$ is at most one. Hence, $\|M-M^{\downarrow\uparrow}\|_F^2\leq 2nl(\lambda)$. Together with~\eqref{eq:risk_upper_bound_WM_1}, we conclude that 
\[
	\|M_{\hat \pi_{\WM-SR}^{-1}} - M_{\pi^{*-1}} \|_F^2 \leq c'l(\lambda)
	 \log^{11}\left(\frac{2nd}{\delta \zeta_-}\right)\left[\cR_F(n,\lfloor d/l(\lambda)\rfloor,\zeta) + n\right]\enspace .
\]
with probability at least $1- c'n\log^{9}(\tfrac{nd}{\delta\zeta_-}) \delta$. Since $\delta$ has been chosen small enough and since  $\|M_{\hat \pi_{\WM-SR}^{-1}} - M_{\pi^{*-1}} \|_F^2\leq nd$, we conclude that 
\[
	\E\left[\|M_{\hat \pi_{\WM-SR}^{-1}} - M_{\pi^{*-1}} \|_F^2\right] \leq c'l(\lambda)
	 \log^{11}\left(\frac{2nd}{\zeta_-}\right)\left[\cR_F(n,\lfloor d/l(\lambda)\rfloor,\zeta) + n\right]\enspace . 
\]
Since $l(\lambda)\leq c'\log^{c''}(nd(\lambda\vee 1)/\zeta_-)/\lambda$, we deduce from this bound that 
\beqn 
	\lefteqn{\E\left[\|M_{\hat \pi_{\WM-SR}^{-1}} - M_{\pi^{*-1}} \|_F^2\right] } &&\\ &\leq& c'
	 \log^{c''}\left(\frac{2nd}{\zeta_-}\right)\left[\left(\frac{\zeta}{\sqrt{\lambda}}\right)^2 \left\{ \frac{nd^{1/6}}{(\tfrac{\zeta}{\sqrt{\lambda}})^{1/3}} \wedge \frac{n^{3/4}d^{1/4}}{(\tfrac{\zeta}{\sqrt{\lambda}})^{1/2}}\wedge n\sqrt{d\lambda}\wedge \frac{n^{2/3}\sqrt{d}\lambda^{1/3}}{(\tfrac{\zeta}{\sqrt{\lambda}})^{1/3}}+  n \right\}	
 	 + \frac{n}{\lambda}\right] \\ 
	 &\leq & c'
	 \log^{c''}\left(\frac{2nd}{\zeta_-}\right)\left[\cR_F[n,d,\zeta/\sqrt{\lambda}]+  \frac{n}{\lambda} \right]\ , 
\eeqn 
since $\lambda_-\leq 1$. Again, since $\lambda_-\leq 1$, we have $\lambda\leq c_3 \log^{c_4}\left(nd(\lambda\vee 1)/\zeta\right)$ for some numerical constant $c_3$ and $c_4$. We conclude that
\[
	\E\left[\|M_{\hat \pi_{\WM-SR}^{-1}} - M_{\pi^{*-1}} \|_F^2\right]\leq 
	c'
	 \log^{c''}\left(\frac{2nd}{\zeta_-}\right)\left[\cR_F[n,d,\zeta/\sqrt{\lambda}]+  \frac{n}{\lambda}e^{-\tfrac{\lambda}{c_3 \log^{c_4}\left(nd(\lambda\vee 1)/\zeta\right)} } \right]\ ,
\]
which concludes the proof.

\section{Permutation estimation in the semi-random model}\label{sec:semi_random}

\subsection{Model and Algorithm}

We now consider a slightly different model with $\Upsilon^*$ samples $Y^{(1)},\ldots ,Y^{(\Upsilon^*-1)}$. The noise matrices $E^{(1)}, \ldots, E^{(\Upsilon^*-1)}$ are sampled independently (as previously) and $Y^{(t)}_{ij}=E^{(t)}_{ij}+ M^{(t)}_{ij}$ where $M^{(t)}_{ij}$ is chosen by an adversary in $[M_{ij},M_{i,j+1}]$. This slightly different model is mainly motivated by the analysis of the partial observation scheme in Section~\ref{sec:partial}. In particular, building upon this model and relying on the corresponding modifications in the algorithm allows us to recover the right dependency with respect to $\zeta$ in Section~\ref{sec:partial}.

\medskip

We consider a slight variant $\hat{\pi}_{WM-SR}$ of the estimator $\hat{\pi}_{WM}$ to handle the adversarial differences. The procedure $\hat{\pi}_{WM-SR}$ is computed exactly as $\hat{\pi}_{WM}$ except that
\begin{itemize}
	\item In $\algoPivot$ (Algorithm~\ref{alg:pivot}), the threshold $\beta_{\tris}\sqrt{\log(\frac{2|P|}{\delta})}$ is replaced by $\beta_{\tris}\sqrt{\log(\frac{2|P|}{\delta})}+4\|\omega\|_{\infty}/\|\omega\|_2$ and $\overline{\beta}_{\tris}\sqrt{\log(\frac{2|P|}{\delta})}$ is replaced by $\overline{\beta}_{\tris}\sqrt{\log(\frac{2|P|}{\delta})}+8\|\omega\|_{\infty}/\|\omega\|_2$ 
	\item In $\algoDimensionReductionWM$ (Algorithm~\ref{algo:dimension_reduction_wm}), we respectively replace the definitions of the  CUSUM and empirical width by
		\begin{align}
			\widehat \bDelta^{(\ext)}_{k, r'}(\cV^+, \cV^-) & = \sum_{k' = k - r'}^{k + r' - 1} \overline y_{k'}(\cV^+) - \overline y_{k'-1}(\cV^-)   \ ;\\ \widehat \bC^{(\ext)}_{k, r'}(\cV)           &    = \sum_{k' = k}^{k + r' - 1} \overline y_{k'}(\cV) - \sum_{k' = k - r'-1}^{k - 2}\overline y_{k'}(\cV)\label{eq:stat_cusum_wm-SR}\enspace . 
		\end{align}
\end{itemize}

\begin{theorem}\label{th:second_estimator_WM-SR}
	There exist three numerical constants $c$, $c'$, and $c_0$ such that the following holds. Fix $\delta>0$ and assume that $\Upsilon \geq c_0 \log^{8}\left(nd/(\delta\zeta_-)\right)$.  For any permutation $\pi^*\in \Pi_n$ and any matrix $M$ such that $M_{\pi^{*-1}}\in \mathbb{C}_{\text{BISO}}$, the hierarchical sorting tree estimator with memory $\hat{\pi}_{WM-SR}$ satisfies 
	\beq\label{eq:risk_upper_bound_WM-SR}
	\|M_{\hat \pi_{\WM-SR}^{-1}} - M_{\pi^{*-1}} \|_F^2 \leq c
	 \log^{11}\left(\frac{2nd}{\delta \zeta_-}\right)\left[\cR_F(n,d,\zeta) + n\right]\enspace ,
	\eeq
    with probability at least $1- c'n\log^{9}(\tfrac{nd}{\delta\zeta_-}) \delta$.
\end{theorem}

\subsection{Proof of Theorem~\ref{th:second_estimator_WM-SR}}

The proof follows the main steps as that of Theorem~\ref{th:second_estimator_WM} and we mainly emphasize here the differences. In the proof of~\Cref{th:second_estimator_WM}, we often work with the aggregated model~\eqref{eq:model_aggregated} $Z=\Theta+\NOISEC$ which is restricted to a subset $P$ of experts and a subset $Q\subset \cQ_r$ of questions aggregated at scale $r$ -- see $\algoEncodeMatrix$ for details. For $t=0,\ldots, \Upsilon-1$, the counterpart of~\eqref{eq:model_aggregated} is the following
\beq\label{eq:model_aggregated-SR}
	Z^{(t)} = \Theta^{(t)} + \NOISEC^{(t)} \ ,
\eeq
where the entries of $\NOISEC^{(t)}$ are independent and $\zeta$-subGaussian and $\Theta^{(t)}$ stands for the corresponding aggregation of the matrix $M^{(t)}$. Since the total variation of each row of $M$ is at most one, one readily checks that 
\beq\label{eq:majoration_upper_SR_fundamental}
\sum_{j\in Q}|\Theta^{(t)}_{ij}-\Theta_{ij}|\leq 1\ . 
\eeq

\medskip

Since $\hat{\pi}_{WM-SR}$ is a hierarchical sorting tree estimator, we are in position to control its loss using \Cref{prop:algo_tree_sorting}. For this purpose, we need to prove that~\Cref{cor:trisection} still holds in the semi-random model which, in turn, would imply that   Corollary~\ref{cor:tree_sorting} is true. In fact, the proof of~\Cref{cor:trisection} is verbatim the same except that \Cref{cor:meta_analysis_of_pivot} is replaced by the following lemma.

We remind that $P'= \overline{P}\setminus (\nL\cup \nU)$, and $\overline{P'}= \overline{P}\setminus (\overline{\nL}\cup \overline{\nU})$.
\begin{lemma}\label{cor:meta_analysis_of_pivot-SR}
	For any non-zero vector $w \in \bbR_+^Q$, any pivot $\gamma \in \{1, \dots, |\overline P| \}$
	, we have $\P[\cP_3]\geq 1-\delta$. Besides, on the same event of probability at least $1-\delta$, we have
	\begin{equation}\label{eq:condition_espacement_pivot-SR}
		\abs{\proscal<\Theta_{i, \cdot} - \Theta_{i_{\gamma}, \cdot},\frac{w}{\|w\|_2}>} \leq (2\zeta\sqrt{2}+ \overline{\beta}_{\tris})\sqrt{\Log{\frac{2|\overline P|}{\delta}}} + 10\frac{\|w\|_\infty}{\|w\|_2}  \quad  \text{ if } i\in P'\enspace .
	\end{equation}
\end{lemma}

\begin{proof}[Proof of Lemma~\ref{cor:meta_analysis_of_pivot-SR}]
Consider any sample $t\in [0,\Upsilon-1]$, any vector $w\in \mathbb{R}^q$, and any $i\in \overline{P}$. As a straightforward consequence of~\eqref{eq:majoration_upper_SR_fundamental}, we deduce that 
\beq\label{eq:upper_SR_important}
\Big|	\langle \Theta^{(t)}_{i, \cdot}- \Theta_{i, \cdot}, \frac{w}{\|w\|_2}\rangle  \Big| \leq \frac{\|w\|_{\infty}}{\|w\|_2}\ . 
\eeq
We then deduce from a union bound, that with probability higher than $1-\delta$, we have 
\[
		\left|  \langle Z^{(t)}_{i, \cdot},\frac{w}{\|w\|_2}\rangle   -   \langle \Theta_{i, \cdot},\frac{w}{\|w\|_2}\rangle  \right|\leq \zeta \sqrt{2\log\left(\frac{2|\overline{P}|}{\delta}\right)}+ \frac{\|w\|_{\infty}}{\|w\|_2} \ .
\]
	simultaneously for all $i$ in $\overline{P}$. The rest of the proof of Lemma~\ref{cor:meta_analysis_of_pivot-SR} is left unchanged provided that we replace $\sqrt{2\log\left(\frac{2|\overline{P}|}{\delta}\right)}$ 	by $\sqrt{2\log\left(\frac{2|\overline{P}|}{\delta}\right)}+ \frac{\|w\|_{\infty}}{\|w\|_2}$. 
\end{proof}

Then, being in position to apply Corollary~\ref{cor:tree_sorting}, we state the counterpart of Proposition~\ref{prop:rate_of_trisection_WEI}. 
\begin{proposition}\label{prop:rate_of_trisection_WEI-SR}
	On the intersection of event $\xi$ (defined in \Cref{cor:tree_sorting}) and an event of probability higher than {$1-5\cdot2^{t}\tau_{\infty}\delta$}, it holds  that
	$$\sum_{\overline P \in \overline \cL_{t}}\|M(\overline P) - \overline M(\overline P)\|^2 \lesssim   \log^{9}\left(\frac{6nd}{\delta \zeta_-}\right)\left[\cR_F(n,d,\zeta)+ n\right]\enspace .$$
\end{proposition}
We conclude the proof of Theorem~\ref{th:second_estimator_WM-SR} by combining~\Cref{prop:rate_of_trisection_WEI-SR} with Corollary~\ref{cor:tree_sorting}. Hence, we only need to prove the last proposition.

\subsection{Proof of Proposition~\ref{prop:rate_of_trisection_WEI-SR}}

Again, we only emphasize the differences with the proof of~\Cref{prop:rate_of_trisection_WEI}. We start with the analysis of $\algoDimensionReductionWM$. Recall that we slightly changed the definition of the CUSUM statistics 
\begin{align*}
	\widehat \bC^{(\ext)}_{k, 2r_{\cp}} & = \frac{1}{2r_{\cp}}\left(\sum_{k' = k}^{k + 2r_{\cp} - 1} \overline y_{k'}(\cV_{2r_{\cp}}) - \sum_{k' = k - 2r_{\cp}-1}^{k - 2}\overline y_{k'}(\cV_{2r_{\cp}})\right)           
\end{align*}
by shifting the second sum by one index. The definition of the population CUSUM statistic $\bC^{*(\ext)}_{k, 2r_{\cp}}$ is left unchanged. Similarly, we slightly changed the definition of $\widehat{\bDelta}^{(\ext)}_{k, r_{\cp}}$ to 
$$\widehat \bDelta^{(\ext)}_{k, r_{\cp}} = \frac{1}{2r_{\cp}}\sum_{k' = k - r_{\cp}}^{k + r_{\cp} - 1} \overline y_{k'}(\cV^+_{r_{\cp}}) - \overline y_{k'-1}(\cV^-_{r_{\cp}})\ ,  $$
by shifting again the right hand-side observation by one. With these simple shifts, $\widehat \bDelta^{(\ext)}_{k, r_{\cp}}$ and $\widehat \bC^{(\ext)}_{k, 2r_{\cp}}$ both overestimates $\bDelta^{*(\ext)}_{k, r_{\cp}}$ and $\bC^{*(ext)}_{k, 2r_{\cp}}$ and arguing as in the proof of Lemma~\ref{lem:analysis_of_residual_dim}, we will prove that 	$Q^*_{WM} \subset \widehat Q_{WM}$ with probability at least $1-\delta$ --see Lemma~\ref{lem:analysis_of_residual_dim_SR} below. However, we need to adapt the definition of 
$\overline D^*_{WM}(\cT, P, h,r, r_{\cp})$ to cope with this possible bias. Define
\begin{align}\label{eq:definition_D_WM_SR_cp_h_1}
	\overline D^*_{WM-SR-1}(\cT, P, h,r, r_{\cp}) & = \left\{ k \in 1, \dots, d ~:~ \bC^{*(\ext)}_{k, 2r_{\cp}} \geq \frac{h}{128} ~ \mathrm{and} ~ \bDelta^{*(\ext)}_{k, r_{\cp}} \geq \frac{h}{128} \right\} \ ;  \\
	\overline D^*_{WM-SR-2}(\cT, P, h,r, r_{\cp}) & = \left\{ k \in 1, \dots, d ~:~ \overline m_{k+r_{\cp}}(\cV^+_{r_{\cp}}) - \overline m_{k-r_{\cp}}(\cV^+_{r_{\cp}}) \geq \frac{h r_{\cp}}{128} \right\} 
	\label{eq:definition_D_WM_SR_cp_h_2}\ ; \\
	\overline D^*_{WM-SR-3}(\cT, P, h,r, r_{\cp}) & = \left\{ k \in 1, \dots, d ~:~ \overline m_{k+r_{\cp}-1}(\cV^-_{r_{\cp}}) - \overline m_{k-r_{\cp}-1}(\cV^-_{r_{\cp}}) \geq \frac{h r_{\cp}}{128} \right\} 	\label{eq:definition_D_WM_SR_cp_h_3}\ ; \\
	\overline D^*_{WM-SR-4}(\cT, P, h,r, r_{\cp}) & = \left\{ k \in 1, \dots, d ~:~ \overline m_{k+2r_{\cp}}(\cV_{2r_{\cp}}) - \overline m_{k-2r_{\cp}-1}(\cV_{2r_{\cp}}) \geq \frac{h r_{\cp}}{128} \right\} 	\label{eq:definition_D_WM_SR_cp_h_4}\ .
\end{align}
Then, we define the corresponding subsets $\overline Q^*_{WM-SR-1}$, 
$\overline Q^*_{WM-SR-2}$, $\overline Q^*_{WM-SR-3}$, and $\overline Q^*_{WM-SR-4}$ of $\cQ_r$. For short, we write $\overline{Q}^*_{WM-SR}= \overline Q^*_{WM-SR-1}\cup \overline Q^*_{WM-SR-2}\cup \overline Q^*_{WM-SR-3}\cup \overline Q^*_{WM-SR-4}$. We have the following counterpart of~\Cref{lem:analysis_of_residual_dim}. 

\begin{lemma}\label{lem:analysis_of_residual_dim_SR}
	Consider any valid hierarchical sorting tree $\cT$, any  subset $P$ of a leaf $G$ of $\cT$, any  $h \in \cH$, and any $r \in \cR$. With probability at least $1- \delta$, it holds that 
	\begin{equation}\label{eq:blocks_residu_SR}
		Q^*_{WM} \subset \widehat Q_{WM} \subset \overline Q^*_{WM-SR} \enspace .
	\end{equation}
\end{lemma}
Obviously, Lemma~\ref{lem:restriction_to_residual_dim} is still true since it does not depend on the data generating process. Then, we adapt Propositions~\ref{prop:rate_analysis_of_pivot} and \ref{prop:rate_analysis_of_trisectionACP} to this adversarial setting.

\begin{proposition}\label{prop:rate_analysis_of_pivot-SR}
	Consider any $\overline P\subset [n]$, any $r\in \cR$, and any subset $Q\subset Q_r$. Also, fix any $\eta>0$ and any  $\phi>0$. Provided that
	$$\|\left[\Theta(\overline P,Q) - \overline \Theta(\overline P,Q)\right]_{\eta}\|^2_F \geq \frac{1}{\phi}\|\Theta(\overline P,Q) - \overline \Theta(\overline P,Q)\|_F^2
		\geq 8\eta|\overline P| \left[\phi_{l_1}\sqrt{\log(\tfrac{2|\overline P|}{\delta})}\sqrt{|Q|}+20\right] \enspace,$$
	then, with probability higher than $1-\delta$, we have
	$$\|\Theta(\overline P',Q) - \overline \Theta(\overline P',Q)\|^2_F \leq \left(1 - \frac{1}{16 \phi }\right) \|\Theta(\overline  P,Q) - \overline \Theta(\overline  P,Q)\|^2_F \enspace .$$
\end{proposition}

Recall the definition~\eqref{eq:definition_phi_l_1} of $\phi_{l_1}$. Henceforth, the matrix  $\Theta(\widetilde{P}, Q)$ is said to be indistinguishable in $l_1$-norm if it satisfies
\begin{equation}\label{eq:undistinguishable_in_L1_norm_SR}
	\max_{i,j\in \overline P}  \|\Theta_{i, \cdot}(\widetilde{P},Q) - \Theta_{j, \cdot}(\widetilde P,Q)\|_1 \leq \phi_{l_1} \sqrt{|Q|\log\left(\tfrac{2|\overline P|}{\delta}\right)}+  20 \enspace .
\end{equation}

\begin{proposition}\label{prop:rate_analysis_of_trisectionACP-SR}
	Let $\overline P\subset[n]$ and $Q\subset[d]$.  If $\Theta(\widetilde{P},Q)$ is indistinguishable in $l_1$-norm and if
	\begin{align}\label{eq:condition_signal_trisection_pca-SR}
		\|\Theta(\widetilde P,Q) - \overline \Theta(\widetilde P,Q)\|_F^2\geq  10^6 \log^3\left(\frac{6nd}{\delta \zeta_-}\right)\left[ \zeta^2 \left(\sqrt{|\widetilde P||Q|} + |\widetilde P|\right)+ |\widetilde{P}|\right] \enspace ,
	\end{align}
	then, with probability higher than $1-3\delta$, we have
	$$\|\Theta(\overline P',Q) - \overline \Theta(\overline P',Q)\|^2_F \leq \left(1 - \frac{1}{200\log^2(nd/\zeta_-)}\right) \|\Theta(\overline P,Q) - \overline \Theta(\overline P,Q)\|^2_F \enspace .$$
\end{proposition}

Equipped with these two propositions, we arrive at the counterpart of Propositions~\ref{prop:rate_analysis_of_block_sorting} and \ref{prop:rate_analysis_of_block_sorting_wm}. Recall Definition~\eqref{eq:definition_psi} of the function $\Psi$ by $\Psi(p,r,h,q)= \frac{hp\sqrt{rq}}{\zeta}\land \sqrt{pq}+p$.
\begin{proposition}\label{prop:rate_analysis_of_block_sorting_wm-SR} Assume that $\cT_t$ is a valid hierarchical sorting tree. 
	Consider a leaf $G$ of $\cT$ of type $\0$ or $\1$ at depth $t$. 
	With probability higher than $1- 5\tau_{\infty}\delta$, there exists a subset $\overline P^{\dagger}$ such that  $\overline{P}\subseteq \overline P^{\dagger}\subseteq G$ and the following property  holds. For some $r_{\cp}^{\dagger} \geq r^{\dagger} \in \cR$ and some $h^{\dagger}\in \cH$, upon writing
	$Q^{\dagger}_{WM} = Q^*_{WM}$ and $\overline Q^{\dagger}_{WM-SR} = \overline Q^*_{WM-SR}$,  we have simultaneously
	\begin{align}\label{eq:block_sorting_wm_1-SR}
		\|\left[\Theta(\overline P^{\dagger}, Q^{\dagger}_{WM-SR}) - \overline \Theta(\overline P^{\dagger}, Q^{\dagger}_{WM})\right]_{\sqrt{r^{\dagger}}h^{\dagger}}\|^2_F
		& \leq  2\cdot 10^6 \log^{3}\left(\frac{6nd}{\delta \zeta_-}\right) \left[\zeta^2\Psi(|\overline P^{\dagger}|,r^{\dagger},h^{\dagger},|\overline Q^{\dagger}_{WM-SR}|)+ |\overline P^{\dagger}|\right] \enspace  ;                     \\
		\label{eq:block_sorting_wm_2-SR}
		\|M(\overline P^{\dagger})-\overline{M}(\overline P^{\dagger})\|^2_F & \leq 16\zeta^2 + 96|\cR||\cH|\|[\Theta(\overline P^{\dagger}, Q^{\dagger}_{WM}) - \overline \Theta(\overline P^{\dagger}, Q^{\dagger}_{WM})]_{\sqrt{r^{\dagger}}h^{\dagger}}\|^2_F \enspace .
	\end{align}
\end{proposition}
The proof is analogous to that of~\Cref{prop:rate_analysis_of_block_sorting_wm}, up to some numerical constants, and is omitted.

Then, we state the counterpart of Lemma~\ref{lem:constraints_on_residual_dim} to control $|\overline Q^{\dagger}_{WM-SR}|$. In comparison to this lemma, we have an additional term $n/(prh)$.

\begin{lemma}\label{lem:constraints_on_residual_dim-SR}
	Assume that $\cT_t$ is a valid hierarchical sorting tree. 
	For any $h\in \cH$, $r\in R$, any integer $p$, any sequence $\overline{P}_v$ of subsets of $G_v$, it holds that 
		\begin{equation}
			\sum_{\overline P \in \cP^*(p)} |\overline Q^*_{WM-SR}(\overline P,h,r)| \lesssim 
			\log(d)	\left[\left\{  \frac{\sqrt{nd(r_0\vee r)}}{rh\sqrt{p}} \wedge\frac{d(r_0\vee r)}{r^2 h}\wedge  \frac{nd}{pr} \wedge \frac{n(r_0\vee r)}{p rh} \right\} + \frac{n}{prh} \right]\ .
		\end{equation}
	\end{lemma}
	Then, we apply~\Cref{prop:rate_analysis_of_block_sorting_wm-SR} to control the loss on an additional event of probability higher than $1-5\cdot2^{t}\tau_{\infty}\delta$.
	\begin{align}\nonumber
	\lefteqn{	\sum_{\overline{P} \in \overline{\cL}_{t}}\|M(\overline{P}) - \overline M(\overline{P})\|_F^2 } &\\ & \nonumber \leq 16\zeta^2|\overline{\cL}_t| + 96 |\cR||\cH|\sum_{p,h,r}\sum_{\overline P^{\dagger} \in \cP^*(p,h,r)}\|\left[\Theta(\overline P^{\dagger}, \overline Q^{\dagger}_{WM}) - \overline \Theta(\overline P^{\dagger}, \overline Q^{\dagger}_{WM})\right]_{\sqrt{r}h}\|_F^2 \\ \nonumber
	& \lesssim  \log^{5}\left(\frac{6nd}{\delta  \zeta_-}\right)\sum_{p,h,r}\sum_{\overline P^{\dagger} \in \cP^*(p,h,r)}\left[ \zeta^2\sqrt{[\frac{h^2pr}{\zeta^2}\wedge 1] p|\overline Q^*_{WM}(\overline P^{\dagger},h,r)|}+ (\zeta^2\vee 1) p\right] \\ \nonumber
		& \lesssim \log^{5.5}\left(\frac{6nd}{\delta  \zeta_-}\right)  \sum_{p,h,r}
		\left[  \zeta^2 \sqrt{ \frac{n r}{r\vee r_0}  \sum_{\overline P^{\dagger} \in \cP^*(p,h,r)}|\overline Q^*_{WM}(\overline{P}^{\dagger},h,r)\big|}+ 	(\zeta^2\vee 1)n\right]\\  
	 \nonumber 
		& \lesssim \log^{6}\left(\frac{6nd}{\delta  \zeta_-}\right)\sum_{p,h,r}\left[\zeta^2\left(n^{3/4}d^{1/4}\left(\frac{1}{(r_0\vee r)ph^2}\right)^{1/4} \wedge  n\sqrt{\frac{d}{p(r_0\vee r)}} \wedge \frac{n}{\sqrt{ph}} \wedge\sqrt{\frac{nd}{rh}}\right)  + \zeta^2n \sqrt{\frac{1}{pr_0 h}}+ (\zeta^2\vee 1)n\right]\\	\nonumber
	 & \stackrel{(a)}{\lesssim}\log^{7}\left(\frac{6nd}{\delta  \zeta_-}\right)\sum_{p,h}\left[\zeta^2\left(\frac{n^{3/4}d^{1/4}}{\zeta^{1/2}}\wedge  n\sqrt{d}\wedge \frac{nh}{\zeta}\sqrt{d} \wedge \frac{n}{\sqrt{h}}\wedge \frac{\sqrt{nd}}{\sqrt{h}} \right)+  (\zeta^2\vee 1)n\right]\\ \nonumber
	 	 &\lesssim  \log^{9}\left(\frac{6nd}{\delta \zeta_-}\right)\left[\zeta^2\left(\frac{n^{3/4}d^{1/4}}{\zeta^{1/2}}\wedge  n\sqrt{d}\wedge \frac{n^{2/3}\sqrt{d}}{\zeta^{1/3}} \wedge \frac{nd^{1/6}}{\zeta^{1/3}} \right)+  (\zeta^2\vee 1)n\right]\ ,
	\end{align}
	where, in (a), we use that $pr_0h\geq pr_0h^2\gtrsim 1$, the rest of the bounds being analogous to the proof of Proposition~\ref{prop:rate_of_trisection_WEI-SR}. This concludes the proof. 

\subsection{Proofs of the  lemmas}

\begin{proof}[Proof of Lemma~\ref{lem:analysis_of_residual_dim_SR}]
	By a union bound and arguing as in the proof of Lemma~\ref{lem:analysis_of_residual_dim}, we deduce that, with probability higher than $1-\delta$, we have simultaneously
		\begin{align}\label{eq:0_1-SR}
			\max_{r_{\cp}\in [4r, \tilde r] \cap \cR}\max_{k\in [d]}\big|\widehat{\bC}^{*(\ext)}_{k,2r_{\cp}}- \E\big[\widehat{\bC}^{*(\ext)}_{k,2r_{\cp}}\big]\big|& \leq h/32\ ; \\
			\max_{r_{\cp}\in [4r, \tilde r] \cap \cR}\max_{k\in [d]}\big|\widehat{\bDelta}^{*(\ext)}_{k,r_{\cp}}- \E\big[\widehat{\bDelta}^{(\ext)}_{k,r_{\cp}}\big]\big|&\leq h/32 \ .	\label{eq:0_2-SR}
		\end{align}
	Because of the adversarial observations, we now have 
	\begin{align*}
		\bC^{*(\ext)}_{k,2r_{\cp}} \leq & \E\left[\widehat{\bC}^{*(\ext)}_{k,2r_{\cp}}\right]\leq 	\bC^{*(\ext)}_{k,2r_{\cp}} + \frac{\overline{m}_{k+2r_{\cp}}(\cV)- \overline{m}_{k-2r_{\cp}-1}(\cV)}{2r_{\cp}} \ ,\\
		\bDelta^{*(\ext)}_{k,r_{\cp}} \leq & \E\left[\widehat{\bDelta}^{*(\ext)}_{k,r_{\cp}}\right]\leq 	\bDelta^{*(\ext)}_{k,r_{\cp}} + \frac{\overline{m}_{k+r_{\cp}}(\cV^+)- \overline{m}_{k-r_{\cp}-1}(\cV^-)}{2r_{\cp}}+\frac{\overline{m}_{k+r_{\cp}-1}(\cV^+)- \overline{m}_{k-r_{\cp}}(\cV^-)}{2r_{\cp}} \enspace .
	\end{align*}
	Combining the above bounds with~\eqref{eq:0_1-SR} and~\eqref{eq:0_2-SR} allows us to conclude. 
	\end{proof}

\begin{proof}[Proof of Proposition~\ref{prop:rate_analysis_of_pivot-SR}]
	With the notation of the proof of Proposition~\ref{prop:rate_analysis_of_pivot}, the condition~\eqref{eq:condition_l_1_pivot} is now replaced by 
	\begin{equation}\label{eq:condition_l_1_pivot-SR}
		\max_{i,j\in \overline P'} \|\Theta(\overline P')_{i, \cdot} - \Theta(\overline P')_{j, \cdot}\|_1 \leq  \Phi_{l_1} \sqrt{|Q|} + 20 \ ,
	\end{equation}
	where we used \Cref{cor:meta_analysis_of_pivot-SR} with $w=\1_{Q}$. The rest of the proof is left unchanged except that we replace $\Phi_{l_1}\sqrt{|Q|}$ by $\Phi_{l_1}\sqrt{|Q|}+20$. 
\end{proof}

\begin{proof}[Proof of Proposition~\ref{prop:rate_analysis_of_trisectionACP-SR}]
 Lemma~\ref{lem:structure_isotonic_matrices} is still true. However, Lemma~\ref{lem:concentration_pca} needs to be updated to 
\begin{lemma}\label{lem:concentration_pca-SR}
	Fix any $\delta\in (0,1)$.
	If
	\begin{equation}\label{eq:condition_pca-SR}
		\|\Theta - \overline{\Theta}\|_{\mathrm{op}}^2 \geq 6400\left[|\widetilde P|+ \zeta^2\left[\sqrt{|Q|(5|\widetilde P|+\log(6/\delta))} +7|\widetilde  P|+2\log(6/\delta)\right]\right] \enspace ,
	\end{equation}
	then, with probability higher than $1-\delta$, we have  
	$$ \|\hat v^T \left(\Theta - \overline \Theta\right) \|_2^2 \geq \frac{1}{2}\| \Theta - \overline \Theta \|_{\mathrm{op}}^2 \enspace .$$
\end{lemma}
In light of Condition~\eqref{eq:condition_signal_trisection_pca-SR}, this assumption is valid. Together with Lemma~\ref{lem:structure_isotonic_matrices}, we deduce that  there exists an event of probability higher than $1-\delta$ such that
\[
	\|\hat v^T \left(\Theta - \overline \Theta\right) \|^2_2 \geq \frac{1}{2}\|\Theta - \overline \Theta\|^2_{\mathrm{op}}\geq 
	\frac{1}{32\log^2(nd/\zeta_-)}\|\Theta - \overline \Theta\|^2_F \enspace .
\]
As the vectors $\hat{z}$ and $\hat{w}$ are defined though $Z^{(3)}$, we rather focus on $\Theta^{(3)}$. By~\eqref{eq:majoration_upper_SR_fundamental}, we have $\|{\Theta}^{(3)}-\Theta\|_{\mathrm{op}}\leq \sqrt{|\widetilde{P}|}$. 
\begin{eqnarray} \nonumber
	\|\hat v^T (\Theta^{(3)} - \overline \Theta^{(3)})\|^2_2&\geq& \|\hat v^T (\Theta- \overline \Theta)\|^2_2  - 4\|\Theta-\Theta^{(3)}\|_{\mathrm{op}} \|\Theta-\overline{\Theta}\|_{\mathrm{op}}\geq  \|\hat v^T (\Theta- \overline \Theta)\|^2_2  - 4\sqrt{|\widetilde{P}|}\|\Theta-\overline{\Theta}\|_{\mathrm{op}}\\ 
&\geq & 	 	\frac{9}{20}\|\Theta - \overline \Theta\|^2_{\mathrm{op}}\geq 	\frac{1}{36\log^2(nd/\zeta_-)}\|\Theta - \overline \Theta\|^2_F\ .	\label{eq:pca_signal_with_frobenius-SR}
\end{eqnarray}
Then, the analysis of $\hat{z}$ and $z^*$ follows the same steps as in the original proofs, - see \Cref{subsec:rate_analysis_of_trisectionACP} - the main difference being that we invoke~\eqref{eq:condition_l_1_pivot-SR} instead of~\eqref{eq:condition_l_1_pivot}. More precisely, we still have 
\begin{equation}\label{eq:lb_estw_truew-SR}
	\abs{\hat v^T (\Theta^{(3)} - \overline \Theta^{(3)}) \frac{\hat w}{\|\hat w\|_2}}^2 \geq \frac{16}{25}\|w^* \|_2^2 \enspace .
\end{equation}
and 
\begin{equation}\label{eq:wstar_on_S-SR}
	\|w^* \|_2^2 = \|z^*\|_2^2 - \sum_{l\in S^{*c}}(z^*_l)^2\ .
\end{equation}
The control of $\sum_{l\in S^{*c}}(z^*_l)^2$ is slightly different. 
\beqn
\left[\sum_{l\in S^{*c}}(z^*_l)^2\right]^2&=& \left[\sum_{l\in S^{*c}}[\hat v^T(\Theta^{(3)} - \overline \Theta^{(3)} )]_lz^*_l\right]^2\\ &\leq&  \|\left(\Theta^{(3)} - \overline \Theta^{(3)} \right)z^*_{S^{*c}}\|^2_2=
\sum_{i \in \widetilde P}\left(\sum_{l \in S^{*c}}(\Theta^{(3)}_{i,l} - \overline \theta^{(3)}_{l})z^*_{l} \right)^2 \\
&\leq &  \frac{18\zeta^2}{|\widetilde  P|^2} \log\left(\frac{2|Q|}{\delta}\right) \sum_{i \in \widetilde  P}\left(\sum_{l \in S^{*c}}\sum_{j\in \widetilde  P}|\Theta^{(3)}_{i,l} - \Theta^{(3)}_{j,l}|\right)^2\\
&\leq & \frac{18\zeta^2}{|\widetilde P|^2} \log\left(\frac{2|Q|}{\delta}\right) \sum_{i \in \widetilde P}\left(\sum_{j\in \widetilde P}\|\Theta^{(3)}_{i,\cdot} - \Theta^{(3)}_{j,\cdot}\|_1\right)^2\\
&\leq & 18\zeta^2  \log\left(\frac{2|Q|}{\delta}\right) |\widetilde{P}|\left[ \phi_{l_1} \log^{1/2}\left(\frac{2|\widetilde P|}{\delta}\right)\sqrt{Q}+22\right]^2\\
&\leq & \left[250 \zeta^2 \log\left(\frac{2|Q||\widetilde P|}{\delta}\right)(\sqrt{|\widetilde{P}||Q|}+1)+400 |\widetilde{P}|\right]^2 \enspace  ,
\eeqn
where we used~\eqref{eq:condition_l_1_pivot-SR} as well as the fact $\|\Theta^{(3)}_{i,\cdot}- \Theta_{i,\cdot}\|_1\leq 1$. 
Recall that  $z^*= \hat{v}^T(\Theta^{(3)}-\overline{\Theta}^{(3)})$.
Combining \Cref{eq:pca_signal_with_frobenius-SR}, \Cref{eq:wstar_on_S-SR}, and Condition~\eqref{eq:condition_signal_trisection_pca-SR}, we deduce that 
\[
	\|w^* \|_2^2 \geq \frac{1}{72\log^2(nd/\zeta_-)}\|\Theta - \overline \Theta\|_F^2\ ,
\]
which, together with \Cref{eq:lb_estw_truew-SR}, yields
\[
	\left\|(\Theta^{(3)} - \overline \Theta^{(3)}) \frac{\hat w}{\|\hat w\|_2} \right\|_2^2 \geq \abs{\hat v^T (\Theta^{(3)} - \overline \Theta^{(3)}) \frac{\hat w}{\|\hat w\|_2}}^2 \geq \frac{1}{120\log^2(nd/\zeta_-)}\|\Theta - \overline \Theta\|_F^2 \enspace .
\]
Then, we come back to the matrix $\Theta - \overline \Theta$ using again~\eqref{eq:majoration_upper_SR_fundamental}.
\[
	\left\|(\Theta - \overline \Theta) \frac{\hat w}{\|\hat w\|_2} \right\|_2^2\geq 
\frac{9}{10}\abs{\hat v^T (\Theta^{(3)} - \overline \Theta^{(3)}) \frac{\hat w}{\|\hat w\|_2}}^2 - 9 |\widetilde{P}|\ .
\]
Then, we apply Harris' inequality as in the original proof of the lemma to conclude that 
\begin{equation}\label{eq:lower_former_signal-SR}
	\left\|(\Theta - \overline \Theta) \frac{\hat w^+}{\|\hat w^+\|_2} \right\|_2^2\geq \left\|(\Theta - \overline \Theta) \frac{\hat w}{\|\hat w\|_2} \right\|_2^2  \geq  \frac{9}{1200\log^2(nd/\zeta_-)}\|\Theta - \overline \Theta\|_F^2 - 9p\geq \frac{1}{150\log^2(nd/\zeta_-)}\|\Theta - \overline \Theta\|_F^2 \enspace .
\end{equation}
Applying the pivot algorithm to $\hat w^+$, we deduce from Lemma~\ref{cor:meta_analysis_of_pivot-SR} that there exists an event of probability higher than $1-\delta$ such that
$$\max_{i,j\in \overline P'}\abs{\proscal<\Theta(\overline P')_{i,\cdot} - \Theta(\overline P')_{j,\cdot},\frac{\hat w^+}{\|\hat w^+\|_2}>}  \leq \phi_{l_1}\sqrt{\log\left(\tfrac{2|\overline P|}{\delta}\right)}+20\ . $$
By convexity, it follows that
\[
	\left\|[\Theta(\overline P') - \overline \Theta(\overline P')]\tfrac{\hat w^+}{\|\hat w^+\|_2} \right\|_2^2 \leq 2\phi_{l_1}^{2}\log\left(\tfrac{2|\overline P|}{\delta}\right)|\overline P'|+800|\overline P'|  \leq|\widetilde P| \left[2\phi_{l_1}^{2}\log\left(\tfrac{2|\overline P|}{\delta}\right)+ 800\right]\ .
\]
In light of Condition~\eqref{eq:condition_signal_trisection_pca}, this quantity is small compared to $\| \Theta - \overline{\Theta}\|_F^2$.
\begin{equation}\label{eq:upper_new_signal-SR}
	\| (\Theta(\overline P') - \overline \Theta(\overline P'))\tfrac{\hat w^+}{\|\hat w^+\|_2} \|_2^2 \leq\frac{1}{200\log^2(nd/\zeta_-)} \| \Theta - \overline{\Theta}\|_F^2\ .
\end{equation}
Then, we conclude from~\eqref{eq:upper_new_signal-SR} as we did from~\eqref{eq:upper_new_signal} in the original proof.

\end{proof}

\begin{proof}[Proof of Lemma~\ref{lem:concentration_pca-SR}]
	For short, we write $p=|\widetilde{P}|$. 
Since $Z^{(t)}=\Theta^{(t)}+N^{(t)}$ for $t=1,2$, the difference wih Lemma~\ref{lem:concentration_pca} is that $\Theta^{(1)}$ and $\Theta^{(2)}$ are involved in the terms $v^{T}(Z^{(1)}-\overline{Z}^{(1)})$ and $v^{T}(Z^{(2)}-\overline{Z}^{(2)})$.  
Hence, arguing as in the proof of Lemma~\ref{lem:concentration_pca}, we derive that, on an event  of probability higher than $1-3\delta$, we have simultaneously for all
	$v \in \mathbb R^{p}$ with $\|v\|_2\leq 1$ that
	\begin{align*}
		 & \Big|\|v^T(Z^{(1)} - \overline{Z}^{(1)})\|_2^2 - \|v^T(\Theta^{(1)} - \overline{\Theta}^{(1)})\|_2^2+\\ &  \frac{1}{2}\|v^T(\Theta^{(1)} - \overline{\Theta}^{(1)}-\Theta^{(2)}+\overline{\Theta}^{(2)})\|_2^2 - \frac{1}{2}\|v^T(Z^{(1)} - \overline{Z}^{(1)} - Z^{(2)} + \overline{Z}^{(2)})\|_2^2 \Big| \\ &\leq 10\zeta\left[\|\Theta^{(1)} - \overline{\Theta}^{(1)}\|_{\mathrm{op}}+ \frac{1}{2}\|\Theta^{(1)} - \overline{\Theta}^{(1)}-\Theta^{(2)}+\overline{\Theta}^{(2)}\|_{\mathrm{op}}  \right]\sqrt{2p+\log(6/\delta)}\\ &\quad  + 192 \zeta^2\left[\sqrt{q(3p+\log(6/\delta))}+ (3p+  \log(6/\delta))\right] \enspace .
	\end{align*}
By~\eqref{eq:majoration_upper_SR_fundamental}, we have $\|{\Theta}^{(t)}-\Theta\|_{\mathrm{op}}\leq \sqrt{p}$ for $t=1,2$. Hence, the above bound simplifies in
\begin{align*}
	& \Big|\|v^T(Z^{(1)} - \overline{Z}^{(1)})\|_2^2 - \|v^T(\Theta - \overline{\Theta})\|_2^2  - \frac{1}{2}\|v^T(Z^{(1)} - \overline{Z}^{(1)} - Z^{(2)} + \overline{Z}^{(2)})\|_2^2 \Big| \\ &\leq 10\zeta\left[\|\Theta- \overline{\Theta}\|_{\mathrm{op}}+4\sqrt{p}  \right]\sqrt{2p+\log(6/\delta)}  + 192 \zeta^2\left[\sqrt{q(3p+\log(6/\delta))}+ (3p+  \log(6/\delta))\right]\\ & \quad \quad + 12 p + 4\sqrt{p}\|\Theta-\overline{\Theta}\|_{\mathrm{op}} \enspace .
\end{align*}
Since we assume that  $\|\Theta - \overline{\Theta}\|_{\mathrm{op}}^2 \geq 6400\Big[p+  \zeta^2[\sqrt{q(5p+\log(6/\delta))} +7p+2\log(6/\delta)]\Big]$, we deduce that, on the same
	event, we have
	\[
		\sup_{v\in \mathbb{R}^p: \ \|v\|_2\leq 1} \Big|\|v^T(Z^{(1)} - \overline{Z}^{(1)})\|_2^2 - \|v^T(\Theta - \overline{\Theta})\|_2^2 - \frac{1}{2}\|v^T(Z^{(1)} - \overline{Z}^{(1)} - Z^{(2)} + \overline{Z}^{(2)})\|_2^2 \Big|\leq \frac{1}{4}\|\Theta - \overline{\Theta}\|_{\mathrm{op}}^2\enspace .\]
The rest of the proof  is left unchanged.

\end{proof}

\begin{proof}[Proof of Lemma~\ref{lem:constraints_on_residual_dim-SR}]
	Recall that $\overline Q^*_{WM-SR}(\overline P,h,r)$ decomposes as the union of  $\overline Q^*_{WM-SR-1}$,  $\overline Q^*_{WM-SR-2}$, and  $\overline Q^*_{WM-SR-3}$, $\overline Q^*_{WM-SR-4}$. Since $\overline Q^*_{WM-SR-1}$ is defined analogously to $\overline Q^*_{WM}$ --but with a different numerical constant--, we can argue as in the proof of Lemma~\ref{lem:constraints_on_residual_dim-SR}, which yields
	\[
		\sum_{\overline P \in \cP^*(p)} |\overline Q^*_{WM-SR-1}(\overline P,h,r)| \lesssim 
		\log(d)	\left[ \frac{\sqrt{nd(r_0\vee r)}}{rh\sqrt{p}} \wedge\frac{d(r_0\vee r)}{r^2 h}\wedge  \frac{nd}{pr} \wedge \frac{n(r_0\vee r)}{p rh} \right] \ . 
	\]
	It remains to consider the three last sets. We only focus on $\overline Q^*_{WM-SR-2}(\overline P,h,r)$, the last ones being analogous. We first focus on a single set  $\overline Q^*_{WM-SR-2}(\overline P, h, r, r_{\cp})$. If $k$ belongs to  $\overline D^*_{WM-SR-2}(\overline P, h, r, r_{\cp})$, this implies that the total variation of  $\overline m(\cV^+_{r_{\cp}})$ between $k-r_{\cp}$ and $k+r_{\cp}$ is at least $hr_{\cp}/128$. Since the total variation of $\overline m(\cV^+_{r_{\cp}})$ is at most one, there are at most $c/(hr_{\cp})$ regions of $\cQ_{r_{\cp}}$ that contain at least a point $\overline D^*_{WM-SR-2}(\overline P, h, r, r_{\cp})$, which entails that there are at most $c'\tfrac{r_{\cp}}{r}\cdot \tfrac{1}{hr_{\cp}}= c'/(hr)$ regions of $\cQ_r$  that contain at least a point $\overline D^*_{WM-SR-2}(\overline P, h, r, r_{\cp})$. Since $r_{\cp}$ takes at most a logarithmic number of values and since $|\cP^*(p)|\leq n/p$, we obtain
	\[
		\sum_{r_{\cp}} 	\sum_{\overline P \in \cP^*(p)}|\overline Q^*_{WM-SR-2}(\overline P, h, r, r_{\cp})| \lesssim \log(d)	\frac{n}{p r h} \ , 	
	\]
	which concludes the proof. 
\end{proof}

\section{Proofs for the $l_{\infty}$ loss}

\begin{proof}[Proof of Lemma~\ref{lem:l_infty}]
	Without loss of generality, we assume that $\pi^*$ is the identity. Fix any $i\in [n]$ and assume that $\hat{\pi}^{-1}(i)\neq i$. Consider for instance the case where $m= \hat{\pi}^{-1}(i) > i$. As a consequence, there are at least $l=m-i$ experts that are below $m$ in the oracle order and above $m$ in the estimated order . Denote $j$ the smallest of those experts. Hence, we have $j\leq i < m$ and  $\hat{\pi}(j)\geq \hat{\pi}(m)$. Besides, since $j\leq i\leq m$, we deduce from the bi-isotonic assumption that 
	\[
		\|M_{i,.}-M_{\hat{\pi}^{-1}(i),.}\|_2^2\leq \|M_{j,.}-M_{\hat{\pi}^{-1}(i),.}\|_2^2= \|M_{j,.}-M_{m,.}\|_2^2\ . 
	\]
	Taking the supremum over all $i$ implies that $l_{\infty}(\hat{\pi},\pi^*)\leq l_{err}(\hat{\pi},\pi^*)$. Let us turn to the second inequality. 
	Consider any $i <j$ such that $\hat{\pi}(i)> \hat{\pi}(j)$. We consider three subcases.
	\begin{enumerate}
	    \item[(i)]  If $\hat{\pi}(i)\geq j$, then we have $\|M_{i,.}-M_{j,.}\|_2^2\leq \|M_{i,.}-M_{\hat{\pi}(i),.}\|_2^2$. 
	    \item[(ii)] If $\hat{\pi}(j)\leq  i$, then $\|M_{i,.}-M_{j,.}\|_2^2\leq \|M_{\hat{\pi}(j),.}-M_{j,.}\|_2^2$. 
	    \item[(iii)] It remains to consider the case where we have $i< \hat{\pi}(j)< \hat{\pi}(i)< j$. As a consequence, for each $k\in [d]$, we have  $M_{j,k}-M_{i,k}\leq M_{j,k}-M_{\hat{\pi}(j),k}+M_{\hat{\pi}(i),k}- M_{i,k}$, which in turn implies that 
	\[
		\|M_{j,.}-M_{i,.}\|_2^2\leq 4l_{\infty}(\hat{\pi},\pi^*)\ . 
		\]
	\end{enumerate}
		Taking the supremum over all $i$ and reminding the definition of $j$ concludes the proof.
	\end{proof}

	\begin{proof}[Proof of Proposition~\ref{prp:risk_minimax_sup}]
		For $n=2$, all the losses are equal. Hence, the minimax lower bound~\eqref{eq:lower_bound_l_infty} is a straightforward consequence of the general minimax lower bound of Theorem~\ref{th:lower_bound_poisson} by a reduction to the case where $n=2$ (recall that $\zeta=1$ here) - This reduction is achieved by putting to $0$ the signal corresponding to all $n-2$ experts that do not corresponds to the $2$ experts of interest that will be most difficult to distinguish so that estimating the permutation amounts to deciphering between these two experts. Hence, we derive that 
		\[
		\inf_{\hat \pi}  \sup_{\pi^*\in \Pi_n}\sup_{M:\,  M_{\pi^{*-1}}\in \mathbb{C}_{\text{BISO}}} \mathbb E_{(\pi^*, M)} [l_{\infty}(\hat{\pi},\pi^*)]
		\geq c''  \left[\left(\frac{ d^{1/6}}{\lambda^{5/6} } \bigwedge \frac{\sqrt{d}}{\lambda} + \frac{1}{\lambda}\right)\bigwedge d\right] \ . 
		\]
		It turns out that the term $1/\lambda$ is higher than $d$ if $\lambda\leq 1/d$ and is smaller than $\frac{ d^{1/6}}{\lambda^{5/6} } \bigwedge \frac{\sqrt{d}}{\lambda}$ for larger $\lambda$'s. Hence, we can conclude that 
		\[
			\inf_{\hat \pi}  \sup_{\pi^*\in \Pi_n}\sup_{M:\,  M_{\pi^{*-1}}\in \mathbb{C}_{\text{BISO}}} \mathbb E_{(\pi^*, M)} [l_{\infty}(\hat{\pi},\pi^*)]
			\geq c''  \left[\frac{ d^{1/6}}{\lambda^{5/6} } \bigwedge \frac{\sqrt{d}}{\lambda} \bigwedge d\right] \ . 		
		\]

		Regarding the upper bound, we build upon the analysis of $\hat{\pi}_{WMP}$ in the specific case of $n=2$. Consider any fixed $i$ and $j$. With probability higher than $1- c\delta\log^{c'}[nd(\lambda \vee 1 )]$, it follows from the proof of Theorems~\ref{th:second_estimator_WM} and~\ref{th:second_estimator_WM_poisson_intro} that (i) $\hat{\pi}_{WMP}$ builds a valid hierarchical sorting tree and (ii) the set $\overline{P}\subset \{i,j\}$ built at the end of $\algoBlockSorting$ satisfies
		\beq\label{eq:upper_P}
		\|M(\overline{P})-\overline{M}(\overline{P})\|_F^2\leq 	c_1   \log^{c_2}\left(nd(\lambda\vee 1)\right)\left[\frac{ d^{1/6}}{\lambda^{5/6} }   \bigwedge \frac{\sqrt{d}}{\lambda} + \frac{1}{\lambda}\right]  \enspace  .
		\eeq
		It follows from (i) that  $(i,j)$ (resp. $(j,i)$) is added to $\cP\cC$ only if $\pi^*(i)<\pi^*(j)$ (resp. $\pi^*(i)<\pi^*(j)$). Besides, if 
		\beq\label{eq:separation_P}
		\|M_{i,.}-M_{j,.}\|_2> 2 	c_1   \log^{c_2}\left(nd(\lambda\vee 1)\right)\left[\frac{ d^{1/6}}{\lambda^{5/6} }   \bigwedge \frac{\sqrt{d}}{\lambda} + \frac{1}{\lambda}\right]  \enspace  ,
		\eeq
		then, this implies that $|\overline{P}|\leq 1$, otherwise this would contradict Equation~\eqref{eq:upper_P}. 
		
		Then, taking a union bound over all possible $(i,j)$, we deduce that there exists an event of probability higher than $1- cn^2 \delta\log^{c'}(nd(\lambda \vee 1 ))$, such that $\cP\cC$ is consistent and contains all $2$-tuples of experts that satisfy~\eqref{eq:separation_P}.

		Turning to the estimated permutation $\hat{\pi}_{PC}$, we consider any two experts such that $\pi^*(i) < \pi^* (j)$ and $\hat{\pi}_{PC}(i)> \hat{\pi}_{PC}(j)$. The latter condition implies that $\phi(i)\geq \phi(j)$. Since $\cP\cC$ is consistent, we have $\pi^*(i)\geq 1+ \phi(i)$. Define $\pi_-^*(j)$ as the number of experts $k$ that are below $j$ and are far apart from $j$ in the sense of Equation~\eqref{eq:separation_P}. We know that, under the above event, we have that $\phi(j)\geq \pi_-^*(j)$. This implies that $\pi^*(i) > \pi_-^*(j)$. As a consequence, $i$ and $j$ are not far apart in the sense of Equation~\eqref{eq:separation_P}. This implies that 
		\[
			l_{err}(\hat{\pi}_{PC},\pi^*)\leq 2 	c_1   \log^{c_2}\left(nd(\lambda\vee 1)\right)\left[\frac{ d^{1/6}}{\lambda^{5/6} }   \bigwedge \frac{\sqrt{d}}{\lambda} + \frac{1}{\lambda}\right]\enspace .
		\]
		Since $l_{err}(\hat{\pi}_{PC},\pi^*)$ is equivalent to $l_{\infty}(\hat{\pi}_{PC},\pi^*)$, this bound also holds  (with a larger constant) for the latter loss. Since $\delta$ has been chosen small enough and since the loss is always smaller than $d$, we arrive at the following risk bound 
		\[
			\E\left[l_{\infty}(\hat{\pi}_{PC},\pi^*)\right]\leq  c'_1   \log^{c'_2}\left(nd(\lambda\vee 1)\right)\left[\frac{ d^{1/6}}{\lambda^{5/6} }   \bigwedge \frac{\sqrt{d}}{\lambda} + \frac{1}{\lambda}\right]\enspace , 	
		\]
		which, in turn, implies that 
		\[
		\E[l_{\infty}(\hat{\pi}_{PC},\pi^*)] \leq c   \log^{c'}\left(nd(\lambda\vee 1)\right)\left[\frac{ d^{1/6}}{\lambda^{5/6} }   \bigwedge \frac{\sqrt{d}}{\lambda} \right]\bigwedge d  \enspace  .
		\]
		
		\end{proof}
		
\section{Proof of the Minimax lower bounds}

\subsection{Proof of Theorem~\ref{th:lower_bound_poisson}}

\subsubsection{Noiseless minimax lower bound}

Here, we shall prove the following minimax lower bound holding in the noiseless case $\zeta=0$. 
\begin{equation}\label{eq:lower_bound_noiseless}
	\cR^*[n,d,\lambda,0] \geq c \left[\frac{n}{\lambda}e^{-2\lambda} \wedge nd \right] 
\end{equation}
Obviously, the bound remains valid for general $\zeta\geq 0$.  Define the positive integer $d_-= 1\vee [\lfloor 1/\lambda \rfloor \wedge d]\leq d$.   We build a prior distribution $\nu$ of $M$ as follows. For each row $i=1,\ldots, n$, we sample $W_i\sim \mathcal{B}(1/2)$. If $\zeta_i=1$, the $i$-th row of $M$ is constant and equal to 1. if $W_i=0$, then the $i$-th row of $M$ has its $d_-$ first entries equal to $0$, while the remaining entries are equal to $1$. 

We write $\mathbf{P}$ and $\mathbf{E}$ for the corresponding marginal probability and expectations of the data $(x_t,y_t)$.  
\[
\cR^*[n,d,\lambda,0] \geq \inf_{\hat{\pi}} \mathbf{E}\left[\|M_{\hat{\pi}^{-1}}-M_{\pi^{*-1}}\|_2^2\right]\ .
\]
For each entry $i=1,\ldots, n$, we write $N_i= \sum_{t} \1_{x_t\in \{i\}\times [d_-]}$ the number of observations on the $d_-$ first columns of the $i$-th row. If $N_i\geq 1$, then the statistician knows the value of $W_i$. Conversely, if $N_i=0$, then she has no information on the value of $W_i$. Given an estimator $\hat{\pi}$, it is always possible to reduce its loss by ranking at the top the experts such that $N_i\geq 1$ and $W_i=1$, ranking below the experts such that $N_i\geq 1$ and $W_i=0$, and putting in between the experts such that $N_i=0$. Conditionally to the observations $(x_t,y_t)$, the values of $W_i$ such that $N_i=0$ are still distributed according to a Bernoulli distribution. As a consequence, for any $\hat{\pi}$ which has been rearranged as explained above, the conditional risk satisfies
\[
	\mathbf{E}\left[\|M_{\hat{\pi}^{-1}}-M_{\pi^{*-1}}\|_2^2\big| (x_t,y_t)\right]\geq d_- \times g(\sum_{i=1}^n \1_{N_i=0})\ ,
\]
where $g(k)$ corresponds to the expected number of error of $\hat{\pi}$ when there are exactly $k$ rows without any observations. Since conditionally to $\hat{\pi}$, the corresponding values of $W$ have been sampled independently as Bernoulli random variables with parameter $1/2$, we arrive at the following expression for $g(k)$:
\[
	g(k)= \sum_{i=1}^k  \P[\{W_i= 1\}\cap \{\sum_{j=1}^k W_j \leq k-i \}]+  \P[\{W_i= 0\}\cap \{\sum_{j=1}^k W_j > k-i \}]	\enspace . 
\]
We have $g(1)=0$, $g(2)= 1/2$, $g(3)= 1$. For $k\geq 4$, we focus on the $\lfloor k/4\rfloor $ first and  $\lfloor k/4\rfloor $ last entries to deduce that  
\beqn 
g(k)&\geq &\E\left[\sum_{i=1}^{\lfloor k/4\rfloor} W_i\right]P\left[\sum_{i=\lfloor k/4\rfloor+1}^k W_i\leq k/2 \right]+ \E\left[\sum_{i=k-\lfloor k/4\rfloor+1}^{k} (1-W_i)\right]P\left[\sum_{i=1}^{k-\lfloor k/4\rfloor} (1- W_i)\leq k/2 \right]\\&\geq& 0.5\lfloor \frac{k}{4}\rfloor\enspace .  \\
\eeqn 
Hence, there exists a universal constant $c>0$ such that we have $g(k)\geq c(k-1)$ for any $k\geq 1$. 
Since $N_i$ follows a Poisson distribution with parameter $\lambda d_-$, $V= \sum_{i=1}^n \1_{N_i=0}$ follows a binomial distribution with parameters $(e^{-\lambda d_-}, n)$.  
We obtain $\cR^*[n,d,\lambda,0]\geq cd_- \E[(V-1)_+]$.
If $\E[V]\geq 2$, then we simply use  $\E[(V-1)_+] \geq \E[V]/2$. If $\E[V]<2$, we use  $\E[(V-1)_+]\geq \P[V=2]= \frac{n(n-1)}{2}e^{-2\lambda d_-}(1- e^{-2\lambda d_-})^{n-2}\geq c' n^2e^{-2\lambda d_-}$. In any case, we conclude that 
\[
	\cR^*[n,d,\lambda,0]\geq c''d_- n e^{-2\lambda d_-}\enspace . 
\]
If $\lambda \leq 1/d$, then $d_-=d$, and the right hand-side is higher than $c''nde^{-2}$. If $\lambda\in [1/d,1]$, then we have $d_-\in [1/(2\lambda), 1/\lambda]$ and the right hand-side risk is higher than $c n/\lambda$. Finally, if $\lambda\geq 1$, we take $d_-=1$ and the right hand-side is higher than $c'ne^{-2\lambda}$. We have proved Equation~\eqref{eq:lower_bound_noiseless}.

\subsubsection{Proof of the remaining regimes}

Since the minimax risk is increasing with $n$ and $d$, we can assume without loss of generality that both $n$ and $d$ express as a power of $2$.

We shall first build a collection of prior distributions $\nu_{{\bf G}}$ indexed by ${\bf G}\in \boldsymbol{\mathcal{G}}$ on $M$. We denote  $\mathbf{P}_{\bf G}^{(\mathbf{full})}$ and  $\mathbf{E}_{\bf G}^{(\mathbf{full})}$ the corresponding marginal probability distributions and expectations on the data $(x_t,y_t)$. Since we aim at proving the lower bound in the Gaussian setting, we assume that the data $y_t$ is a normal random variable with mean $M_{x_t}$ and variance $\zeta^2$ conditionally on $M$ and $x_t$.
The minimax risk~\eqref{eq:minimax_risk_perm} is higher than the worst Bayesian risk. 
\begin{equation}\label{eq:starting_point }
\cR^*[n,d,\lambda, \zeta]\geq \inf_{\hat{\pi}}\sup_{{\bf G}\in  \boldsymbol{\mathcal{G}}}\mathbf{E}_{\bf G}^{\mathbf{full}}\left[\|M_{\hat{\pi}^{-1}}- M_{\pi^{*-1}}\|_F^2\right]\ . 
\end{equation}
We first spend some time defining the corresponding prior distributions before applying a sequence of reduction arguments.

\subsubsection{Construction of the Prior distribution on $M$}

Let $\tilde{n}\in [n]$ be an a power of 2 so that $n/\tilde{n}$ is an integer. From a broad perspective, the general purpose of this prior construction is to break down the permutation estimation problem into $n/\tilde{n}$ independent bisection problems of size $\tilde{n}$. We will fix the value of $\tilde{n}$ at the end of the proof. The permuted matrix $M_{\pi^{*-1}}$ will turn out to be block constant and we introduce $\tilde{d}\in [d]$ the number of blocks of questions, each of them being of size $d/\tilde{d}$. Here we assume that $\tilde{d}$ is a power of $2$ so that $d/\tilde{d}$ is an integer. 
$\tilde{d}$ will be also fixed at the end of the proof.

We introduce the staircase matrix $C$ of dimension $(n/\tilde{n})\times \tilde{d}$ such that $C_{\iota,\kappa}= \iota\tilde n/(4n) + \kappa/(4\tilde d)$. Also write $U$ for the constant $\tilde{n}\times d/\tilde{d}$ matrix whose entries are all equal to one. With this notation, the Kronecker product matrix $C\otimes U$ of size $n\times d$ is a bi-isotonic staircase matrix with blocks of size $\tilde{n}\times (d/\tilde{d})$. 

Then, we shall perturb the matrix $C\otimes U$ in order to simultaneously craft $n/\tilde{n}$ independent clustering problems of size $\tilde{n}$ each. Set $\tilde \lambda = \lambda \frac{d}{\tilde d}$ and $\lambda_0 = \tilde n \tilde \lambda$.  Let $\upsilon$ be a positive number and  let also $q$ be an integer smaller than or equal to $\tilde{d}$  and
\begin{equation} \label{eq:definition_M_lower_bound}
M = C \otimes U + \upsilon \frac{\zeta}{\sqrt{\lambda_0}} B^{(\mathbf{full})}\ , 
\end{equation}
where the random matrix $B^{(\mathbf{full})}\in \{0,1\}^{n\times d}$ is defined below.

For this purpose, we consider a collection $\mathcal{G}$ of subsets of $[\tilde{n}]$ with size $\tilde{n}/2$ that are well-separated in symmetric difference as defined by the following lemma.
\begin{lemma}\label{lem:packing}
There exists a numerical constant $c_0$ such that the following holds for any even integer $\tilde{n}$. 
There exists a collection $\mathcal G$ of subsets of $[\tilde{n}]$ with size $\tilde{n}/2$ whose satisfies $\log(|\cG|)\geq c_0 |\tilde n|$ and whose elements are $\tilde{n}/4$-separated, that is 
$|G_1 \Delta G_2| \geq \tilde n/4$ for any $G_1\neq G_2$. 
\end{lemma}
The above result is a straightforward consequence of Varshamov-Gilbert's lemma -- see e.g.~\cite{tsybakov}.

\medskip 

For each block $\iota \in [n/\tilde{n}]$, we fix a subset $G^{(\iota)}$ from $\mathcal{G}$. Then, we consider its 'translation' $G^{t(\iota)}=\{x+(\iota-1)\tilde n: x\in G^{\iota}\}$. The experts of $G^{t(\iota)}$ will correspond to the subgroup of 'higher' experts in the group $\iota$. We write ${\bf G}= (G^{t(1)}, \ldots, G^{t(n/\tilde{n})})$ and $\boldsymbol{\mathcal{G}}$ the corresponding collection of all possible ${\bf G}$. Given any such ${\bf G}$, we  shall define a prior distribution $\nu_{{\bf G}}$ on $M$.

 For $\iota\in [n/\tilde n]$, we sample uniformly a subset $Q^{(\iota)}$ of $q$ block of questions among the $\tilde{d}$ blocks. In each of these $q$ blocks, the corresponding rows of $B^{(\mathbf{full})}$ are equal to one. More formally, upon writing $\mathbf{1}_{d/\tilde{d}}$ for the constant vector of size $d/\tilde{d}$, we have  
\begin{equation}\label{eq:definition_B_Full}
B^{(\mathbf{full})} = \sum_{\iota=1}^{n/\tilde n}\mathbf 1_{G^{t(\iota)}} (Q^{(\iota)} \otimes \mathbf{1}_{d/\tilde{d}})^T \enspace .
\end{equation}
To sum up, we define a prior distribution $\nu_{\bf G}$ on $B^{(\mathbf{full})}$ (and equivalently on $M$) such that, under  $\nu_{\bf G}$, 
all the rows of $B^{(\mathbf{full})}$ that do not belong to any $G^{t(\iota)}$ are zero. All the rows belonging to the same set $G^{t (\iota)}$ are equal and block constants with $\tilde{d}$ blocks of size $d/\tilde{d}$, among which $q$ blocks are exactly equal to one. 

\begin{figure}\center
	\includegraphics[scale=0.4]{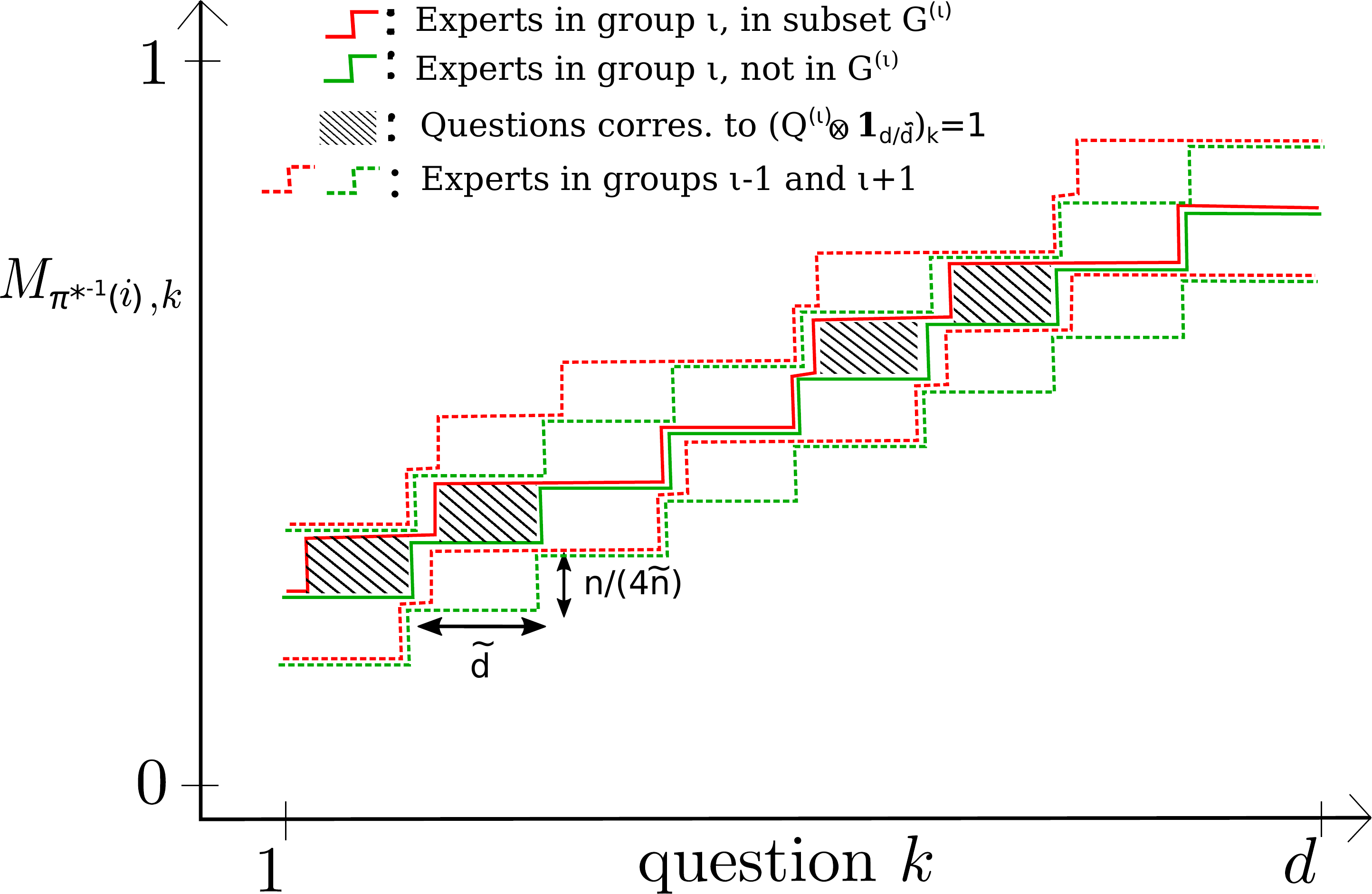}
	\caption{ Example of a matrix $M$ sampled from $\nu_{\bf G}$.
	\label{fig:LB}}
\end{figure}

\medskip 
Coming back to the matrix $M$ defined in~\eqref{eq:definition_B_Full}, we see that as soon as 
\begin{equation}\label{eq:condition_v}
2 \upsilon \zeta/\sqrt{\lambda_0} \leq \tilde n/(4n) \land 1/(4\tilde d) \ , 
\end{equation}
then, almost surely, the matrix $M$, is up to a (non-unique) permutation, bi-isotonic and its coefficients are in $[0, 1]$. Defining the subset $\overline{G}^{(\iota)}= \{i+ (\iota-1)\tilde{n}:\  i \in [\tilde{n}]\}$, we see that, under $\nu_{\bf G}$, recovering a suitable permutation $\pi^*$ is exactly equivalent to estimating the subgroup $G^{t(\iota)}\subset \overline{G}^{(\iota)}$ for each $\iota = 1,\ldots, n/\tilde{n}$.  
This construction of $M$ is illustrated in Figure~\ref{fig:LB}. To sum up, the prior distribution distribution $\nu_{\bf G}$ on $M$ requires the choice of the parameters $\tilde{n}\in [n]$, $\tilde{d}\in [d]$, the sparsity $q\in [\tilde{d}]$, and some signal level $\upsilon>0$ satisfying~\eqref{eq:condition_v}.

\medskip 

As we shall use several reduction arguments, we need to introduce some new notation. First, we respectively denote  $\mathbf{P}_{\bf G}^{(\mathbf{full})}$ and  $\mathbf{E}_{\bG}^{(\mathbf{full})}$ for the marginal probability and expectation with respect to the data when $M$ is sampled according to  $\nu_{\bG}$.

The distribution of the rows $\overline{G}^{t(\iota)}$ in $M$ under $\nu_{\bf G}$ only depends on $G^{t(\iota)}$. In what follows, we write $\nu_{G^{t(\iota)}}$ for this distribution. Similarly, we write $\mathbf{P}^{(\mathbf{full})}_{G^{t(\iota)}}$ for the corresponding marginal distribution of the observations $(x_t,y_t)$ such that $(x_{t})_1\in \overline{G}^{t(\iota)}$. By the poissonization trick, the distribution $\mathbf{P}^{(\mathbf{full})}_{\bf G}$ is a product measure of $\mathbf{P}^{(\mathbf{full})}_{G^{t(\iota)}}$ for $\iota=1,\ldots, n/\tilde{n}$. We write $\mathbf{E}^{(\mathbf{full})}_{G^{t(\iota)}}$ for the corresponding expectation.

\paragraph{Step 2: Problem Reduction} 
We start with prior distributions $\nu_{\bf G}$.  
\[
\cR^*[n,d,\lambda, \zeta]  \geq \inf_{\hat \pi} \sup_{{\bf G}\in \boldsymbol{\mathcal{G}}}\E^{(\mathbf{full})}_{\bf G}\|M_{\pi^{*-1}}-M_{\hat{\pi}^{-1}}\|_2^2
\]
For each of these matrices $M$ sampled from a distribution $\nu_{\bf G}$, it turns out that $\pi^{*}(\overline{G}^{(\iota)})=\overline{G}^{(\iota)}$. Hence, to estimate $\pi^*$, we only need to estimate each $G^{t(\iota)}\subset \overline{G}^{(\iota)}$ from the data. Intuitively, we therefore can restrict ourselves to estimators $\hat{\pi}$ satisfying  $\hat{\pi}(\overline{G}^{(\iota)})=\overline{G}^{(\iota)}$. More precisely, if an estimator $\tilde{\pi}$ does not satisfy this condition, then we can modify $\tilde{\pi}$ in $\hat{\pi}$ in order to enforce the $\overline{G}^{(\iota)}$'s to be be stable. Since, by Condition~\eqref{eq:condition_v} experts in different  $\overline{G}^{(\iota)}$ are far from each other, it turns out that the loss of $\hat{\pi}$ is smaller than that of $\tilde{\pi}$. 
\begin{align*}
	\cR^*[n,d,\lambda, \zeta]&\geq \inf_{\hat \pi:\  \hat{\pi}(\overline{G}^{(\iota)})=\overline{G}^{(\iota)} }\sup_{{\bf G}\in \boldsymbol{\mathcal{G}}}  \sum_{\iota=1}^{ n/\tilde n }  \mathbf{E}^{(\mathbf{full})}_{G^{t(\iota)}} \left[\|\big(M_{\hat \pi^{-1}} - M_{\pi^{*-1}}\big)_{\overline{G}^{(\iota)}}\|_F^2\right]\\
	& \geq \inf_{\hat \pi:\  \hat{\pi}(\overline{G}^{(\iota)})=\overline{G}^{(\iota)} } \sum_{\iota=1}^{ n/\tilde n }   \sup_{G^{t(\iota)}}\mathbf{E}^{(\mathbf{full})}_{G^{t(\iota)}} \left[\|\big(M_{\hat \pi^{-1}} - M_{\pi^{*-1}}\big)_{\overline{G}^{(\iota)}}\|_F^2\right]\\
&\geq  \sum_{\iota=1}^{ n/\tilde n } \inf_{\hat{\pi}^{(\iota)}}  \sup_{G^{t(\iota)}}\mathbf{E}^{(\mathbf{full})}_{G^{t(\iota)}} \left[\|\big(M_{\hat \pi^{(\iota)-1}} - M_{\pi^{*-1}}\big)_{\overline{G}^{(\iota)}}\|_F^2\right]\ , 
\end{align*}
where, in the last line, $\hat \pi^{(\iota)}$ stands for any estimator of the restriction $\pi^*$ to $\overline{G}^{(\iota)}$. By symmetry, we arrive at 
\begin{equation}\label{eq:lower_minimax_1}
	\cR^*[n,d,\lambda, \zeta]\geq  \frac{n}{\tilde{n}}\inf_{\hat{\pi}^{(1)}}  \sup_{G^{t(1)}}\mathbf{E}^{(\mathbf{full})}_{G^{t(1)}} \left[\|\big(M_{\hat \pi^{(1)-1}} - M_{\pi^{*-1}}\big)_{\overline{G}^{(1)}}\|_F^2\right]
\end{equation}
In summary, we have reduced the problem of estimating $\pi^*$ into the sum of $n/\tilde{n}$ problems of size $\tilde{n}$. Under $\nu_{G^{t(\iota)}}$, the restriction of $M$ to $\overline{G}^{(\iota)}$ contains $\tilde{n}/2$ good experts (those in $G^{t(\iota)}$) and $\tilde{n}/2$ bad experts.  The square Euclidean distance between these two types of experts is  $\frac{q\upsilon^2 d \zeta^2}{\lambda_0 \tilde{d}}$. If we denote $\hat{G}^{t(\iota)}$ the set of the $\tilde{n}/2$ best experts according to $\hat \pi^{(\iota)}$, then the loss writes as 
\[
	\|\big(M_{\hat \pi^{(\iota)-1}} - M_{\pi^{*-1}}\big)_{\overline{G}^{(\iota)}}\|_F^2= 	\frac{q\upsilon^2 d \zeta^2}{\lambda_0 \tilde{d}} \big|\hat{G}^{(\iota)}\Delta G^{t(\iota)}|\enspace . 
\]
Coming back to~\eqref{eq:lower_minimax_1}, we obtain 
\[
	\cR^*[n,d,\lambda, \zeta]\geq  \frac{nq\upsilon^2 d \zeta^2}{\tilde{n}\lambda_0 \tilde{d}}\inf_{\hat{G}^{(1)}}\sup_{G^{t(1)}}\mathbf{E}^{(\mathbf{full})}_{G^{t(1)}}\left[\big|\hat{G}^{(1)}\Delta G^{t(1)}|\right]\enspace . 
\]
Since all possible values of $G^{t(1)}$ are $\tilde{n}/4$-apart by definition of the collection $\cG$, we deduce that 
\[	\cR^*[n,d,\lambda, \zeta]\geq  \frac{nq\upsilon^2 d \zeta^2}{8\lambda_0 \tilde{d}}\inf_{\hat{G}^{(1)}}\sup_{G^{t(1)}}\mathbf{P}^{(\mathbf{full})}_{G^{t(1)}}\left[\hat{G}^{(1)}\neq  G^{t(1)}\right]\enspace . 	
\]
For any group $G^{t(1)}$, under $\nu_{G^{t(1)}}$, the rows of the restrictions of $M$ to $\overline{G}^{t(1)}$ are block-constant with $\tilde{d}$ blocks of $d/\tilde{d}$ questions. Consider the $\tilde{n}\times \tilde{d}$ matrices $N$ and $Y^{\downarrow}$ defined by 
\[
	N_{i,j}= \sum_{t} \1_{x_t\in \{i\}\times [(j-1)(d/\tilde{d})+1, j(d/\tilde{d})+1]} \ ;\quad \quad  Y^{\downarrow}_{i,j} =  \sum_{t} \1_{x_t\in \{i\}\times [(j-1)(d/\tilde{d})+1, j(d/\tilde{d})+1]} \left(y_t - \frac{\tilde{n}}{4n}- \frac{j}{4\tilde{d}}\right) \enspace . 
\]
To simplify the notation, we write henceforth $G$ and $\hat{G}$ for  $G^{t(1)}$ and $\hat{G}^{(1)}$ respectively. We also write $\mathbf{P}_{G}$ for the corresponding marginal distribution of $N$ and $Y^{\downarrow}$. By a sufficiency argument, it turns out that
\[
\inf_{\hat{G}} \sup_{G}	\mathbf{P}^{(\mathbf{full})}_{G}\left[\hat{G}\neq  G\right]=\inf_{\hat{G}} \sup_{G}	\mathbf{P}_{G}\left[\hat{G}\neq  G\right]\ .  
\]
Hence, we arrive at the following conclusion
\begin{equation}\label{eq:lower_minimax_G_T1}
	\cR^*[n,d,\lambda, \zeta]\geq  \frac{nq\upsilon^2 d \zeta^2}{8\lambda_0 \tilde{d}}\inf_{\hat{G}} \sup_{G}	\mathbf{P}_{G}\left[\hat{G}\neq  G\right]\enspace . 
\end{equation}
Let us introduce a third-part distribution $\mathbf{P}_0$ on $N$ and  $Y^{\downarrow}$ corresponding to the case $\upsilon=0$. Each of the entry of $N$ therefore follows an independent  Poisson distribution with parameter $\tilde{\lambda}$ and, given $N_{i,j}$, we have $Y^{\downarrow}\sim \cN(0,N_{i,j}\zeta^2)$. We then deduce from Fano's inequality~\cite{tsybakov} that  
\begin{align}\label{eq:KLFano}
	\inf_{\hat{G}} \sup_{G\in \cG} \mathbf P_{G}(\hat{G} \neq G) &\geq 1 - \frac{1 + \max_{G \in \cG} \mathrm{KL}(\mathbf P_{G}||\mathbf P_0)}{\log(|\cG|)}\ ,  
	\end{align}
	where $\mathrm{KL}(.||.)$ stands for the Kullback-Leibler divergence. Then, the following lemma bounds these Kullback-Leibler divergences.
\begin{lemma}\label{lem:KLborn}
Assume that $\lambda_0 = \tilde n \lambda d/\tilde d \geq 1$ and that $8\upsilon^2\leq 1$. For any $G\in \cG$, we have 
\begin{align*}
\mathrm{KL}(\mathbf P_G||\mathbf P_0)&\leq \frac{4\upsilon^2q^2}{\tilde d}\ .
\end{align*}
In the specific case where $\tilde{d}=q=1$, we have $\mathrm{KL}(\mathbf P_G||\mathbf P_0)=\upsilon^2/2$ for any $G\in \cG$, any $\lambda_0>0$, and any $\upsilon>0$. 
\end{lemma}

Let us summarize our findings by combining~\eqref{eq:lower_minimax_G_T1}, \eqref{eq:KLFano}, with Lemma~\ref{lem:KLborn} and the different constraints on the parameters \Cref{eq:condition_v}.
\begin{proposition}\label{prp:lower_minimax_abstract}
Provided that $\tilde{n}$, $\tilde{d}$, $q$, and $\upsilon$ satisfy the two following conditions
\begin{eqnarray}\label{eq:cond_1}
\lambda &\geq& \frac{\tilde{d}}{\tilde{n} d}\enspace ; \\
\upsilon & \leq& 2^{-3/2} \bigwedge \left[c_0 \frac{\sqrt{\tilde{d}\tilde{n}}}{q} \bigwedge  c_1 \frac{\sqrt{\lambda}}{\zeta}\left[ \frac{\tilde n^{3/2}d^{1/2}}{n \tilde{d}^{1/2} } \land \frac{\sqrt{\tilde n d}}{\tilde d^{3/2}}\right] \right]\enspace , \label{eq:cond_2}
\end{eqnarray}
then, we have 
\begin{equation}
\cR^*[n,d,\lambda,\zeta]\geq   c'' \frac{n q  \upsilon^2 \zeta^2}{\tilde{n}\lambda } \ .  
\end{equation}
\end{proposition}

In the specific case where we fix $\tilde{d}=q= 1$ and $\tilde{n}=n$, we can deduce from 
combining~\eqref{eq:lower_minimax_G_T1}, \eqref{eq:KLFano}, and the second part of~Lemma~\ref{lem:KLborn} that
\[
	\cR^*[n,d,\lambda,\zeta]\geq   c'' \frac{n\upsilon^2 \zeta^2}{\tilde{n}\lambda } \ ,   	
\]
provided that 
$\upsilon^2 \leq c'\frac{\lambda nd} {\zeta^2}\wedge n$.
By choosing $\upsilon^2$ of the order of the right-hand side, we then deduce that 
\begin{equation}\label{eq:lower_bound_n_lambda}
\cR^*[n,d,\zeta] \geq c \left[\frac{n\zeta^2}{\lambda}\wedge nd \right] \enspace . 
\end{equation}

\subsubsection{Step 3. Choice of the parameters and conclusion}

Writing $\lambda'= \lambda/\zeta^2$, recall that we aim  at proving that 
\begin{equation}\label{eq:objective_lower}
R[n,d,\lambda,\zeta]\geq c \left[	\left[ \frac{nd^{1/6}}{\lambda'^{5/6}} \bigwedge \frac{n^{3/4}d^{1/4}}{\lambda'^{3/4}}\bigwedge \frac{n^{2/3}\sqrt{d}}{\lambda'^{5/6}} \bigwedge \frac{n\sqrt{d}}{\lambda'}\right]	+\frac{n}{\lambda'}+
  \frac{n}{\lambda}e^{-2\lambda}\right]\bigwedge nd\enspace . 
\end{equation}
Since we  have proved the lower bound~\eqref{eq:lower_bound_noiseless} and \eqref{eq:lower_bound_n_lambda}, we only have to prove the corresponding minimax lower bound for the remaining four rates. For this purpose, we shall fix the values of $\tilde{n}$, $\tilde{d}$, $q$, and $\upsilon$ and apply from Proposition~\ref{prp:lower_minimax_abstract}. In the sequel we write $\lfloor x\rfloor_{dya}$ for $2^{\lfloor \log_2(x)\rfloor }$.

\medskip 

\noindent 
{\bf Case 1}: Rate $\frac{nd^{1/6}}{\lambda'^{5/6}}$. This rate can only occur if $n\leq d$,  $\lambda' \in [n^3/d, d^2]$ and $\lambda \geq 1\wedge [\lambda'^{5/6}/d^{1/6}]$. In this case, we take $\tilde{n}=2$, $\tilde{d}= \lfloor (\lambda'd)^{1/3}\rfloor_{dya}$, and $q=\lfloor \sqrt{\tilde{d}}\rfloor$. One readily checks that the conditions~\eqref{eq:cond_1} and~\eqref{eq:cond_2} are satisfied for a universal numerical value of $\upsilon$. Then, Proposition \ref{prp:lower_minimax_abstract} leads to the desired rate. 

\medskip

\noindent 
{\bf Case 2}: Rate $\frac{n^{3/4}d^{1/4}}{\lambda'^{3/4}}$. This rate can only occur if $\lambda\geq [1\wedge (n\lambda'^3/d)^{1/4}]$ and  (a) either  $n\leq d$ and $\lambda'\in [\frac{n}{d}, \frac{n^3}{d}]$ or (b) $n\in [d; d^2]$ and $\lambda'\in [\frac{n}{d},\frac{d^3}{n}]$. In this case, we take $\tilde{d}= \lfloor (\lambda' nd)^{1/4}\rfloor_{dya}$,  $\tilde{n}=\lfloor n/\tilde{d}\rfloor_{dya}$ , and $q=\lfloor \sqrt{\tilde{n}\tilde{d}}\rfloor$. One readily checks that the conditions~\eqref{eq:cond_1} and~\eqref{eq:cond_2} are satisfied for an universal numerical value of $\upsilon$. Then, Proposition~\ref{prp:lower_minimax_abstract} leads to the desired rate.

\medskip 

\noindent 
{\bf Case 3}: Rate $\frac{n^{2/3}\sqrt{d}}{\lambda'^{5/6}}$. This rate can only occur if $\lambda\geq [1\wedge \frac{\lambda'^{5/6}n^{1/3}}{\sqrt{d}}]$ and  (a) either  $n\in [d,d^2]$ and $\lambda'\in  [\frac{d^3}{n}, n^2]$ or (b) $n\geq d^2$ and $\lambda'\in [\frac{n^2}{d^3},n^2]$. In this case, we take $\tilde{n}=\lfloor (n^2/\lambda')^{1/3}\rfloor_{dya}$,  $\tilde{d}=d$, and $q=\lfloor \sqrt{\tilde{n}\tilde{d}}\rfloor$. One readily checks that the conditions~\eqref{eq:cond_1} and~\eqref{eq:cond_2} are satisfied for an universal numerical value of $\upsilon$. Then, Proposition~\ref{prp:lower_minimax_abstract} leads to the desired rate. 

\medskip 

\noindent 
{\bf Case 4}: Rate $\frac{n\sqrt{d}}{\lambda'}$. This rate can only occur if $\lambda\geq 1$ and $\lambda'\geq (n\vee d)^{2}$. In this case, we take $\tilde{n}=2$,  $\tilde{d}=d$, and $q=\lfloor \sqrt{d}\rfloor$. One readily checks that the conditions~\eqref{eq:cond_1} and~\eqref{eq:cond_2} are satisfied for a universal numerical value of $\upsilon$. Then, Proposition~\ref{prp:lower_minimax_abstract} leads to the desired rate. 
This concludes the proof.

\subsubsection{Proof of Lemma~\ref{lem:KLborn}}

\begin{proof}[Proof of Lemma~\ref{lem:KLborn}]
In order to bound the Kullback-Leibler discrepancy $\mathrm{KL}(\mathbf P_G||\mathbf P_0)$, we first observe that the rows of $N$ and $Y^{\downarrow}$ outside $G$ have the same distribution on $\mathbf P_G$ and $\mathbf P_0$. Besides, all the rows of $N$ and $Y^{\downarrow}$ in $G$ are identically distributed on $\mathbf P_G$ and on $\mathbf P_0$.
Define the vectors $\overline{N}$ and $\overline{Y}^{\downarrow}$ by $\overline{N}_j= \zeta^{-1}\sum_{i\in G}N_{i,j}$ and ${Y}_j^{\downarrow}=\zeta^{-1}\sum_{i\in G}{Y}^{\downarrow}_{i,j}$  are a sufficient statistic for deciphering $\mathbf P_G$ and $\mathbf P_0$, we have $\mathrm{KL}(\mathbf P_G||\mathbf P_0)= \mathrm{KL}(\mathbf P'||\mathbf P)$ where $\mathbf P'$ and $\mathbf P$ stand for the corresponding marginal distributions of $\overline{N}$ and $\overline{Y}^{\downarrow}$.

\medskip 

Set $u = \upsilon/\sqrt{\lambda_0}$. Under $\mathbf{P}$, given $\overline{N}$, the $\overline{Y}^{\downarrow}_j$'s are independent and satisfy $\overline{Y}^{\downarrow}_j\sim \cN(0,\overline{N}_j)$. Under $\mathbf{P}'$, conditionally to the subset $Q$ of size $q$ and conditionally to $\overline{N}$, the  $\overline{Y}^{\downarrow}_j$'s are independent and satisfy $\overline{Y}^{\downarrow}_j\sim \cN(u \overline{N}_j  \1\{j\in Q\},\overline{N}_j)$.

\medskip

In the specific case of $q=\tilde{d}=1$, we can explicitely compute the Kullback Leibler divergence. Conditionally to $\overline{N}_1=x$, $\overline{Y}^{\downarrow}$ is either distributed $\cN(0,x)$ under $\mathbf{P}$ and $\cN(ux,x)$ under $\mathbf{P}'$. Hence, their conditional Kullback-divergence is $u^2x/2$. Integrating with respect to $x$, we conclude that 
\[
	\mathrm{KL}(\mathbf P'||\mathbf P) = \mathbf{E}\left[\frac{u^2}{2} \overline{N}\right] = \frac{u^2 \lambda_0}{2} = \frac{\upsilon^2}{2}\enspace .  	
\]
We have shown the second result.

\medskip 
Let us come back to the general case. 
For $z=1,0$, define 
$$\alpha_{z}(x,y) = \frac{\lambda_0^x e^{-\lambda_0}}{x!} \frac{1}{\sqrt{2\alpha x}} \exp(- \frac{(y-uxz)^2}{2x})\enspace .$$
Then, the density of $\mathbf P$ with respect to $\mu\otimes \lambda$ where $\mu$ is the discrete measure  and $\lambda$ is the Lebesgues measure is
$\prod_j \alpha_0(\overline{N}_j, \overline{Y}^{\downarrow}_j)$. Besides, the density of $\mathbf P'$ is
$$ \int \Big[\prod_j  \alpha_{\1_{j\in Q}}(\overline{N}_j, \overline{Y}^{\downarrow}_j)\Big]d\eta(Q)\enspace ,$$
where $\eta$ stands for the uniform distribution over $\{Q \in \{0,1\}^{\tilde d}: \|Q\|_0 = q\}$. It is more convenient to first control the $\chi^2$ distance $\mathbf P$ and $\mathbf P'$. Since this distance is, up to an additive term of order 1, the second moment of the likelihood ratio between $\mathbf P$ and $\mathbf P'$, we arrive at the following
\begin{align*}
\lefteqn{\chi^2(\mathbf P',\mathbf P)+1} &\\ &= \int  \Big[\prod_{j\in Q\cap Q'}\frac{[\alpha_1(x_j,y_j)]^2}{\alpha_0(x_j,y_j)}   d\mu(x_j)dy_j\Big]\Big[\prod_{j\in Q\Delta Q'}\alpha_1(x_j,y_j) d\mu(x_j)dy_j\Big]  d\eta(Q)d\eta(Q') \\
&= \int  \Big[\prod_{j\in Q\cap Q'}\frac{[\alpha_1(x_j,y_j)]^2}{\alpha_0(x_j,y_j)}   d\mu(x_j)dy_j\Big] d\eta(Q)d\eta(Q')  \enspace  ,
\end{align*}
since $\alpha_1$ is a density. Let us work out each of these ratios.  
\begin{align*}
\int  \frac{ \alpha^2_{1}(x, y)}{\alpha_0(x, y)} dxdy &= \int \alpha_0(x,y) \exp\big[ \frac{2yux - u^2x^2}{ x}\big] d\mu(x)dy\\
&= \sum_{x=0}^{\infty}  \frac{\lambda_0^x e^{-\lambda_0}}{x!} e^{u^2x} = \exp(\lambda_0 (e^{u^2} - 1)) := \exp(\mathcal I)\enspace .
\end{align*}
Coming back to the $\chi^2$ distance, we arrive at the following equality
\begin{align*}
\chi^2(\mathbf P',\mathbf P) &=\int \exp\left( \mathcal I |Q\cap Q'|\right)d\eta(Q)d\eta(Q') - 1\enspace .
\end{align*}
Here, $|Q\cap Q'|$ is distributed as an Hypergeometric distribution with parameters  $\tilde{d}$, $q$, and $q/\tilde{d}$.
We know from 
Aldous (p.173) \cite{aldous85} that $|Q\cap Q'|$ follows the same distribution as 
the random variable $\mathbb{E}(W|\mathcal{B})$ where $W$ is a binomial random 
variable of parameters $q$, $q/\tilde{d}$ and $\mathcal{B}$ is some suitable 
$\sigma$-algebra. By Jensen's inequality, we deduce that
\begin{align*}
\chi^2(\mathbf P',\mathbf P) &\leq \mathbb E[\exp( \mathcal I W)] - 1 =  \left[1 + \frac{q}{\tilde d} (\exp( \mathcal I ) - 1)\right]^{q} - 1 \enspace .
\end{align*}
Recall that $\lambda_0 u^2 = \upsilon^2 \leq 1/8$. Hence, provided that  $\lambda_0 = \tilde n \lambda d/\tilde d \geq 1$, we have
$\mathcal I \leq 2\lambda_0  u^2 = 2\upsilon^2$. 
It then follows that 
\begin{align*}
\chi^2(\mathbf P',\mathbf P) &\leq  \exp\left(q^2/\tilde d (\exp( \mathcal I ) - 1)\right) - 1 \leq  \exp\left( 4\upsilon^2q^2/\tilde d  \right) - 1.
\end{align*}
To conclude, we use the classical bound $\mathrm{KL}(\mathbf P'||\mathbf P)\leq \log\left(1+ \chi^2(\mathbf P',\mathbf P)\right)$ --see e.g.~\cite{tsybakov}. This leads us to
\[
	\mathrm{KL}(\mathbf P'||\mathbf P)\leq \frac{4\upsilon^2q^2}{\tilde d }	\ . 
\]

\end{proof}

\subsection{Proof of Theorem~\ref{thm:lower_bound_full_observation}}

Fix $n$, $d$, $\zeta$, and $\kappa\geq 2$, and assume that, for some $c'$, there exists an estimator $\hat{\pi}$ satisfying
\beq\label{eq:contradiction_lower_full}
	  \sup_{\pi^*\in \Pi_n}\sup_{M:\,  M_{\pi^{*-1}}\in \mathbb{C}_{\text{BISO}}} \mathbb E_{(\pi^*, M)} \|M_{\hat \pi^{-1}} - M_{\pi^{*-1}}\|_F^2
	\leq c'\left[\log^{-\kappa}(nd/\zeta_-) \cR_F[n,d,\zeta]\bigwedge nd\right]\ ,
\eeq
with $\Upsilon= \lfloor \log^{\kappa}(nd/\zeta_-)\rfloor $ samples.  

Let us show that this bound would contradict the minimax lower bound in the Poisson setting. Fix $\lambda = \tfrac{112}{3} \log^{\kappa}(nd/\zeta_-)$ and consider the model~\eqref{eq:model_partial}. 
Define the estimator $\tilde{\pi}$ such that $\tilde{\pi}=\hat{\pi}$ under the event $\cA$ such that there are at least $\Upsilon$ observations on each entry and $\tilde{\pi}$ is computed arbitrarily otherwise. By~\eqref{eq:contradiction_lower_full},  $\tilde\pi$ satisfies 
\beq\label{eq:contradict_2}
	\mathbb E_{(\pi^*, M)} \|M_{\hat \pi^{-1}} - M_{\pi^{*-1}}\|_F^2\leq nd\mathbb{P}[\cA^c]+ c'\left[\log^{-\kappa}(nd/\zeta_-) \cR_F[n,d,\zeta]\bigwedge nd\right]\ . 
\eeq
By Chernoff'inequality for Poisson random variable, we deduce that 
\beqn
	\P[\cA^c]&\leq& nd \exp\left[-\frac{3}{28}\lambda \right]\leq nd e^{-4\log^{\kappa}(nd/\zeta_-)}\leq \frac{\zeta_-^2}{nd} e^{-4\log^{\kappa}(nd/\zeta_-)+ 2\log(nd/\zeta_-)}
\eeqn 	
There exists a constant $c_0$ such that for any $\kappa\geq 2$, $e^{-4x^\kappa+2x}\leq \frac{c_0}{x^{2\kappa}}$. We deduce that
\beqn 
\P[\cA^c] 
	&\leq & \frac{\zeta_-^2}{nd} \frac{c_0}{\log^{2\kappa}(nd/\zeta_-)}\ ,
\eeqn
where we used that $e^{x}\geq 1+ x^{\beta}/\beta$ for any $x\geq 0$ and any $\beta>0$ and that $\kappa\geq 2$. 
We then deduce from~\eqref{eq:contradict_2} that 
\beq\label{eq:contradict_3}
	\mathbb E_{(\pi^*, M)} \|M_{\hat \pi^{-1}} - M_{\pi^{*-1}}\|_F^2\leq  \left(c'+ \frac{c_0}{\log^{\kappa}(nd/\zeta_-)}\right)\left[\log^{-\kappa}(nd/\zeta_-) \cR_F[n,d,\zeta]\bigwedge nd\right]\ . 
\eeq
For $\lambda \geq 1$, $\cR_F[n,d,\zeta/\sqrt{\lambda}]\geq \frac{\cR_F[n,d,\zeta]}{\lambda}$. We deduce that 
\[
	\mathbb E_{(\pi^*, M)} \|M_{\hat \pi^{-1}} - M_{\pi^{*-1}}\|_F^2\leq  \frac{112}{3}\left(c'+ \frac{c_0}{\log^{\kappa}(nd/\zeta_-)}\right)\left[\cR_F[n,d,\zeta/\sqrt{\lambda}]\bigwedge nd\right]\ . 	
\]
Taking $c'$ small enough compared to the numerical constant $c$ in Theorem~\ref{th:lower_bound_poisson} contradicts this last theorem provided that $nd/\zeta_-$ is larger than some some numerical constant. Hence, no estimator can achieve~\eqref{eq:contradiction_lower_full} for this constant $c'$ when $(nd/\zeta_-)$ is large enough. 

 \medskip 

It remains to consider the case where $nd/\zeta_-$ is smaller than some constant $c''\geq 2$. We only need to prove that the minimax risk is lower bounded by $\frac{c_0}{\Upsilon}$ where $\Upsilon$ is the sample size. Since the minimax risk is non-decreasing with respect to $n$, $d$, and $\zeta$, we only have to consider the case $n=2$, $d=1$, $\zeta=2/c''$. Define $a= \zeta/\sqrt{\Upsilon}$. Consider a problem where either  $M= (a,0)^T$ or $M= (0,a)^T$. Then, with positive probability, no test is able to distinguish both hypotheses and the risk of any estimator is at most of the order $a^2=\zeta^2/\Upsilon$. The result follows.

\bibliographystyle{abbrv}
\bibliography{biblio}

\end{document}

%% file: style.tex
\usepackage{bm}
\usepackage[utf8]{inputenc}
\usepackage{mathrsfs}
\usepackage[T1]{fontenc}
\usepackage[french,english]{babel}

\usepackage{comment}
\usepackage{listings}
\usepackage{algorithmicx}
\usepackage[ruled]{algorithm}
\usepackage{algpseudocode}
\usepackage{algpascal}
\usepackage{algc}

\usepackage{bigints}
\usepackage{amsmath,amsfonts,amssymb}
\usepackage{MnSymbol}

\usepackage{graphicx}
\usepackage[left=3cm,right=3cm,top=2cm,bottom=2cm]{geometry}

\usepackage{hyperref} 
\hypersetup{
    colorlinks=true,
    linkcolor=blue,
    filecolor=magenta,
    urlcolor=cyan,
}

\usepackage{ulem}

\usepackage{amsthm}

\usepackage{cleveref} 
\crefformat{equation}{(#2#1#3)}

\theoremstyle{definition}
\newtheorem{definition}{Definition}
\theoremstyle{plain}
\newtheorem{theorem}{Theorem}[section]
\newtheorem{corollary}[theorem]{Corollary}

\newtheorem{proposition}[theorem]{Proposition}
\newtheorem{property}{Property}
\newtheorem{lemma}[theorem]{Lemma}

\theoremstyle{remark}

\usepackage[french]{babel}
\usepackage[hang,small]{caption}
\usepackage{boxedminipage}
\usepackage{dsfont}

\usepackage[usenames,dvipsnames]{xcolor}
\usepackage[textwidth=1.5cm, textsize=scriptsize]{todonotes}
\usepackage{pdfpages}

\definecolor{darkgreen}{rgb}{0,0.6,0}

\DeclareMathOperator{\Prob}{\mathbb{P}}

\def\beq{\begin{equation}}
  \def\eeq{\end{equation}}
  \def\beqn{\begin{eqnarray*}}
  \def\eeqn{\end{eqnarray*}}
  \def\bitem{\begin{itemize}}
  \def\eitem{\end{itemize}}
  \def\benum{\begin{enumerate}}
  \def\eenum{\end{enumerate}}
  \def\bmult{\begin{multline*}}
  \def\emult{\end{multline*}}
  \def\bcenter{\begin{center}}
  \def\ecenter{\end{center}}

\DeclareMathOperator*{\argmax}{arg\, max}
\newcommand{\argmin}{\mathop{\mathrm{arg\,min}}}

\def\cA{\mathcal{A}}

\def\cC{\mathcal{C}}
\def\cD{\mathcal{D}}

\def\cG{\mathcal{G}}
\def\cH{\mathcal{H}}

\def\cL{\mathcal{L}}

\def\cN{\mathcal{N}}

\def\cP{\mathcal{P}}
\def\cQ{\mathcal{Q}}
\def\cR{\mathcal{R}}
\def\cS{\mathcal{S}}
\def\cT{\mathcal{T}}
\def\cU{\mathcal{U}}
\def\cV{\mathcal{V}}

\def\cY{\mathcal{Y}}
\def\cZ{\mathcal{Z}}

\def\nw{w}

\def\nF{F}

\def\nI{I}

\def\nL{L}

\def\nO{O}

\def\nQ{Q}

\def\nU{U}

\def\bC{\mathbf{C}}

\def\bG{\mathbf{G}}

\def\bp{\mathbf{p}}

\def\1{{\mathbf 1}}
\def\0{{\mathbf 0}}

\newcommand{\bDelta}{{\boldsymbol\Delta}}

\def\bbP{\mathbb{P}}

\def\bbR{\mathbb{R}}

\def\bbZ{\mathbb{Z}}

\newcommand{\E}{\operatorname{\mathbb{E}}}
\renewcommand{\P}{\operatorname{\mathbb{P}}}

\newcommand{\Log}[1]{\log \left( #1 \right) }

\newcommand{\floor}[1]{\left\lfloor#1\right\rfloor} 

\newcommand{\abs}[1]{\left\lvert#1\right\rvert}

\def\proscal<#1,#2>{\langle #1,#2\rangle}
\newcommand{\poubelle}[1]{}

\newcommand{\cp}{{\mathrm{cp}}}

\newcommand{\PCA}{{\mathrm{pca}}}
\newcommand{\pca}{{\mathrm{pca}}}
\newcommand{\ext}{{\mathrm{ext}}}

\newcommand{\loc}{{\mathrm{loc}}}

\newcommand{\tris}{{\mathrm{tris}}}

\newcommand{\WM}{\mathrm{WM}}

\newcommand{\peel}{{\mathrm{thres}}}
\newcommand{\Capa}{\mathrm{W}_{\infty, 1}}

\newcommand{\env}{\mathrm{env}}

\newcommand{\algoPivot}{\mathbf{Pivot}}
\newcommand{\algoEncodeMatrix}{\mathbf{Encode-Matrix}}

\newcommand{\algoEncodeSet}{\mathbf{Encode-Set}}

\newcommand{\algoDimensionReduction}{\mathbf{DimensionReduction}}
\newcommand{\algoDimensionReductionWM}{\mathbf{DimensionReduction-WM}}

\newcommand{\algoDoubleTrisectionLocal}{\mathbf{DoubleTrisection-Local}}
\newcommand{\algoDoubleTrisection}{\mathbf{DoubleTrisection}}
\newcommand{\algoDoubleTrisectionPCA}{\mathbf{DoubleTrisection-PCA}}
\newcommand{\algoDoubleTrisectionWEI}{\mathbf{DoubleTrisection-WM}}
\newcommand{\algoBlockSorting}{\mathbf{BlockSort}}

\newcommand{\algoOrder}{\mathbf{Order}}
\newcommand{\algoTreeSorting}{\mathbf{TreeSort}}
\newcommand{\algoAddChild}{\mathbf{AddChild}}

\newcommand{\algoLeaf}{\mathbf{Leaf}}
\newcommand{\algoDepth}{\mathbf{Depth}}

\newcommand{\algoType}{\mathbf{Type}}
\newcommand{\Nodes}{\mathbf{Nodes}}

\newcommand{\Leaves}{\mathbf{Leaves}}

\newcommand{\Noise}{E}
\newcommand{\NOISEC}{N}
\newcommand{\spaceAnd}{\quad \text{ and } \quad}